\documentclass[authoryear, 11pt]{imsart}
\RequirePackage[OT1]{fontenc}
\RequirePackage{amsthm,amsmath,natbib}
\RequirePackage{hypernat}
\startlocaldefs
\numberwithin{equation}{section}
\theoremstyle{plain}

\endlocaldefs
\usepackage{psfrag, epsfig, graphicx}
\usepackage{amsfonts,amsmath,latexsym,amssymb}
\usepackage[active]{srcltx}
\usepackage{fullpage}
\newtheorem{lemma}{Lemma}
\newtheorem{theorem}{Theorem}
\newtheorem{proposition}{Proposition}

\newtheorem{assumption}{Assumption}

\newtheorem*{ass4.1}{Assumption $4^\prime$}
\newtheorem*{ass4.2}{Assumption $4^{\prime\prime}$}

\newtheorem*{assGauss}{Assumption 3 [Gaussian case]}

\newtheorem{remark}{Remark}
\newtheorem{definition}{Definition}
\newtheorem{corollary}{Corollary}

\renewcommand{\kappa}{\varkappa}

\newcommand{\te}{\theta}
\newcommand{\vth}{\vartheta}

\newcommand{\ve}{\epsilon}
\newcommand{\e}{\varepsilon}
\newcommand{\rf}{\mathrm{f}}

\newcommand{\rd}{\mathrm{d}}
\newcommand{\ra}{\mathrm{a}}
\newcommand{\rb}{\mathrm{b}}
\newcommand{\ri}{\mathrm{i}}

\newcommand{\rk}{\mathrm{k}}

\newcommand{\cA}{{\cal A}}
\newcommand{\cB}{{\cal B}}
\newcommand{\cC}{{\cal C}}

\newcommand{\cE}{{\cal E}}
\newcommand{\cF}{{\cal F}}
\newcommand{\cG}{{\cal G}}
\newcommand{\cH}{{\cal H}}
\newcommand{\cI}{{\cal I}}

\newcommand{\cL}{\mathcal{L}}
\newcommand{\cM}{{\cal M}}
\newcommand{\cN}{{\cal N}}

\newcommand{\cQ}{{\cal Q}}
\newcommand{\cR}{{\cal R}}
\newcommand{\cS}{{\cal S}}

\newcommand{\cU}{{\cal U}}
\newcommand{\cV}{{\cal V}}

\newcommand{\cX}{{\cal X}}

\newcommand{\cZ}{{\cal Z}}
\newcommand{\bB}{\mathbb B}
\newcommand{\bC}{\mathbb C}

\newcommand{\bE}{\mathbb E}

\newcommand{\bH}{\mathbb H}

\newcommand{\bK}{{\mathbb K}}
\newcommand{\bL}{{\mathbb L}}

\newcommand{\bN}{{\mathbb N}}

\newcommand{\bP}{{\mathbb P}}

\newcommand{\bR}{{\mathbb R}}
\newcommand{\bS}{{\mathbb S}}
\newcommand{\bT}{{\mathbb T}}

\newcommand{\bX}{{\mathbb X}}

\newcommand{\bZ}{{\mathbb Z}}
\newcommand{\mA}{\mathfrak{A}}
\newcommand{\mB}{\mathfrak{B}}

\newcommand{\mH}{\mathfrak{H}}
\newcommand{\mS}{\mathfrak{S}}
\newcommand{\mL}{\mathfrak{L}}
\newcommand{\ms}{\mathfrak{s}}

\newcommand{\mh}{\mathfrak{h}}
\newcommand{\md}{\mathfrak{d}}

\newcommand{\mw}{\mathfrak{w}}

\newcommand{\mm}{\mathfrak{m}}
\newcommand{\mn}{\mathfrak{n}}

\newcommand{\mt}{\mathfrak{t}}
\newcommand{\mz}{\mathfrak{z}}
\newcommand{\mT}{\mathfrak{T}}
\newcommand{\mX}{\mathfrak{X}}
\newcommand{\mZ}{\mathfrak{Z}}
\newcommand{\mR}{\mathfrak{R}}
\newcommand{\ma}{\mathfrak{a}}
\newcommand{\mE}{\mathfrak{E}}
\newcommand{\mb}{\mathfrak{b}}

\newcommand{\vs}{\varsigma}
\newcommand{\z}{\zeta}
\newcommand{\D}{\Delta}
\newcommand{\epr}{\hfill\hbox{\hskip 4pt
                \vrule width 5pt height 6pt depth 1.5pt}\vspace{0.5cm}\par}

\begin{document}
\begin{frontmatter}
\title{Upper functions
for positive random
functionals}
\runtitle{Upper functions}
%\thankstext{T1}{}
\begin{aug}
\author{\fnms{Oleg} \snm{Lepski}
%\thanksref{t3}
\ead[label=e2]{lepski@cmi.univ-mrs.fr}}
%\and
%\author{\fnms{Third} \snm{Author}\thanksref{t1}
%\ead[label=e3]{third@somewhere.com}
%\ead[label=u1,url]{http://www.foo.com}}
%\thankstext{t1}{Some comment}
%\thankstext{t2}{Supported by the ISF grant 389/07}
%\thankstext{t3}{Supported by the ANR grant 0000}
%\thankstext{t3}{Second supporter of the project}
\runauthor{ O. Lepski}

%\affiliation{University of Haifa and Universit\'e Aix--Marseille I}

\address{Laboratoire d'Analyse, Topologie, Probabilit\'es\\
 Universit\'e Aix-Marseille  \\
 39, rue F. Joliot-Curie \\
13453 Marseille, France\\
\printead{e2}\\ }
\phantom{E-mail:\ }
\end{aug}

\begin{abstract}
In this paper we
are interested in finding upper functions for a collection real-valued random variables
$\big\{\Psi\big(\chi_\theta\big), \theta\in \Theta\big\}$.
Here
$\{\chi_\theta, \theta\in \Theta\}$ is the family of continuous  random mappings,
$\Psi$ is a given sub-additive positive functional and $\Theta$ is a totally bounded subset of a metric space.
We seek a non-random function $U:\Theta\to\bR_+$ such that
$
\sup_{\theta\in \Theta} \big\{\Psi\big(\xi_\theta\big)-U(\theta)\big\}_+
$
is "small" with prescribed  probability. We apply the developed  machinery to the variety of  problems
related to gaussian random functions and empirical processes.

\end{abstract}

\begin{keyword}[class=AMS]
\kwd[Primary ]{60E15}
\kwd[; secondary ]{62G07, 62G08}
\end{keyword}

\begin{keyword}
\kwd{upper function}
\kwd{empirical processes}
\kwd{gaussian random function}
\kwd{metric entropy}
\kwd{doubling measure}
\end{keyword}
\end{frontmatter}
%\date{}
%\author{
%Alexander Goldenshluger  \and Oleg Lepski}
%
%
%\maketitle

%%%%%%%%%%%%%%%%%%%%
%\setlength{\evensidemargin}{-0.5cm}
\def\huh{\hbox{\vrule width 2pt height 8pt depth 2pt}}
\def\blacksquare{\hbox{\vrule width 4pt height 4pt depth 0pt}}
\def\square{\hbox{\vrule\vbox{\hrule\phantom{o}\hrule}\vrule}}
\def\inter{\mathop{{\rm int}}}
\tableofcontents
%\ref{sec:introduction}
%
\section{Introduction}
\label{sec:introduction}

The main objective of this paper is to look from the unique point of view at some phenomena  arising in different areas of  probability theory
and mathematical statistics. We will try to understand what is common between   classical probabilistic results, such as the law of iterated logarithm for example, and
 well-known problem in adaptive estimation called price to pay for adaptation. Why exists two different kinds of this price?  What relates exponential inequalities for $M$-estimators, so-called uniform-in-bandwidth
consistency in density or regression model  and the bounds for modulus of continuity of gaussian random functions defined on a metric space equipped with doubling measure?

It turned out that  all these  and many others problems  can be reduced to  following one. Let $\mT$ be  a set  and let  $\left(\Omega,\mB,\mathrm{P}\right)$ be a complete  probability space.
Let  $\chi$ defined on $\mT\times\Omega$ be a given  $\mB$-measurable  map into linear metric space $\mS$  and
let  $\Psi:\mS\to\bR_+$ be a given {\sl continuous  sub-additive} functional. Let $\Theta\subset\mT$ and suppose that we for any
$\theta\in\Theta$ and  $z>0$ one can find $U(\theta,z)$ such that
\begin{eqnarray}
\label{eq1:introduction}
&&\mathrm{P}\left\{\left[\Psi(\chi_\theta)-U(\theta,z)\right]>0\right\}\leq ce^{-z},\;\;c>0.
\end{eqnarray}
Assuming additionally that $\lambda U(\cdot,z)\geq  U(\cdot,\lambda z)$ for any $z>0,\;\lambda\geq 1$ we also have for any $q\geq 1$
\begin{eqnarray}
\label{eq2:introduction}
&&\mathrm{E}\left\{\left[\Psi(\chi_\theta)-U(\theta,z)\right]_+\right\}^q\leq c \Gamma(q+1)\big[U(\theta,1)\big]^qe^{-z},\quad \forall z\geq 1,
\end{eqnarray}
where $\Gamma$ is gamma-function and $[a]_+$ is a positive part of $a$.

The problem which we address now consists in a finding of $\mathbf{U}(\theta,z)$ and  $\boldsymbol{U}(\theta,z)$ satisfying
\begin{eqnarray}
\label{eq3:introduction}
&&\mathrm{P}\left\{\sup_{\theta\in\Theta}\left[\Psi(\chi_\theta)-\mathbf{U}(\theta,z)\right]>0\right\}\leq \mathbf{c}e^{-z};\quad\forall z\geq 1
\\
\label{eq4:introduction}
&&\mathrm{E}\bigg\{\sup_{\theta\in\Theta}\left[\Psi(\chi_\theta)-\boldsymbol{U}_q(\theta,z)\right]_+\bigg\}^q\leq \mathbf{c}_q\left[\inf_{\theta\in\Theta}\boldsymbol{U}_q(\theta,1)\right]^qe^{-z},\quad\forall z\geq 1,
\end{eqnarray}
where $\mathbf{c}$ and $\mathbf{c}_q$ are numerical constants. If (\ref{eq3:introduction}) and (\ref{eq4:introduction}) hold we will say that
$\mathbf{U}(\cdot,\cdot)$ and $\boldsymbol{U}_q(\cdot,\cdot)$ are {\it upper functions } for the collection of random variables $\left\{\Psi(\chi_\theta),\;\theta\in\Theta\right\}$.

The main questions on which we would like to answer are the following.
\begin{itemize}
\item
Do $\mathbf{U}(\cdot,\cdot)$ and $\boldsymbol{U}_q(\cdot,\cdot)$ coincide with
$U(\cdot,\cdot)$ up to numerical constants or there is a  "price to pay" for passing from pointwise (
%for fixed $\theta\in\Theta$)
results  \ref{eq1:introduction})--(\ref{eq2:introduction}) to uniform ones given in (\ref{eq3:introduction})--(\ref{eq4:introduction})?
\item
 Do $\mathbf{U}(\cdot,\cdot)$ and $\boldsymbol{U}_q(\cdot,\cdot)$ coincide up to numerical constants? In other words should  one  to pay the same price for the probability and moment's bounds?

\end{itemize}

We will show that a payment exists and in general $\boldsymbol{U}_q(\cdot,\cdot)\gg \mathbf{U}(\cdot,\cdot)\gg U(\cdot,\cdot)$. Thus,
we will seek  $\mathbf{U}(\cdot,\cdot)$ and $\boldsymbol{U}_q(\cdot,\cdot)$  satisfying   (\ref{eq3:introduction}) and (\ref{eq4:introduction}) and  "minimally" separated away from $U(\cdot,\cdot)$.
We will realize this program under the following condition.
\begin{assumption}
\label{ass:fixed_theta_local}
\begin{enumerate}
\item
There exist $A:\mT\to \bR_+$,    $B:\mT\to \bR_+$ and
$\textsf{c}>0$
such that $\forall z>0$
\begin{equation}
\label{eq1:ass_fixed_theta_local}
\mathrm{P}\left\{\Psi(\chi_\theta)\geq z\right\}\leq \mathrm{c}\exp{\left\{-\frac{z^{2}}{A^{2}(\theta)+ B(\theta)z}\right\}},\;\;
\forall\theta\in\Theta.
\end{equation}
\item
There exist $\mathrm{a}:\mT\times\mT\to \bR_+$ and  $\mathrm{b}:\mT\times\mT\to \bR_+$ such that $\forall z>0$
\begin{equation}
\label{eq2:ass_fixed_theta_local}
\mathrm{P}\left\{\Psi(\chi_{\theta_1}-\chi_{\theta_2})\geq z\right\}\leq \mathrm{c}\exp{\left\{-\frac{z^{2}}{\mathrm{a}^{2}(\theta_1,\theta_2)+ \mathrm{b}(\theta_1,\theta_2)z}\right\}},\;\;
\forall\theta_1,\theta_2\in\Theta.
\end{equation}
\end{enumerate}
\end{assumption}
\begin{remark}
\label{rem:linearity_after_ass1_local}
If Assumption \ref{ass:fixed_theta_local} ({\it 1}) holds on $\mT$ (not only on $\Theta$), $\mT$ is linear space and if, additionally, the map $\chi_t$ is linear on $\mT$ , then the Assumption \ref{ass:fixed_theta_local} ({\it 2}) is automatically fulfilled  since  one can take $\mathrm{a}(\mt_1,\mt_2)=A(\mt_1-\mt_2)$ and $\mathrm{b}(\mt_1,\mt_2)=B(\mt_1-\mt_2),\;\mt_1,\mt_2\in\mT$.

\end{remark}
\begin{remark}
We can easily deduce from   (\ref{eq1:ass_fixed_theta_local})  that for any $\theta\in\Theta$
\begin{eqnarray}
\label{eq100:discussion-price-to-pay}
&&\mathrm{P}\Big\{\Psi(\chi_\theta)\geq A(\theta)\sqrt{z}+B(\theta)z\Big\}\leq \mathrm{c}\exp{\left\{-z\right\}},\;\;
\forall z\geq 0;
%\nonumber
\\*[2mm]
\label{eq2:discussion-price-to-pay}
&&\mathrm{E}\Big\{\Psi(\chi_\theta)- \Big[A(\theta)\sqrt{z}+B(\theta)z\Big]\Big\}_+^{q}\leq \mathrm{c}\Gamma(q+1)\Big[A(\theta)+ B(\theta)\Big]^{q}\exp{\left\{-z\right\}},\;\;\forall z\geq 1.
\end{eqnarray}
Therefore, (\ref{eq1:introduction})--(\ref{eq2:introduction}) hold with $U(\theta,z)=A(\theta)\sqrt{z}+B(\theta)z$.
\end{remark}

 Assumption \ref{ass:fixed_theta_local} is not  new. In particular,  it can be found in slightly different form  in \cite{wellner},  \cite{Tal-book}, where this assumption is used for deriving the bound for
$\mathrm{E}\left[\sup_{\theta\in\Theta}\Psi(\chi_\theta)\right]$. The usual technique is based on the {\em chaining argument}
available in view of (\ref{eq2:ass_fixed_theta_local}).
It is worth mentioning that
uniform probability and moment bounds for $\left[\sup_{\theta\in\Theta}\Psi(\chi_\theta)\right]$ in the case where
$\chi_\theta$ is
 empirical or gaussian  process are a subject of vast literature, see, e.g.,
\cite{alexander},  \cite{Talagrand}, \cite{Lif}, \cite{wellner}, \cite{van-de-Geer},
\cite{massart}, \cite{bousquet},
\cite{gine-kolt} among many others.
Such bounds play an important role in establishing
the laws of iterative logarithm and central limit
theorems [see, e.g., \cite{alexander} and
\cite{gine-zinn}].

However much less attention was paid to finding of  upper functions.  The majority of the papers, where such problems are considered,
contains asymptotical results, see, i.e. \cite{Kalina}, \cite{Qua&Vatan}, \cite{Bobkov}, \cite{Shiryaev} and references therein.
We would like especially mention the paper \cite{Ostrovskii}, where upper function satisfying the inequalities similar to
 (\ref{eq3:introduction}), was obtained  for the modulus of continuity of random fields satisfying the  Cramer
 condition.

The researches carried out in the present paper complete the investigations done in  \cite{GL3}, where the upper functions as well as inequalities (\ref{eq3:introduction})--(\ref{eq4:introduction}) were obtained under following condition:
 $\chi_t$ is linear and
there are  $A:\mT\to\bR_+$, $\;B:\mT\to\bR_+$, $\;V:\mT\to\bR_+$
such that
\begin{equation*}
\label{eq:individual-prob}
\mathrm{P}\left\{\Psi(\chi_\mt)-V(\mt)\geq z\right\}\leq g\bigg(\frac{z^2}{A^{2}(\mt)+B(\mt)z}\bigg),\quad \forall \mt\in\mT,
\end{equation*}
where $g:\bR_+\to\bR_+$ is a strictly  decreasing to zero function. We note that if $g(x)=e^{-x}$ and $V\equiv 0$ this assumption coincides with
(\ref{eq1:ass_fixed_theta_local}) and, since $\chi_t$ is linear (\ref{eq2:ass_fixed_theta_local}) is automatically fulfilled, see Remark \ref{rem:linearity_after_ass1_local}.  In \cite{GL3} under additional assumption imposed on $A,B, V$ and $\Theta\subset\mT$ the upper functions %$\mathbf{V}$ and $\boldsymbol{V}_q$
for the collection  $\left\{\Psi(\chi_\theta),\;\theta\in\Theta\right\}$ were found. As it was shown that  they distinguish from the function $V(\cdot,\cdot)$ only by numerical constants!  The imposed assumptions do not admit the case $V\equiv 0$ that, as it was said above,  leads to completely different solution of the problem at hand.

To derive upper functions satisfying (\ref{eq3:introduction})--(\ref{eq4:introduction}) we complete Assumption \ref{ass:fixed_theta_local}
by the following conditions.
\begin{assumption}
\label{ass:metric-case-local}
%\begin{enumerate}
%\item
$\chi_{\bullet}:\mT\to\mS$ is continuous $\mathrm{P}$-a.s.
%\item

Mappings $\mathrm{a}$ and $\mathrm{b}$ are  semi-metrics on $\mT$ and
  $\Theta$ is totally bounded with respect to  $\ra\vee\rb$.

\vskip0.1cm

$\overline{A}_\Theta:=\sup_{\theta\in\Theta}A(\theta)<\infty,\; \overline{B}_\Theta:=\sup_{\theta\in\Theta}B(\theta)<\infty.$

%\end{enumerate}
\end{assumption}
Denote by $\bS$ the following set of real functions:
$$
\bS=\bigg\{s:\bR\to\bR_+\setminus \{0\}:\;\sum_{k=0}^{\infty}s\big(2^{k/2}\big)\leq 1\bigg\}.
$$
For any $\widetilde{\Theta}\subseteq\Theta$ and any  semi-metric $\rd$ on $\mT$ let
$\mathfrak{E}_{\widetilde{\Theta},\;\mathrm{d}}(\delta),\;\delta>0,$ denote the entropy of $\widetilde{\Theta}$ measured in $\mathrm{d}$.
 For any $x>0$,  $\widetilde{\Theta}\subseteq\Theta$ and  $s\in\bS$  define the quantities
\begin{eqnarray}
\label{eq3:def_local}
&&e^{(\ra)}_s\big(x,\widetilde{\Theta}\big)=\sup_{\delta>0}\delta^{-2}
\mathfrak{E}_{\widetilde{\Theta},\;\mathrm{a}}\left(x(48\delta)^{-1}s(\delta)\right),\quad
e^{(\rb)}_s\big(x,\widetilde{\Theta}\big)=\sup_{\delta>0}\delta^{-1}\mathfrak{E}_{\widetilde{\Theta},\;\mathrm{b}}\left(x(48\delta)^{-1}s(\delta)\right).
\end{eqnarray}

\begin{assumption}
\label{ass:parameter_local}
There exist $s_1,s_2\in\bS$  such that $\forall x>0$
$$
e^{(\ra)}_{s_1}\big(x,\Theta\big)<\infty,\quad e^{(\rb)}_{s_2}\big(x,\Theta\big)<\infty.
$$

\end{assumption}

\paragraph{Organization of the paper} In Section \ref{sec:Key proposition_} we construct  upper functions for
 $\left\{\Psi(\chi_\theta),\;\theta\in\Theta\right\}$ and prove for them the inequalities (\ref{eq3:introduction})--(\ref{eq4:introduction})
 under Assumptions \ref{ass:fixed_theta_local}--\ref{ass:parameter_local}. In fact we present two different constructions which will be refereed to upper functions of the first and second type (Propositions \ref{prop_uniform_local2} and \ref{prop_uniform_local3}). We also derive some consequences related to the upper functions for modulus of continuity of random real-valued mappings (Propositions \ref{prop_uniform_local4} and \ref{prop_uniform_local5} ). In Section \ref{sec:gauss} we apply Propositions \ref{prop_uniform_local3} and \ref{prop_uniform_local4}
to gaussian random functions. In Section \ref{sec:subsection-norm-gauss} we derive upper functions for  $L_p$-norm of some Wiener integrals (Theorem \ref{th:gauss-norm}) and in Section \ref{sec:subsection-gauss-local-continuity} we study the local modulus of continuity of gaussian functions defined on a metric space satisfying doubling condition (Theorem \ref{th:gauss-local-modulus}). Section \ref{sec:empirical-processes}
is devoted to the detailed  consideration of generalized empirical processes. We provide with rather general assumption (Assumption \ref{ass:bounded_case}) under which
the upper functions admit the  explicit expression, Section \ref{sec:totally-bounded-case} (Theorem \ref{th:empiric_totaly_bounded_case}) and Section
\ref{sec:partially_totally-bounded-case} (Theorem \ref{th:empiric-partially_totaly_bounded_case}). We also establish non-asymptotical versions
of the law of iterated logarithm (Theorem \ref{th:LIL-nonasym}) and the law of logarithm (Theorem \ref{th:LL-nonasym}). Section \ref{sec:localized-empiric} is devoted to the application of Theorems \ref{th:empiric_totaly_bounded_case} and \ref{th:empiric-partially_totaly_bounded_case} to empirical processes possessing some special structure, Theorems \ref{th:empirical-product-general}--\ref{th:LL-product}. Proofs of main results are given in  Sections \ref{sec:proofs-of-propositions}--\ref{sec:proofs-of-empiric-theorems} and technical lemmas are proven in Appendix.

\section{General setting}
\label{sec:Key proposition_}

 Denote by $\cS_{\ra,\rb}$ the subset of $\bS\times\bS$ for which Assumption \ref{ass:parameter_local} holds and
let  $A, B, \mathrm{a}$  and $\mathrm{b}$ be any mappings for which Assumption \ref{ass:fixed_theta_local} is
fulfilled. For any
$\vec{s}=(s_1,s_2)\in\cS_{\ra,\rb}$, any $\kappa=(\kappa_1,\kappa_2),\;\kappa_1>0,\kappa_2>0,$ and any $\widetilde{\Theta}\subseteq\Theta$ put
\begin{eqnarray}
\label{eq:def-function-e-vec-s}
e_{\vec{s}}\big(\kappa,\widetilde{\Theta}\big)=e^{(\ra)}_{s_1}\big(\kappa_1,\widetilde{\Theta}\big)+e^{(\rb)}_{s_2}\big(\kappa_2,\widetilde{\Theta}\big).
\end{eqnarray}

\subsection{Inequalities for the suprema}
Put  for any  $\widetilde{\Theta}\subseteq\Theta$, any $\e>0$ and any $y\geq 0$
\begin{eqnarray*}
%\label{eq4:def_local}
U^{(\e)}_{\vec{s}}\big(y,\kappa,\widetilde{\Theta}\big)=\kappa_1\sqrt{2\big[1+\e^{-1}\big]^{2} e_{\vec{s}}\big(\kappa,\widetilde{\Theta}\big)+y}
+\kappa_2\Big(2\big[1+\e^{-1}\big]^{2}e_{\vec{s}}\big(\kappa,\widetilde{\Theta}\big)+y\Big).
\end{eqnarray*}

\begin{proposition}
\label{prop_uniform_local1}
 Let   Assumptions \ref{ass:fixed_theta_local}-\ref{ass:parameter_local}  hold and let $\widetilde{\Theta}\subseteq\Theta$ be fixed.
 Then for any $\widetilde{\kappa}$ satisfying  $\widetilde{\kappa}_1\geq \sup_{\theta\in\widetilde{\Theta}}A(\theta)$ and $\widetilde{\kappa}_2\geq \sup_{\theta\in\widetilde{\Theta}}B(\theta)$, any
  $ \vec{s}\in\cS_{\mathrm{a},\mathrm{b}}$, $\e\in \big(0,\sqrt{2}-1\big]$ and $ y\geq 1$,
 \begin{eqnarray*}
\mathrm{P}\left\{\sup_{\theta\in\widetilde{\Theta}}\Psi\left(\chi_{\theta}\right)\geq  U^{(\e)}_{\vec{s}}\big(y,\widetilde{\kappa},\widetilde{\Theta}\big)\right\} \leq
2\mathrm{c}\exp{\left\{-y/(1+\e)^{2}\right\}}.
\end{eqnarray*}
Moreover, for any $q\geq 1$
$$
\mathrm{E}\left\{\sup_{\theta\in\widetilde{\Theta}}\Psi\left(\chi_{\theta}\right)-U^{(\e)}_{\vec{s}}\big(y,\widetilde{\kappa},\widetilde{\Theta}\big)\right\}^{q}_+ \leq
2\mathrm{c}\Gamma(q+1)\left[(1+\e)^{2}y^{-1}U^{(\e)}_{\vec{s}}
\big(y,\widetilde{\kappa},\widetilde{\Theta}\big)\right]^q\exp{\left\{-y/(1+\e)^{2}\right\}}.
$$
%where  $\Gamma$ is gamma-function.
\end{proposition}
We remark that  $\sup_{\theta\in\widetilde{\Theta}}\Psi\left(\chi_{\theta}\right)$ is $\mB$-measurable for any $\widetilde{\Theta}\subseteq\Theta$ since $\Psi$ is continuous, the mapping $\theta\mapsto\chi_{\theta}$ is
continuous $\mathrm{P}$-a.s., $\Theta$ is a totally bounded set and considered probability space is complete (see, e.g. Lemma \ref{lem:measurability} below).

%\smallskip

\paragraph{Discussion}

We will see that the Proposition \ref{prop_uniform_local1} is crucial technical tool for deriving upper functions. It contains the main ingredient
of our future construction  the quantity $e_{\vec{s}}(\cdot,\cdot)$.
The  important issue in this context  is the choice of  $\vec{s}\in\cS_{\ra,\rb}$. For many particular problems it is sufficient to choose $\vec{s}=(s^*,s^*)$, where
\begin{equation}
\label{eq:function-s-star}
s^*(x)=(6/\pi^2)\big(1+[\ln{x}]^{2}\big)^{-1},\;x\geq 0.
\end{equation}
This choice is explained by two simple reasons: its explicit description allowing to compute the quantity $e_{\vec{s}}$ in particular problem
%constants appeared in Propositions \ref{prop_uniform_local1}
and the logarithmical  decay of this function when $x\to\infty$. In view of the latter remark  we can  consider the set $\Theta$ whose entropy obeys the restriction which is closer to the minimal one (c.f. Sudakov lower bound for gaussian random functions \cite{Lif}).
We note, however, that there exist examples where $\vec{s}$ has to be chosen on more special way (see Theorem \ref{th:gauss-norm}).

Let us now discuss the  role of parameter $\e$. In most particular problems considered in the paper we will not be interested in optimization of the numerical constants involved in the description of upper functions. If so,  the choice of this parameter can be done in arbitrary way and we will put  $\e=\sqrt{2}-1$ to simplify the notations and computations. Note, however, that there are some problems (see, for instance Section \ref{sec:subsection-gauss-local-continuity}), where $\e$ must be chosen carefully. The typical requirements to this choice is $\e=\e(y)$ and
$$
\e(y)\to 0,\quad y\e^{2}(y)\to 0,\;\; y\to\infty.
$$

The bounds similar to whose presented in Proposition  \ref{prop_uniform_local1} are the subject of vast literature see, for instance, the books  \cite{Lif},  \cite{wellner} or \cite{van-de-Geer}.
Note, however, that the results presented in the proposition may have an independent interest, at least, for the problems where the quantity
 $e_{\vec{s}}(\cdot,\cdot)$ can be expressed explicitly.
In this case under rather general conditions it is possible, putting  $\widetilde{\Theta}=\Theta$ and $\widetilde{\kappa}=\big(\overline{A}_\Theta,\overline{B}_\Theta\big)$, to compute the tail probability as well as the expected value of the suprema of random mappings.
Note also that Assumptions \ref{ass:fixed_theta_local}-\ref{ass:parameter_local} guaranty that $\mathrm{E}\left\{\sup_{\theta\in\Theta}\Psi\left(\chi_{\theta}\right)\right\}^{q}$ is finite for any $q\geq 1$.

%\medskip

\subsection{Upper functions of the first and second type} We will now use Proposition \ref{prop_uniform_local1} in order to derive the upper functions for $\Psi\left(\chi_{\theta}\right)$ on $\Theta$.
%We enforce Assumption \ref{ass:metric-case-local} ({\it 3}) supposing that $\Theta$ is a \textsf{compact} set.
%Throughout of this section we will suppose that
%$$
% \sup_{\theta\in\Theta} A(\theta)=:\overline{A}<\infty,\;\; \sup_{\theta\in\Theta} B(\theta)=:\overline{B}<\infty,
%$$
Denote $\underline{A}=\inf_{\theta\in\Theta} A(\theta)$ and $\underline{B}=\inf_{\theta\in\Theta} B(\theta)$.

We present two kinds  of upper functions for $\Psi\left(\chi_{\theta}\right)$ on $\Theta$ which we will refer to upper functions of
the first and second type.
The first construction is completely determined by
the functions $A$, $B$ and by the  semi-metrics $\ra$ and $\rb$. It requires however the additional condition
$\underline{A}>0, \underline{B}>0$. We will use corresponding results for the particular  problems studied in Section \ref{sec:empirical-processes}.

The second construction is related to some special structure imposed on the set $\Theta$. Namely we will suppose that $\Theta=\cup_\alpha\Theta_{\alpha\in\mA}$, where $\big\{\Theta_\alpha,\; \alpha\in\mA\big\}$ is a given collection of sets. Here we will be interested in a finding of upper function for
$
\sup_{\theta\in \Theta_\alpha}\Psi\left(\chi_{\theta}\right)
$
on $\mA$, which can be also viewed as an upper function  for  $\Psi\left(\chi_{\theta}\right)$ on $\Theta$. The corresponding results are  used in order to obtain rather precise inequalities for the modulus of continuity of random functions, Section \ref{sec:subsection-gauss-local-continuity}.
Moreover we apply this bound for deriving of an  upper function for the $\bL_p$-norms of Wiener integrals,  Section \ref{sec:subsection-norm-gauss}.  We deduce the corresponding inequality directly from Proposition \ref{prop_uniform_local3} below  without passing to the concentration inequalities. %It, in its  turn, shows that the use of standard  chaining arguments may leads to rather sharp and deep results.

\vskip0.1cm

We finish this short introduction with the following remark. In order to establish  the inequalities (\ref{eq3:introduction})--(\ref{eq4:introduction}) for the upper functions presented below we will need to prove that corresponding supremum is a random variable. The result below is sufficient for all problems considered in the paper and until proofs
we will not discuss  the measurability issue.

\begin{lemma}
\label{lem:measurability}
Let $\boldsymbol{\mT}$ be the set equipped with the metric $\md$, $\left(\boldsymbol{\Omega},\boldsymbol{\mB},\boldsymbol{\mathrm{P}}\right)$ be a complete probability space and  $\zeta:\boldsymbol{\Omega}\times\boldsymbol{\mT}\to \bR$ be $\boldsymbol{\mathrm{P}}$-a.s. continuous. Let $\mZ$ be a set,  $g:\mZ\to\bR$  be a given function and $\left\{\mT_\mz\subseteq\mT,\;\mz\in\mZ\right\}$ be an arbitrary sequence of sets.
If $\boldsymbol{\mT}$ is  totally bounded then  $\sup_{\mz\in\mZ}\Big[\sup_{\mathfrak{t}\in\boldsymbol{\mT}_\mz}\zeta(\mt,\cdot)-g(\mz)\big]$ is $\boldsymbol{\mB}$-measurable.
\end{lemma}

The proof of the lemma is given in Appendix. We would like to emphasize  that there is no any assumption imposed on the function $g$, index set $\mZ$ and on the collection  $\left\{\mT_\mz\subseteq\mT,\;\mz\in\mZ\right\}$.

Putting $\mZ=\mT$ and $\mT_\mt=\{\mt\}$ we come to the following consequence of Lemma \ref{lem:measurability}.

\begin{corollary}
\label{cor:after-lem-measurability}
 Under assumptions of Lemma \ref{lem:measurability}  $\sup_{\mathfrak{t}\in\boldsymbol{\mT}}\big[\zeta(\mt,\cdot)-g(\mt)\big]$ is $\boldsymbol{\mB}$-measurable.
 %for any $g:\mT\to\bR$.
\end{corollary}

\paragraph{Upper functions of the first type} As it was said above throughout  this section we will suppose that $\underline{A}>0, \underline{B}>0$.
Put for any $t>0$
\begin{eqnarray*}
&&\Theta_{A}(t)=\Big\{\theta\in\Theta:\; A(\theta)\leq t\Big\},\quad \Theta_{B}(t)=\Big\{\theta\in\Theta:\; B(\theta)\leq t\Big\}.
 %\\*[2mm]
%&&\bZ^{*}_{\cA}(a)=\left\{\zeta\in\bZ:\; \cA\big([\zeta]\big)\leq a\right\},\qquad\qquad\quad\;\; \bZ^{*}_{\cB}(b)=\left\{\zeta\in\bZ:\; \cB\big([\zeta]\big)\leq b\right\}.
\end{eqnarray*}
 For any $\vec{s}\in\cS_{\ra,\rb}$ introduce the  function
\begin{eqnarray}
\label{eq:def-of-function-e}
&&\mathcal{E}_{\vec{s}}(u,v)=
e^{(\ra)}_{s_1}\Big(\underline{A}u,\Theta_{A}\big(\underline{A}u\big)\Big)+ e^{(\rb)}_{s_2}\Big(\underline{B}v,\Theta_{B}\big( \underline{B}v\big)\Big),\quad u,v\geq 1.
\end{eqnarray}
Denote also
%\begin{eqnarray*}
$
\ell(u)=\ln{\left\{1+\ln{(u)}\right\}}+2\ln{\left\{1+\ln{\left\{1+\ln{(u)}\right\}}\right\}}
$
and set for any $\theta$ and   $\e>0,r\geq 0$
\begin{eqnarray}
\label{eq:def-price-to-pay-proba}
&&\quad P_{\e}(\theta)=2\big[1+\e^{-1}\big]^{2}\cE\big(\cA_\e(\theta),
\cB_\e(\theta)\big)+(1+\e)^{2}\big[
\ell\big(\cA_\e(\theta)\big)+
\ell\big(\cB_\e(\theta)\big)\big];
\\*[2mm]
\label{eq:def-price-to-pay-moments}
 &&\quad M_{\e,r}(\theta)=(1+\e)^{2}\left\{2\big[1+\e^{-1}\big]^{2}\cE\big(\cA_\e(\theta),
\cB_\e(\theta)\big)+
(\e+r)\ln{\big[\cA_\e(\theta)
\cB_\e(\theta)\big]}\right\},
\end{eqnarray}
where $\cA_\e(\theta)=(1+\e)\big[A(\theta)\big/\underline{A}\big]$ and $\cB_\e(\theta)=(1+\e)\big[B(\theta)\big/\underline{B}\big].$
\;Define for any $z\geq 0$
\begin{eqnarray}
\label{eq:def-of-upper-function-proba}
\mathrm{V}^{(z,\e)}(\theta)&=&(1+\e)^{2}\left(A(\theta)\sqrt{P_{\e}(\theta)+(1+\e)^2z}+
B(\theta)\Big[P_{\e}(\theta)+(1+\e)^2z\Big]\right);
\\*[2mm]
\label{eq:def-of-upper-function-moments}
\mathrm{U}^{(z,\e,r)}(\theta)&=&(1+\e)^2\left(A(\theta)\sqrt{M_{\e,r}(\theta)+(1+\e)^2z}+
B(\theta)\Big[M_{\e,r}(\theta)+(1+\e)^2z\Big]\right).
 \end{eqnarray}
In the proposition below we prove that the functions defined in (\ref{eq:def-of-upper-function-proba}) and (\ref{eq:def-of-upper-function-moments}) are upper functions for $\Psi\left(\chi_{\theta}\right)$ on $\Theta$.
%We remark that all constants involved in their descriptions depend only on
We remark that the behavior of  $\Psi\left(\chi_{\theta}\right)$ on $\Theta$ is completely determined by the functions $A$ and $B$ and by the entropies of their level sets measured in  semi-metrics $\ra$ and $\rb$. The number $\e$ and the couple of  functions $\vec{s}$ can be viewed as tuning  parameters allowing either to weaken assumptions or to obtain sharper bounds but they are not  related to the random functional $\Psi\left(\chi_{\theta}\right)$ itself.
%In particular, the choice $\e=\sqrt{2}-1$ and $\vec{s}=(s^*,s^*)$, where $s^*$ is defined in (\ref{eq:function-s-star}), is convenient for the majority of considered problems.

\begin{proposition}
\label{prop_uniform_local2}
 Let   Assumptions \ref{ass:fixed_theta_local}-\ref{ass:parameter_local} be fulfilled.
Then $\forall \vec{s}\in\cS_{\ra,\rb}$, $\forall \e\in \big(0,\sqrt{2}-1\big]$ and
$\forall z\geq 1$
\begin{gather*}
\mathrm{P}\left\{\sup_{\theta\in\Theta}\Big[\Psi\left(\chi_{\theta}\right)- \mathrm{V}^{(z,\ve)}(\theta)\Big]\geq 0\right\}
\leq 2\mathrm{c}\left[1+\Big[\ln{\left\{1+\ln{(1+\e)}\right\}}\Big]^{-2}\right]^{2}\exp{\left\{-z\right\}};
\\*[2mm]
\mathrm{E}\left\{\sup_{\theta\in\Theta}\Big[\Psi\left(\chi_{\theta}\right)- \mathrm{U}^{(z,\ve,q)}(\theta)\Big]\right\}^q_+
%\\*[2mm]&&
\leq \mathrm{c}2^{(5q/2)+2}\Gamma(q+1)\;\e^{-q-4}\;\big[\underline{A}\vee\underline{B}\big]^q\exp{\left\{-z\right\}}.
\end{gather*}
%where $C_{\e,q}=$
\end{proposition}

It is obvious that the assertions of the proposition remain valid if one replaces the function $\cE$ by any its upper bound. It is important since
the exact computation of this function is too complicated in general.
%\begin{remark}
%\label{rem:after-proposition4-local}
 %As we already mentioned  both upper functions   are  determined by the function $\cE$.
 We note that the role of the latter function in our construction  is similar to whose which Dudley integral plays in the computations of the expectation of the  suprema of gaussian or sub-gaussian   processes \cite{Lif}, \cite{Tal-book}.

%\end{remark}

\paragraph{Price to pay for "uniformity"}

 We remark that in view of (\ref{eq100:discussion-price-to-pay}) and (\ref{eq2:discussion-price-to-pay}), the function $U^{(z)}(\theta):=A(\theta)\sqrt{z}+B(\theta)z$ can be viewed as "pointwise upper function" for $\Psi(\chi_\theta)$, i.e. for fixed $\theta$. Comparing the inequalities (\ref{eq100:discussion-price-to-pay}) and (\ref{eq2:discussion-price-to-pay})  with  whose given  in Proposition \ref{prop_uniform_local2} we conclude that they differ from each other by numerical constants only. In this context, the functions $P(\cdot)$ and $M_q(\cdot)$ given by (\ref{eq:def-price-to-pay-proba}) and (\ref{eq:def-price-to-pay-moments})
can be viewed as price to pay for "uniformity". That means that in order to pass from "pointwise" result to the "uniform" one we need, roughly speaking, to multiply $A(\cdot)$ by $\sqrt{P(\cdot)}$ or $\sqrt{M_q(\cdot)}$ and $B(\cdot)$ by $P(\cdot)$ or $M_q(\cdot)$. The question, arising naturally: is  such payment necessary or minimal? We do not think that the answer can be done under "abstract considerations", i.e. under Assumptions \ref{ass:fixed_theta_local}-\ref{ass:parameter_local}.  However, for particular problems it is seemed possible.
Unexpectedly the answer on the formulated above question can come from the solution of the problem studied in mathematical statistics. In this context it worths to mention the relation between well-known  phenomenon in adaptive estimation, called {\it price to pay for adaptation} \cite{lepski}, \cite{LepSpok} and \cite{Spok}, and
what we call here {\it price to pay for uniformity}. We have no place here to describe this relation in detail and mention only  several facts.

First let us remark that Proposition \ref{prop_uniform_local2} contains the results which can be directly used for the construction of adaptive procedures. Indeed, almost all constructions of adaptive estimators (model selection \cite{Barron-Birge}, risk hull minimization \cite{golubev}, Lepski method \cite{lepski}, or recently developed universal estimation routine \cite{GL1})  involve the upper functions for stochastic objects of different kinds.
Next, it is known that there are two types of {\it price to pay for adaptation}: $(\ln)$-price,  \cite{lepski} and $(\ln\ln)$-price, \cite{Spok}. The $(\ln)$-price appears in the problems where the risk of  estimation procedures is described by a power loss-functions and it corresponds to the function $M_q(\cdot)$, where the parameter $q$ is a power. The  $(\ln\ln)$-price appears in the case of bounded losses that corresponds to the function $P(\cdot)$. Since the theory of adaptive estimation is equipped with very developed criteria of optimality, \cite{lepski}, \cite{Tsyb}, \cite{Kluch}, we might assert that the  payment for uniformity is optimal if the use of corresponding upper function leads to  optimal adaptive estimators.

We finish the discussion concerning the statements of Proposition \ref{prop_uniform_local2}  with the following remark. Comparing the result given in (\ref{eq2:discussion-price-to-pay}) with the second assertion of Proposition \ref{prop_uniform_local2} we can state that the inequality obtained there is very precise since, remind, $\underline{A}=\inf_{\theta\in\Theta}A(\theta)$ and $\underline{B}=\inf_{\theta\in\Theta}B(\theta)$.

\paragraph{Upper functions of the second type} Suppose that we are given by the collection $\big\{\Theta_\alpha,\; \alpha\in\mA\big\}$, satisfying $\Theta=\cup_{\alpha\in\mA}\Theta_\alpha$, and by two mappings $\tau_1:\mA\to \big(0,\overline{\tau}_1\big], \;\tau_2:\mA\to\big(0,\overline{\tau}_2\big]$, where $\overline{\tau}_1,\overline{\tau}_2<\infty$.\;
For any $u>0$ put
\begin{eqnarray*}
\Theta^\prime_1(u)=\bigcup_{\alpha:\;\tau_1(\alpha)\leq u}\Theta_\alpha,\quad g^*_A(u)=\sup_{\theta\in\Theta^\prime_1(u)}A(\theta);
\\
\Theta^\prime_2(u)=\bigcup_{\alpha:\;\tau_2(\alpha)\leq u}\Theta_\alpha, \quad g^*_B(u)=\sup_{\theta\in\Theta^\prime_2(u)}B(\theta),
\end{eqnarray*}
and let  $g_A$ and $g_B$ be arbitrary chosen  increasing functions, satisfying $g_A\geq g^*_A$ and $g_B\geq g^*_B$ (we note that obviously $ g^*_A$ and $ g^*_B$ are increasing).

\smallskip

Since $\Theta^\prime_1(\cdot),\Theta^\prime_2(\cdot)\subseteq\Theta$, in view of Assumption \ref{ass:parameter_local} for any $u,v>0$ one can find the functions $s_1(u,\cdot)$ and $s_2(v,\cdot)$ for which the latter assumption is fulfilled on $\Theta^\prime_1(u)$ and $\Theta^\prime_2(v)$ respectively. Let us suppose additionally that
\begin{eqnarray}
\label{eq:def-of-quantities-lambda}
\lambda_1:=\sup_{t\in[1,\sqrt{2}]}\;\sup_{x>\underline{\tau}_1}\;\sup_{\delta>0}\frac{s_1(xt,\delta)}{s_1(x,\delta)}<\infty,\quad \lambda_2:=\sup_{t\in[1,\sqrt{2}]}\;\sup_{x>\underline{\tau}_2}\;\sup_{\delta>0}\frac{s_2(xt,\delta)}{s_2(x,\delta)}<\infty,
\end{eqnarray}
where $\underline{\tau}_1=\inf_\alpha \tau_1(\alpha)$ and $\underline{\tau}_2=\inf_\alpha \tau_2(\alpha)$.

We remark that if the functions $s_1(u,\cdot)$ and $s_2(v,\cdot)$ are chosen independently of $u,v$ then  $\lambda_1=\lambda_2=1$.
It is also obvious that  $\lambda_1,\lambda_2\geq 1$.

The condition (\ref{eq:def-of-quantities-lambda}) allows us to define the function:
\begin{eqnarray}
\label{eq:def-of-function-tilde-e-alpha}
&&\mathcal{E}^{\prime}(u,v)=e^{(\ra)}_{s_1(u,\cdot)}\Big(\lambda_1^{-1}g_A(u),\Theta^\prime_1(u)\Big)+e^{(\rb)}_{s_2(v,\cdot)}
\Big(\lambda_2^{-1}g_B(v),\Theta^\prime_2(v)\Big),\;\; u,v>0.
\end{eqnarray}
We note that the function $\mathcal{E}^{\prime}$ is constructed similarly to the function $\cE$ used in the previous section, but now the  functions $s_1$ and $s_2$ can be chosen in accordance with  considered level sets.

At last, for any $\alpha\in\mA$ and any  $\e>0$ set
\begin{gather*}
\widehat{\cE}^{(\e)}(\alpha)=\cE^\prime\Big((1+\e)\tau_1(\alpha),(1+\e)\tau_2(\alpha)\Big).
%\\ \widehat{\ell}_{V}(\alpha)=\ell\Big(\overline{\tau}_1\big/\tau_1(\alpha)\Big)+\ell\Big(\overline{\tau}_2\big/\tau_2(\alpha)\Big),
%\quad\widehat{\ell}_{U}(\alpha)=\ln{\Big(\overline{\tau}_1\big/\tau_1(\alpha)\Big)}+\ln{\Big(\overline{\tau}_2\big/\tau_2(\alpha)\Big)},
\end{gather*}
Put $\delta_j=(1+\e)^{-j}, j\geq 0,$ and let  $R_r:\bR_+\times\bR_+\to \bR_+ ,r\geq 0,$ be an arbitrary family of  {\it increasing} in both arguments functions, satisfying for any $\e\in\big(0,\sqrt{2}-1\big]$
\begin{gather}
\label{eq:condition-on the functions-R}
%\sum_{j=0}^{\infty}\exp{\left\{-R_V\big(\overline{\tau}_\ri \delta_{j}\big)\right\}}=:\cR^{(\ri)}(\e)<\infty;\;\;\ri=1,2;
%\\
\sum_{j=0}^{J}\sum_{k=0}^{K}\left[g_A\big(\overline{\tau}_1 \delta_{j}\big)\vee g_B\big(\overline{\tau}_2 \delta_{k}\big)\right]^{r}\exp{\left\{-R_r\big(\overline{\tau}_1 \delta_{j},\overline{\tau}_2 \delta_{k}\big)\right\}}=:\cR^{(\e,r)}<\infty.
\end{gather}
Here integers $J,K$ are defined as follows.
$$
J=\left\lfloor\ln_{1+\ve}\big(\overline{\tau}_1\big/\underline{\tau}_1\big)\right\rfloor+1,\quad K=\left\lfloor\ln_{1+\ve}\big(\overline{\tau}_2\big/\underline{\tau}_2\big)\right\rfloor+1.
%\;\;\underline{\tau}_\ri=\inf_{\alpha\in\mA}\tau_\ri(\alpha),\;\ri=1,2.
$$
If $\underline{\tau}_\ri=0,\;\ri=1,2, $ the corresponding quantity is put equal to infinity.

Set $\widehat{R}^{(\e)}_r(\alpha)=R_r\Big((1+\e)\tau_1(\alpha),(1+\e)\tau_2(\alpha)\Big)$
and define
\begin{eqnarray*}
%\widehat{V}^{(z,\e)}(\alpha)&=&g_A\Big([1+\e]^2\tau_1(\alpha)\Big)\sqrt{2\big[1+\e^{-1}\big]^{2}\;
%\widehat{\mathcal{E}}^{(\e)}(\alpha)+(1+\e)^{2}\Big(\widehat{R}^{(\e)}_{0}(\alpha)+z\Big)}\\*[1mm]
%\;&+&
%g_B\Big([1+\e]^2\tau_2(\alpha)\Big)\Big[2\big[1+\e^{-1}\big]^{2}\;
%\widehat{\mathcal{E}}^{(\e)}(\alpha)+(1+\e)^{2}\Big(\widehat{R}^{(\e)}_0(\alpha)+z\Big)\Big].
%\\*[2mm]
\widehat{\mathrm{U}}^{(z,\e,r)}(\alpha)&=&(1+\e)g_A\Big([1+\e]^2\tau_1(\alpha)\Big)\sqrt{2\big[1+\e^{-1}\big]^{2}\;
\widehat{\mathcal{E}}^{(\e)}(\alpha)+\widehat{R}^{(\e)}_r(\alpha)+z}\\*[1mm]
\;&+&
(1+\e)^{2}g_B\Big([1+\e]^2\tau_2(\alpha)\Big)\Big(2\big[1+\e^{-1}\big]^{2}\;\widehat{\mathcal{E}}^{(\e)}(\alpha)+
\widehat{R}^{(\e)}_r(\alpha)+z\Big).
\end{eqnarray*}
Below we assert that  $\widehat{U}^{(z,\e,r)}, r=0, r=q,$ are  upper functions for $\Big[\sup_{\theta\in\Theta_\alpha}\Psi\left(\chi_{\theta}\right)\Big]$ on $\mA$. However, before to present
 exact statements, let us briefly discuss  some possible choices of the functions $R_r$. We would like to emphasize that the opportunity to select these functions allows to obtain quite different and  precise results.
First possible choice is given by
\begin{equation}
\label{eq:choice-of-function-R}
R_{0}(u,v)=\ell\Big(\overline{\tau}_1u^{-1}\Big)+\ell\Big(\overline{\tau}_2v^{-1}\Big),\quad R_{r}(u,v)=\e\Big[\ln{\Big(\overline{\tau}_1u^{-1}\Big)}+\ln{\Big(\overline{\tau}_2v^{-1}\Big)}\Big],\; r>0.
\end{equation}
These functions are used in the problems in which $\widehat{\cE}^{(\e)}(\cdot)$ is bounded by some absolute constant independent of all quantities involved in the description of the problem, assumptions etc.

This choice leads to the following values of the constants in (\ref{eq:condition-on the functions-R}):
\begin{equation}
\label{eq:constants-choice-of-function-R}
\cR^{(\e,0)}\leq\left[2+\Big[\ln{\left\{1+\ln{(1+\e)}\right\}}\Big]^{-2}\right]^{2},\quad \cR^{(\e,r)}<4\Big[g_A\big(\overline{\tau}_1\big)\vee g_A\big(\overline{\tau}_2\big)\Big]^{r}\e^{-4}.
\end{equation}
Another important choice is given by $R_r=\cE^{\prime}$ independently of $r$, see, for instance, Theorem \ref{th:gauss-norm}. In view of (\ref{eq:condition-on the functions-R}), this choice corresponds to the case when the function  $\cE^{\prime}$ increases to infinity.
\begin{proposition}
\label{prop_uniform_local3}
 Let   Assumptions \ref{ass:fixed_theta_local}-\ref{ass:parameter_local} be fulfilled.
\;Then for any $s_1,s_2$ satisfying (\ref{eq:def-of-quantities-lambda}) and any $R_r,\;r\geq 0,$ satisfying (\ref{eq:condition-on the functions-R}), for any $ \e\in\big(0,\sqrt{2}\big]$ and any
$z\geq 1, q\geq 1$
\begin{eqnarray*}
&&
\mathrm{P}\left\{\sup_{\alpha\in\mA}\bigg[\sup_{\theta\in\Theta_\alpha}\Psi\left(\chi_{\theta}\right)- \widehat{\mathrm{U}}^{(z,\e,0)}(\alpha)\bigg]\geq 0\right\}
\leq 2\mathrm{c}\cR^{(\e,0)}\exp{\left\{-z\right\}};
\\*[2mm]
&&\mathrm{E}\left\{\sup_{\alpha\in\mA}\bigg[\sup_{\theta\in\Theta_\alpha}\Psi\left(\chi_{\theta}\right)- \widehat{\mathrm{U}}^{(z,\e,q)}(\alpha)\bigg]\right\}^{q}_+\leq \mathrm{c}2^{(5q/2)+1}\Gamma(q+1)\cR^{(\e,q)}\;\e^{-q}\exp{\left\{-z\right\}}.
\end{eqnarray*}
\end{proposition}

\begin{remark}
We note that the results of the proposition is very general. Indeed, there are no assumptions imposed on the collection $\Theta_\alpha,\;\alpha\in\mA,$ and  the functions $\tau_1,\tau_2$ can be chosen  arbitrary.  Moreover, the condition (\ref{eq:condition-on the functions-R}) is very mild, so the choice of functions $R_r$ is quite flexible.
\end{remark}

\subsection{Upper functions for the modulus of continuity of random mappings} In this section we apply Proposition \ref{prop_uniform_local3} in order to derive upper functions for the local and global modulus of continuity of real-valued random mappings. It is worth mentioning that in this circle of problems the upper functions are actively  exploited, see e.g. \cite{Ostrovskii} and the references therein. We will suppose that Assumption \ref{ass:fixed_theta_local} ({\it 2}), Assumption \ref{ass:metric-case-local} and
Assumption \ref{ass:parameter_local} are verified, $\chi_\mt$ is real-valued random mapping defined on the metric space $\mT$,  $\rd$ is a semi-metric on $\mT$ and $\Psi(\cdot)=|\cdot|$.

\paragraph{Upper function for local modulus of continuity} Let $\theta_0$ be a fixed element of $\Theta$ and set for any $\D\in \big(0, D_\rd(\Theta)\big]$,  where
$D_\rd(\Theta)$ is the diameter of $\Theta$ measured in the  semi-metric $\rd$,
$$
\mm_\D(\theta_0)=\sup_{\theta\in\Theta_\D}\big|\chi_\theta-\chi_{\theta_0}\big|,\quad \Theta_\D=\Big\{\te\in\Theta:\;\;\rd\big(\theta,\te_0\big)\leq \D\Big\}.
$$
Thus, $\mm_\D(\theta_0),\;\D\in \big(0, D_\rd(\Theta)\big],$ is the local modulus of continuity of $\chi_{\theta}$ in $\theta_0$ measured in $\rd$.

\smallskip

If we put $\widetilde{\chi}_\theta=\chi_\theta-\chi_{\theta_0},\;\te\in\Theta,$ we assert first that  Assumption \ref{ass:fixed_theta_local} ({\it 2}) can be viewed as  Assumption \ref{ass:fixed_theta_local} ({\it 1}) for  $\widetilde{\chi}_\theta$ on $\Theta$  with  $A(\cdot)=\ra(\cdot,\te_0)$ and $B(\cdot)=\rb(\cdot,\te_0)$.
Next, noting that $\widetilde{\chi}_{\theta_1}-\widetilde{\chi}_{\theta_2}=\chi_{\theta_1}-\chi_{\theta_2}$ for any $\theta_1,\theta_2\in\Theta$ we conclude that  Assumption \ref{ass:fixed_theta_local} ({\it 2}) is verified for $\widetilde{\chi}_{\theta}$ on $\Theta$ with $\ra$ and $\rb$.
%In view of the latter remark, since we suppose that   Assumption \ref{ass:parameter_local} is fulfilled for $\chi_\theta$, this assumption is verified for $\widetilde{\chi}_{\theta}$ as well.

Thus, we can apply Proposition \ref{prop_uniform_local3} with $\alpha=\Delta,\;\Theta_\alpha=\Theta_\D,\;\mA=\big(0,D_\rd(\Theta)\big]$
and we choose $\tau_1(\Delta)=\tau_2(\Delta)=\Delta$. This choice implies obviously for any $u\leq D_\rd(\Theta) $
$$
\Theta^\prime_1(u)=\Theta^\prime_2(u)=\Theta_u,\quad g_A(u)=\sup_{\te:\; \rd\big(\theta,\te_0\big)\leq u}\ra\big(\theta,\te_0\big),\;\; g_B(u)=\sup_{\te:\; \rd\big(\theta,\te_0\big)\leq u}\rb\big(\theta,\te_0\big).
$$
Fix $\vec{s}\in\cS_{\ra,\rb}$ and  put for any $\D\in \big(0, D_\rd(\Theta)\big]$ and any $\e\in\big(0,\sqrt{2}-1\big]$
\begin{equation*}
\widehat{\cE}^{(\e)}(\D,\te_0)=e^{(\ra)}_{s_1}\Big(g_A\big([1+\e]\D\big),\Theta_{[1+\e]\D}\Big)+e^{(\rb)}_{s_2}
\Big(g_B\big([1+\e]\D\big),\Theta_{[1+\e]\D}\Big).
\end{equation*}
Here  $e^{(\ra)}_{s_1}$ and $e^{(\rb)}_{s_2}$ are defined by (\ref{eq3:def_local}).  We also set $\lambda_1=\lambda_2=1$ since the functions $s_1,s_2$ are chosen independently of the collection $\left\{\Theta_\D,\;\D\in \big(0, D_\rd(\Theta)\big]\right\}$.

 Choose also
$
R_{0}(u,v)=\ell\Big(\overline{\tau}_1u^{-1}\Big)+\ell\Big(\overline{\tau}_2v^{-1}\Big)
$
and define
\begin{eqnarray*}
\widehat{V}_{\vec{s}}^{(z,\e)}(\D,\theta_0)&=&(1+\e)g_\cA\Big([1+\e]^2\D\Big)\sqrt{2\big[1+\e^{-1}\big]^{2}\widehat{\cE}^{(\e)}(\D,\theta_0)
+(1+\e)^{2}\Big[2\ell\Big(D_\rd(\Theta)\big/\D\Big)+z\Big]}
\\
&+&
(1+\e)^{2}g_\cB\Big([1+\e]^2\D\Big)\Big\{2\big[1+\e^{-1}\big]^{2}\widehat{\cE}^{(\e)}(\D,\theta_0)+(1+\e)^2\Big[2\ell\Big(D_\rd(\Theta)\big/\D\Big)+z\Big]\Big\}.
\end{eqnarray*}
Then, applying Proposition \ref{prop_uniform_local3} and taking into account (\ref{eq:constants-choice-of-function-R}) we come to the following result.
\begin{proposition}
\label{prop_uniform_local4}
Let   Assumptions \ref{ass:fixed_theta_local}-\ref{ass:parameter_local} be fulfilled.
Then $\forall s\in\cS_{\ra,\rb}$, $\forall \e>0$ and
$\forall z\geq 1$
\begin{eqnarray*}
&&
\mathrm{P}\left\{\sup_{\D\in \big(0, D_\rd(\Theta)\big]}\Big[\mm_\D- \widehat{V}_s^{(z,\e)}(\D,\theta_0)\Big]\geq 0\right\}
\leq 4\mathrm{c}\left[2+\Big[\ln{\left\{1+\ln{(1+\e)}\right\}}\Big]^{-2}\right]^{2}\exp{\left\{-z\right\}}.
\end{eqnarray*}
\end{proposition}
In Section \ref{sec:gauss} we apply Proposition \ref{prop_uniform_local4} to  gaussian random functions defined on a metric space satisfying so-called {\it doubling condition}.
\begin{remark}
\label{rem1:after-proposition4-local} If $\rb\equiv 0$,  $\rd=\ra$ and
$
\sup_{\D\in \big(0, D_\rd(\Theta)\big]}\widehat{\cE}^{(\e)}(\D,\theta_0)=:\widehat{\cE}^{(\e)}(\theta_0)<\infty,
$
the upper function $\widehat{V}_s^{(z,\e)}$ has very simple form
\begin{equation}
\label{eq:rem-after-proposition4}
\widehat{V}_s^{(z,\e)}(\D,\theta_0)=(1+\e)^3\D\sqrt{4\big[1+\e^{-1}\big]^{2}\widehat{\cE}^{(\e)}(\theta_0)
+(1+\e)^{2}\Big[\ell\Big(D_\rd(\Theta)\big/\D\Big)+z\Big]}.
\end{equation}
Hence, the result of Proposition \ref{prop_uniform_local4} can be viewed as the non-asymptotical version of the law of iterated logarithm for sub-gaussian processes defined on some totaly bounded subset of metric space.
In this context it is worth mentioning the paper \cite{Ostrovskii} where the upper functions for local and global modulus of continuity  were found for the stochastic processes satisfying Cramer's condition.
\end{remark}
\begin{remark}
\label{rem2:after-proposition4-local}
We also note that  we replaced in (\ref{eq:rem-after-proposition4}) the factor $2\ell\big(D_\rd(\Theta)\big/\D\big)$ appeared in the upper function used in Proposition
\ref{prop_uniform_local4} by $\ell\big(D_\rd(\Theta)\big/\D\big)$. It is explained by the fact that $\overline{\tau}_2=0$ in (\ref{eq:choice-of-function-R}) since $B,\rb\equiv 0$. By the same reason, the probability bound in this case is given by
$
4\mathrm{c}\left[2+\Big[\ln{\left\{1+\ln{(1+\e)}\right\}}\Big]^{-2}\right]\exp{\left\{-z\right\}}.
$
\end{remark}

\paragraph{Upper function for global modulus of continuity}

Set $\Theta^{(2)}=\Theta\times\Theta$ and let for any $\vth=(\theta_1,\te_2)\in\Theta^{(2)}$ and any $\D\in \big(0, D_\rd(\Theta)\big]$,
$$
\z(\vth)=\chi_{\theta_1}-\chi_{\theta_2},\quad \mm_\D=\sup_{\vth\in\Theta^{(2)}_\D}\big|\z_\vth\big|,\quad \Theta^{(2)}_\D=\Big\{\vth\in\Theta^{(2)}:\;\;\rd\big(\theta_1,\te_2\big)\leq \D\Big\}.
$$
Thus, $\mm_\D,\;\D\in \big(0, D_\rd(\Theta)\big],$ is the global modulus of continuity of $\chi_{\theta}$ on $\Theta$ measured in  $\rd$.

Put $\mathbf{A}(\vth)=\ra(\theta_1,\te_2),\; \mathbf{B}(\vth)=\rb(\theta_1,\te_2)$,  $\vth=(\theta_1,\te_2)\in\Theta^{(2)}$, and equip $\Theta^{(2)}$
with the following  semi-metrics: $\vth=(\te_1,\te_2),\vs=(\vs_1,\vs_2)\in\Theta^{(2)}$
$$
\ra^{(2)}(\vth,\vs)=2\left[\ra(\te_1,\vs_1)\vee\ra(\te_2,\vs_2)\right],\;\; \rb^{(2)}(\vth,\vs)=2\left[\rb(\te_1,\vs_1)\vee\rb(\te_2,\vs_2)\right].
$$
Some remarks are in order. We note first that Assumption \ref{ass:fixed_theta_local} ({\it 2}) can be viewed as  Assumption \ref{ass:fixed_theta_local} ({\it 1}) for  $\z(\vth)$ on $\Theta^{(2)}$  with  $A=\mathbf{A}$ and $B=\mathbf{B}$.

Next we obtain  in view of  Assumption \ref{ass:fixed_theta_local} ({\it 2}) $\forall \vth,\vs\in\Theta^{(2)}$ and  $\forall z>0$
\begin{eqnarray*}
&&\mathrm{P}\Big\{|\z(\vth)-\z(\vs)|\geq z\Big\}\leq \mathrm{P}\Big\{\big|\chi_{\theta_1}-\chi_{\vs_1}\big|\geq z/2\Big\}+\mathrm{P}\Big\{\big|\chi_{\theta_2}-\chi_{\vs_2}\big|\geq z/2\Big\}
\\
&&\leq
\mathrm{c}\exp{\left\{-\frac{z^{2}}{4\big[\mathrm{a}(\theta_1,\vs_1)\big]^{2}+ 2\mathrm{b}(\theta_1,\vs_1)z}\right\}}+
\mathrm{c}\exp{\left\{-\frac{z^{2}}{4\big[\mathrm{a}(\theta_2,\vs_2)\big]^{2}+ 2\mathrm{b}(\theta_2,\vs_2)z}\right\}}
\\
&&\leq
\mathrm{c}^{(2)}\exp{\left\{-\frac{z^{2}}{\big[\mathrm{a}^{(2)}(\vth,\vs)\big]^{2}+ \mathrm{b}^{(2)}(\vth,\vs)z}\right\}}
\end{eqnarray*}
We conclude that Assumption \ref{ass:fixed_theta_local} ({\it 2}) holds for  $\z(\vth)$ on $\Theta^{(2)}$  with $\ra=\ra^{(2)}$, $\rb=\rb^{(2)}$
and $\mathrm{c}^{(2)}=2\mathrm{c}$.

\smallskip

Since obviously
\begin{eqnarray}
\label{eq:def-entropy-in-modulus-continuity}
\mE_{\ra^{(2)},\Theta^{(2)}}(\varsigma)\leq 2\mE_{\ra,\Theta}(\varsigma/2),\quad \mE_{\rb^{(2)},\Theta^{(2)}}(\varsigma)\leq 2\mE_{\rb,\Theta}(\varsigma/2),\;\;\varsigma>0,
\end{eqnarray}
we assert that Assumptions \ref{ass:metric-case-local} and \ref{ass:parameter_local} are fulfilled on $\Theta^{(2)}$  with $\ra=\ra^{(2)}$ and  $\rb=\rb^{(2)}$.

\smallskip

Put $\overline{\Theta}^{(2)}=\cup_{\D>0}\Theta^{(2)}_\D$. Since $\overline{\Theta}^{(2)}\subset\Theta^{(2)}$ we can apply Proposition \ref{prop_uniform_local3} with $\alpha=\Delta,\;\Theta_\alpha=\Theta^{(2)}_\D,\;\mA=\big(0,D_\rd(\Theta)\big]$
%We remark that in our case $\overline{\cA}=D_\ra(\Theta),\overline{\cB}=D_\rb(\Theta)$
and we choose $\tau_1(\Delta)=\tau_2(\Delta)=\Delta$.

\smallskip

The latter choice implies obviously for any $u\leq D_\rd(\Theta) $
$$
\Theta^\prime_1(u)=\Theta^\prime_2(u)=\Theta^{(2)}_u,\quad g_\mathbf{A}(u)=\sup_{\vth\in\Theta^{(2)}_u}\mathbf{A}(\vth),\;\; g_\mathbf{B}(u)=\sup_{\vth\in\Theta^{(2)}_u}\mathbf{B}(\vth).
$$
Fix $\vec{s}\in\cS_{\ra,\rb}$ and  put for any $\D\in \big(0, D_\rd(\Theta)\big]$  and any $\e\in\big(0,\sqrt{2}-1\big]$
\begin{equation*}
\widehat{\cE}^{(\e)}(\D)=e^{(\ra^{(2)})}_{s_1}\Big(g_A\big([1+\e]\D\big),\Theta^{(2)}_{[1+\e]\D}\Big)+e^{(\rb^{(2)})}_{s_2}
\Big(g_B\big([1+\e]\D\big),\Theta^{(2)}_{[1+\e]\D}\Big).
\end{equation*}
Here  $e^{(\ra^{(2)})}_{s_1}$ and $e^{(\rb^{(2)})}_{s_2}$ are defined by (\ref{eq3:def_local}), where $\ra,\rb$ are replaced by  $\ra^{(2)}$ and $\rb^{(2)}$ respectively. We also set $\lambda_1=\lambda_2=1$ since the functions $s_1,s_2$ are chosen independently of the collection $\left\{\Theta^{(2)}_\D,\;\D\in \big(0, D_\rd(\Theta)\big]\right\}$.
%\begin{equation*}
%\widehat{\ell}(\D)=\ell\left(\frac{D_\ra(\Theta)}{g_\cA(\D)}\right)+\ell\left(\frac{D_\rb(\Theta)}{g_\cB(\D)}\right).
%\end{equation*}
Choose
$
R_{0}(u,v)=\ell\Big(\overline{\tau}_1u^{-1}\Big)+\ell\Big(\overline{\tau}_2v^{-1}\Big)
$ and define
\begin{eqnarray*}
\widehat{V}_{\vec{s}}^{(z,\e)}(\D)&=&(1+\e)g_\cA\Big([1+\e]^2\D\Big)\sqrt{2\big[1+\e^{-1}\big]^{2}\widehat{\cE}^{(\e)}(\D)
+(1+\e)^{2}\Big[2\ell\Big(D_\rd(\Theta)\big/\D\Big)+z\Big]}
\\
&+&
(1+\e)^{2}g_\cB\Big([1+\e]^2\D\Big)\Big\{2\big[1+\e^{-1}\big]^{2}\widehat{\cE}^{(\e)}(\D)+(1+\e)^2\Big[2\ell\Big(D_\rd(\Theta)\big/\D\Big)+z\Big]\Big\}.
\end{eqnarray*}
Then, applying Proposition \ref{prop_uniform_local3} and taking into account (\ref{eq:constants-choice-of-function-R}) we come to the following result.
\begin{proposition}
\label{prop_uniform_local5}
Let   Assumptions \ref{ass:fixed_theta_local}-\ref{ass:parameter_local} be fulfilled.
Then $\forall s\in\cS_{\ra,\rb}$, $\forall \e>0$ and
$\forall z\geq 1$
\begin{eqnarray*}
&&
\mathrm{P}\left\{\sup_{\D\in \big(0, D_\rd(\Theta)\big]}\Big[\mm_\D- \widehat{V}_s^{(z,\e)}(\D)\Big]\geq 0\right\}
\leq 4\mathrm{c}\left[2+\Big[\ln{\left\{1+\ln{(1+\e)}\right\}}\Big]^{-2}\right]^{2}\exp{\left\{-z\right\}}.
\end{eqnarray*}
\end{proposition}
The obtained inequality allows, in particular, to prove that the families of probabilities measures generated by $\chi_\theta$ is dense.
This, in its turn, is crucial step in proving of the weak convergence of probabilities measures.

\section{Gaussian random functions}
\label{sec:gauss}

In this section we apply Propositions \ref{prop_uniform_local2}-\ref{prop_uniform_local4} to the family of zero-mean gaussian random functions. Thus, let
$\chi_\theta,\theta\in\Theta,$ is a real valued continuous gaussian  random function such that $\mathrm{E}\chi_\theta=0,\;\forall\theta\in\Theta$. We are interested first in finding an upper function for $\big|\chi_\theta\big|,\theta\in\Theta$.
Let
$$
V(\theta)=\sqrt{\mathrm{E}\left|\chi_\theta\right|^{2}},\quad \rho(\theta_1,\theta_2)=\sqrt{\mathrm{E}\left|\chi_{\theta_1}-\chi_{\theta_2}\right|^{2}}
$$
We remark that Assumption \ref{ass:fixed_theta_local}
holds with $\mathrm{c}=2$,  $B\equiv 0$ and $\mathrm{b}\equiv 0$ and $\forall A\geq\sqrt{2}V$, $\forall \mathrm{a}\geq\sqrt{2}\rho$.
%It is  obvious  that this assumption  is  fulfilled if, for instance $\mathrm{a}=\sqrt{2}\rho$ and $ A\geq\sqrt{2}V$.
Since $\mathrm{b}\equiv 0$ Assumption \ref{ass:parameter_local}
is reduced to
\begin{assGauss}
\label{ass:Gauss-case-local}
%\label{ass:parameter_local}
There exist $s\in\bS$ such that for any $x>0$
$$
\sup_{\delta>0}\delta^{-2}
\mathfrak{E}_{\widetilde{\Theta},\;\mathrm{a}}\left(x(48\delta)^{-1}s(\delta)\right)<\infty.
$$
\end{assGauss}
\noindent Thus, if the latter assumption holds, Propositions \ref{prop_uniform_local2}-\ref{prop_uniform_local4} can be applied.

The aim of this section is to find uppers functions for quite different functionals of various gaussian processes. We would like to emphasize that the original problem is not always related to the consideration of $\big|\chi_\theta\big|,\theta\in\Theta,$ although such problems are also studied.   The idea is to reduce it (if necessary)   to whose for which one of Propositions \ref{prop_uniform_local2}-\ref{prop_uniform_local4} can be used.
%related to finding an upper function for, say, $\Phi(\xi_w),\; w\inW\$
Without special  mentionning we will always consider a separable modification of $\chi_\theta,\theta\in\Theta.$

\subsection{Upper functions for $\bL_p$-norms of Wiener integrals}
\label{sec:subsection-norm-gauss}

Let $K:\bR^d\to\bR$ be a continuous compactly supported function such that $\|K\|_\infty<\infty $. Without loss of generality we will assume that the support of $K$ is $[-1/2,1/2]^{d}$.
Let $0<h^{(\min)}\leq h^{(\max)}\leq 1$ be given numbers.
%with strictly positive coordinates and such that $h^{(\min)}_j\leq h^{(\max)}_j,\;j=\overline{1,d}$. Without loss of generality we will assume that $h^{(\max)}\in (0,1]^{d}$.
Put $\cH=\Big[h^{(\min)}, h^{(\max)}\Big]^{d}$    and let $K_h({\boldsymbol{\cdot}})=h^{-d}K\left({\boldsymbol{\cdot}}/h\right),\;h\in\cH$. Here $u/v,\; u,v\in\bR^d$ denotes coordinate-wise division.
Let $b(\rd t)$ is white noise on $\bR^d$ and consider the family of gaussian random fields
$$
\xi_h(t)=h^{-d}\int_{\bR^d}K_h\left(t-u\right)b(\rd u),\quad h\in\cH.
$$
%\subsubsection{Upper functions for $\bL_p$-norms of Wiener integrals}
%\label{sec:subsection-norm-gauss}
%Furthermore we will  assume that  $\big|\int K\big|\geq 1/\sqrt{2}$ although our results are valid whenever $\int K \neq 0$.
Let $\bK_\mu=[-\mu/2,\mu/2]^{d},\; \mu\geq 1,$ be a given cube and let for any $1\leq p<\infty$
$$
\big\|\xi_h\big\|_p=\left(\int_{\bK_\mu}\big|\xi_h(t)|^p\rd t\right)^{\frac{1}{p}}.
$$
The objective is to find an upper function for $\big\|\xi_h\big\|_p$ on $\cH$ and later on $C_1,C_2\ldots,$ denote the constants completely determined by $d,p,\mu$, $\gamma$ and $K$. It is  worth mentioning  that the explicit values of these constants can be found and some of them are given in the proof of the theorem.

We will be interested only the case $2\leq p<\infty$, since for any $p\in [1,2)$ we obviously have
$$
\big\|\xi_h\big\|_p\leq (\mu)^{\frac{d(2-p)}{2p}}\big\|\xi_h\big\|_2
$$
and, therefore, we can use the upper function found for $p=2$ for any  $p\in [1,2)$.

Let $\bB^{s}_{q,r}, s>0,\;1\leq q,r\leq \infty,$ denote the Besov space on $\bR^d$, see e.g. \cite{EdmTri}, and later on $\bH_q(s,L)$ denote the the ball of the radius $L>0$ in $\bB^{s}_{q,\infty}$.

Suppose that $K\in\bH_\infty(\gamma,L)$ and  without loss of generality we assume that $L=1$ that implies in particular that $\|K\|_\infty\leq 1.$
\begin{theorem}
\label{th:gauss-norm}
Assume  that $\gamma>d/2$. Then for any  $2\leq p<\infty$,  $h^{(\min)}, h^{(\max)}\in (0, 1)$,
 and   $ q\geq 1$
\begin{eqnarray*}
\mathrm{P}\bigg\{\sup_{h\in\cH}\left[\big\|\xi_h\big\|_p - C_1 h^{-d/2}\right]\geq 0\bigg\}
&\leq& C_2\exp{\left\{-2^{-3/2}\left(h^{(\max)}\right)^{-2d/p}\right\}}
\\
%*[2mm]
\mathrm{E}\bigg\{\sup_{h\in\cH}\left[\big\|\xi_h\big\|_p - C_1 h^{-d/2}\right]\bigg\}_+^q
&\leq& C_3(q)\;\left(h^{(\max)}\right)^{\frac{qd(2-p)}{2p}}
\exp{\left\{-2^{-3/2}\left(h^{(\max)}\right)^{-2d/p}\right\}}.
\end{eqnarray*}

\end{theorem}
The proof of the theorem is given in Section \ref{sec:proofs-of-gauss-theorems}.
The constant $C_1$ involved in the description of found upper functions is bounded function of $p$ on any bounded interval. Thus,
the upper functions are independent of $p$ if $p\in [2,p_0]$ for any given $p_0\geq 2$.

Also it is important to mention that the obtained upper functions are sharp. Indeed, it is not difficult to prove that for any $h>0$
$$
C_4h^{-d/2}\leq \mathrm{E}\big\|\xi_h\big\|_p\leq C_5 h^{-d/2}.
$$
This, together with the concentration inequality for gaussian processes, \cite{Talagrand},  yields in particular for any given $h>0$
$$
\mathrm{P}\left\{\big\|\xi_h\big\|_p \geq \widetilde{C}_1 h^{-d/2}\right\}
\leq \widetilde{C}_2\exp{\left\{-\widetilde{C}_3 h^{-2d/p}\right\}}.
$$
This inequality coincides, up to numerical constants, with the first inequality in Theorem \ref{th:gauss-norm} in particular case when
$h^{(\min)}=h^{(\max)}=h$.

\subsection{Upper functions for  local modulus of continuity under doubling condition}
\label{sec:subsection-gauss-local-continuity}

Let metric space $(\mT,\rd)$ be equipped with Borel measure $\kappa$. This measure is doubling if $\exists Q\geq 1$ such that
$$
\kappa\big\{\bB_\rd(\mt,2r)\big\}\leq Q\kappa\big\{\bB_\rd(\mt,r)\big\},\;\;\forall \mt\in\mT,\;\;\forall r>0,
$$
where $\bB_\rd(\mt,r)$ is the closed ball with center $\mt$ and radius $r$. For example, if $\mT=\bR^d$ and $\kappa$ is Lebesgue measure then
$Q=2^d$.

As it was proved in \cite{CoifmanWeiss} the existence of a doubling measure on $\mT$ implies that the space $\mT$ is doubling. It means that there exists $N_\rd\in\bN^*$ depending only on $Q$ such that for any $r>0$ each closed ball in $\mT$ of radius $r$ can be covered by at most $N_\rd$ closed balls of radius $r/2$.  This yields that $\bB(\mt,r)$ is totally bounded for any $\mt\in\mT$, $r>0$ and, moreover,
$$
\mE_{\bB_\rd(\mt,r),\;\rd}(\delta)\leq \ln{(N_\rd)}\left(\left[\log_2{\big\{r/\delta\big\}}\right]_++1\right),\;\;\forall\delta>0.
$$
Let $\mt\in\mT$ and $r>0$ be fixed. In this section, using Proposition \ref{prop_uniform_local4} we establish the upper function for local modulus of continuity of $\chi_\theta$  on $\Theta:=\bB_\rd(\mt,r)$. The simplest consequence of this result will be  the law of iterated logarithm (LIL) for $|\chi_\theta-\chi_\mt|$  as well as its non-asymptotical version. Studying the local modulus of continuity we are obviously interested in the case when $r$ is small even $r\to 0$. Thus, without loss of generality we will assume that $r\leq 1$.

To apply Proposition \ref{prop_uniform_local4}  we need to define the function $g_A$, $A(\cdot)=\sqrt{2}\rho(\cdot-\mt)$,  to choose the function $s_1\in\cS_{\ra,\mathrm{0}}$, and to compute the function $\widehat{\cE}^{(\e)}(\D,\mt)$, $\D\in (0,r]$, $\e\in\big(0,\sqrt{2}-1\big]$, given by
\begin{equation*}
\widehat{\cE}^{(\e)}(\D,\mt)=e^{(\ra)}_{s_1}\Big(g_A\big([1+\e]\D\big),\bB_\rd\big(\mt,[1+\e]\D\big)\Big),\quad \ra=\sqrt{2}\rho.
\end{equation*}
Introduce the function
$$
\psi(x)=\sqrt{2}\sup_{\substack{\te_1,\te_2\in\bB_\rd\big(\mt,1\big):
\\\rd\big(\te_1,\te_2\big)\leq x}}\rho\big(\te_1,\te_2\big),\;\;x\in(0,2],
$$
and suppose that $\psi(2)<\infty$. Note that obviously $\psi(0)=0$, since $\rd$ is a metric, and $\psi$ is increasing. Moreover, for any $u\in (0,r]$
$$
g^*_A(u):=\sqrt{2}\sup_{\theta:\;\rd(\theta,\mt)\leq u}\rho(\theta,\mt)\leq \psi(u),
$$
that allows us to put $g_A=\psi$.
%The following assumption is supposed to be held:
%\begin{equation}
%\label{eq:asumption-on-psi}
%\sup_{\D\in (0,r]}\D\big[\psi(\D)]^{-1}=:C_\psi(r)<\infty\;\;\forall r>0.
%\end{equation}
Denoting $\psi^{-1}$ the inverse function of $\psi$ we assert that  $\forall u\in(0,r]$
\begin{equation}
\label{eq:formula-for-entropy}
\mE_{\bB_\rd(\mt,u),\;\sqrt{2}\rho}(\delta)\leq \mE_{\bB_\rd(\mt,u),\;\rd}\Big(\psi^{-1}(\delta)\Big)\leq \ln{(N_\rd)}\left(\left[\log_2{\left\{u\big/\psi^{-1}(\delta)\right\}}\right]_++1\right),\;\;\forall\delta>0.
\end{equation}
Hence, if the function $\psi$ is such that Assumption \ref{ass:Gauss-case-local} is fulfilled then Proposition \ref{prop_uniform_local4} is applicable and that provides us with  the upper function for $|\chi_\theta-\chi_\mt|$  on $\bB_\rd(\mt,r)$. However, this upper function does not admit an explicit expression, in particular its dependence on the variable $\D$ cannot be analyzed in general. So, we prefer to impose an additional assumption on the function $\psi$ that allows us to obtain the explicit expression of the upper function for $|\chi_\theta-\chi_\mt|$ and analyze it as well as the corresponding probability bound for small values of the radius $r$. We will not be tending here to the maximal generality and
 suppose that there exist $0<\underline{c}\leq\overline{c}<\infty$ and $\beta>0$ such that
\begin{equation}
\label{eq:asumption-on-psi}
\underline{c}u^{\beta}\leq\psi(u)\leq \overline{c}u^{\beta},\quad \forall u\in (0,1].
\end{equation}
For example, if $\rd=\rho$ one has $\psi(u)=\sqrt{2}u$, and, therefore, (\ref{eq:asumption-on-psi}) holds.
Under (\ref{eq:asumption-on-psi}) obviously $\big(\delta/\overline{c}\big)^{\frac{1}{\beta}}\leq\psi^{-1}(\delta)\leq\big(\delta/\underline{c}\big)^{\frac{1}{\beta}}$
and we get from (\ref{eq:formula-for-entropy})
\begin{equation}
\label{eq:formula-for-entropy-continuation}
\mE_{\bB_\rd(\mt,u),\;\sqrt{2}\rho}(\delta)\leq \ln{(N_\rd)}\left(\left[\log_2{\left\{u\big(\delta/\overline{c}\big)^{-\frac{1}{\beta}}\right\}}\right]_++1\right),\;\;\forall\delta>0.
\end{equation}
Taking into account that $g_A(u)=\psi(u)\geq \underline{c}u^{\beta} $
and choosing $
s(x)=s^*(x)=(6/\pi^2)\big(1+[\ln{x}]^{2}\big)^{-1},$
 we obtain from (\ref{eq:formula-for-entropy-continuation}) for any $\D\in (0,r]$
\begin{eqnarray*}
\widehat{\cE}^{(\e)}(\D,\mt)&=&\sup_{\delta>0}\delta^{-2}\mE_{\bB_\rd\big(\mt,[1+\e]\D\big),\;\sqrt{2}\rho}
\left(g_A\big([1+\e]\D\big)(48\delta)^{-1}s(\delta)\right)
\\*[2mm]
&\leq& \sup_{\delta>0}\delta^{-2}\mE_{\bB_\rd\big(\mt,[1+\e]\D\big),\;\sqrt{2}\rho}
\left(\underline{c}\big([1+\e]\D\big)^{\beta}(48\delta)^{-1}s(\delta)\right)
\\*[2mm]
&\leq& \ln{(N_\rd)}\sup_{\delta>0}\delta^{-2}\left(\beta^{-1}\left[\log_2{\left\{\frac{48\overline{c}}{\underline{c}}\right\}}+
\log_2{\left\{\frac{\delta}{s^{*}(\delta)}\right\}}\right]_++1\right)=:\mathrm{C}(\beta,\underline{c},\overline{c},\rd).
\end{eqnarray*}
As we see $\widehat{\cE}^{(\e)}(\D,\mt)$ is independent of $\D,\mt$ and bounded from above by the constant which is completely determined by the triplet  $(\mT,\rd,\kappa)$ and by the quantities $\underline{c}, \overline{c}$ and $\beta$. Thus, in view of Proposition \ref{prop_uniform_local4} the upper function for  has the following form (see also Remark \ref{rem2:after-proposition4-local}).
\begin{eqnarray*}
\widehat{V}_{s^*}^{(z,\e)}(\D,\mt)=\overline{c}(1+\e)^{1+2\beta}\D^{\beta}\sqrt{2\big[1+\e^{-1}\big]^{2}
\mathrm{C}(\beta,\underline{c},\overline{c},\rd)
+(1+\e)^{2}\Big[\ell\Big(2r\big/\D\Big)+z\Big]},
\end{eqnarray*}
where, remind, $
\ell(y)=\ln{\left\{1+\ln{(y)}\right\}}+2\ln{\left\{1+\ln{\left\{1+\ln{(y)}\right\}}\right\}}, y>0.
$

\smallskip

Choose  $z=z(r)=\ln{\left\{1+\ln{\left\{1+\big|\ln{(r)}\big|\right\}}\right\}}$ and $\e=\e(r):=z^{-1}(r)$ and define
\begin{eqnarray*}
\mathfrak{a}(r)&=&(1+\e(r))^{1+2\beta}\left(\sup_{\D\in (0,r] }\sqrt{\frac{2\big[1+\e^{-1}(r)\big]^{2}\mathrm{C}(\beta,\underline{c},\overline{c},\rd)
+(1+\e(r))^{2}\Big[\ell\Big(2r\big/\D\Big)+z(r)\Big]}{\ln{\left\{1+\big|\ln{(\D)}\big|\right\}}}}\right);
\\
\mathfrak{p}(r)&=&\frac{2+\Big[\ln{\left\{1+\ln{(1+\e(r))}\right\}}\Big]^{-2}}{1+\ln{\left\{1+\big|\ln{(r)}\big|\right\}}}.
\end{eqnarray*}
We note that if $r\to 0$ then
\begin{eqnarray}
\label{eq:asymp-property-of-ma(r)}
\mathfrak{a}(r)=1+\mathcal{O}\left(\frac{\ln{\left\{1+\ln{\left\{1+\big|\ln{(r)}\big|\right\}}\right\}}}{\ln{\left\{1+\big|\ln{(r)}\big|
\right\}}}\right),\quad
\mathfrak{p}(r)=\mathcal{O}\left(\frac{\left[\ln{\left\{1+\ln{\left\{1+\big|\ln{(r)}\big|\right\}}\right\}}\right]^2}
{\ln{\left\{1+\big|\ln{(r)}\big|
\right\}}}\right).
\end{eqnarray}
The following result is immediate consequence of  Proposition \ref{prop_uniform_local4}. Put  $\mm(\D)=\sup_{\te\in B_\rd(\mt,\D)}|\chi_\theta-\chi_\mt|$.
\begin{theorem}
\label{th:gauss-local-modulus}
Let $\mT$ be doubling space and suppose that (\ref{eq:asumption-on-psi}) holds.
Then,  we have for any $\mt\in\mT$ and any $r\in (0,1)$
\begin{eqnarray*}
&&
\mathrm{P}\left\{\sup_{\D\in \big(0,r\big]}\left[\frac{\mm(\D)}{\overline{c}\D^{\beta}\sqrt{\ln{\left\{1+\big|\ln{(\D)}\big|\right\}}}}\right]\geq \mathfrak{a}(r)\right\}
\leq 8\mathfrak{p}(r).
\end{eqnarray*}

\end{theorem}

The first consequence of Theorem \ref{th:gauss-local-modulus} is the law of iterated logarithm.  Indeed,  taking into account
that $\mathfrak{p}(r)\to 0, \;\;\mathfrak{a}(r)\to 1,\;r\to 0, $ we come to the following assertion.
\begin{corollary}
\label{cor1:after-theorem-gauss-modulus} Let $\mT$ be doubling space and suppose that (\ref{eq:asumption-on-psi}) holds. Then for any $\mt\in\mT$ $\mathrm{P}-\text{a.s.}$
$$
\limsup_{\D\to 0+}\left[\frac{\mm(\D)}{\D^{\beta}\sqrt{\ln{\left\{1+\big|\ln{(\D)}\big|\right\}}}}\right]\leq \overline{c}.
$$

\end{corollary}
We note, that although the statement of Corollary \ref{cor1:after-theorem-gauss-modulus} is traditional in probability theory, the non-asymptotical statement of
Theorem \ref{th:gauss-local-modulus} is much more informative.

The next consequence of Theorem \ref{th:gauss-local-modulus} seems more curious. We remark that if $\mT$ is doubling with respect to the intrinsic  semi-metric $\rho$, then the normalizing sequence appeared in the theorem  is  independent of $\rho$. Moreover, the function $\ma(\cdot)$
depends only on $\rho$.

Indeed, if $\rd=\rho$ then $\phi(u)=\sqrt{2}u$, and therefore, $\beta=1$ and $\underline{c}=\overline{c}=\sqrt{2}$. It yields, in particular,
that $\mathrm{C}(\beta,\underline{c},\overline{c},\rho)=\mathrm{C}_{N_\rho}$, where
$$
\mathrm{C}_{N_\rho}=\ln{(N_\rho)}\sup_{\delta>0}\delta^{-2}\left(\left[3+2\log_2{\left(\pi\right)}+
\log_2{\left\{\delta\big(1+[\ln{\delta}]^{2}\big)\right\}}\right]_++1\right),
$$
and, therefore, $\ma(\cdot)=\ma_{N_\rho}(\cdot)$, where $\ma_{N_\rho}(\cdot)$ is completely determined by $\rho$ via the quantity $N_\rho$.
\begin{corollary}
\label{cor2:after-theorem-gauss-modulus}
Let $\mT$ be doubling space with respect to $\rd=\rho$.
Then,  for any $\mt\in\mT$ and $r\in (0,1)$
\begin{eqnarray*}
&&
\mathrm{P}\left\{\sup_{\D\in \big(0,r\big]}\left[\frac{\mm(\D)}{\D\sqrt{2\ln{\left\{1+\big|\ln{(\D)}\big|\right\}}}}\right]\geq \mathfrak{a}_{N_\rho}(r)\right\}
\leq 8\mathfrak{p}(r).
\end{eqnarray*}

\end{corollary}
We note that if $\mathfrak{R}^*$ be the set of metrics $\rho$ for which $\mT$ is doubling and such that $N_\rho\leq N^*$ for some fixed
$N^*\in\bN^*$ then the function $\ma_{N_\rho}(\cdot)$ in the assertion of Corollary \ref{cor2:after-theorem-gauss-modulus} can be replaced by the universal on $\mathfrak{R}^*$ function $\ma_{N^*}(\cdot)$. The  corresponding inequality becomes "metric free".

We finish this section by the consideration of several examples, where the condition (\ref{eq:asumption-on-psi}) is verified. In these examples
$\mT$ is either $\bR^d$ or $\bR^d_+, d\geq 1,$ $\kappa$ is Lebesgue measure and $\rd$ is the euclidian distance.

\paragraph{Example 1. L\'evy function} Here $\chi_\mt,\;\mt\in\bR^d,$ is zero-mean gaussian random field with $\rho=\sqrt{d}$. Hence,  (\ref{eq:asumption-on-psi}) holds with $\underline{c}=\overline{c}=\sqrt{2}$ and $\beta=1/2$. We deduce from Theorem \ref{th:gauss-local-modulus} that
\begin{eqnarray*}
&&
\mathrm{P}\left\{\sup_{\D\in \big(0,r\big]}\left[\frac{\mm(\D)}{\sqrt{2\D\ln{\left\{1+\big|\ln{(\D)}\big|\right\}}}}\right]\geq \mathfrak{a}(r)\right\}
\leq 8\mathfrak{p}(r),\;\;r\in(0,1).
\end{eqnarray*}
\paragraph{Example 2. Fractional  brownian motion} Here $\chi_\mt,\; \mt\in\bR_+,$ is zero-mean gaussian random process with $\rho=d^{\alpha/2},\alpha\in (0,2]$. Hence,  (\ref{eq:asumption-on-psi}) holds with $\underline{c}=\overline{c}=\sqrt{2}$ and $\beta=\alpha/2$. We get  from Theorem \ref{th:gauss-local-modulus}
\begin{eqnarray*}
&&
\mathrm{P}\left\{\sup_{\D\in \big(0,r\big]}\left[\frac{\mm(\D)}{\sqrt{2\D^{\alpha}\ln{\left\{1+\big|\ln{(\D)}\big|\right\}}}}\right]\geq \mathfrak{a}(r)\right\}
\leq 8\mathfrak{p}(r),\;\;r\in(0,1).
\end{eqnarray*}
\paragraph{Example 3. Ornstein-Uhlenbeck process} Here $\chi_\mt,\;\mt\in\bR_+,$ is given by
$$
\chi_t=\frac{\sigma}{\sqrt{2\lambda}}b\left(e^{2\theta t}\right)e^{-\lambda t},\;\;\lambda,\sigma>0,
$$
where $b$ is the standard Wiener process.
In this case
$\rho=\lambda^{-1/2}\sigma\sqrt{1-\exp{\{-\lambda d\}}}$.
Since we consider $r\in (0,1)$ then
(\ref{eq:asumption-on-psi}) holds with $\underline{c}=\sigma\sqrt{2}e^{-1}e^{-\lambda/2}$, $\overline{c}=\sigma\sqrt{2}$, $\beta=1/2$ and we have
\begin{eqnarray*}
&&
\mathrm{P}\left\{\sup_{\D\in \big(0,r\big]}\left[\frac{\mm(\D)}{\sigma\sqrt{2\D\ln{\left\{1+\big|\ln{(\D)}\big|\right\}}}}\right]\geq \mathfrak{a}(r)\right\}
\leq 8\mathfrak{p}(r),\;\;r\in(0,1).
\end{eqnarray*}

\section{Generalized  empirical processes}
\label{sec:empirical-processes}

Let $\big(\cX,\mX,\nu\big)$ be $\sigma$-finite space  and let $\left(\Omega,\mA,\textsf{P}\right)$ be a   probability space.
Let $X_i, \; i\geq 1,$ be a the collection of $\cX$-valued {\it independent} random variables  defined on $\left(\Omega,\mA, \textsf{P}\right)$ and  having the densities $f_i$ with respect to measure $\nu$. Furthermore, $\bP_{\rf},\;\rf=(f_1,f_2,\ldots),$ denotes the probability law   of  $(X_1,X_2,\ldots)$ and $\bE_{\rf}$ is  mathematical expectation with respect to $\bP_\rf$.

Let $G:\boldsymbol{\mH}\times\cX\to\bR$ be a given mapping, where  $\boldsymbol{\mH}$ is a set.
Put  $\forall n\in\bN^*$
\begin{eqnarray}
\label{eq:generilized-empirical-process}
\xi_\mh(n)=n^{-1}\sum_{i=1}^{n}\Big[G\big(\mh,X_i\big)-\bE_{\rf} G(\mh,X_i)\Big],\quad \mh\in\boldsymbol{\mH}.
\end{eqnarray}
 We will say that $\xi_\mh(n),\;\mh\in\boldsymbol{\mH},$ is generalized empirical process. Note that if $\mh:\cX\to\bR$ and  $G(\mh,x)=\mh(x),
 \;\mh\in\boldsymbol{\mH}, x\in\cX,$ then $\xi_\mh(n)$ is the standard empirical process parameterized by $\boldsymbol{\mH}$.

Throughout this section we will suppose that
\begin{eqnarray}
\label{eq:boundness-empirical-proc}
\overline{\mathbf{G}}_\infty(\mh):=\sup_{x\in\cX}\big|G(\mh,x)\big|<\infty,\quad \forall \mh\in\boldsymbol{\mH},
\end{eqnarray}
and it will be referred to {\it bounded case}. Some generalizations  concerning the situations where this assumption fails are discussed in Section
\ref{sec:some-generalizations}.

The condition (\ref{eq:boundness-empirical-proc}) implies that the random variables
$G(\mh,X_i),\mh\in\boldsymbol{\mH}$, and $G(\mh_1,X_i)-G(\mh_2,X_j)$, $\mh_1,\mh_2\in\boldsymbol{\mH},$ $i=\overline{1,n}$, are bounded, and we obtain in view of Bernstein inequality $\forall z>0$
\begin{eqnarray}
\label{eq:empircal-bernstein1}
\bP_\rf\left\{\big|\xi_\mh(n)\big|>z\right\}&\leq& 2\exp{\left\{-\frac{z^{2}}{A_\rf^{2}(\mh)+z B_\infty(\mh) }\right\}};
\\
\label{eq:empircal-bernstein2}
\bP_\rf\left\{\Big|\xi_{\mh_1}(n)-\xi_{\mh_2}(n)\Big|>z\right\}
%\\*[2mm]&&
&\leq& 2\exp{\left\{-\frac{z^{2}}{\ra_\rf^{2}(\mh_1,\mh_2)+z \rb_\infty(\mh_1,\mh_2) }\right\}},
\end{eqnarray}
where
\begin{eqnarray}
\label{eq1:def-A-B-a-b-empirical}
&\displaystyle{A^{2}_{\rf}(\mh)=2n^{-2}\sum_{i=1}^{n}\bE_fG^{2}(\mh,X_i)},\quad &\ra^{2}_{\rf}(\mh_1,\mh_2)=2n^{-2}\sum_{i=1}^{n}\bE_f\big(G(\mh_1,X_i)-G(\mh_2,X_i)\big)^{2};
\\
\label{eq2:def-A-B-a-b-empirical}
&\;\;\;\;\;B_\infty(\mh)=(4/3)n^{-1}\sup_{x\in\cX}\big|G(\mh,x)\big|,\quad &\rb_\infty(\mh_1,\mh_2)=(4/3)n^{-1}\sup_{x\in\cX}\big|G(\mh_1,x)-G(\mh_2,x)\big|.
\end{eqnarray}
We conclude that {\it Assumption \ref{ass:fixed_theta_local} is fulfilled} with $\Psi(\cdot)=|\cdot|$,\; $A=A_\rf,\; B=B_\infty,\; \ra=\ra_\rf,\;\rb=\rb_\infty$ and $\mathrm{c}=2$.

 It is easily seen that $\ra_\rf$ and $\rb_\infty$ are  semi-metrics on $\boldsymbol{\mH}$.
 We note also that   $\xi_\bullet:\boldsymbol{\mH}\to\bR$  is $\mathrm{P}$-a.s  continuous  in the topology generated by  $\rb_\infty$.
 Thus, if $\mH\subseteq \boldsymbol{\mH}$ is  totally bounded with respect to $\ra_f\vee\rb_\infty$  and such that $\overline{A}_\rf:=\sup_{\mh\in\mH}A_\rf(\mh)<\infty$,
 $\overline{B}_\infty:=\sup_{\mh\in\mH}B_\infty(\mh)<\infty$, then
  we conclude that {\it Assumption \ref{ass:metric-case-local} is verified}.

 Thus, in the problems for which Assumption \ref{ass:parameter_local} is verified the machinery developed in Propositions
\ref{prop_uniform_local2}-\ref{prop_uniform_local3} can be applied for $\big|\xi_\mh(n)\big|,\mh\in\mH.$
We would like to emphasize, however, that   problems studied below are not always related to the consideration of
$\big|\xi_\mh(n)\big|,\mh\in\mH$ with $\mH$ being totally bounded, although such problems are also studied.
The idea is to reduce them (if necessary)   to whose for which one of Propositions \ref{prop_uniform_local2}-\ref{prop_uniform_local3}
can be used. For instance, we will be interested in finding upper functions for $\big|\xi_\mh(n)\big|$ on $\mh\in\mH$ not only for given $n$
but mostly on
$\mathbf{N}\times\mH$, where $\mathbf{N}$ is a given subset of $\bN^*$.
It will allow, in particular, to study generalized empirical processes with random number of summands.

However
the application of Propositions
\ref{prop_uniform_local2}-\ref{prop_uniform_local3} requires to compute
the functions $\cE$ or $\widehat{\cE}$ and there is no a general receipt how to do it.
The main goal of  this section is to provide with rather general assumptions under which the
latter quantities can be computed explicitly. As it was already mentioned
in Introduction upper functions for random objects appear in various
areas of probability theory and mathematical statistics.
As the consequence the different nature of problems requires to specify the imposed assumptions.
The assumptions presented below are oriented mostly to the problems arisen in mathematical statistics
 that definitely reflects author's scientific interests. However, some pure probabilistic results like the law of iterated logarithm
 and the law of logarithm will be established as well.

\subsection{Problem formulation and examples. Main condition}
\label{sec:subsec-prob-form}

In this section we find upper functions  for several functionals of the generalized empirical process  $\xi_\mh(n)$ defined in (\ref{eq:generilized-empirical-process}) under condition (\ref{eq:boundness-empirical-proc}).
We remark that the parameter $\mh$ may possess a composite structure and its components may have very different nature.
In order to treat such situations it will be convenient for us to assume that for some $m\geq 1$
\begin{eqnarray}
\label{eq:mH-space-decomposition}
\boldsymbol{\mH}=\mH_1\times\cdots\times\mH_m,
\end{eqnarray}
where $\mH_j,\,j=\overline{1,m},$ be given sets.
We will use the following notations. For any given $k=\overline{0,m}$ put
$$
\mH_{1}^{k}=\mH_1\times\cdots\times\mH_k,\quad \mH_{k+1}^{m}=\mH_{k+1}\times\cdots\times\mH_m,
$$
with the agreement that $\mH_{1}^{0}=\emptyset, \;\mH_{m+1}^{m}=\emptyset$. The elements of $\mH_{1}^{k}$ and $\mH_{k+1}^{m}$
will be denoted by $\mh^{(k)}$ and $\mh_{(k)}$ respectively.
%%%%%%%%%%%%%%%%%%%%%%%%%%%%%%%%%%%%%%%%%%%%
\iffalse
We will be interested in finding an upper function for
$$
\zeta_{\mh^{(k)}}:=\sup_{\mh_{(k)}\in\mH_{k+1}^{m}}\big|\xi_\mh\big|,
$$
on $\mH_{1}^{k}$. It will include as particular cases the upper function for $\xi_\mh$ on $\mH$ ($k=m$) and the upper functions over cylindric sets $\mH_{1}^{k}\times\big\{\bar{\mh}_{(k)}\big\}$, where $\bar{\mh}_{(k)}:=\big(\bar{\mh}_{k+1},\ldots,\bar{\mh}_{m}\big)$ is fixed.
\fi
%%%%%%%%%%%%%%%%%%%%%%%%%%%%%%%%%%%%%%%%
We will suppose that for any $ j=\overline{k+1,m}$ the set $\mH_j$ is
 endowed with the   semi-metric $\varrho_j$ and the Borel measure $\kappa_j$.

 In the next two sections we  find  upper functions for $|\xi_\mh(n)|$ on
 on some subsets of  $\mH$ (possibly  depending on $n$!)
 and we will consider two cases.

%\smallskip

{\it Totally bounded case.}  In this case we will suppose that  $\mH_j$ is totally bounded with respect to
$\varrho_j$ for any $j=\overline{k+1,m}$.

\smallskip

{\it Partially totally bounded case.}  Here we  first suppose that for some $p\geq 1$
\begin{eqnarray}
\label{eq:cX-space-decomposition}
(\cX,\nu)=\big(\cX_1\times\cdots\times\cX_p,\;\nu_1\times\cdots\times\nu_p\big),
\end{eqnarray}
where $(\cX_l,\nu_l)\;l=\overline{1,p},$ are of measurable spaces and $\nu$ is the product measure.

Next we will assume that $\mH_{m}=\cX_1$. As the consequence, the assumption, that $\mH_m$ is totally bounded, is too restrictive. In particular,  it does not verified in the case $\cX=\cX_1=\bR^d$ which appears in many examples. Before to start with the presentation of the results let us consider several examples.
%Moreover, we will suppose that for some $\mu>0$
%\begin{eqnarray}
%\label{eq:localization}
%G(\mh,x)=0, \quad \forall x=(x_1,\dots,x_p)\in\cX,\;\forall\mh=(\mh_1,\dots,\mh_m)\in\mH:\;\; \varrho_1(x_1,\mh_{k+1})\leq\mu
%\end{eqnarray}
%and this explains why the studied case is called localized.

\smallskip

{\bf Example 1. Density model.}
Let  $K:\bR^d\to\bR$ be a given function compactly supported on $[-1/2,1/2]^{d}$ and
denote for any $h=(h_1,\ldots,h_d)\in (0,1]^d$
$$
K_h(\cdot)=\left[\prod_{i=1}^{d}h_i\right]^{-1}K\left(\cdot/h_1,\ldots, \cdot/h_d\right),
$$
where, as previously, for two vectors $u,v\in\bR^d$ the notation $u/v$ denotes the coordinate-vice division.

%Let $h^{\min}=\big(h^{\min}_1,\ldots h^{\min}_d\big),\;h^{\max}=\big(h^{\max}_1,\ldots h^{\max}_d\big)\in (0,1]^d$ be given vectors and let
%$$
%\cH^d=\cH_1\times\cdots\cH_d,\quad \cH_i:=\big[h^{\min}_i,h^{\max}_i\big],\;i=\overline{1,d}.
%$$
Put $p=1$, $m=d+1$, $k=d$, $\cX_1=\mH_{d+1}=\bR^d$, $\mH_i=(0,1],\;i=\overline{1,d}$ and consider for any $\mh=(h,x)\in\boldsymbol{\mH}:=(0,1]^d\times\bR^{d}$
$$
\xi_\mh(n)=\widehat{\xi}_{h,x}(n):=n^{-1}\sum_{i=1}^{n}\bigg[K_{h}\left(X_i-x\right)-\bE_{\rf}\left\{
K_h\left(X_i-x\right)\right\} \bigg].
$$
We have come to the well-known in nonparametric statistics kernel density estimation process.
Here the function $K$ is a kernel and the vector $h$ is a multi-bandwidth.
%Note that, in particular, in this case $\zeta_{\mh^{(d)}}=\|\widehat{\xi}_h\|_\infty$. In Section \ref{sec:integrability_property} we will study also the case  $\mH=\cH_d\times \{t\}$, where $t\in\bR^d$ is fixed.

\smallskip

{\bf Example 2. Regression model.}
Let $\e_i,i=\overline{1,n},$ be independent real random variables distributed on  $\cI\subseteq\bR$ and such that $\bE\e_i=0$ for any
$i=\overline{1,n}$. Let $Y_i,i=\overline{1,n},$ be independent $d$-dimension random vectors. The sequences
$\left\{\e_i,i=\overline{1,n}\right\}$ and $\left\{Y_i,i=\overline{1,n}\right\}$ are assumed independent.
%\smallskip
Let $\cM$ be a given set of $d\times d$ invertible  matrices and let $\cI\subseteq \bR$ and $\cX_1\subseteq\bR^d$  be  given interval.

Put  $p=2$, $m=d+2$, $k=d$, $\cX_1=\mH_{d+2}=\bR^d$, $\cX_2=\cI$, $\mH_j=(0,1],\;j=\overline{1,d}$ and  $\mH_{d+1}=\cM$.
Consider for any $\mh=(x,h,M)\in\boldsymbol{\mH}:=(0,1]^d\times\cM\times\bR^d$
$$
\xi_\mh(n)=\widetilde{\xi}_{h,M,x}(n):=n^{-1}|\text{det}(M)|\sum_{j=i}^{n}K_{h}\Big[M(Y_i-x)\Big]\e_i.
$$
The family of random fields $\left\{\widetilde{\xi}_{x,h,M}(n),\;x,h,M\in(0,1]^d\times\cM\times \bR^d\right\}$
appears in non-parametric regression under single index hypothesis, \cite{Stone}.
%Note also that $\zeta_{\mh^{(d)}}=\|\widetilde{\xi}_{h,M}\|_\infty$.

If $\cI$ is bounded interval, i.e. $\e_i$ are bounded random variables, then (\ref{eq1:def-A-B-a-b-empirical}) and (\ref{eq2:def-A-B-a-b-empirical}) hold and the results from Section \ref{sec:Key proposition_} are applicable.
However this assumption is too restrictive and it does not satisfied even in the classical gaussian regression. At the first glance it is seemed that if $\cI=\bR$  Propositions \ref{prop_uniform_local2}-\ref{prop_uniform_local3} are not applicable here. Although the aforementioned  problem  lies beyond of the scope of the paper, let us briefly discuss how to reduce it to the problem in which the machinery developed in Propositions \ref{prop_uniform_local2}-\ref{prop_uniform_local3} can be applied.

\paragraph{Some generalizations.}
\label{sec:some-generalizations}
Let $\left(\e_i,i=\overline{1,n}\right)$ be the sequence of independent real-valued  random variables such that $\bE\e_i=0$ (later on for simplicity we assume that $\e_i$ has symmetric distribution) and $\bE\e_i^{2}=:\sigma^2_i<\infty$. Let $\bar{X}_i,\;i=\overline{1,n},$ be a the collection of $\bar{\cX}$-valued {\it independent} random elements and
suppose also that $\left(\bar{X}_i,i=\overline{1,n}\right)$ and $\left(\e_i,i=\overline{1,n}\right)$ are independent.
Consider the generalized empirical process
$$
\bar{\xi}_\mh(n)=n^{-1}\sum_{i=1}^{n}\bar{G}\big(\mh,\bar{X}_i\big)\e_i,\quad \mh\in\boldsymbol{\mH},
$$
where, as previously, $\bar{G}:\boldsymbol{\mH}\times\cX\to\bR$ be a given mapping satisfying (\ref{eq:boundness-empirical-proc}).
For any $y>0$ define
$$
\bar{\xi}_\mh(n,y)=n^{-1}\sum_{i=1}^{n}\bar{G}\big(\mh,\bar{X}_i\big)\e_i\mathrm{1}_{[-y,y]}(\e_i),
\quad\eta_n(y)=\sup_{i=\overline{1,n}}\big|\e_i\big|\big[1-\mathrm{1}_{[-y,y]}(\e_i)\big].
$$
Obviously,  for any $y>0$
$$
\bar{\xi}_\mh(n,y)=n^{-1}\sum_{i=1}^{n}\Big[G_y\big(\mh,X_i\big)-\bE_{\rf} G_y(\mh,X_i)\Big],\quad X_i=\big(\bar{X}_i,\e_i\big),
$$
where $G_y(\mh,x)=\bar{G}(\mh,\bar{x})\mathrm{1}_{[-y,y]}(u), \;x=(\bar{x},u)\in \cX:=\bar{\cX}\times\bR,\;\mh\in\boldsymbol{\mH}$.
Since $G_y$ is bounded for any $y>0$ the inequalities (\ref{eq:empircal-bernstein1}) and (\ref{eq:empircal-bernstein2}) hold and, analogously to (\ref{eq1:def-A-B-a-b-empirical}) and (\ref{eq1:def-A-B-a-b-empirical}), we have
\begin{eqnarray*}
&\displaystyle{A^{2}_{\rf}(\mh)=2n^{-2}\sum_{i=1}^{n}\sigma^2_i\bE_f\bar{G}^{2}(\mh,\bar{X}_i)},\; &\ra^{2}_{\rf}(\mh_1,\mh_2)=2n^{-2}\sum_{i=1}^{n}\sigma^2_i\bE_f\big(\bar{G}(\mh_1,\bar{X}_i)-G(\mh_2,\bar{X}_i)\big)^{2};
\\
&\;\;\;\;\;B_\infty(\mh)=(4y/3) n^{-1}\sup_{x\in\cX}\big|\bar{G}(\mh,\bar{x})\big|,\; &\rb_\infty(\mh_1,\mh_2)=(4/3)y n^{-1}\sup_{x\in\cX}\big|\bar{G}(\mh_1,\bar{x})-\bar{G}(\mh_2,\bar{x})\big|.
\end{eqnarray*}
Let also $\mH\subseteq\boldsymbol{\mH}$ be such that  the results obtained in Propositions \ref{prop_uniform_local2}-\ref{prop_uniform_local3} are applicable to $\left|\bar{\xi}_\mh(n,y)\right|$ on $\mH$ for any  $y>0$. It is extremely important to emphasize  that neither $A_{\rf}(\cdot)$ nor $\ra_{\rf}(\cdot,\cdot)$ depend on $y$.

This yields, in view of Theorems \ref{th:empiric_totaly_bounded_case} and \ref{th:empiric-partially_totaly_bounded_case} below, that upper functions for $\left|\bar{\xi}_\mh(y)\right|,\;\mh\in\mH$ (for brevity $V(\mh,y)$ and $U_q(\mh,y), q\geq 1$) can be found in the form:
$$
V(\mh,y)=V_1(\mh)+yV_2(\mh),\quad U_q(\mh,y)=U_{q,1}(\mh)+yU_{q,2}(\mh).
$$
It means that  we are able to bound from above any $y>0$
$$
\bP_f\left\{\sup_{\mh\in\mH}\left[\bar{\xi}_\mh(n,y)-V(\mh,y\right]>0\right\},\quad \bE_f\left\{\sup_{\mh\in\mH}\left[\bar{\xi}_\mh(n,y)-U_q(\mh,y\right]\right\}^q_+
$$
Moreover, we obviously have  for any $y>0$
\begin{eqnarray*}
&&\bP_f\bigg\{\sup_{\mh\in\mH}\left[\bar{\xi}_\mh(n)-V(\mh,y)\right]>0\bigg\}\leq
\bP_f\bigg\{\sup_{\mh\in\mH}\left[\bar{\xi}_\mh(n,y)-V(\mh,y\right]>0\bigg\}+
\bP_f\left\{\eta_n(y)>0\right\};
\\
&&
\bE_f\bigg\{\sup_{\mh\in\mH}\left[\bar{\xi}_\mh(n)-U_q(\mh,y\right]\bigg\}^q_+\leq \bE_f\bigg\{\sup_{\mh\in\mH}\left[\bar{\xi}_\mh(n,y)-U_q(\mh,y\right]\bigg\}^q_++\bigg(\sup_{\mh\in\mH}\overline{G}_\infty(\mh)\bigg)^q
\bE\left(\eta_n\right)^q.
\end{eqnarray*}
Typically, $V(\cdot,y)=V^{(n)}(\cdot,y)$ and $U_q(\cdot,y)=U^{(n)}_q(\cdot,y)$ and $V^{(n)}_2(\cdot)\ll V^{(n)}_1$ and
$U_{q,2}^{(n)}(\cdot)\ll U_{q,1}^{(n)}$ for all $n$ large enough. It allows to choose $y=y_n$ in optimal way, i.e. to balance both terms in latter inequalities, that usually leads to  sharp upper functions $V^{(n)}_1(\cdot)+y_nV^{(n)}_2(\cdot)$ and $U^{(n)}_{q,1}(\cdot)+y_nU^{(n)}_{q,2}(\cdot)$.

%\subsection{Upper functions in bounded case}
%\label{sec:upper functions in bounded case}

%Now let us come back to the consideration of generalized empirical processes obeying (\ref{eq:boundness-empirical-proc}).

\paragraph{Main Assumption}

Now let us come back to the consideration of generalized empirical processes obeying (\ref{eq:boundness-empirical-proc}).
Assumption \ref{ass:bounded_case}  below is the main tool allowing us to compute \textsf{explicitly}  upper functions.
Introduce the following notation: for any $\mh^{(k)}\in\mH_{1}^k$
%for any function $g:\cX\to \bR$ put $\|g\|_\infty=\sup_{x\in\cX}|g(x)|$ and set
%any $y\in\mH_j$ we put $\mh(j,y)=(\mh_1,\ldots,\mh_{j-1},y,\mh_{j+1},\ldots,\mh_m)$.
$$
\mathbf{G_\infty}\big(\mh^{(k)}\big)=\sup_{\mh_{(k)}\in\mH_{k+1}^m}\sup_{x\in\cX}|G(\mh,x)|,
$$
and let $G_\infty:\mH_{1}^k\to\bR_+$ be any mapping satisfying
\begin{equation}
\label{eq:def-of-function-G-infty}
\mathbf{G_\infty}\big(\mh^{(k)}\big)\leq G_\infty\big(\mh^{(k)}\big),\;\;\forall\mh^{(k)}\in\mH_{1}^k.
\end{equation}
Let $\left\{\mH_j(n)\subset\mH_j,\; n\geq 1\right\}, j=\overline{1,k},$ be a sequence of sets
and denote $\mH_{1}^k(n)=\mH_1(n)\times\cdots\mH_k(n)$.
Set for any $n\geq 1$
 $$
 \underline{G}_n=\inf_{\mh^{(k)}\in\mH_{1}^k(n)}G_\infty\big(\mh^{(k)}\big), \quad \overline{G}_n=\sup_{\mh^{(k)}\in\mH_{1}^k(n)}G_\infty\big(\mh^{(k)}\big).
$$
For  any $n\geq 1$, $j=\overline{1,k}$ and any $\mh_j\in\mH_j(n)$ define
$$
G_{j,n}(\mh_j)=\sup_{\mh_1\in\mH_1(n),\ldots,\mh_{j-1}\in\mH_{j-1}(n),\mh_{j+1}\in\mH_{j+1}(n),
\ldots,\mh_{k}\in\mH_{k}(n)}G_\infty\big(\mh^{(k)}\big),
\qquad \underline{G}_{j,n}=\inf_{\mh_j\in\mH_j(n)}G_{j,\infty}(\mh_j).
$$
Noting that $\big|\ln{(t_1)}-\ln{(t_2)}\big|$ is a metric on $\bR_+\setminus \{0\}$,
we equip  $\mH_{1}^k(n)$ with the following semi-metric.
For any $n\geq 1$ and any  $\hat{\mh}^{(k)},\bar{\mh}^{(k)}\in\mH_{1}^{k}(n)$ set
$$
\varrho_n^{(k)}\Big(\hat{\mh}^{(k)},\bar{\mh}^{(k)}\Big)=\max_{j=\overline{1,k}}\left|
\ln\big\{G_{j,n}(\hat{\mh}_j)\big\}-\ln\big\{G_{j,n}(\bar{\mh}_j)\big\}\right|,
$$
where $\hat{\mh}_j,\bar{\mh}_j,\; j=\overline{1,k},$ are the coordinates of
 $\hat{\mh}^{(k)}$ and $\bar{\mh}^{(k)}$ respectively.

\begin{assumption}
\label{ass:bounded_case}
\begin{itemize}

\item[$\mathbf{(i)}$]
$0<\underline{G}_n\leq\overline{G}_n<\infty $  for any  $n\geq 1$ and for any  $j=\overline{1,k}$
 $$
 \frac{G_\infty\big(\mh^{(k)}\big)}{\underline{G}_n}\geq \frac{G_{j,n}(\mh_j)}{\underline{G}_{j,n}},\quad \forall \mh^{(k)}=(\mh_1,\ldots,\mh_k)\in\mH_{1}^k(n),\;\;\forall n\geq 1;
 $$
 \item[$\mathbf{(ii)}$]
 There exist  functions $L_j:\bR_+\to \bR_+,\;D_j:\bR_+\to\bR_+,\;j=0,k+1,\ldots,m,$
  satisfying  $L_j$ non-decreasing and bounded on each bounded interval, $D_j\in\bC^{1}\big(\bR\big),\;D(0)=0,$  and such that
\begin{eqnarray*}
\left\|G(\mh,\cdot)-G(\overline{\mh},\cdot)\right\|_\infty&\leq&
\Big\{G_\infty\Big(\mh^{(k)}\Big)\vee G_\infty\Big(\overline{\mh}^{(k)}\Big)\Big\}D_0\Big\{\varrho_n^{(k)}\Big(\mh^{(k)},\overline{\mh}^{(k)}\Big)\Big\}
\\
&+&\sum_{j=k+1}^{m}L_j\Big\{G_\infty\Big(\mh^{(k)}\Big)\vee G_\infty\Big(\overline{\mh}^{(k)}\Big)\Big\}D_j
\Big(\varrho_j\big(\mh_j,\mh^\prime_j\big)\Big),
\end{eqnarray*}
for any $\mh,\mh^\prime\in\mH_{1}^k(n)\times\mH_{k+1}^m$ and $n\geq 1$.
\end{itemize}
\end{assumption}
\noindent We remark that Assumption \ref{ass:bounded_case} ($\mathbf{i}$) is automatically fulfilled  if $k=1$.

\begin{remark}
\label{rem:after-ass-bounded-case}
If $n\geq 1$ is fixed or $\mH_j(n),\;j=\overline{1,k},$ are independent on $n$, for example $\mH_j(n)=\mH_j,\;j=\overline{1,k},$ for all $n\geq 1$
then  upper functions for $|\xi_h(n)|$ can be derived under Assumption \ref{ass:bounded_case}.
%and Assumption \ref{ass:sec:totaly_bounded_case} (or Assumptions \ref{ass:partially-bounded-case} and \ref{ass:as-assumption-totally-bounded})  below
 However, if we are interested in finding of upper functions  for $|\xi_h(n)|$ when $n$ is varying,  we cannot do it in general without specifying the dependence of $\mH_j(n),\;j=\overline{1,k},$ on
$n$.

\end{remark}

In view of latter remark we will seek upper functions for $|\xi_h(n)|$ when $\mh\in\widetilde{\mH}(n):=\widetilde{\mH}_{1}^k(n)\times\mH_{k+1}^m$.
Here $\widetilde{\mH}_{1}^k(n)=\widetilde{\mH}_1(n)\times\cdots\widetilde{\mH}_k(n)$ and  $\left\{\widetilde{\mH}_j(n)\subset\mH_j(n),\; n\geq 1\right\}, j=\overline{1,k},$ be a sequence of sets satisfying additional restriction. We will not be tending here to the maximal generality and  complete Assumption \ref{ass:bounded_case} by the following condition.
%For any $\mathbf{m}\geq 1$ set $\mathbf{N_m}=$.

\begin{assumption}
\label{ass:dependence-on-n}

For any $\mathbf{m}\in\bN^*$ there exists $n[\mathbf{m}]\in \{\mathbf{m},\mathbf{m}+1,\ldots,2\mathbf{m}\}$ such that
\begin{equation*}
%\label{eq:order}
\bigcup_{n\in \{\mathbf{m},\mathbf{m}+1,\ldots,2\mathbf{m}\}}\widetilde{\mH}_1^k(n)\subseteq\mH_1^k\big(n[\mathbf{m}]\big).
%\quad \forall j=\overline{1,k}.
\end{equation*}
\end{assumption}

We note that Assumption \ref{ass:dependence-on-n} obviously  holds if  for any $j=\overline{1,k}$ the sequence $\left\{\widetilde{\mH}_j(n),\;n\geq 1\right\}$ is increasing/decreasing sequence of sets.

\subsection{Totally bounded case}
\label{sec:totally-bounded-case}

The objective is to find upper functions  for $|\xi_\mh(n)|$ under Assumption \ref{ass:bounded_case} enforced, if necessary, by Assumption \ref{ass:dependence-on-n} and the condition imposed on the entropies of the sets $\mH_j,\;j=\overline{k+1,m}$.

\subsubsection{Assumptions and main result}

 The following condition will be additionally imposed  in this section.
\begin{assumption}
\label{ass:sec:totaly_bounded_case}
Suppose that (\ref{eq:mH-space-decomposition}) holds and there exist $N,R<\infty$ such that  for any $\varsigma>0$  and any  $j=\overline{k+1,m}$
\begin{eqnarray*}
%\label{eq:entropy-sec-bounded-case}
\mathfrak{E}_{\mH_j,\varrho_j}(\varsigma)\leq N\left(\left[\log_2{\big\{R/\varsigma\big\}}\right]_++1\right),
\end{eqnarray*}
where, as previously, $\mathfrak{E}_{\mH_j,\varrho_j}$ denotes the entropy of $\mH_j$ measured in $\varrho_j$.
\end{assumption}

We remark that Assumption \ref{ass:sec:totaly_bounded_case}  is fulfilled, in particular, when
$\big(\mH_j,\varrho_j,\kappa_j\big),\;j=\overline{k+1,m},$ are bounded and
satisfy  doubling condition. Note also that this assumption can be considerably weakened, see discussion after Theorem \ref{th:empiric_totaly_bounded_case}.

\paragraph{Notations}

Let  $3\leq\mathbf{n_1}\leq\mathbf{n_2}<2\mathbf{n_1}$ be fixed and set $\widetilde{\mathbf{N}}=\{\mathbf{n_1},\ldots,\mathbf{n_2}\}$.
For any $\mh\in\boldsymbol{\mH}$ set
$$
{F}_{\mathbf{n_2}}(\mh)=\left\{\begin{array}{ll}
\sup_{i=\overline{1,\mathbf{n_2}}}\bE_\rf\big|G(\mh,X_i)\big|,\quad& \mathbf{n_1}\neq \mathbf{n_2};
\\*[2mm]
(\mathbf{n_2})^{-1}\sum_{i=1}^{\mathbf{n_2}}\bE_\rf\big|G(\mh,X_i)\big|,\quad& \mathbf{n_1}=\mathbf{n_2},
\end{array}
\right.
$$
and  remark that if additionally
$X_i,\;i\geq 1,$ are identically distributed then we have the same definition of ${F}_{\mathbf{n_2}}(\cdot)$ in both cases.
We note that
$$
\displaystyle{{F_{\mathbf{n_2}}}:=\sup_{n\in\widetilde{\mathbf{N}}}
\sup_{\mh\in\widetilde{\mH}(n)}{F_{\mathbf{n_2}}}(\mh)\leq\sup_{n\in\widetilde{\mathbf{N}}}\overline{G}_n<\infty}
$$
in view of Assumption \ref{ass:bounded_case} ($\mathbf{i}$).
\;Let $\boldsymbol{b}>1$  be fixed and put
$$
\mathbf{n}=\left\{
\begin{array}{ll}
\mathbf{n_1},\;& \mathbf{n_1}=\mathbf{n_2};
\\
n[\mathbf{n_1}],\;& \mathbf{n_1}\neq\mathbf{n_2},
\end{array}
\right.
\qquad
\beta=\left\{
\begin{array}{ll}
0,\;& \mathbf{n_1}=\mathbf{n_2};
\\
\boldsymbol{b},\;& \mathbf{n_1}\neq\mathbf{n_2},
\end{array}
\right.
$$
where, remind, that $n[\cdot]$ is defined in Assumption \ref{ass:dependence-on-n}.

Define $\widehat{L}_j(z)=\max\big\{z^{-1}L_j(z),1\big\}$ and  $\cL^{(k)}(z)=\sum_{j=k+1}^{m}\log_2{\left\{\widehat{L}_j\left(2z\right)\right\}}$
 and introduce the following quantities: for any $\mh^{(k)}\in\mH_1^k$ and any $q>0$
\begin{gather*}
P\left(\mh^{(k)}\right)=(36k\delta^{-2}_*+6)\ln{\left(1+\ln{\left\{2\underline{G}^{-1}_{\mathbf{n}} G_\infty\big(\mh^{(k)}\big)\right\}}\right)}
+36N\delta^{-2}_*\cL^{(k)}\Big( G_\infty\big(\mh^{(k)}\big)\Big)+18C_{N,R,m,k};
\\*[2mm]
M_q\left(\mh^{(k)}\right)=\big(72k\delta^{-2}_*+2.5q+1.5\big)\ln{\left(2\underline{G}^{-1}_{\mathbf{n}} G_\infty\big(\mh^{(k)}\big)\right)}
+72N\delta^{-2}_*\cL^{(k)}\Big( G_\infty\big(\mh^{(k)}\big)\Big)+36C_{N,R,m,k}.
\end{gather*}
Here $\delta_*$ it is  the smallest solution of the equation  $(48\delta)^{-1}s^{*}(\delta)=1$, where, remind, $s^*(\delta)=(6/\pi^2)\big(1+[\ln{\delta}]^{2}\big)^{-1},\;\delta\geq 0$.
%It is obvious that $\delta_*>0$.
The quantities $N, R$ are defined in Assumption \ref{ass:sec:totaly_bounded_case}.

The \textsf{explicit} expression of the constant $C_{N,R,m,k}$, as well as \textsf{explicit} expressions of the constants $\lambda_1,\;\lambda_2$
  and $C_{D,\boldsymbol{b}}$ used in the description of the results below, are given in Section \ref{sec:constants-th-totally-bounded} which precedes the proof of Theorem \ref{th:empiric_totaly_bounded_case}.

\paragraph{Result}
For any $\mathbf{r}\in \overline{\bN}$ put $F_{\mathbf{n_2},\mathbf{r}}(\mh)=\max\left[F_{\mathbf{n_2}}(\mh),e^{-\mathbf{r}}\right]$ and
define  for any $\mh\in\boldsymbol{\mH}$, $u\geq 0$ and $q>0$
\begin{eqnarray*}
\cV_{\mathbf{r}}^{(u)}(n,\mh)&=&\lambda_1\sqrt{G_\infty\big(\mh^{(k)}\big)\Big(F_{\mathbf{n_2},\mathbf{r}}(\mh)n^{-1}\Big) \Big(P\big(\mh^{(k)}\big)+
2\ln{\left\{1+\left|\ln{\left(F_{\mathbf{n_2},\mathbf{r}}(\mh)\right)}\right|\right\}}
+u\Big)}
\\
&&\hskip-0.3cm +
\lambda_2G_\infty\big(\mh^{(k)}\big)\Big(n^{-1}\ln^{\beta}{(n)}\Big)
\Big(P\big(\mh^{(k)}\big)+2\ln{\left\{1+\left|\ln{\left(F_{\mathbf{n_2},\mathbf{r}}(\mh)\right)}\right|\right\}}+u\Big);
%\label{eq:def-of-upper-function-moments}
\\*[2mm]
\cU_{\mathbf{r}}^{(u,q)}(n,\mh)&=&\lambda_1\sqrt{G_\infty\big(\mh^{(k)}\big)\Big(F_{\mathbf{n_2},\mathbf{r}}(\mh)n^{-1}\Big)
\Big(M_q\big(\mh^{(k)}\big)+2\ln{\left\{1+\left|\ln{\left(F_{\mathbf{n_2},\mathbf{r}}(\mh)\right)}\right|\right\}}
+u\Big)}
\\
&&\hskip-0.3cm+
\lambda_2G_\infty\big(\mh^{(k)}\big)\Big(n^{-1}\ln^{\beta}{(n)}\Big)
\Big(M_q\big(\mh^{(k)}\big)+2\ln{\left\{1+\left|\ln{\left(F_{\mathbf{n_2},\mathbf{r}}(\mh)\right)}\right|\right\}}+u\Big).
\end{eqnarray*}

\begin{theorem}
\label{th:empiric_totaly_bounded_case}
Let Assumptions \ref{ass:bounded_case} and \ref{ass:sec:totaly_bounded_case} be fulfilled. If $\mathbf{n_1}\neq \mathbf{n_2}$ suppose additionally that Assumption \ref{ass:dependence-on-n} holds.
Then for any $\mathbf{r}\in\bN$, $\boldsymbol{b}>1$ $ u\geq 1$  and $q\geq 1$
\begin{gather*}
\bP_\rf\left\{\sup_{n\in\widetilde{\mathbf{N}}}\sup_{\mh\in\widetilde{\mH}(n)}\Big[\big|\xi_\mh(n)\big|-\cV_{\mathbf{r}}^{(u)}(n,\mh)\Big]\geq 0\right\}
\leq 2419\; e^{-u}
;
\\
%*[2mm]
\bE_\rf\left\{\sup_{n\in\widetilde{\mathbf{N}}}\sup_{\mh\in\widetilde{\mH}(n)}
\Big[\big|\xi_\mh(n)\big|-\cU^{(u,q)}_{\mathbf{r}}(n,\mh)\Big]\right\}^q_+
%\\*[2mm]&&
\leq c_q\left[\sqrt{(\mathbf{n_1})^{-1}F_{\mathbf{n_2}}\underline{G}_{\mathbf{n}} }
\vee\left( (\mathbf{n_1})^{-1}\ln^{\beta}{(\mathbf{n_2})}\underline{G}_{\mathbf{n}} \right) \right]^q e^{-u},
\end{gather*}
where
$c_q= 2^{(7q/2)+5}3^{q+4}\Gamma(q+1)(C_{D,\boldsymbol{b}})^{q}$.
\end{theorem}

\begin{remark}
\label{rem:after-theorem3}
The inspection of  the proof of the theorem allows us to  assert that Assumption \ref{ass:sec:totaly_bounded_case} can be weakened.
The condition that is needed in view of the used technique: for some  $\alpha\in (0,1),\;L<\infty$
\begin{eqnarray}
\label{eq:rem-after-theorem3}
\sup_{\varsigma>0}\varsigma^{-\alpha}\mathfrak{E}_{\mH_j,\varrho_j}(\varsigma)\leq L,\;\; j=\overline{k+1,m}.
\end{eqnarray}
In particular, it allows to consider the generalized empirical processes indexed by the sets of smooth functions.
However the latter assumption does not permit to express upper functions explicitly as it is done in Theorem \ref{th:empiric_totaly_bounded_case}.
This explains why we prefer to state our results under Assumption \ref{ass:sec:totaly_bounded_case}.

\end{remark}

Several  other remarks are in order.

$1^0.$  First we note that the results presented in the theorem are obtained without any assumption imposed on the densities
 $f_i,\; i\geq 1$. In particular,  found upper functions remain finite even if the densities $f_i,\; i\geq 1$ are unbounded.

$2^0.$ Next, putting $\mathbf{r}=+\infty$ we get the results of the theorem with $F_{\mathbf{n_2},\mathbf{r}}(\cdot)=F_{\mathbf{n_2}}(\cdot)$.
It improves the first terms in the expressions of $\cV_{\mathbf{r}}^{(u)}(\cdot,\cdot)$ and $\cU^{(u,q)}_{\mathbf{r}}(\cdot,\cdot)$,
however the second terms may explode if
$
F_{\mathbf{n_2}}(\mh)=0
$
for some $\mh\in\mathbf{\mH}$.
The latter fact explains the necessity to "truncate" $F_{\mathbf{n_2}}(\cdot)$ from below, i.e. to consider $F_{\mathbf{n_2},\mathbf{r}}(\cdot)$ instead of $F_{\mathbf{n_2}}(\cdot)$.

\subsubsection{Law of iterated logarithm}

Our goal here is to use the first assertion of Theorem \ref{th:empiric_totaly_bounded_case} in order to establish a non-asymptotical version of the law of iterated logarithm for
$$
 \eta_{\mh^{(k)}}(n):=\sup_{\mh_{(k)}\in\mH_{k+1}^m}\big|\xi_{\mh}(n)\big|.
$$
Let us suppose that for some $\mathfrak{c}>0, \;\mb> 0$
\begin{equation}
\label{eq1:LIL}
\mathfrak{c}\leq \underline{G}_n\leq
\overline{G}_n\leq \mathfrak{c} n^{\mb},\quad \forall n \geq 1.
\end{equation}
We would like to emphasize that the restriction $ \underline{G}_n\geq \mathfrak{c}$ is imposed for the simplicity of the notations and the results presented below are valid if $\underline{G}_n$ decreases to zero polynomially in $n$.

Moreover we will assume that
\begin{equation}
\label{eq2:LIL}
\sup_{n\geq 1}\sup_{\mh\in\widetilde{\mH}(n)}\sup_{i\geq 1}\bE_\rf\big|G(\mh,X_i)\big|=:\mathbf{F}<\infty.
\end{equation}
We will see that the latter condition is checked in various particular problems if the densities $f_i,\; i\geq 1$ are uniformly bounded.
Suppose finally that for some $\ma>0$
\begin{equation}
\label{eq22:LIL}
\cL^{(k)}(z)\leq \ma\ln{\big\{1+\ln{(z)}\big\}},\quad\forall z\geq 3.
\end{equation}
For any $a>0$ and $n\geq 3$ define
$$
\overline{\mH}_1^k(n,a)=\widetilde{\mH}_1^k(n)\cap\left\{\mh^{(k)}: \;\;G_\infty\big(\mh^{(k)}\big)\leq n\big[\ln(n)\big]^{-a}\right\}.
$$

\begin{theorem}
\label{th:LIL-nonasym}
Let Assumptions \ref{ass:bounded_case}, \ref{ass:dependence-on-n} and  \ref{ass:sec:totaly_bounded_case}  be fulfilled  and suppose additionally that (\ref{eq1:LIL}), (\ref{eq2:LIL}) and (\ref{eq22:LIL}) hold.
Then there exists $\Upsilon>0$ such that for any $\mathbf{j}\geq 3$ and any $a>2$
\begin{gather*}
\bP_\rf\left\{\sup_{n\geq \mathbf{j}}\;\sup_{\mh^{(k)}\in\overline{\mH}_1^k(n,a)}\Bigg[\frac{\sqrt{n}\;\eta_{\mh^{(k)}}(n)}
{\sqrt{G_\infty\big(\mh^{(k)}\big)\ln{\big(1+\ln{(n)}\big)}}}\Bigg]\geq \Upsilon\right\}
\leq  \frac{2419}{\ln(\mathbf{j})}.
\end{gather*}

\end{theorem}
The explicit expression of the constant $\Upsilon$ can be easily derived but it is quite cumbersome and we omit its derivation.
\begin{remark}
\label{rem:after-th-LIL}
The inspection of the proof of the theorem shows that for any $y\geq 0$ one can find $0<\Upsilon(y)<\infty$ such that the assertion of the theorem remains true if one replaces $\Upsilon$ by   $\Upsilon(y)$ and the right hand side of the obtained inequality by
$2419\big[\ln(\mathbf{j})\big]^{-(1+y)}$. It makes reasonable the consideration of small values of $\mathbf{j}$.

\end{remark}

The simple corollary of Theorem \ref{th:LIL-nonasym}  is the law of iterated logarithm:
\begin{equation}
\label{eq3:LIL}
\limsup_{n\to\infty}\sup_{\mh^{(k)}\in\overline{\mH}_1^k(n,a)}\Bigg[\frac{\sqrt{n}\;\eta_{\mh^{(k)}}(n)}
{\sqrt{G_\infty\big(\mh^{(k)}\big)\ln{\ln{(n)}}}}\Bigg]\leq \Upsilon, \quad \bP_\rf-\text{a.s.}
\end{equation}

\subsection{Partially totally bounded  case}
\label{sec:partially_totally-bounded-case}

We begin this section with the following definition used in the sequel. Let $\bT$ be a set equipped with  a semi-metric $\md$ and let
$\mathfrak{n}\in\bN^*$ be fixed.

\begin{definition}
\label{def:locally-finite-space}
We say that $\left\{\bT_{\mathbf{i}}\subset\bT,\; \mathbf{i}\in \mathbf{I}\right\}$ is $\mathfrak{n}$-totally bounded cover of $\bT$ if
\begin{itemize}
\item
\quad $\displaystyle{\bT=\cup_{\mathbf{i}\in \mathbf{I}}\bT_{\mathbf{i}}}$ and $\mathbf{I}$ is countable;

\smallskip

\item
\quad $\bT_{\mathbf{i}}$ is totally bounded  for any $\mathbf{i}\in\mathbf{I}$;

\smallskip

\item
\quad
$
 \text{card}\Big(\left\{\mathbf{k}\in\mathbf{I}:\;\; \bT_{\mathbf{i}}\cap\bT_{\mathbf{k}}\neq \emptyset\right\}\Big)\leq \mn
$
for any $\mathbf{i}\in\mathbf{I}$.
\end{itemize}
\end{definition}
Let us illustrate the above definition by some examples.

Let $\bT=\bR^d,\;d\geq 1$. Then any countable partition of $\bR^d$ consisted of bounded sets forms $1$-totally bounded cover of $\bR^d$.
Note, however, that the partitions will not be suitable choice for particular problems studied later. We will be mostly interested in $\mathfrak{n}$-totally bounded covers satisfying the following \textsf{separation property}: there exists $\mathfrak{r}>0$ such that  for all $\mathbf{i},\mathbf{k}\in\mathbf{I}$ satisfying $\bT_{\mathrm{i}}\cap\bT_{\mathbf{k}}=\emptyset$
\begin{equation}
\label{eq:separation-property}
\inf_{t_1\in\bT_{\mathrm{i}},t_2\in\bT_{\mathrm{k}}}\md(t_1,t_2)\geq\mathfrak{r}.
\end{equation}
Let us return to $\bR^d$ that we equip with the metric  generated by the supremum norm. Denote by $\bB_r(t),\;t\in\bR^d, r>0,$
the  closed ball in this metric
 with the radius $r$ and the center $t$. For given $\mathfrak{r}>0$  consider the collection $\left\{\bB_{\frac{\mathfrak{r}}{2}}(\mathfrak{r}\mathbf{i}),\; \mathbf{i}\in \bZ^d\right\}$, where we understand $\mathfrak{r}\mathbf{i}$ as coordinate-wise multiplication. It is easy to check that this collection is
  $3^d$-totally bounded cover of $\bR^d$ satisfying (\ref{eq:separation-property}).

 \par

We would  like to emphasize that $\mathfrak{n}$-totally bounded covers satisfying the separation property can be often constructed when
 $\bT$ is a homogenous  metric space  endowed with the Borel measure obeying doubling condition. Some useful results for this construction can be found in the recent paper \cite{Kerk}, where such spaces were scrutinized.

We finish the discussion about $\mathfrak{n}$-totally bounded covers with the following notation: for any $t\in\bT$ put
$$
\bT(t)=\bigcup_{
\mathbf{i}\in\mathbf{I}:\;t\in\bT_{\mathbf{i}}
}\;\bigcup_{
\mathbf{k}\in\mathbf{I}:\;\bT_{\mathbf{i}}\cap\bT_{\mathbf{k}}\neq\emptyset}\bT_{\mathbf{k}}.
$$

\subsubsection{Assumptions and main result}

 Throughout this section we will assume that the representation (\ref{eq:cX-space-decomposition}) holds and the elements of $\cX_l,\; l=\overline{1,p},$  will be denoted by $x_l$ .
We keep all notations from previous section and  replace Assumption \ref{ass:sec:totaly_bounded_case}  by the following  conditions.

 \begin{assumption}
\label{ass:partially-bounded-case}
\begin{itemize}
\item[$\mathbf{(i)}$]
Let (\ref{eq:mH-space-decomposition}) and (\ref{eq:cX-space-decomposition}) hold with $\cX_1=\mH_m$ and for some $\mathfrak{n}\in\bN^*$ there exists
a collection $\Big\{\mathrm{H}_{m,\mathbf{i}},\;\mathbf{i}\in \mathbf{I}\Big\}$ being the $\mn$-totally bounded cover of $\mH_m$ satisfying
for some
$N,R<\infty$
%$\forall $
 %any $$
%and any  $j=\overline{k+1,m-1}$
\begin{eqnarray*}
%\label{eq:entropy-sec-bounded-case}
&&\mathfrak{E}_{\mathrm{H}_{m,\mathbf{i}},\varrho_m}(\varsigma)\leq N\left(\left[\log_2{\big\{R/\varsigma\big\}}\right]_++1\right),\quad\forall \mathbf{i}\in\mathbf{I},\;\;\forall \varsigma>0.
\end{eqnarray*}

\item[$\mathbf{(ii)}$]
For
%$\forall $
 any $\varsigma>0$
%and any  $j=\overline{k+1,m-1}$
\begin{eqnarray*}
&&\mathfrak{E}_{\mH_j,\varrho_j}(\varsigma)\leq N\left(\left[\log_2{\big\{R/\varsigma\big\}}\right]_++1\right),\quad \forall j=\overline{k+1,m-1}.
\end{eqnarray*}
\end{itemize}
\end{assumption}
Usually one can construct  many $\mn$-totally bounded covers satisfying Assumption \ref{ass:partially-bounded-case} $(\mathbf{i})$.
The condition below restricts this choice and relates it to properties of the mapping  $G(\cdot,\cdot)$ describing generalized empirical process.

\begin{assumption}
\label{ass:as-assumption-totally-bounded}
 For any $n\geq 1$ and any $\mh=(\mh_1,\ldots,\mh_m)\in\mH(n)$
$$
%\qquad\sup_{\substack{\mathbf{i},\mathbf{k}\in\mathbf{I}:\\
%\mathrm{H}_{m,\mathbf{i}}\cap\mathrm{H}_{m,\mathbf{k}}=\emptyset}}
%.
\sup_{x\in\cX:\; x_1\notin \mH_m(\mh_m)}|G(\mh,x)|\leq n^{-1}G_\infty\big(\mh^{(k)}\big).
$$

\end{assumption}
We would like to emphasize that in order to satisfy   Assumption \ref{ass:as-assumption-totally-bounded} in particular examples, the
$\mn$-totally bounded cover $\Big\{\mathrm{H}_{m,\mathbf{i}},\;\mathbf{i}\in \mathbf{I}\Big\}$ should usually possess the separation property. Indeed, one of the typical examples, where Assumption \ref{ass:as-assumption-totally-bounded} is fulfilled, is the following: there exist $\gamma>0$ such that for $G(x,\mh)=0$ for any $x\in\cX,\;\mh\in\mH,$ satisfying $\rho_m(x_1,\mh_m)\geq \gamma.$

\paragraph{Result}  For any $i=\overline{1,n}$ we denote $X_i=\big(X_{1,i},\ldots,X_{p,i}\big)$,
$$
f_{1,i}(x_1)=\int_{\cX_2\times\cdots\times\cX_p} f_i(x_1,\ldots,x_p)\prod_{l=2}^{p}\nu_l\big(\rd x_l\big).
$$
and if  $\cX=\cX_1\; (p=1)$ then we put $X_{1,i}=X_i$ and  $f_{1,i}=f_{i}$.

Put for any $n\geq 1$, $v> 0$ and any $\mh_m\in\mH_m$
\begin{equation*}
\mL_{n,v}(\mh_m)=-\ln{\bigg(\bigg[n^{-1}\sum_{i=1}^{n}\int_{\mH_m(\mh_m)}f_{1,i}(x)\nu_1\big(\rd x\big)\bigg]\vee n^{-v}\bigg)}.
\end{equation*}
Note that obviously
$
0\leq \mL_{n,v}(\mh_m)\leq v\ln{(n)},\;\forall\mh_m\in\mH_m.
$
Put for any $\mh\in\boldsymbol{\mH}$
\begin{eqnarray*}
\widetilde{P}(\mh)=P\big(\mh^{(k)}\big)+\mL_{n,v}\big(\mh_m\big)+
2\ln{\left\{1+\left|\ln{\left(F_{\mathbf{n_2},\mathbf{r}}(\mh)\right)}\right|\right\}};
\\*[2mm]
\widetilde{M}_q(\mh)=M_q\big(\mh^{(k)}\big)+\mL_{n,v}\big(\mh_m\big)+
2\ln{\left\{1+\left|\ln{\left(F_{\mathbf{n_2},\mathbf{r}}(\mh)\right)}\right|\right\}}.
\end{eqnarray*}
Define for any $\mh\in\boldsymbol{\mH}$, $\mathbf{r}\in \overline{\bN}$, $z\geq 0$ and $q>0$
\begin{gather*}
\widetilde{\cV}_{\mathbf{r}}^{(v,z)}(n,\mh)=\lambda_1\sqrt{G_\infty\big(\mh^{(k)}\big)\Big(F_{\mathbf{n_2},\mathbf{r}}(\mh)n^{-1}\Big) \big(\widetilde{P}\big(\mh\big)+z\big)}+
%\\
%&&\hskip-0.3cm +
\lambda_2G_\infty\big(\mh^{(k)}\big)\Big(n^{-1}\ln^{\beta}{(n)}\Big)
\big(\widetilde{P}\big(\mh\big)+z\big);
%\label{eq:def-of-upper-function-moments}
\\*[2mm]
\widetilde{\cU}_{\mathbf{r}}^{(v,z,q)}(n,\mh)=\lambda_1\sqrt{G_\infty\big(\mh^{(k)}\big)\Big(F_{\mathbf{n_2},\mathbf{r}}(\mh)n^{-1}\Big)
\big(\widetilde{M}_q(\mh) +z\big)}
+
\lambda_2G_\infty\big(\mh^{(k)}\big)\Big(n^{-1}\ln^{\beta}{(n)}\Big)
\big(\widetilde{M}_q(\mh) +z\big).
\end{gather*}
\begin{theorem}
\label{th:empiric-partially_totaly_bounded_case}
Let Assumptions \ref{ass:bounded_case}, \ref{ass:partially-bounded-case} and \ref{ass:as-assumption-totally-bounded} hold. If $\mathbf{n_1}\neq \mathbf{n_2}$ suppose additionally that Assumption \ref{ass:dependence-on-n} holds.
Then for any $\mathbf{r}\in\bN$, $v\geq 1$, $ z\geq 1$  and $q\geq 1$
\begin{eqnarray*}
&&\hskip-0.6cm\bP_\rf\left\{\sup_{n\in\widetilde{\mathbf{N}}}\sup_{\mh\in\widetilde{\mH}(n)}
\Big[\big|\xi_\mh(n)\big|-\widetilde{\cV}_{\mathbf{r}}^{(v,z)}(n,\mh)\Big]\geq 0\right\}
\leq  \mn^5\Big\{4838e^{-z}+2\mathbf{n_1}^{2-v}\Big\};
\\*[0mm]
&&\hskip-0.6cm\bE_\rf\left\{\sup_{n\in\widetilde{\mathbf{N}}}\sup_{\mh\in\widetilde{\mH}(n)}
\Big[\big|\xi_\mh(n)\big|-\widetilde{\cU}^{(v,z,q)}_{\mathbf{r}}(n,\mh)\Big]\right\}^q_+
%\\*[2mm]&&
\leq 2\mn^5 c_q\left[\sqrt{(\mathbf{n_1})^{-1}F_{\mathbf{n_2}}\underline{G}_\mathbf{n} }
\vee\left( (\mathbf{n_1})^{-1}\ln^{\beta}{(\mathbf{n_2})}\underline{G}_\mathbf{n} \right) \right]^qe^{-z}
\\
&&\hskip6.9cm+
2^{q+1}\mn^5(\overline{G}_\mathbf{n})^{q}\;\mathbf{n_1}^{2-v}.
\end{eqnarray*}

\end{theorem}
Although  the assertions of the theorem are true whenever $v\geq 1$ the presented results are obviously reasonable
only if $v>2$. For example (as we will see later) the typical choice of this parameter for the "moment bound" is $v=q+2$.

\smallskip

In spite of the fact that upper functions presented in Theorem \ref{th:empiric-partially_totaly_bounded_case}
are found explicitly their expressions are quite cumbersome. In particular, it is unclear how to compute the function
$\mL_{n,v}(\cdot)$. Of course, since  $\mL_{n,v}(\mh_m)\leq v\ln{(n)},\;\forall\mh_m\in\mH_m,$ one can replace it by
$v\ln{(n)}$ in the definition of $\widetilde{P}(\cdot)$ and $\widetilde{M}_q(\cdot)$, but the corresponding upper functions are
not always sufficiently tight.

Our goal now is to simplify the expressions for upper functions given  in Theorem \ref{th:empiric-partially_totaly_bounded_case}. Surprisingly, that if $n$ is fixed, i.e. $\mathbf{n_1}=\mathbf{n_2}$, it can be done without any additional assumption.

Set for any $v>0$ and $\mh\in\boldsymbol{\mH}$
$$
\widehat{P}_v\big(\mh^{(k)}\big)=P\big(\mh^{(k)}\big)+
2v\left|\ln{\left(2G_\infty\big(\mh^{(k)}\big)\right)}\right|,\quad \widehat{M}_{q,v}\big(\mh^{(k)}\big)=M_q\big(\mh^{(k)}\big)+
2v\left|\ln{\left(2G_\infty\big(\mh^{(k)}\big)\right)}\right|,
$$
and let $\widehat{{F}}_{\mathbf{n_2}}(\mh)=\max[{F}_{\mathbf{n_2}}(\mh),\mathbf{n_2}^{-1}]$.

\begin{corollary}
\label{cor:after-th:empiric-partially_totaly_bounded_case}
Let the assumptions of Theorem \ref{th:empiric-partially_totaly_bounded_case} hold. If $\mathbf{n_1}\neq\mathbf{n_2}$ suppose additionally
that $X_{i,1},\;i\geq 1,$ are identically distributed.

\noindent Then, the results of Theorem \ref{th:empiric-partially_totaly_bounded_case} remain valid if one replaces $\widetilde{\cV}_{\mathbf{r}}^{(v,z)}(n,\mh)$ and $\widetilde{\cU}^{(v,z,q)}_{\mathbf{r}}(n,\mh)$  by
\begin{eqnarray*}
\widehat{\cV}^{(v,z)}(n,\mh)&=&\lambda_1\sqrt{G_\infty\big(\mh^{(k)}\big)\Big(\widehat{{F}}_{\mathbf{n_2}}(\mh)n^{-1}\Big) \Big(\widehat{P}_v\big(\mh^{(k)}\big)+2(v+1)\big|\ln{\big\{\widehat{{F}}_{\mathbf{n_2}}(\mh)\big\}}\big|
+z\Big)}
\\
&&\hskip-0.3cm +
\lambda_2G_\infty\big(\mh^{(k)}\big)\Big(n^{-1}\ln^{\beta}{(n)}\Big)
\Big(\widehat{P}_v\big(\mh^{(k)}\big)+2(v+1)\big|\ln{\big\{\widehat{{F}}_{\mathbf{n_2}}(\mh)\big\}}\big|+z\Big);
%\label{eq:def-of-upper-function-moments}
\\*[2mm]
\widehat{\cU}^{(v,z,q)}(n,\mh)&=&\lambda_1\sqrt{ G_\infty\big(\mh^{(k)}\big)\Big(\widehat{{F}}_{\mathbf{n_2}}(\mh)n^{-1}\Big)
\Big(\widehat{M}_{q,v}\big(\mh^{(k)}\big)+2(v+1)\big|\ln{\big\{\widehat{{F}}_{\mathbf{n_2}}(\mh)\big\}}\big|
+z\Big)}
\\
&&\hskip-0.3cm+
\lambda_2G_\infty\big(\mh^{(k)}\big)\Big(n^{-1}\ln^{\beta}{(n)}\Big)
\Big(\widehat{M}_{q,v}\big(\mh^{(k)}\big)+2(v+1)\big|\ln{\big\{\widehat{{F}}_{\mathbf{n_2}}(\mh)\big\}}\big|
+z\Big).
\end{eqnarray*}

\end{corollary}

We would like to emphasize that we do not require that $X_i,\;i\geq 1,$ would be identically distributed. In particular, coming back to the generalized empirical process considered in Example 2, Section \ref{sec:subsec-prob-form}, where $X_i=(Y_i,\e_i)$,  the design points $Y_i, \; i\geq 1,$ are often supposed to be uniformly  distributed on some bounded domain of $\bR^d$.
As to the noise variables $\varepsilon_i, \; i\geq 1,$ the restriction that they are identically distributed cannot be justified in general.

\subsubsection{Law of  logarithm}

Our goal here is to use the first assertion of Corollary \ref{cor:after-th:empiric-partially_totaly_bounded_case} in order to establish the result referred later to   {\it the law of  logarithm}. Namely we show that for some $\boldsymbol{\Upsilon}>0$
\begin{equation}
\label{eq1:LL}
\limsup_{n\to\infty}\sup_{\mh^{(k)}\in\overline{\mH}_1^k(n,a)}\frac{\sqrt{n}\;\eta_{\mh^{(k)}}(n)}
{\sqrt{G_\infty\big(\mh^{(k)}\big)\Big[\ln{\left\{G_\infty\big(\mh^{(k)}\big)\right\}}\vee\ln{\ln{(n)}}\Big]}}\leq \boldsymbol{\Upsilon} \quad \bP_\rf-\text{a.s.}
\end{equation}
As previously we will first provide with
the non-asymptotical version of (\ref{eq1:LL}).

We will suppose that (\ref{eq1:LIL}) and (\ref{eq2:LIL}) are fulfilled and replace (\ref{eq22:LIL}) by  the following assumption.
For some $\ma>0$
\begin{equation}
\label{eq2:LL}
\cL^{(k)}(z)\leq \ma\ln{(z)},\quad\forall z\geq 2.
\end{equation}

\begin{theorem}
\label{th:LL-nonasym}
Let Assumptions \ref{ass:bounded_case}, \ref{ass:dependence-on-n}, \ref{ass:partially-bounded-case} and \ref{ass:as-assumption-totally-bounded} be fulfilled. Suppose also that (\ref{eq1:LIL}), (\ref{eq2:LIL}) and (\ref{eq2:LL}) hold and assume that $X_{i,1},\;i\geq 1,$ are identically distributed.

Then  there exits $\boldsymbol{\Upsilon}>0$ such that for any $\mathbf{j}\geq 3$ and any $a>4$
\begin{gather*}
\bP_\rf\left\{\sup_{n\geq \mathbf{j}}\;\
\sup_{\mh^{(k)}\in\overline{\mH}_1^k(n,a)}\frac{\sqrt{n}\;\eta_{\mh^{(k)}}(n)}
{\sqrt{G_\infty\big(\mh^{(k)}\big)\Big[\ln{\left\{G_\infty\big(\mh^{(k)}\big)\right\}}\vee\ln{\ln{(n)}}\Big]}}\geq \boldsymbol{\Upsilon}
\right\}
\leq  \frac{4840\mn^5}{\ln{(\mathbf{j})}}.
\end{gather*}

\end{theorem}

Some remarks are in order. The explicit expression of the constant $\boldsymbol{\Upsilon}$ is available and  the generalization , similar to one announced in Remark \ref{rem:after-th-LIL}, is possible.  Also, (\ref{eq1:LL}) is an obvious consequence of Theorem \ref{th:LL-nonasym}.
At last, we note that in view of (\ref{eq1:LIL}) the factor  $\Big[\ln{\left\{G_\infty\big(\mh^{(k)}\big)\right\}}\vee\ln{\ln{(n)}}\Big]$
can be replaced by $\ln(n)$ which is, up to a constant, its upper estimate. The corresponding result is, of course, rougher than one presented in the theorem, but its derivation does not require $X_{i,1},\; i\geq 1,$ to be identically distributed. This result is deduced directly from
Theorem \ref{th:empiric-partially_totaly_bounded_case}. Its proof is almost the same as the proof of Theorem \ref{th:LL-nonasym} and based on the trivial  bound $\mL_{n,v}(\mh_m)\leq v\ln{(n)},\;\forall\mh_m\in\mH_m$.

%The explicit expression of the function  $\Upsilon(\cdot)$ can be easily derived but it is quite cumbersome and we omit its derivation.
%Analyzing  the proof of the theorem we can assert that $\Upsilon(\cdot)$ is continuous  on $\bR_+$ function. Hence, choosing $y=y(\mathbf{j})=\big[\ln(\mathbf{j})]^{-1/2}$ we come to (\ref{eq1:LL}) with $\boldsymbol{\Upsilon}=\Upsilon(0)$.

\subsection{Application to localized processes}
\label{sec:localized-empiric}

Let $\big(\bX_l,\mu_l,\rho_l\big),\; l=\overline{1,d+1},\;d\in\bN,$ be the collection of measurable metric spaces.
Throughout this section we will suppose that (\ref{eq:cX-space-decomposition}) holds with $p=2$,
$$
\cX=\cX_1\times\cX_2,\quad\big(\cX_1,\nu_1\big)=\big(\bX_1\times\cdots\times\bX_d,\mu_1\times\cdots\times\mu_d\big)=:\big(\bX_1^d,\mu^{(d)}\big),\quad
\big(\cX_2,\nu_2\big)=\big(\bX_{d+1},\mu_{d+1}\big),
$$
$x_j$ denotes the element of $\bX_j,\; j=\overline{1,d+1},$ and   $x^{(d)}$ will denotes the element of $\bX_1^d$. We equip  the space $\bX_1^d$ with the semi-metric $\rho^{(d)}=\max_{l=\overline{1,d}}\rho_l$.

%\smallskip

\paragraph{Problem formulation} This section is devoted to the application of Theorems \ref{th:empiric_totaly_bounded_case} and \ref{th:empiric-partially_totaly_bounded_case} in the following case:
\begin{itemize}
%\item
%$\cX_l,l=\overline{1,d}$, are equipped with semi-metric $\rho_l$;
%\smallskip

\item
$\mH_1^d:=\mH_{1}\times\cdots\times\mH_d=(0,1]\times\cdots\times(0,1]=(0,1]^d$, (i.e. $k=d$);

\smallskip

\item
$\mH_{d+1}^{d+2}=\mH_{d+1}\times\mH_{d+2}:=\cZ\times\bar{\bX}_1^d\;$, i.e. $m=d+2$, where $\bar{\bX}_1^d:=\bar{\bX}_1\times\cdots\times\bar{\bX}_d$ be a given subset of $\bX_1^d$ and $(\cZ,\md)$ is a given metric space.

\smallskip

\item
The function $G(\cdot,\cdot)$ obeys some structural assumption described below and for any $\mh:=\big(r,\mz,\bar{x}^{(d)}\big)\in(0,1]^d\times\cZ\times\bar{\bX}_1^d$ the function $G(\mh,\cdot)$ "decrease rapidly " outside of the set $\Big\{x_1\in\bX_1:\;\;\rho_1\big(x_1,\bar{x}_1\big)\leq r_1\Big\}\times\cdots\times\Big\{x_d\in\bX_d:\;\;\rho_d\big(x_d,\bar{x}_d\big)\leq r_d\Big\}\times\bX_{d+1}$.
% where $\bar{x}^{(d)}\in\bar{\cX}_1^d$.
\end{itemize}

Let  $K:\bR^d\to\bR$ be a given function, $\big(\gamma_1,\ldots,\gamma_d\big)\in\bR^d_+$ be given vector and set for any $r\in (0,1]^d$
$$
K_{r}(\cdot)=V_r^{-1}K\left(\cdot/r_1,\ldots, \cdot/r_d\right),\quad V_r=\prod_{l=1}^d r_l^{\gamma_l}.
$$
where, as previously, for  $u,v\in\bR^d$ the notation $u/v$ denotes the coordinate-wise division.
%$\rmm\in\bM_\cR(\gamma,L)$.
Let
\begin{equation}
\label{eq:structural-assumption-empiric-general}
G(\mh,x)=g\big(\mz,x\big)K_{r}\left(\vec{\rho}\big(x^{(d)},\bar{x}^{(d)}\big)\right),\quad \mh=\Big(r,\mz,\bar{x}^{(d)}\Big)\in(0,1]^d\times\cZ\times\bar{\bX}_1^d=:\boldsymbol{\mH},
\end{equation}
where $g:\cZ\times\cX\to \bR$ is a given function whose properties will be described later and
$$
\vec{\rho}\big(x^{(d)},\bar{x}^{(d)}\big)=\left(\rho_1\big(x_1,\bar{x}_1\big),\ldots,\rho_d\big(x_d,\bar{x}_d\big)\right).
$$
The corresponding generalized empirical process is given by
$$
\xi_{\mh}(n)=n^{-1}\sum_{i=1}^{n}\bigg[g\big(\mz,X_i\big)K_{r}\left(\vec{\rho}
\Big(\left[X_i\right]^{(d)},\bar{x}^{(d)}\Big)\right)-
\bE_{\rf}\left\{g\big(\mz,X_i\big)K_{r}\left(\vec{\rho}
\Big(\left[X_i\right]^{(d)},\bar{x}^{(d)}\Big)\right)\right\} \bigg].
$$
We will seek  upper functions for the random field
$
\displaystyle{\zeta_r\big(n,\bar{x}^{(d)}\big):=\sup_{x_{d+1}\in\bX_{d+1}}\big|\xi_{r,\mz,\bar{x}^{(d)}}(n)\big|}
$
%on $\cR(n)\times\bar{\bX}_1^d$
in two cases: $\bar{\bX}_1^d=\bX_1^d$ and $\bar{\bX}_1^d=\left\{\bar{x}^{(d)}\right\}$ for a fixed $\bar{x}^{(d)}\in\bX_1^d$.

 To realize this program we will apply Theorems \ref{th:empiric_totaly_bounded_case} and \ref{th:empiric-partially_totaly_bounded_case}
 to $\xi_{\mh}(n),\; \mh=\Big(r,\mz,\bar{x}^{(d)}\Big)$.
  It is worth mentioning that  corresponding upper functions can be used for constructing of estimation procedures in different areas of mathematical statistics: $M$-estimation with locally polynomial fitting (non-parametric regression),  kernel density estimation and  many others.

Moreover, we apply Theorem \ref{th:LIL-nonasym}  for establishing a non-asymptotical version of the law of iterated logarithm for $\zeta_r\big(\bar{x}^{(d)},n\big)$  in the case where $\bar{\bX}_1^d=\left\{\bar{x}^{(d)}\right\}$ .
We also apply  Theorem  \ref{th:LL-nonasym} for deriving a non-asymptotical version of the law of  logarithm for $\left\|\zeta_r(n)\right\|_\infty:=\sup_{\bar{x}^{(d)}\in\bX_1^d}\left|\zeta_r\big(\bar{x}^{(d)},n\big)\right|$. Our study here generalizes  in several directions the existing results \cite{mason-ein}, \cite{gine}, \cite{mason-ein2}, \cite{Dony1}, \cite{Dony2}.

\paragraph{Assumptions and notations}

\begin{assumption}
\label{ass:on-fucnctions-g-and-K}
$(\mathbf{i})$\;  $\|K\|_\infty<\infty$ and  for some  $L_1>0$
$$
\hskip0cm \left|K(t)-K(s)\right|\leq \frac{L_1|t-s|}{1+|t|\wedge|s|},\quad\forall t,s\in\bR^d,
%\sup_{u\in\bR^d}\bigg(\sup_{|t_1|\geq |u_1|,\ldots |t_d|\geq |u_d|}\;\sup_{l=\overline{1,d}}\left|\frac{\partial K(t)}{\partial t_l}\right|\bigg)\bigg(1+\sum_{l=1}^d |u_l|\bigg)\leq L_1.
$$
where $|\cdot|$ denotes supremum norm on $\bR^d$.

\smallskip

\qquad\qquad\qquad\quad\;\;$(\mathbf{ii})$\;\;$\|g\|_\infty:=\displaystyle{\sup_{\mz\in\cZ,\;x\in\cX}}\;g\big(\mz,x\big)<\infty$, and for some $\alpha\in (0,1]$, $L_\alpha>0$,
$$
\hskip0.9cm\sup_{x\in\cX}\left|g\big(\mz,x\big)-g\big(\mz^\prime,x\big)\right|\leq L_\alpha \left[\md(\mz,\mz^\prime)\right]^\alpha,\;\forall \mz,\mz^\prime\in\cZ;
$$

\end{assumption}

The conditions ($\mathbf{i}$) and ($\mathbf{ii}$) are quite standards. In particular ($\mathbf{i}$) holds if $K$ is compactly supported and lipschitz continuous.
If  $g\big(\mz,\cdot\big)=\bar{g}(\cdot)$, for any $\mz\in\cZ$, then ($\mathbf{ii}$) is verified for any bounded $\bar{g}$.

Let $0<r^{(\min)}_l(n)\leq r^{(\max)}_l(n)\leq 1, \;l=\overline{1,d},\; n\geq 1,$ be given decreasing sequences and let
\begin{eqnarray*}
&&\mH(n)=\cR(n)\times\cZ\times\bar{X}_1^d,\qquad\cR(n)=\prod_{l=1}^d \big[r^{(\min)}_l(2n),r^{(\max)}_l(n)\big];
\\
&&\widetilde{\mH}(n)=\widetilde{\cR}(n)\times\cZ\times\bar{X}_1^d,\qquad\widetilde{\cR}(n)=\prod_{l=1}^d \big[r^{(\min)}_l(n),r^{(\max)}_l(n)\big].
\end{eqnarray*}
We note that $\widetilde{\mH}(n)\subseteq\mH(n)$ for any $n\geq 1$ since $r^{(\min)}_l(\cdot),r^{(\max)}_l(\cdot), \;l=\overline{1,d},$ are decreasing, and obviously $\widetilde{\mH}(n)\subseteq \mH(\mathbf{m})$ for any $n\in \{\mathbf{m},\ldots,\mathbf{2m}\}$ and any $\mathbf{m}\geq 1$.

\begin{remark}
\label{rem:verification_ass-dependence-on-n}
Assumption \ref{ass:dependence-on-n} is  fulfilled with $n[\mathbf{m}]=\mathbf{m}$.
\end{remark}

%Later on we will always assume that $\lim_{n\to\infty}r^{(\min)}_j(n)=0,\quad j=\overline{1,d}.$

\begin{lemma}
\label{lem:loc-empiric-pointwise}
Suppose that Assumption \ref{ass:on-fucnctions-g-and-K} is fulfilled  and let $\bar{\bX}_1^d\subseteq\bX_1^d$ be an arbitrary subset.
Then, for  arbitrary sequences  $0<r^{(\min)}_l(n)\leq r^{(\max)}_l(n)\leq 1, \;l=\overline{1,d},\; n\geq 1,$  Assumption \ref{ass:bounded_case} holds with
\begin{gather*}
 \varrho^{(d)}_n\left(r,r^\prime\right)=\max_{l=\overline{1,d}}\big|\gamma_l\ln{\left(r_l/r^\prime_l\right)}\big|,
\quad \varrho_{d+1}=\big[\md\big]^\alpha,\quad \varrho_{d+2}=\max_{l=\overline{1,d}}\rho_l;
\\
D_0(z)=\exp{\{dz\}}-1+(L_1/\|K\|_\infty) \Big(\exp{\left\{\gamma^{-1}z\right\}}-1\Big),\quad \gamma=\min_{l=\overline{1,l}}\gamma_l;
\\
D_{d+1}(z)=\big(L_\alpha/ \|g\|_\infty\big) z,\quad D_{d+2}(z)=\big(L_1/ \|K\|_\infty\big) z\quad L_{d+1}(z)=z,\quad L_{d+2}(z)=z^{2}.
\end{gather*}
Additionally, if $\bar{\bX}_1^d$ consists of a single point $\bar{x}^{(d)}\in \bX_1^d$  then $L_{d+2}\equiv 0.$

\end{lemma}
The proof of the lemma is postponed to Appendix. We remark that $\varrho_{d+1}$ is a semi-metric, since $\alpha\in (0,1]$, and    the semi-metric $\varrho^{(d)}_n$ is independent on $n$. In view of latter remark
all quantities involved in Assumption \ref{ass:bounded_case} are independent on the choice of $r^{(\min)}_l(\cdot),r^{(\max)}_l(\cdot), \;l=\overline{1,d},$. We want to emphasize nevertheless that the assertion of the lemma is true for an arbitrary but {\it a priory}
chosen $r^{(\min)}_j(\cdot), r^{(\max)}_l(\cdot), \;l=\overline{1,d}$.

Thus,  Lemma \ref{lem:loc-empiric-pointwise} guarantees the verification of the main assumption of Section \ref{sec:subsec-prob-form}, that
makes possible the application of Theorems \ref{th:empiric_totaly_bounded_case}--\ref{th:LL-nonasym}.
Hence, we have to match the notations of these theorems to the notations used in the present section.

\smallskip

Since $k=d$ and $\mH_1^k=(0,1]^d$ we have $\mh^{(k)}=r$ and, therefore, in view of Assumption \ref{ass:on-fucnctions-g-and-K}
\begin{eqnarray*}
%\label{eq1:proof-theorem55-empirical}
&&\mathbf{G_\infty}(r):=\sup_{(\mz,\bar{x}^{(d)})\in\cZ\times\bar{\bX}_1^d}\;\sup_{x\in\bX_1^{d+1}}
\left|G\left(\Big\{r,\mz,\bar{x}^{(d)}\Big\},x\right)\right|\leq V^{-1}_r\|g\|_\infty \|K\|_\infty=:G_\infty(r).
\\*[1mm]
%\label{eq2:proof-theorem55-empirical}
&& \underline{G}_n:=\inf_{r\in\cR(n)}G_\infty(r)=V^{-1}_{r^{(\max)}(n)}\|g\|_\infty \|K\|_\infty,\quad\forall n\geq 1.
\end{eqnarray*}
We remark that the function $G_\infty(\cdot)$ is  independent of the choice of $\bar{\bX}_1^d$.
Define
$$
f^{(d)}_i\big(x^{(d)}\big)=\int_{\bX_{d+1}}f_i(x)\left[\mu^{(d)}\big(\rd x^{(d)}\big)\times\mu_{d+1}\big(\rd x_{d+1}\big)\right], \;\; i\geq 1,
$$
 and
let $3\leq \mathbf{n_1}\leq \mathbf{n_2\leq 2\mathbf{n_1}}$ be fixed. Set for any $\left(r,\bar{x}^{(d)}\right)\in (0,1]^d\times\bar{\bX}_1^d$
$$
F_{\mathbf{n_2}}\left(r,\bar{x}^{(d)}\right)=
\left\{
\begin{array}{ll}
\|g\|_\infty(\mathbf{n_2})^{-1}\sum_{i=1}^{\mathbf{n_2}} \int_{\bX_1^{d}}\left|K_{r}\left(\rho\big(x^{(d)},\bar{x}^{(d)}\big)\right)\right|f^{(d)}_{i}\big(x^{(d)}\big)\mu^{(d)}(\rd x^{(d)}),\;&\mathbf{n_1}=\mathbf{n_2};
\\*[2mm]
\|g\|_\infty\displaystyle{\sup_{i=\overline{1,\mathbf{n_2}}}} \int_{\bX_1^{d}}\left|K_{r}\left(\rho\big(x^{(d)},\bar{x}^{(d)}\big)\right)\right|f^{(d)}_{i}\big(x^{(d)}\big)\mu^{(d)}(\rd x^{(d)}),\;&\mathbf{n_1}\neq\mathbf{n_2},
\end{array}
\right.
$$
and note that in view of Assumption \ref{ass:on-fucnctions-g-and-K} $(\mathbf{ii})$
$$
F_{\mathbf{n_2}}(\mh)\leq F_{\mathbf{n_2}}\left(r,\bar{x}^{(d)}\right),\quad\forall \mh\in(0,1]^d\times\cZ\times\bar{\bX}_1^d.
$$
We remark that the function $F_{\mathbf{n_2}}(\cdot,\cdot)$ is  independent of the choice of $\cZ$. Put also
$$ F_{\mathbf{n_2}}:=\sup_{n\in\widetilde{\mathbf{N}}}\sup_{\big(r,\bar{x}^{d}\big)\in\cR(n)\times\bar{\bX}_1^d}
F_{\mathbf{n_2}}\big(r,\bar{x}^{d}\big)<\infty,
$$
where, remind, $\widetilde{\mathbf{N}}=\{\mathbf{n_1},\ldots,\mathbf{n_2}\}$.
 Finally for any $\mathbf{r}\in\overline{\bN}$ set  $F_{\mathbf{n_2},\mathbf{r}}(\cdot,\cdot)=\max\left[F_{\mathbf{n_2}}(\cdot,\cdot),e^{-\mathbf{r}}\right]$.

\subsubsection{Pointwise results} Here we will consider the case, where  $\bar{\bX}_1^d=\left\{\bar{x}^{(d)}\right\}$ and   $\bar{x}^{(d)}$ is a fixed element of $\bX_1^d$. Note that in view of Lemma \ref{lem:loc-empiric-pointwise} $L_{d+1}(z)=z$ and $L_{d+2}\equiv 0$ that implies $\cL^{(k)}\equiv 0$.

We will suppose that Assumption \ref{ass:sec:totaly_bounded_case} holds with $k=d,m=d+1$ and $\left(\mH_{d+1},\varrho_{d+1}\right)=\left(\cZ,[\md]^\alpha\right)$. It is equivalent obviously  to assume that
Assumption \ref{ass:sec:totaly_bounded_case} holds with $\left(\mH_{d+1},\varrho_{d+1}\right)=\left(\cZ,\md\right)$ and with the constants
$\tilde{N}=\alpha N$ and $\tilde{R}=R^{1/\alpha}$.

Let $\beta$ and $C_{N,R,m,k}$ be the constants defined in Theorem \ref{th:empiric_totaly_bounded_case}. Set for any $r\in (0,1]^d$ and $q>0$
\begin{gather*}
P(r)=(36d\delta^{-2}_*+6)\ln{\left(1+\sum_{l=1}^d \gamma_l\ln{\left\{\frac{2r^{(\max)}_l(n)}{r_l}\right\}}\right)}
+18C_{N,R,d+1,d};
\\*[2mm]
M_q(r)=\big(72d\delta^{-2}_*+2.5q+1.5\big)\sum_{l=1}^d \gamma_l\ln{\left(\frac{2r^{(\max)}_l(n)}{r_l}\right)}+36C_{N,R,d+1,d}.
\end{gather*}
and  define for $\mathbf{r}\in \overline{\bN}$ and $u>0$
\begin{eqnarray*}
\cV_{\mathbf{r}}^{(u)}\big(n,r,\bar{x}^{d}\big)&=&\boldsymbol{\lambda_1}\sqrt{\Big[F_{\mathbf{n_2},\mathbf{r}}\big(r,\bar{x}^{d}\big)(nV_r)^{-1}\Big]\Big[P(r)+
2\ln{\left\{1+\left|\ln{\left\{F_{\mathbf{n_2},\mathbf{r}}\big(r,\bar{x}^{d}\big)\right\}}\right|\right\}}
+u\Big]}
\\
&&\hskip-0.3cm +
\boldsymbol{\lambda_2}\Big[(nV_r)^{-1}\ln^{\beta}{(n)}\Big]
\Big[P(r)+
2\ln{\left\{1+\left|\ln{\left\{F_{\mathbf{n_2},\mathbf{r}}\big(r,\bar{x}^{d}\big)\right\}}\right|\right\}}
+u\Big];
\\*[2mm]
\cU_{\mathbf{r}}^{(u,q)}\big(n,r,\bar{x}^{d}\big)&=&\boldsymbol{\lambda_1}\sqrt{\Big[F_{\mathbf{n_2},\mathbf{r}}\big(r,\bar{x}^{d}\big)(nV_r)^{-1}\Big]\Big[M_q(r)+
2\ln{\left\{1+\left|\ln{\left\{F_{\mathbf{n_2},\mathbf{r}}\big(r,\bar{x}^{d}\big)\right\}}\right|\right\}}
+u\Big]}
\\
&&\hskip-0.3cm +
\boldsymbol{\lambda_2}\Big[(nV_r)^{-1}\ln^{\beta}{(n)}\Big]
\Big[M_q(r)+
2\ln{\left\{1+\left|\ln{\left\{F_{\mathbf{n_2},\mathbf{r}}\big(r,\bar{x}^{d}\big)\right\}}\right|\right\}}
+u\Big],
\end{eqnarray*}
where $\boldsymbol{\lambda_1}=\sqrt{\|g\|_\infty \|K\|_\infty}\lambda_1$, $\boldsymbol{\lambda_2}=\|g\|_\infty \|K\|_\infty\lambda_2$ and $\lambda_1,\;\lambda_2$ are defined in Theorem \ref{th:empiric_totaly_bounded_case}.

The result below is the direct consequence of Theorem \ref{th:empiric_totaly_bounded_case} and Lemma \ref{lem:loc-empiric-pointwise}.
We remark that  defined above quantities are   functions of $r$ and $n$ since $\bar{x}^{d}$ is fixed. Since they do not depend on the variable $\mz$, these quantities will be  automatically upper functions for
$$
\displaystyle{\zeta_r\big(n,\bar{x}^{(d)}\big):=\sup_{x_{d+1}\in\bX_{d+1}}\Big|\xi_{r,\mz,\bar{x}^{(d)}}\big(\bar{x}^{(d)}\big)\Big|}.
$$

\begin{theorem}
\label{th:empirical-product-general}
Let Assumption \ref{ass:on-fucnctions-g-and-K}  be fulfilled and suppose that Assumption \ref{ass:sec:totaly_bounded_case} holds with $k=d,m=d+1$ and $\left(\mH_{d+1},\varrho_{d+1}\right)=\left(\cZ,[\md]^\alpha\right)$.

Then for any given  decreasing sequences $0<r^{(\min)}_l(n)\leq r^{(\max)}_l(n)\leq 1, \;l=\overline{1,d},\; n\geq 1,$ any $\bar{x}^{d}\in\bar{\bX}_1^d$
 any $\mathbf{r}\in\bN$, $\boldsymbol{b}>1$ $ u\geq 1$  and $q\geq 1$
\begin{gather*}
\bP_\rf\left\{\sup_{n\in\widetilde{\mathbf{N}}}\sup_{r\in\widetilde{\cR}(n)}
\Big[\zeta_r\big(n,\bar{x}^{(d)}\big)-\cV_{\mathbf{r}}^{(u)}\big(n,r,\bar{x}^{d}\big)\Big]\geq 0\right\}
\leq 2419\; e^{-u}
;
\\
%*[2mm]
\bE_\rf\left\{\sup_{n\in\widetilde{\mathbf{N}}}\sup_{r\in\widetilde{\cR}(n)}
\Big[\zeta_r\big(n,\bar{x}^{(d)}\big)-\cU_{\mathbf{r}}^{(u,q)}\big(n,r,\bar{x}^{d}\big)\Big]\right\}^q_+
%\\*[2mm]&&
\leq c^{\prime}_q\left[\sqrt{\frac{F_{\mathbf{n_2}}}{\mathbf{n_1}V_{r^{(\max)}(\mathbf{n_1})}}}
\vee\left(\frac{\ln^{\beta}{(\mathbf{n_2})}}{V_{r^{(\max)}(\mathbf{n_1})}\mathbf{n_1}} \right) \right]^q e^{-u},
\end{gather*}
where
$c^\prime_q= 2^{(7q/2)+5}3^{q+4}\Gamma(q+1)\left(C_{D,\boldsymbol{b}}\max\left[\sqrt{\|g\|_\infty \|K\|_\infty},\|g\|_\infty \|K\|_\infty\right]\right)^{q}$.
\end{theorem}

The explicit expression for $C_{D,\boldsymbol{b}}$ can be also found in Theorem \ref{th:empiric_totaly_bounded_case}. In the case considered here it is completely determined by  $(\gamma_1,\ldots\gamma_d)$,  $L_1$,  $L_\alpha$ and $\boldsymbol{b}$.

As well as  the assertions of   Theorem \ref{th:empiric_totaly_bounded_case} the latter theorem is proved without any assumption imposed on the densities
$f_i,\; i=\overline{1,n}$. The choice of $r^{(\min)}_l(n), r^{(\max)}_l(n), \;l=\overline{1,d},\; n\geq 1,$ is also assumption free.
Additionally, Assumption \ref{ass:sec:totaly_bounded_case} can be replaced by (\ref{eq:rem-after-theorem3}), see Remark \ref{rem:after-theorem3}.

Note also that if   $g\big(\mz,\cdot\big)=\bar{g}(\cdot)$, for any $\mz\in\cZ$, then Assumption \ref{ass:sec:totaly_bounded_case} in not needed anymore and, moreover, Assumption \ref{ass:on-fucnctions-g-and-K}
($\mathbf{ii}$) is verified for an arbitrary bounded $\bar{g}$. Hence, in this case the assertions of Theorem \ref{th:empirical-product-general} are established  under very mild  Assumption \ref{ass:on-fucnctions-g-and-K}
($\mathbf{i}$) imposed on the function $K$.
\begin{remark}
\label{rem:after-cor-after-theorem-empirical-product-general}
We note that the discussed in Introduction so-called price to pay for uniformity disappears if  $r=r^{(\max)}$. Indeed, $P\left(r^{(\max)}\right)$ and $M_q\left(r^{(\max)}\right)$ are absolute constants. This property is crucial, in particular, for constructing statistical procedures used in the estimation of functions possessing inhomogeneous smoothness, see  \cite{LepMam}, \cite{KerkLep}.

\end{remark}

\paragraph{Some additional assumptions and their consequences} To apply Theorem \ref{th:empirical-product-general}  to specific problems one needs to find an  efficient upper bound for the quantity $F_{\mathbf{n_2}}(\cdot,\cdot)$.
Below we provide with sufficient condition allowing to solve this problem under general consideration and we will not be tending here to the maximal generality. We impose some additional restrictions on the densities $f_i,\; i=\overline{1,n},$ and on the measures $\mu_l$ of $\rho_l$-balls in the spaces $\bX_l,\; l=\overline{1,d}$. Moreover, we should precise  the behavior of  the function $K$ at infinity.
Then,  we will use these assumptions   for establishing of the law of iterated logarithm.

Introduce the following notations.
For any $t\in\bR^d_+$ define
$$
\check{K}(t)=\sup_{|u|\notin \Pi_t}|K(u)|,\quad \Pi_t=[0,t_1]\times\cdots\times[0,t_d].
$$
For any $l=\overline{1,d}$,  $\; x_l\in\bX_l$  and  $\mathrm{r}>0$ set $\bB_l\big(\mathrm{r},x_l\big)=\big\{y\in\bX_l:\;\;\rho_l\big(y,x_l\big)\leq \mathrm{r}\big\}$.
\begin{assumption}
\label{ass:after-th-empirical-product-general}
\noindent  There exists $L_2>0$ such that
\begin{eqnarray}
\label{eq2:ass-normlization-using-gammas}
&& \sup_{t\in\bR^d_+}\bigg[\Big(\prod_{l=1}^dt_l^{1+\gamma_l}\Big)\check{K}(t)\bigg]\leq L_2;
\end{eqnarray}
For any $l=\overline{1,d}$ and any $x_l\in\bX_l$ one has
$
\bX_l=\cup_{\mathrm{r}>0}\left(\bB_l\big(\mathrm{r},x_l\big)\right)
$
 and there exist $L^{(l)}>0$
\begin{eqnarray}
\label{eq1:ass-normlization-using-gammas}
&& \mu_l\Big(\bB_l\big(\mathrm{r},x_l\big)\Big)\leq L^{(l)}\mathrm{r}^{\gamma_l},\quad \forall \mathrm{r}>0;
\end{eqnarray}
Moreover,
\begin{equation}
\label{eq:ass-boundness-of densities}
\sup_{i\geq 1}\sup_{x^{(d)}\in\bX_1^d}f^{(d)}_i\big(x^{(d)}\big)=: \rf_\infty<\infty.
\end{equation}

\end{assumption}

The condition (\ref{eq2:ass-normlization-using-gammas}) is obviously fulfilled if $K$ is compactly supported on $[0,1]^d$. It is also satisfied
in the case of Gaussian or Laplace kernel.

The condition (\ref{eq1:ass-normlization-using-gammas}) can be  easily checked if $\bX_l,\;l=\overline{1,d}$ are doubling metric spaces.
In particular, if $\bX_l=\bR$ and $\mu_l, \;l=\overline{1,d}, $ are the Lebesgue measures than (\ref{eq1:ass-normlization-using-gammas}) holds with
$L^{(l)}=1, \;\gamma_l=1,\;l=\overline{1,d}.$ If $\bX_l=\bR^{d_l}, \;\;l=\overline{1,d}, $ then (\ref{eq1:ass-normlization-using-gammas}) holds with $\gamma_l=d_l$ and the constants $L^{(l)}$ depending on the choice of the distances $\rho_l$.

% Theorem \ref{th:empiric_totaly_bounded_case} is applicable on $\cR\times\mH_{d+1}\times\bar{\cX}_1^d$.

 As to condition (\ref{eq:ass-boundness-of densities}) we remark that the boundedness of the entire density $f_i$ is not required. For example, under independence structure, i.e. $f_i(x)=f^{(d)}_i\big(x^{(d)}\big)p_i\big(x_{d+1}\big)$, the densities $p_i$ may be unbounded.
\begin{lemma}
\label{lem:bound-for-F_n2}
The following bound holds under Assumption \ref{ass:after-th-empirical-product-general}:
\begin{equation*}
\sup_{\mathbf{n_2}\geq 1}\sup_{r\in (0,1]}\sup_{\bar{x}^{(d)}\in\bX_1^d}F_{\mathbf{n_2}}\left(r,\bar{x}^{(d)}\right)\leq
2^{d}\rf_\infty \|g\|_\infty L_2\prod_{l=1}^d2^{\gamma_l}L^{(l)}.
\end{equation*}

\end{lemma}
The proof of lemma is postponed to Appendix.
Our goal now is to deduce the law of iterated logarithm for $\zeta_r\big(n,\bar{x}^{(d)}\big)$ from Theorem \ref{th:LIL-nonasym}.
Set for  any $n\in\bN^*$ and $a>0$
$$
\overline{\cR}_a(n)=\Big\{r\in (0,1]^d:\;\;V_{r}\geq \; n^{-1}(\ln n)^a \Big\}.
$$
and choose $h^{(\max)}=(1,\ldots,1)$ and $h^{(\min)}=\left(1/n,\ldots,1/n\right)$.
\begin{remark}
\label{rem:precedes-LIL-product}
$\quad 1^0.$ Note that
$
\overline{\cR}_a(n)\subset [ n^{-1},1]^d=:\widetilde{\cR}(n)
$
for any $n\geq 3$ and any $a>0$
and, therefore, the assertion of Lemma \ref{lem:loc-empiric-pointwise} holds.

%$ 2^0.$ We note that $\overline{\cR}_a(n)\subset\overline{\cR}_a(n+1)$ for any $n\geq 3$ since the sequence $n^{-1}(\ln n)^a,\;n\geq 3,$ is strictly decreasing for any $a>0$. Therefore, Assumption \ref{ass:dependence-on-n} is fulfilled for $\overline{\cR}_a(n)\times\cZ,\; n\geq 3.$
\vskip0.1cm

$ 2^0.$ We have
$
\underline{G}_n=\|K\|_\infty \|g\|_\infty,\quad \overline{G}_n=\|K\|_\infty \|g\|_\infty n^{-d}
$
 for any $n\geq 1$ and, therefore, (\ref{eq1:LIL}) is verified with $\mathfrak{c}=\|K\|_\infty \|g\|_\infty$ and $\mb=d$.

\vskip0.1cm

$ 3^0.$ Lemma \ref{lem:bound-for-F_n2} implies that the condition (\ref{eq2:LIL}) holds with $\mathbf{F}\leq 2^{d}\rf_\infty \|g\|_\infty L_2\prod_{l=1}^d2^{\gamma_l}L^{(l)}.$

\vskip0.1cm

$ 4^0.$ In view of Lemma \ref{lem:loc-empiric-pointwise} $L_{d+1}(z)=z$  and $L_{d+2}\equiv 0$, that implies $\cL^{(k)}\equiv 0$.
Hence, the condition  (\ref{eq22:LIL}) is fulfilled for any $\ma>0$.

\end{remark}

Thus, all assumptions of Theorem \ref{th:LIL-nonasym} are checked and we come to the following statement.

\begin{theorem}
\label{th:LIL-nonasym-product}
Let Assumptions \ref{ass:on-fucnctions-g-and-K} and \ref{ass:after-th-empirical-product-general} be fulfilled and suppose that Assumption \ref{ass:sec:totaly_bounded_case} holds with $k=d,m=d+1$ and $\left(\mH_{d+1},\varrho_{d+1}\right)=\left(\cZ,[\md]^\alpha\right)$.
Then there exists $\Upsilon>0$ such that  for any $\bar{x}^{d}\in\bar{\bX}_1^d$ and any $a>2$
\begin{gather*}
\bP_\rf\left\{\sup_{n\geq \mathbf{j}}\;\sup_{r:\;n^{-1}(\ln n)^a \leq V_{r}\leq 1}\Bigg[\frac{\sqrt{nV_r}\;\zeta_r\big(n,\bar{x}^{(d)}\big)}
{\sqrt{\ln{\big(1+\ln{(n)}\big)}}}\Bigg]\geq \Upsilon\right\}
\leq  \frac{2419}{\ln(\mathbf{j})}.
\end{gather*}

\end{theorem}

\begin{remark}
\label{rem:after-theorem-LIL-product}
The inspection of the proof of Theorem  \ref{th:LIL-nonasym}  together
with  Lemma \ref{lem:bound-for-F_n2} allows us to assert  that the statement of Theorem \ref{th:LIL-nonasym-product} is \textsf{uniform} over the set of bounded
densities.

More precisely, for any $\mathfrak{f}>0$ there exists $\Upsilon(\mathfrak{f})$ such that
\begin{gather}
\label{eq:rem-after-th-LIL-product}
\sup_{\mathrm{f}\in\cF_\mathfrak{f}}\bP_\rf\left\{\sup_{n\geq \mathbf{j}}\;\sup_{r:\;n^{-1}(\ln n)^a \leq V_{r}\leq 1}\Bigg[\frac{\sqrt{nV_r}\;\zeta_r\big(n,\bar{x}^{(d)}\big)}
{\sqrt{\ln{\big(1+\ln{(n)}\big)}}}\Bigg]\geq \Upsilon(\mathfrak{f})\right\}
\leq  \frac{2419}{\ln(\mathbf{j})},
\end{gather}
where $\cF_\mathfrak{f}=\big\{(f_i,i\geq 1):\;\; \mathrm{f}_\infty\leq \mathfrak{f}\big\}$. As before the explicit expression of $\Upsilon(\cdot)$ is available.

\end{remark}

The following consequence of Theorem \ref{th:LIL-nonasym-product} is straightforward.
\begin{equation}
\label{eq:LIL-product}
\limsup_{n\to\infty}\sup_{r:\;n^{-1}(\ln n)^a \leq V_{r}\leq 1}\Bigg[\frac{\sqrt{nV_r}\;\zeta_r\big(n,\bar{x}^{(d)}\big)}
{\sqrt{\ln{\big(1+\ln{(n)}\big)}}}\Bigg]
\leq\Upsilon
\quad \bP_f-\mathbf{a.s.}
\end{equation}
Theorem \ref{th:LIL-nonasym-product} generalizes the existing results, see for example \cite{Dony2},
in the following directions.
\begin{enumerate}
\item
{\it Structural assumption.} The  structural condition (\ref{eq:structural-assumption-empiric-general}) is imposed in cited papers but with additional restriction: either $g(\mz,x)\equiv \text{const} $ ("density case") or $g(\mz,x)=\bar{g}(x)$ ("regression case"). It excludes, for instance, the problems appearing in robust estimation. We note that  Assumption \ref{ass:on-fucnctions-g-and-K} ($\mathbf{ii}$) is fulfilled here if $\bar{g}$ is bounded function and Assumption \ref{ass:sec:totaly_bounded_case} is not needed anymore, since $\bar{g}$ is independent of $\mz$.
\item
{\it Anisotropy.} All known to the author results treat the case where $\bX_l=\bR,\; l=\overline{1,d},$ and $\cR(n)=\left\{(r_1,\ldots,r_d)\in (0,1]^d: \;\; r_l=\mathrm{r},\;\forall l=\overline{1,d},\;\;\mathrm{r}\in\left[\mathrm{r}^{(min)}(n),\mathrm{r}^{(max)}(n)\right]\right\}$ (isotropic case). We remark that  (\ref{eq1:ass-normlization-using-gammas}) is automatically fulfilled  with $\gamma_l=1, L^{(l)}=1,\;l=\overline{1,d}$,  and $V_{r}=\mathrm{r}^{d}$.
Note also that we consider independent but not necessarily identically distributed random variables. This is important, in particular,
for various estimation problems arising  in nonparametric regression model.
\item
{\it Kernel.} We do not suppose that the function $K$ is compactly supported. For instance, one can use the gaussian or laplace kernel.
It allows, for instaince, to consider the problems where $X_1^d$ is not linear space. In particular, it can be some  manifold satisfying doubling condition.
\item
{\it Non-asymptotic nature.} The existing results are presented as in (\ref{eq:LIL-product}). Note, however, that the random field $\zeta_r\big(n,\bar{x}^{d}\big)$ appears  in various  areas of nonparametric estimation (density estimation, regression). As the consequence a.s. convergence has no much sense since there is no a unique probability measure (see, also Remark \ref{rem:after-theorem-LIL-product}).
\end{enumerate}

\subsubsection{Sup-norm results} Here we  consider  $\bar{\bX}_1^d=\bX_1^d$. We  assume that there exists
$\left\{\mathrm{X}_{\mathbf{i}},\;\mathbf{i}\in\mathbf{I}\right\}$ which is $\mn$-totally bounded cover of $\left(\bX_1^d,\rho^{(d)}\right)$
satisfying Assumption \ref{ass:partially-bounded-case} ($\mathbf{i}$) and  possessing  the separation property.

\begin{assumption}
\label{ass:sec-sup-norm-case-general} There exists $\mathfrak{t}>0$ such that for any $\mathbf{i},\mathbf{k}\in\mathbf{I}$ satisfying
$\mathrm{X}_\mathbf{i}\cap\mathrm{X}_{\mathbf{k}}=\emptyset$
$$
\inf_{x^{(d)}\in\mathrm{X}_{\mathbf{i}}}\inf_{\;y^{(d)}\in\mathrm{X}_{\mathbf{k}}}\rho^{(d)}\left(x^{(d)},y^{(d)}\right)>\mt.
$$
\end{assumption}

 Also we  suppose that  Assumption \ref{ass:partially-bounded-case} ($\mathbf{ii}$) holds with $k=d,m=d+1$ and $\left(\mH_{d+1},\varrho_{d+1}\right)=\left(\cZ,[\md]^\alpha\right)$. We remark that in the considered case this assumption coincides with
Assumption \ref{ass:sec:totaly_bounded_case}.

Let, as previously,   $0<r^{(\min)}_l(n)\leq r^{(\max)}_l(n)\leq 1, \;l=\overline{1,d},\; n\geq 1,$ be given decreasing sequences,
\begin{eqnarray*}
&&\mH(n)=\cR(n)\times\cZ\times X_1^d,\qquad\cR(n)=\prod_{l=1}^d \big[r^{(\min)}_l(2n),r^{(\max)}_l(n)\big];
\\
&&\widetilde{\mH}(n)=\widetilde{\cR}(n)\times\cZ\times X_1^d,\qquad\widetilde{\cR}(n)=\prod_{l=1}^d \big[r^{(\min)}_l(n),r^{(\max)}_l(n)\big].
\end{eqnarray*}
Our last condition relates the choice of the vector $r^{(\max)}(n),\; n\geq 1$ and the kernel $K$ with the parameter $\mt$ appearing in Assumption \ref{ass:sec-sup-norm-case-general}. Let us assume that for any $n\geq 1$
 \begin{equation}
\label{eq:condition-on-r-max}
\sup_{r\in\cR(n)}\sup_{|u|\notin (0,\mt]^d}|K(u/r|\leq \|K\|_\infty n^{-1}.
\end{equation}
Note that  (\ref{eq:condition-on-r-max}) holds if $K$ is compactly supported on $[-\mt,\mt]^d$ and  $r^{(\max)}(n)\in (0,\mt)^{d}$ for any $n\geq 1$.
\begin{lemma}
\label{lem:verification-ass8}
Assumption \ref{ass:sec-sup-norm-case-general} and (\ref{eq:condition-on-r-max}) imply Assumption \ref{ass:as-assumption-totally-bounded}.
\end{lemma}

The proof of lemma is given in Appendix.
%To simplify the presentation of our results later on we will assume, without loss of generality, that$\|g\|_\infty\|K\|_\infty\geq 1$.
Set for any $r\in (0,1]^d$ and $v>0$
\begin{gather*}
%\widehat{P}_{v}(r)=\Big([36d+54N]\delta^{-2}_*+2v+6\Big)\ln{\left(2V_r^{-1}\right)}
%+\mathbf{C};
%\\
\widehat{M}_{q,v}(r)=\big([72d+108N]\delta^{-2}_*+2.5q+2v+1.5\big)\ln{\left(2V_r^{-1}\right)}
+\mathbf{C},
\end{gather*}
where we have put $\mathbf{C}=72N\delta^{-2}_*\left|\log_2{\left(\|g\|_\infty\|K\|_\infty\right)}\right|+36C_{N,R,d+1,d}$.

\smallskip

Let $3\leq \mathbf{n_1}\leq \mathbf{n_2}\leq 2\mathbf{n_1}$ be fixed. Set
$\widehat{{F}}_{\mathbf{n_2}}\left(r,\bar{x}^{(d)}\right)=
\max\left[{F}_{\mathbf{n_2}}\left(r,\bar{x}^{(d)}\right),\mathbf{n_2}^{-1}\right]$ and define

\begin{eqnarray*}
\widehat{\cU}^{(v,z,q)}\big(n,r,\bar{x}^{(d)}\big)&=&\boldsymbol{\lambda_1}\sqrt{\Big[\widehat{{F}}_{\mathbf{n_2}}
\left(r,\bar{x}^{(d)}\right)(nV_r)^{-1}\Big] \Big[\widehat{M}_{q,v}(r)+2(v+1)\Big|\ln{\big\{\widehat{{F}}_{\mathbf{n_2}}\left(r,\bar{x}^{(d)}\right)\big\}}\Big|
+z\Big]}
\\
&&\hskip-0.3cm +
\boldsymbol{\lambda_2}\Big[(nV_r)^{-1}\ln^{\beta}{(n)}\Big]
\Big[\widehat{M}_{q,v}(r)+2(v+1)\Big|\ln{\big\{\widehat{{F}}_{\mathbf{n_2}}\left(r,\bar{x}^{(d)}\right)\big\}}\Big|
+z\Big].
\end{eqnarray*}
 Theorem \ref{th:empiric-sup-norm-product} below is the direct consequence of Lemma \ref{lem:loc-empiric-pointwise}, Lemma \ref{lem:verification-ass8} and Corollary \ref{cor:after-th:empiric-partially_totaly_bounded_case}. Remind that $
\displaystyle{\zeta_r\big(n,\bar{x}^{(d)}\big):=\sup_{x_{d+1}\in\bX_{d+1}}\Big|\xi_{r,\mz,\bar{x}^{(d)}}\big(\bar{x}^{(d)}\big)\Big|}
$
and $\widetilde{\mathbf{N}}=\left\{\mathbf{n_1},\ldots,\mathbf{n_2}\right\}$.
\begin{theorem}
\label{th:empiric-sup-norm-product}
Let Assumption \ref{ass:on-fucnctions-g-and-K}  be verified  and suppose that Assumption \ref{ass:partially-bounded-case} ($\mathbf{ii}$) holds with $k=d+1,m=d+2$ and $\left(\mH_{d+1},\varrho_{d+1}\right)=\left(\cZ,[\md]^\alpha\right)$. Suppose also that Assumption \ref{ass:partially-bounded-case} ($\mathbf{i}$) is fulfilled with  $\left(\mH_{d+2},\varrho_{d+2}\right)=\left(\bX_1^d,\rho^{(d)}\right)$ and
$\mathrm{H}_{d+2,\mathbf{i}}=\mathrm{X}_\mathbf{i},\;\mathbf{i}\in\mathbf{I},$ satisfying Assumption \ref{ass:sec-sup-norm-case-general}.
Assume  that (\ref{eq:condition-on-r-max}) holds as well and if $\mathbf{n_1}\neq\mathbf{n_2}$ let
 $\big(X_{i}\big)^d,\;i\geq 1,$ be identically distributed.

Then for any given  decreasing sequences $0<r^{(\min)}_l(n)\leq r^{(\max)}_l(n)\leq 1, \;l=\overline{1,d},\; n\geq 1,$  any $\boldsymbol{b}>1$, $q\geq 1$, $v\geq 1$ and $z\geq 1$
\begin{eqnarray*}
&&\bP_\rf\Bigg\{\sup_{n\in\widetilde{\mathbf{N}}}\sup_{\left(r,\bar{x}^{(d)}\right)\in\widetilde{\cR}(n)\times\bX_1^d}
\Big[\zeta_r\big(n,\bar{x}^{(d)}\big)-\widehat{\cU}^{(v,z,q)}\big(n,r,\bar{x}^{(d)}\big)\Big]\geq 0\Bigg\}
\leq  \mn^5\Big\{4838e^{-z}+2\mathbf{n_1}^{2-v}\Big\};
\\*[0mm]
&&\bE_\rf\Bigg\{\sup_{n\in\widetilde{\mathbf{N}}}\sup_{\left(r,\bar{x}^{(d)}\right)\in\widetilde{\cR}(n)\times\bX_1^d}
\Big[\zeta_r\big(n,\bar{x}^{(d)}\big)-\widehat{\cU}^{(v,z,q)}\big(n,r,\bar{x}^{(d)}\big)\Big]\Bigg\}^q_+
%\\*[2mm]&&
\\
&&\qquad\quad\leq 2\mn^5 c^{\prime}_q\left[\sqrt{\frac{\widehat{{F}}_{\mathbf{n_2}}}{\mathbf{n_1}V_{r^{(\max)}(\mathbf{n_1})}}}
\vee\left(\frac{\ln^{\beta}{(\mathbf{n_2})}}{V_{r^{(\max)}(\mathbf{n_1})}\mathbf{n_1}} \right) \right]^q+
2^{q+1}\mn^5\left(V_{r^{(\min)}(\mathbf{n_1})}\right)^{-q}\;\mathbf{n_1}^{2-v}.
\end{eqnarray*}
\end{theorem}

Remind that
$\widehat{{F}}_{\mathbf{n_2}}=\displaystyle{
\sup_{n\in\widetilde{\mathbf{N}}}\sup_{\left(r,\bar{x}^{(d)}\right)\in\widetilde{\cR}(n)\times\bX_1^d}}
\widehat{{F}}_{\mathbf{n_2}}\left(r,\bar{x}^{(d)}\right)$ and the expression for the constant $ c^{\prime}_q$ can be found in
Theorem \ref{th:empirical-product-general}. We also note that the first assertion of the theorem remains valid if one replaces the quantity
$\widehat{M}_{q,v}(r)$ by the smaller quantity  $\Big([36d+54N]\delta^{-2}_*+2v+6\Big)\ln{\left(2V_r^{-1}\right)}+\mathbf{C}/2$. But the corresponding upper function
will differ from $\widehat{\cU}^{(v,z,q)}$ only by numerical constant.

\vskip0.1cm

We also remark that  $\widehat{{F}}_{\mathbf{n_2}}\leq 2^{d}\rf_\infty \|g\|_\infty L_2\prod_{l=1}^d2^{\gamma_l}L^{(l)} $ for any
$\mathbf{n_2}\geq 3$
under Assumption \ref{ass:after-th-empirical-product-general} in view of Lemma \ref{lem:bound-for-F_n2}. Moreover, if $V_{r^{(\min)}(n)}\geq n^{-p}$ for some $p>0$ then $\widehat{M}_{q,v}(r)$ can be bound from above by
$
\big([72d+108N]\delta^{-2}_*+2.5q+2v+1.5\big)p\ln{\left(2n\right)}
$
which is independent on $r$. Hence, if both restrictions are fulfilled the upper function $\widehat{\cU}^{(v,z,q)}$ in Theorem \ref{th:empiric-sup-norm-product} takes rather simple form, namely
$$
\boldsymbol{\lambda_1}(q)\sqrt{\frac{\ln(n)+z}{nV_r}} +
\frac{\boldsymbol{\lambda_2}(q)\big[\ln^{\beta+1}{(n)}+z\big]}{nV_r},
$$
where the constant $\boldsymbol{\lambda_1}(q)$ and $\boldsymbol{\lambda_2}(q)$ can be easily computed.

\paragraph{Law of logarithm} In this paragraph we will additionally suppose that Assumption \ref{ass:after-th-empirical-product-general} holds.
Then, we remark first that statements $1^{0}-3^{0}$ of Remark \ref{rem:precedes-LIL-product} hold. Next, we note that  $L_{d+1}(z)=z$  and $L_{d+2}(z)=z^2$ in view of Lemma \ref{lem:loc-empiric-pointwise}  that implies $\cL^{(k)}(z)=\ln(z)$ for any $z\geq 1$.
Hence, the condition  (\ref{eq2:LL}) is fulfilled with $\ma=1$.

Thus, all assumptions of Theorem \ref{th:LL-nonasym} are checked and, taking into account that in our case
$$
\eta_{\mh^{(k)}}(n)=\big\|\zeta_r(n)\big\|_\infty:=\sup_{\bar{x}^{(d)}\in\bX_1^d}\zeta_r\big(n,\bar{x}^{(d)}\big),
$$
 we come to the following statement.
\begin{theorem}
\label{th:LL-product}
Let assumptions of Theorem \ref{th:empiric-sup-norm-product} be fulfilled and suppose additionally that that Assumption \ref{ass:after-th-empirical-product-general} holds. Then there exists $\boldsymbol{\Upsilon}$ for any $a>4$

\begin{gather*}
\bP_\rf\Bigg\{\sup_{n\geq \mathbf{j}}\;\
\sup_{r:\;n^{-1}(\ln n)^a \leq V_{r}\leq 1}\frac{\sqrt{nV_r}\;\big\|\zeta_r(n)\big\|_\infty}
{\sqrt{\ln{\left(V^{-1}_r\right)}\vee\ln{\ln{(n)}}}}\geq \boldsymbol{\Upsilon}
\Bigg\}
\leq  \frac{4840\mn^5}{\ln{(\mathbf{j})}}.
\end{gather*}

\end{theorem}
The uniform version over the set of bounded densities, similar to (\ref{eq:rem-after-th-LIL-product}), holds as well.

The immediate consequence of the latter theorem is so-called  "{\it uniform-in-bandwidth consistency}":
\begin{gather}
\label{eq:LL-product}
\limsup_{n\to\infty}\sup_{r:\;n^{-1}(\ln n)^a \leq V_{r}\leq 1}\frac{\sqrt{nV_r}\;\big\|\zeta_r(n)\big\|_\infty}
{\sqrt{\ln{\left(V^{-1}_r\right)}\vee\ln{\ln{(n)}}}}\leq \boldsymbol{\Upsilon}\quad\bP_\rf-\mathbf{a.s}
\end{gather}

The assertion of Theorem \ref{th:LL-product} and its corollary (\ref{eq:LL-product}) generalizes in several directions  the existing results \cite{mason-ein}, \cite{gine}, \cite{mason-ein2}, \cite{Dony1} (see, the discussion after Theorem  \ref{th:LIL-nonasym-product}).

\smallskip

We would like to conclude this section with the following remark. If  $K$ is compactly supported and $g\big(\mz,\cdot\big)=\bar{g}(\cdot)$ for any $\mz\in\cZ$, where $\bar{g}$ is a bounded function,  then all results of this section remain true under Assumptions  \ref{ass:partially-bounded-case}  ($\mathbf{i}$), \ref{ass:on-fucnctions-g-and-K} ($\mathbf{i}$), \ref{ass:sec-sup-norm-case-general}, (\ref{eq1:ass-normlization-using-gammas})
and  (\ref{eq:ass-boundness-of densities}).

\section{Proof of Propositions \ref{prop_uniform_local1}-\ref{prop_uniform_local3}}
\label{sec:proofs-of-propositions}

We start this section with establishing an  auxiliary result.
Let $\mL$ be a set, $\rd_1$ and $\rd_2$ be  semi-metrics on $\mL$ and let $\cL$  be a totally bounded subset of $\mL$ with respect to
$\rd_1$ and $\rd_2$ simultaneously. Let $N_\ri(\delta),\;\delta>0,$  denote the minimal number of balls
of the radius $\delta$ in the metric $\rd_\ri,\;\ri=1,2$ needed to cover the set $\cL$.

\begin{lemma}
\label{lem:two-metrics-local}
Let $l\in\bN^*$ and  $\delta_{1,j}, \delta_{2,j}>0,\; j=\overline{1,l}$ be an arbitrary numbers. One can construct the finite  subset
$L\Big(\big[\delta_{1,j},\delta_{2,j}\big], j=\overline{1,l}\Big):=\left\{\ell_1,\ldots,\ell_N\right\}\subset \cL$ with $N\leq \prod_{j=1}^{l}N_1\Big(\delta_{1,j}/2\Big)N_2\Big(\delta_{2,j}/2\Big)$ and such that
$$
\forall \ell\in\cL\;\;\exists \tilde{\ell}\in L: \;\; \rd_1\big(\ell,\tilde{\ell}\big)\leq\min_{j=\overline{1,l}}\delta_{1,j},\;\quad \rd_2\big(\ell,\tilde{\ell}\big)\leq \min_{j=\overline{1,l}}\delta_{2,j}.
$$

\end{lemma}

\paragraph{Proof of Lemma \ref{lem:two-metrics-local}} Set $N_{\ri,j}= N_\ri\big(\delta_{\ri,j}/2\big),\;\ri=1,2,\; j=\overline{1,l}$. Since $\cL$ is totally bounded in
$\rd_i, i=1,2$ there exist $L^{(\ri,j)}=\left\{\ell^{(\ri,j)}_1,\ldots,\ell^{(\ri,j)}_{N_{\ri,j}}\right\}\subset \cL$ such that
$$
\forall \ell\in\cL\;\;\exists\; \tilde{\ell}\in L^{(\ri,j)}: \;\; \rd_\ri\big(\ell,\tilde{\ell}\big)\leq 2^{-1}\delta_{\ri,j},\;\ri=1,2,\;j=\overline{1,l}.
$$
For any $k_{\ri}=1,\ldots N_{\ri,j},\;\ri=1,2,$ put
$
\cL_{k_\ri}^{(\ri,j)}=\left\{\ell\in\cL:\;\; \rd_\ri\Big(\ell,\ell^{(\ri,j)}_{k_\ri}\Big)\leq 2^{-1}\delta_{\ri,j}\right\}
$
and let $\cL^{(j)}_{k_1,k_2}=\cL^{(1,j)}_{k_1}\cap\cL^{(2,j)}_{k_2}$.
First we note that for any $j=\overline{1,l}$
\begin{equation}
\label{eq1:proof-lemma-metric-local}
\cL=\bigcup_{k_1=1}^{N_{1,j}}\bigcup_{k_2=1}^{N_{2,j}}\cL^{(j)}_{k_1,k_2}
\end{equation}
Moreover, the construction of $\cL_{k_1}^{(1,j)}$ and $\cL_{k_2}^{(2,j)}$ implies that
\begin{equation}
\label{eq2:proof-lemma-metric-local}
\rd_\ri\big(l_1,l_2\big)\leq \delta_{\ri,j},\;\ri=1,2,\;\;\forall l_1,l_2\in\cL^{(j)}_{k_1,k_2}.
\end{equation}
Put $\cN=\otimes_{j=1}^{l}\Big[\left\{1,\ldots N_{1,j}\right\}\times\left\{1,\ldots N_{2,j}\right\}\Big]$  and define for any
$\big(\rk^{(1)},\ldots,\rk^{(l)}\big)\in\cN$
$$
\cL_{\rk^{(1)},\ldots,\rk^{(l)}}=\bigcap_{j=1}^{l}\cL^{(j)}_{\rk^{(j)}}.
$$
%$\vec{k}_j\in\left\{1,\ldots N_{1,j}\right\}\times\left\{1,\ldots N_{2,j}\right\},\; j=\overline{1,l},$ be fixed
The choice of  an arbitrary point in each $\cL_{\rk^{(1)},\ldots,\rk^{(l)}}$ leads to the  construction of  $L\Big(\big[\delta^{(j)}_1,\delta^{(j)}_2\big], j=\overline{1,l}\Big)$ in view of (\ref{eq1:proof-lemma-metric-local}) and (\ref{eq2:proof-lemma-metric-local}).

It remains to note that the cardinality of $\cN$ is equal to $\prod_{j=1}^{l}N_1\Big(\delta_{1,j}/2\Big)N_2\Big(\delta_{2,j}/2\Big)$.

\epr

\subsection{Proof of Proposition \ref{prop_uniform_local1}}

 \paragraph{{\bf I}. Probability bound}
 Fix $\vec{s}\in\cS_{\ra,\rb} $ and put   $s_{1,k}=s_1\Big(2^{k/2}\Big)$  and  $s_{2,k}=s_2\left(2^{k}\right)$, $k\geq 0$.
For any $k\geq 0$  put $\delta_{1}(k)=(24)^{-1}\widetilde{\kappa}_12^{-k/2}s_{1,k},\; \delta_{2}(k)=(24)^{-1}\widetilde{\kappa}_2 2^{-k}s_{2,k}$ and note that
$\delta_{1}(k),\delta_{2}(k)\to 0,\;k\to\infty$ since $s_1,s_2\in\bS$.

 Let $Z_{k}=L\Big(\big[\delta_{1}(k),\delta_{2}(k)\big]\Big),\;k\geq0,$ be the set constructed in Lemma \ref{lem:two-metrics-local} with $\rd_1=\ra$, $\rd_2=\rb$,  $\cL=\widetilde{\Theta}$ and $l=1$. By $N_{k},\; k\geq 0,$ we denote the cardinality of $Z_{k}$.

 \par

Fix $\e>0$  and put $\ve=\e/(1+\e)$,  $k_0=\big\lfloor2\ln_2{(1/\ve)}\big\rfloor+1$.
Let $\theta_{m},\;m=1,\ldots,N_{k_0},$ be the elements of $Z_{k_0}$.
For any $m=1,\ldots,N_{k_0}$ define
$$
\Theta^{(m)}=\Big\{\theta\in\widetilde{\Theta}:\; \mathrm{a}(\theta,\theta_{m})\leq \delta_{1}(k_0),\;\; \mathrm{b}(\theta,\theta_{m})\leq \delta_{2}(k_0)\Big\} ,
$$
and remark that the definition of the sets $Z_{k_0}$ implies that $\widetilde{\Theta}=\bigcup_{m=1}^{N_{k_0}}\Theta^{(m)}$.

In view of the last remark we get
\begin{eqnarray}
\label{p00}
%\qquad\qquad
&&\mathrm{P}\left\{\sup_{\theta\in\widetilde{\Theta}}\Psi\left(\chi_{\theta}\right)\geq U^{(\e)}_{\vec{s}}(y,\widetilde{\kappa},\widetilde{\Theta})\right\} \leq \sum_{m=1}^{N_{k_0}}\mathrm{P}\left\{\sup_{\theta\in\Theta^{(m)}}\Psi\left(\chi_\theta\right)\geq U^{(\e)}_{\vec{s}}(y,\widetilde{\kappa},\widetilde{\Theta})\right\}.
\end{eqnarray}
For any $\theta\in\widetilde{\Theta}$ let  $z_k(\theta)$ be an arbitrary element of $Z_k$ satisfying
\begin{eqnarray}
\label{p00000}
&&\ra\big(\theta,z_k(\theta)\big)\leq \delta_{1}(k),\quad \rb\big(\theta,z_k(\theta)\big)\leq  \delta_{2}(k).
\end{eqnarray}
Fix $m=1,\ldots,N_{k_0}$. The continuity of the mapping  $\theta\mapsto\chi_{\theta}$  guarantees that $\mathrm{P}$-a.s.
%the following relation holds
\begin{equation}
\label{chaining}
\chi_{\theta}=\chi_{\theta_m}+\sum_{k=1}^{\infty}\left[
\chi_{z_{k}(\theta)}-\chi_{z_{k-1}(\theta)}\right],\quad \forall \theta\in\Theta^{(m)},
\end{equation}
where
$z_{k_0}(\theta)=\theta_m,\;\forall \theta\in\Theta^{(m)}$.
Note also that independently of $\theta$ for all $k\geq k_0+1$
\begin{eqnarray}
\label{p20}
&&\ra\big(z_{k}(\theta),z_{k-1}(\theta)\big)\leq \delta_{1}(k)+\delta_{1}(k-1)=:\widetilde{\delta}_{1}(k),
\\
\label{p2000}
&&\rb\big(z_{k}(\theta),z_{k-1}(\theta)\big)\leq\delta_{2}(k)+\delta_{2}(k-1)=:\widetilde{\delta}_{2}(k).
\end{eqnarray}
This is the simplest consequence of triangle inequality and (\ref{p00000}). Introduce the sequence $c_k, k\geq1,$:
$$
c_k=4^{-1}\max\left\{s_{1,k},s_{1,k-1},s_{2,k},s_{2,k-1}\right\}
$$
and remark that $\sum_{k\geq 1}c_k\leq 1$ that follows from the assumption $s_1,s_2\in\bS$.

 We get from sub-additivity of $\Psi$, (\ref{chaining}), (\ref{p20}) and (\ref{p2000}) for any $\theta\in\Theta^{(m)}$
\begin{eqnarray}
\label{p30}
\Psi\left(\chi_{\theta}\right)\leq \Psi\left(\chi_{\theta_m}\right)+
\sup_{k\geq k_0+1}\sup_{\substack{
(u,v)\in Z_{k}\times Z_{k-1}:\\
\mathrm{a}(u,v)\leq \widetilde{\delta}_{1}(k),\;\mathrm{b}(u,v)\leq \widetilde{\delta}_{2}(k)}
}c_k^{-1}\Psi\left(\chi_{u}-\chi_{v}\right),
\end{eqnarray}
%where $Z^{\prime}_k=Z_k, \; k\geq k_0$ and $Z^{\prime}_{k_0-1}=\{\theta_m\}$.
 To simplify the notations we will write $U$ instead of $U^{(\e)}_{\vec{s}}(y,\widetilde{\kappa},\widetilde{\Theta})$ and
 $\cE$ instead of $e_{\vec{s}}(\widetilde{\kappa},\widetilde{\Theta})$.

We obtain from (\ref{p30})
%(to simplify the notation we will write $\kappa(\cZ)$ instead of $\kappa_U(\cZ)$)
\begin{eqnarray}
\label{p40}
&&\mathrm{P}\left\{\sup_{\theta\in\Theta^{(m)}}\Psi\left(\chi_{\theta}\right)\geq U \right\}\leq
\mathrm{P}\Big\{\Psi\left(\chi_{\theta_m}\right)\geq  U(1+\e)^{-1}\Big\}
\\*[2mm]
&&
+\sum_{k=k_0+1}^{\infty}\;\sum_{\substack{(u,v)\in Z_{k}\times Z_{k-1}:\\
\mathrm{a}(u,v)\leq \widetilde{\delta}_{1}(k),\;\mathrm{b}(u,v)\leq \widetilde{\delta}_{2}(k)}}
\mathrm{P}\Big\{\Psi\left(\chi_{u}-\chi_{v}\right)\geq \ve U c_k \Big\}=:I_1+I_2.\nonumber
\end{eqnarray}
We have in view of  Assumption \ref{ass:fixed_theta_local} ({\it 1})
%and \ref{ass:parameter_local}
\begin{eqnarray}
\label{p50_local}
\mathrm{P}\left\{\Psi\left(\chi_{\theta_m}\right)\geq  U(1+\e)^{-1} \right\}\nonumber
%\\*[2mm]
%&&
&\leq&
\mathrm{c}\exp{\left\{
-\frac{(1+\e)^{-2}U^{2}}{\big\{A\left(\theta_m\right)\big\}^{2}+ (1+\e)^{-1}U B\left(\theta_m\right)}
\right\}}
\nonumber
\\*[2mm]
%\nonumber\\*[3mm]
&\leq&\mathrm{c}\exp{\left\{
-\frac{(1+\e)^{-2} U^{2}}{\widetilde{\kappa}_1^{2}+   U \widetilde{\kappa}_2}
\right\}}\leq\mathrm{c}\exp{\left\{-(1+\e)^{-2}\big(y+2\ve^{-2}\cE\big)\right\}}
\nonumber\\*[2mm]
&\leq&\mathrm{c}\exp{\left\{-(1+\e)^{-2}
y-\ve^{-2}\cE\right\}}.
%\leq \mathrm{c}\exp{\left\{-y\right\}}.
\end{eqnarray}
In order to get (\ref{p50_local})  we have first used  that $\widetilde{\kappa}_1\geq \sup_{\theta\in\widetilde{\Theta}}A(\theta)$, $\widetilde{\kappa}_2\geq \sup_{\theta\in\widetilde{\Theta}}B(\theta)$.
Next, we have used that $U\geq v$, where $v$ is the maximal root of the equation
\begin{equation}
\label{eq:root-local}
\frac{u^{2}}{\widetilde{\kappa}_1^{2}+  u\widetilde{\kappa}_2}=y+2\ve^{-2}\cE.
\end{equation}
We also have used that $(1+\e)^{-2}\geq 1/2$.

In view of Assumptions  \ref{ass:fixed_theta_local} ({\it 2}),  \ref{ass:metric-case-local}, (\ref{p20}) and (\ref{p2000}) for any  $u,v\in Z_{k}\times Z_{k-1}$ satisfying $\mathrm{a}(u,v)\leq \widetilde{\delta}_{1}(k),\;\mathrm{b}(u,v)\leq \widetilde{\delta}_{2}(k) $ we have
\begin{eqnarray}
\label{p6_local}
\mathrm{P}\left\{\Psi\left(\chi_{\phi[u]}-\chi_{\phi[v]}\right)
\geq \ve U c_k \right\}
%\nonumber
%\\*[2mm]&&
&\leq& \mathrm{c}\exp\left\{-
\frac{\big(\ve U c_k \big)^{2}}{\big\{\mathrm{a}(u,v)\big\}^{2}+
\big(\ve U c_k\big)\;\mathrm{b}(u,v)}
\right\}
\nonumber
\\*[2mm]
&\leq& \mathrm{c}\exp\left\{-
\frac{(\ve U)^{2}}{\Big\{\widetilde{\delta}_{1}(k) c_k^{-1}\Big\}^{2}+
 U\Big\{\widetilde{\delta}_{2}(k)c_k^{-1}\Big\}}.
\right\}
%\nonumber\\*[2mm]
\end{eqnarray}
Here we have used that $\ve<1$. Let us remark that
\begin{eqnarray*}
\widetilde{\delta}_{1}(k)c_k^{-1}&\leq&
4(24)^{-1}\widetilde{\kappa}_1\left(2^{-k/2}s_{1,k}+2^{-(k-1)/2}s_{1,k-1}\right)\min\left\{s^{-1}_{1,k},s^{-1}_{1,k-1}\right\}\leq \widetilde{\kappa}_1 2^{-k/2-1};
\\
\widetilde{\delta}_{2}(k)c_k^{-1}&\leq&
4(24)^{-1}\widetilde{\kappa}_2\left(2^{-k}s_{2,k}+2^{-k+1}s_{2,k-1}\right)\min\left\{s^{-1}_{2,k},s^{-1}_{2,k-1}\right\}\leq \widetilde{\kappa}_2 2^{-k-1}.
\end{eqnarray*}
Thus, continuing (\ref{p6_local}) we obtain
\begin{eqnarray}
\label{p6000_local}
&&\mathrm{P}\left\{\Psi\left(\chi_{\phi[u]}-\chi_{\phi[v]}\right)
\geq \ve U c_k \right\}
\leq \mathrm{c}\exp\left\{-
\frac{2^{k+1}\ve^{2}\;U^{2}}{\kappa_1^{2}+
 U\kappa_2}
\right\}\leq \mathrm{c}\exp{\left\{-2^{k+1}\ve^{2}\big(y+2\ve^{-2}\cE\big)\right\}}.
%\nonumber\\*[2mm]
%\leq
%\leq \textsf{c}\exp{\left\{-\frac{e_\cZ2^{k}}{(k+1)^{2}}-\frac{4(y-e_\cZ)}{9}\right\}}.
\end{eqnarray}
Here we have used (\ref{eq:root-local}).
Noting that the right hand side of (\ref{p6000_local}) does not depend on $u,v$ we get
\begin{eqnarray}
\label{p7_local}
I_2&\leq&\mathrm{c}
%\exp{\left\{-\frac{4(y-e_\cZ)}{9}\right\}}
\sum_{k=k_0+1}^{\infty}
N_{k}N_{k-1}\exp{\left\{
-2^{k+1}\ve^{2} \big(y+2\ve^{-2}\cE\big)
\right\}}\nonumber
\\
&\leq&\mathrm{c}\exp{(-y)}\sum_{k=k_0+1}^{\infty}
N_{k}N_{k-1}\exp{\left\{
-2^{k+2}\cE-2^{k-k_0}
\right\}}.
\end{eqnarray}
Here we have used the definition of $k_0$ and that $y\geq 1$.

Let us make several remarks. First we note that in view of Lemma \ref{lem:two-metrics-local}
\begin{eqnarray*}
&&\ln{\big(N_k N_{k-1}\big)}\leq 2\ln{\big(N_{k}\big)}
%\\
\leq 2\Big[\mathfrak{E}_{\widetilde{\Theta},\;\mathrm{a}}\Big((24)^{-1}\widetilde{\kappa}_12^{-1-k/2}s_{1,k}\Big)
+\mathfrak{E}_{\widetilde{\Theta},\;\mathrm{b}}\Big((24)^{-1}\widetilde{\kappa}_22^{-k-1}s_{2,k}\Big)\Big].
\end{eqnarray*}
Taking into account that $s_{1,k}=s_1\big(2^{k/2}\big)$
and denoting $\delta_1=2^{k/2}$ we obtain
from (\ref{eq3:def_local})
\begin{eqnarray*}
\mathfrak{E}_{\widetilde{\Theta},\;\mathrm{a}}\Big(\widetilde{\kappa}_12^{-1-k/2}s_{1,k}\Big)=
\mathfrak{E}_{\widetilde{\Theta},\;\mathrm{a}}\Big(\widetilde{\kappa}_1(48\delta_1)^{-1}s_{1}(\delta_1)\Big)\leq \delta_1^{2}e_{s_1}^{(a)}\big(\widetilde{\kappa}_1,\widetilde{\Theta}\big)=
2^{k}e_{s_1}^{(a)}\big(\widetilde{\kappa}_1,\widetilde{\Theta}\big).
\end{eqnarray*}
Taking into account that $s_{2,k}=s_1\big(2^{k}\big)$
and denoting $\delta_2=2^{k}$ we obtain from (\ref{eq3:def_local})
%from (\ref{eq:inequalities-for-e_s})
\begin{eqnarray*}
\mathfrak{E}_{\widetilde{\Theta},\;\mathrm{b}}\Big(\widetilde{\kappa}_22^{-1-k}s_{2,k}\Big)=
\mathfrak{E}_{\widetilde{\Theta},\;\mathrm{b}}\Big(\widetilde{\kappa}_2(48\delta_2)^{-1}s_{2}(\delta_2)\Big)\leq \delta_2e_{s_2}^{(b)}\big(\widetilde{\kappa}_2,\widetilde{\Theta}\big)=
2^{k}e_{s_2}^{(b)}\big(\widetilde{\kappa}_2,\widetilde{\Theta}\big).
\end{eqnarray*}
Thus, we have for any $k\geq 1 $
\begin{eqnarray}
\label{p777_local}
\ln{\big(N_{k}N_{k-1}\big)}\leq 2^{k+1}\left[e_{s_1}^{(a)}(\widetilde{\Theta})+e_{s_2}^{(b)}(\widetilde{\Theta})\right]=2^{k+1}\cE
\end{eqnarray}
%Taking into account that
%$$
%\cE=e_{\vec{s}}\big(\widetilde{\kappa},\widetilde{\Theta}\big):=\big(1+\widetilde{\kappa}_1^{-2}\big)
%e^{(\ra)}_{s_1}\big(\widetilde{\Theta}\big)+\big(1+\widetilde{\kappa}_2^{-1}\big)e^{(\rb)}_{s_2}\big(\widetilde{\Theta}\big).
%$$
and, therefore, $\forall k\geq k_0+1$
$$
\ln{\big(N_{k}N_{k-1}\big)}-2^{k+2}\cE\leq -2^{k_0+2}\cE\leq 2\ve^{-2}\cE.
$$
It yields together with (\ref{p7_local})
\begin{eqnarray}
\label{p8_local}
I_2&\leq&   \mathrm{c}\exp{\left\{-y-2\ve^{-2}\cE\right\}}.
\end{eqnarray}
 We get from
   (\ref{p40}), (\ref{p50_local}) and (\ref{p8_local}) for any $m=1,\ldots,N_{k_0}$
 \begin{eqnarray}
 \label{p80_local}
 \mathrm{P}\left\{\sup_{\theta\in\Theta^{(m)}}\Psi\left(\chi_{\theta}\right)\geq U\right\} \leq 2\mathrm{c}
\exp{\left\{-y(1+\e)^{-2}-\ve^{-2}\cE\right\}}.
\end{eqnarray}
The last bound is independent of $m$ and we have from (\ref{p00}) and  (\ref{p80_local})
\begin{eqnarray*}
% \label{p81_local}
 \mathrm{P}\left\{\sup_{\theta\in\widetilde{\Theta}}\Psi\left(\chi_{\theta}\right)\geq U_s\right\} \leq 2\mathrm{c}
N_{k_0}\exp{\left\{-y(1+\e)^{-2}-\ve^{-2}\cE\right\}}.
\end{eqnarray*}
It remains to note that  $(1-3\ve)^{2}=(1+\e)^{-1}$ and that, similarly to (\ref{p777_local}),
$$
\ln{\big(N_{k_0}\big)}\leq 2^{k_0}\cE\leq \ve^{2}\cE,
$$
 and we come to the first assertion of the proposition.

\paragraph{{\bf II}. Moment's bound}
We get for any $y\geq 1$
\begin{eqnarray}
\label{p001_local}
&&E:=\left(\sup_{\theta\in \widetilde{\Theta}}
\Psi\left(\chi_{\phi[\theta]}\right)-U\right)^q_+
%\nonumber\\*[2mm]&&
=q\int_{0}^{\infty}
x^{q-1}\bP\left\{\sup_{\theta\in\widetilde{\Theta}}\Psi\left(\chi_{\theta}\right)\geq U+x
\right\}{\rd x}\nonumber\\*[2mm]
&&=q\left[U\right]^q\int_{0}^{
\infty}v^{q-1}
\bP\left\{\sup_{\theta\in \widetilde{\Theta}}\Psi\left(\chi_{\theta}\right)\geq
(1+v) U \right\}{\rd v}.
\end{eqnarray}
 Note that $(1+v) U \geq  U_{\vec{s}}^{(\e)}\Big((1+v)y,\widetilde{\kappa},\widetilde{\Theta}\Big)$. Therefore, applying the first statement of the proposition, where
 $y$ is replaced by $vy$ we obtain from (\ref{p001_local})
\begin{eqnarray*}
%\label{p001_local}
E\leq 2\mathrm{c}
\Gamma(q+1)\left[(1+\e)^{2}y^{-1}U\right]^q\exp{\{-y/(1+\e)^{2}\}}.
\end{eqnarray*}
%It remains to note that for any $y\geq 1$
%$$
%y^{-1}U_s\leq \left(\overline{A}_{\Theta_1}+\overline{B}_{\Theta_2}\right)\Big(\tau_s\big[e_{s}(\Theta_1)+ e_{s}(\Theta_2)\big]+1\Big)
%$$
%and we come to the second assertion of the proposition.
\epr

\subsection{Proof of Proposition \ref{prop_uniform_local2}} We start with establishing some technical results used in the sequel.

\paragraph{Preliminaries}
$1^{0}.$ First we formulate  the simple consequence of Proposition \ref{prop_uniform_local1}.

Let $\Theta_1,\Theta_2$ be given subsets of $\Theta$.
For any  $\vec{s}=(s_1,s_2)\in\cS_{\ra,\rb}$ and  any $\kappa=(\kappa_1,\kappa_2)\in\bR^{2}_+\setminus \{0\}$ introduce the following quantity
$$
e_{\vec{s}}\big(\kappa,\Theta_1,\Theta_2\big)=e^{(\ra)}_{s_1}\big(\kappa_1,\Theta_1\big)+e^{(\rb)}_{s_2}\big(\kappa_2,\Theta_2\big).
$$
Put   any $\e>0$ and any $y\geq 0$
\begin{eqnarray*}
%\label{eq4:def_local}
U^{(\e)}_{\vec{s}}\big(y,\kappa,\Theta_1,\Theta_2\big)=\kappa_1\sqrt{2\big[1+\e^{-1}\big]^{2} e_{\vec{s}}\big(\kappa,\Theta_1,\Theta_2\big)+y}
+\kappa_2\Big(2\big[1+\e^{-1}\big]^{2}e_{\vec{s}}\big(\kappa,\Theta_1,\Theta_2\big)+y\Big).
\end{eqnarray*}

\begin{lemma}
\label{lemma_intersection-of-2-sets-local1}
Let   Assumptions \ref{ass:fixed_theta_local}-\ref{ass:parameter_local}  hold and let $\Theta_1,\Theta_2$ be given subsets of $\Theta$.
Let $\kappa$ be chosen such that $\kappa_1\geq \sup_{\theta\in\Theta_1}A(\theta)$ and $\kappa_2\geq \sup_{\theta\in\Theta_2}B(\theta)$.
  Then $\forall \vec{s}\in\cS_{\mathrm{a},\mathrm{b}}$, $\forall \e\in \big(0,\sqrt{2}-1\big]$ and $\forall y\geq 1$,
 \begin{eqnarray*}
\mathrm{P}\left\{\sup_{\theta\in\Theta_1\cap\Theta_2}\Psi\left(\chi_{\theta}\right)\geq
 U^{(\e)}_{\vec{s}}\big(y,\kappa,\Theta_1,\Theta_2\big)\right\} \leq
2\mathrm{c}\exp{\left\{-y/(1+\e)^{2}\right\}}.
\end{eqnarray*}
Moreover, for any $q\geq 1$, putting $C_{\e,q}=2\mathrm{c}\Gamma(q+1)(1+\e)^{2q}$, one has
$$
\mathrm{E}\left\{\sup_{\theta\in\Theta_1\cap\Theta_2}\Psi\left(\chi_{\theta}\right)-
U^{(\e)}_{\vec{s}}\big(y,\kappa,\Theta_1,\Theta_2\big)\right\}^{q}_+ \leq
C_{\e,q}\left[y^{-1}U^{(\e)}_{\vec{s}}\big(y,\kappa,\Theta_1,\Theta_2\big)\right]^q\exp{\left\{-\frac{y}{(1+\e)^{2}}\right\}}.
$$

\end{lemma}
To prove the lemma it suffices to note the following simple facts. In view of the assumptions imposed on $\kappa$ and obvious inclusions $\Theta_1\cap\Theta_2\subseteq\Theta_1$, $\Theta_1\cap\Theta_2\subseteq\Theta_2$ we have
 \begin{eqnarray*}
\kappa_1\geq  \sup_{\theta\in\Theta_1\cap\Theta_2}A(\theta),\quad \kappa_2\geq  \sup_{\theta\in\Theta_1\cap\Theta_2}B(\theta),\quad e_{\vec{s}}\big(\kappa,\Theta_1,\Theta_2\big)\geq e_{\vec{s}}\big(\kappa,\Theta_1\cap\Theta_2\big).
\end{eqnarray*}
It yields
$U^{(\e)}_{\vec{s}}\big(y,\kappa,\Theta_1,\Theta_2\big)\geq U^{(\e)}_{\vec{s}}\big(y,\kappa,\Theta_1\cap\Theta_2\big)$ and to get the assertion of the lemma we apply Proposition \ref{prop_uniform_local1} with $\widetilde{\Theta}=\Theta_1\cap\Theta_2$.

\smallskip

$2^{0}.$  Note that $\Psi\big(\chi_\bullet\big):\Omega\times\Theta\to\bR_+$ is obviously $\mathrm{P}$-a.s. continuous in $\mathrm{a}\vee\mb$ as a composition of two continuous mappings between metric spaces. Hence Corollary  \ref{cor:after-lem-measurability} is applicable with $\boldsymbol{\mT}=\Theta$,
$\md=\mathrm{a}\vee\mb$, $\left(\boldsymbol{\Omega},\boldsymbol{\mB},\boldsymbol{\mathrm{P}}\right)=\left(\Omega,\mB,\mathrm{P}\right)$
$\zeta(\mt,\cdot)=\Psi\big(\chi_\theta(\cdot)\big)$ and $g(\mt)$ is either $\mathrm{V}^{(z,\e)}_{\vec{s}}(\theta)$ or
$\mathrm{U}^{(z,\e,q)}_{\vec{s}}(\theta),\;\mt=\theta.$

\paragraph{Proof of the proposition}
Put $\delta_l=(1+\e)^{l},\; l\geq0,$ and introduce the following sets
\begin{eqnarray*}
&\Theta_{A}^{(l)}=\left\{\theta\in\Theta:\; \underline{A}\delta_{l-1}\leq A(\theta)\leq \underline{A}\delta_{l}\right\},\quad
 \Theta_{B}^{(l)}=\left\{\theta\in\Theta:\; \underline{B}\delta_{l-1}\leq B(\theta)\leq \underline{B}\delta_{l}\right\},\;l\in\bN^*.
\end{eqnarray*}
The idea is to apply Lemma \ref{lemma_intersection-of-2-sets-local1} with $\Theta_1=\Theta_{A}^{(j)}$ and $\Theta_2=\Theta_{B}^{(k)}$ for any given $ j,k\geq 1$. To do that we will need to bound from below
$\mathrm{V}^{(z,\e)}_{\vec{s}}(\theta)$ and
$\mathrm{U}^{(z,\e,q)}_{\vec{s}}(\theta)
$
on  $\Theta_{A}^{(j)}\cap\Theta_{B}^{(k)}$.
%Later on we will use the following notations.
%$$
%\Theta^{(j,k)}=\Theta_{A}^{(j)}\cap \Theta_{B}^{(k)},\quad  \overline{A}_j=\overline{A}_{\Theta_{A}^{(j)}},\quad \overline{B}_k=\overline{B}_{\Theta_{B}^{(k)}}.
%$$
We will consider only $j,k$ such that $\Theta_{A}^{(j)}\cap\Theta_{B}^{(k)}\neq\emptyset$ and supremum over empty set is assumed to be zero. Also we will accept the following agreement: if  $B\equiv 0,\; b\equiv 0$ then $k\equiv 0$ and $\Theta_{B}^{(0)}=\widehat{\Theta}$.

\paragraph{Probability bound} Let $ \theta\in \Theta_{A}^{(j)}\cap\Theta_{B}^{(k)}$ be fixed and  put  $u=\cA_\e(\theta)$, $v=\cB_\e(\theta)$.
Note that
\begin{eqnarray}
\label{eq1-01:proof-prop2-local}
&&e^{(\ra)}_{s_1}\Big(\underline{A}u,\Theta_A\big(\underline{A}u\big)\Big)\geq e^{(\ra)}_{s_1}\Big(\underline{A}u,\Theta_A\big( \underline{A}\delta_{j}\big)\Big)\geq e^{(\ra)}_{s_1}\Big(\underline{A}u,\Theta_A^{(j)}\Big)
\geq  e^{(\ra)}_{s_1}\Big(\underline{A}\delta_{j+1},\Theta_A^{(j)}\Big).
\end{eqnarray}
To get the first two inequalities in (\ref{eq1-01:proof-prop2-local}) we have used that  $\Theta_A^{(j)}\subseteq\Theta_A\big( \underline{A}\delta_{j}\big)\subseteq\Theta_A\big( \underline{A}u\big)$  in view of  $\delta_{j}\leq u$ since $\theta\in\Theta_A^{(j)}$. To get the last inequality in (\ref{eq1-01:proof-prop2-local}) we have used that  the entropy  is decreasing function of its argument and that $\delta_{j+1}\geq u$ since $\theta\in\Theta_A^{(j)}$.

By the same reasons
\begin{eqnarray}
\label{eq1-001:proof-prop2-local}
&& e^{(\rb)}_{s_2}\Big(\underline{B}v,\Theta_B\big(\underline{B}v\big)\Big)\geq e^{(\rb)}_{s_2}\Big(\underline{B}v,\Theta_B\big( \underline{B}\delta_{k}\big)\Big)\geq e^{(\rb)}_{s_2}\Big(\underline{B}v,\Theta_B^{(k)}\Big) \geq  e^{(\rb)}_{s_2}\Big(\underline{B}\delta_{k+1},\Theta_B^{(k)}\Big).
\end{eqnarray}
Taking into account that left hand sides in (\ref{eq1-01:proof-prop2-local})
and (\ref{eq1-001:proof-prop2-local}) are independent of $\theta$, whenever $ \theta\in \Theta_{A}^{(j)}\cap\Theta_{B}^{(k)}$, we deduce from   (\ref{eq:def-of-function-e}), (\ref{eq1-01:proof-prop2-local})
and (\ref{eq1-001:proof-prop2-local})
\begin{eqnarray}
\label{eq1-04:proof-prop2-local}
\cE\Big(\cA_\e(\theta),\cB_\e(\theta)\Big)
&:=&e^{(\ra)}_{s_1}\Big(\underline{A}u,\Theta_A\big(\underline{A}u\big)\Big)+e^{(\rb)}_{s_2}\Big(\underline{B}v,\Theta_B\big(\underline{B}v\big)\Big)
\nonumber
\\
&\geq& e^{(\ra)}_{s_1}\Big(\underline{A}\delta_{j+1},\Theta_A^{(j)}\big)\Big)+e^{(\rb)}_{s_2}\Big(\underline{B}\delta_{k+1},\Theta_B^{(k)}\Big)
%\nonumber\\
= e_{\vec{s}}\left(\kappa,\Theta_A^{(j)},\Theta_B^{(k)}\right),
\end{eqnarray}
for any  $ \theta\in \Theta^{(j,k)}$, where we put  $\kappa=\big(\underline{A}\delta_{j+1},\underline{B}\delta_{k+1}\big)$.

We obtain from (\ref{eq1-04:proof-prop2-local}), putting $y=(1+\e)^{2} \Big[z+\ell\big(\delta_{j}\big)+\ell\big(\delta_{k}\big)\Big]$,
\begin{eqnarray}
\label{p2_local}
&&\mathrm{V}^{(z,\e)}_{\vec{s}}(\theta)\geq U^{(\e)}_{\vec{s}}\Big(y,\kappa,\Theta_A^{(j)},\Theta_B^{(k)}\Big),\;\; \forall \theta\in\Theta_{A}^{(j)}\cap\Theta_{B}^{(k)}.
\end{eqnarray}
Here we have also used that obviously $(1+\e)^{2}A(\theta)\geq \underline{A}\delta_{j+1}=:\kappa_1$ and $(1+\e)^{2}B(\theta)\geq \underline{B}\delta_{k+1}=:\kappa_2$ for any $ \theta\in \Theta_{A}^{(j)}\cap\Theta_{B}^{(k)}$. Moreover, we have used that $2\big[1+\e^{-1}\big]^{2}\geq 4\e^{-2}$ for any $\e\in \big(0,\sqrt{2}-1\big]$.

Therefore,  we obtain $\forall  j,k\geq 1$  in view of (\ref{p2_local})
%denoting $ y=(1+\e)\Big[z+\ell\big(\delta_{j}\vee \delta_{k}\big)\Big]$,
\begin{eqnarray*}
\Psi_{j,k}^*(z):=\sup_{\theta\in\Theta_{A}^{(j)}\cap\Theta_{B}^{(k)}}\left\{\Psi\left(\chi_{\theta}\right)-
\widetilde{\mathrm{V}}^{(z,\ve)}_{\vec{s}}(\theta)\right\}
\leq\sup_{\theta\in\Theta_{A}^{(j)}\cap\Theta_{B}^{(k)}}\Psi\left(\chi_{\theta}\right)-
U^{(\e)}_{\vec{s}}\Big(y,\kappa,\Theta_A^{(j)},\Theta_B^{(k)}\Big).
\end{eqnarray*}
Since  $\kappa_1:=\underline{A}\delta_{j+1}\geq \sup_{\theta\in\Theta_A^{(j)}}A(\theta)$ and $\kappa_2:=\underline{B}\delta_{k+1}\geq \sup_{\theta\in\Theta_B^{(k)}}B(\theta)$,  Lemma  \ref{lemma_intersection-of-2-sets-local1} is applicable
 with $\Theta_1=\Theta_A^{(j)}$ and $\Theta_2=\Theta_B^{(k)}$. Thus, applying it we obtain  $\forall  j,k\geq 1$  and $\forall z\geq 1$
\begin{equation}
\label{p1234_local}
\mathrm{P}\left\{\Psi_{j,k}^*(z)\geq 0\right\}\leq
%\nonumber\\
2\mathrm{c}\exp{\left\{-z\right\}}w_jw_k,
\end{equation}
where we put
$
w_m=\Big[1+m\ln{(1+\ve)}\Big]^{-1}\Big[1+\ln{\left\{1+m\ln{(1+\ve)}\right\}}\Big]^{-2}.
$
Noting that
$$
\sum_{m=1}^{\infty}w_m\leq 1+\Big[\ln{\left\{1+\ln{(1+\e)}\right\}}\Big]^{-2},
$$
taking into account that the union of $\Big\{\Theta_{A}^{(j)}\cap\Theta_{B}^{(k)},\; j,k\geq 1,\Big\}$ covers
 $\widehat{\Theta}$ and    summing up the right hand side in (\ref{p1234_local})
%(\ref{p1_local})
over $j,k$, we arrive at the first assertion of the proposition.

\paragraph{Moment's bound}

Using the same arguments having led to (\ref{p2_local}) and taking into account (\ref{eq1-04:proof-prop2-local}) we obtain for any
 $ \theta\in \Theta^{(j,k)}$ and any $ q\geq 1$
\begin{equation*}
%\label{p99_local}
\mathrm{U}^{(z,\e,q)}_{\vec{s}}(\theta)\geq U^{(\e)}_{\vec{s}}\Big(y,\kappa,\Theta^{(j)}_{A},\Theta^{(k)}_{B}\Big),
\end{equation*}
where $\; y=(1+\e)^{2}\Big[
z+\tau_\e e_{\vec{s}}\left(\kappa,\Theta_A^{(j)},\Theta_B^{(k)}\right)+
(\e+q)\ln{\big(\delta_{j} \delta_{k}\big)}\Big],\;\tau_\e=2(2+\e)\e^{-1}$.

Thus,  applying the second assertion of Lemma \ref{lemma_intersection-of-2-sets-local1}, $\forall  j,k\geq 1$  and $\forall z\geq 1$
\begin{eqnarray}
\label{p11_local}
E_{j,k}&:=&\mathrm{E}\left(\sup_{\theta\in \Theta^{(j,k)}}
\Psi\left(\chi_{\theta}\right)-U^{(\e)}_{\vec{s}}\Big(y,\kappa,\Theta^{(j)}_{A},\Theta^{(k)}_{B}\Big)\right)^q_+
\nonumber\\*[2mm]
&\leq& 2\mathrm{c}\Gamma(q+1)\left[\frac{U^{(\e)}_{\vec{s}}\Big(y,\kappa,\Theta^{(j)}_{A},\Theta^{(k)}_{B}\Big)}{
z+\tau_\e e_{\vec{s}}\left(\kappa,\Theta_A^{(j)},\Theta_B^{(k)}\right)+(\e+q)\ln{\big(\delta_{j} \delta_{k}\big)}}\right]^q\big(\delta_{j}\delta_{k}\big)^{-\e-q}\exp{\left\{-z\right\}}.
\end{eqnarray}
Putting for brevity $e_{j,k}=e_{\vec{s}}\left(\kappa,\Theta_A^{(j)},\Theta_B^{(k)}\right)$ and noting that for any $j,k\geq 1$ and $z\geq 1$
\begin{eqnarray*}
U^{(\e)}_{\vec{s}}\Big(y,\kappa,\Theta^{(j)}_{A},\Theta^{(k)}_{B}\Big)
\leq 2\big[\kappa_1\vee\kappa_2\big]
\left[4\e^{-2}e_{j,k}+y\right]
\leq 2\big[\underline{A}\vee\underline{B}\big]\big[\delta_{j+1}\vee\delta_{k+1}\big]\big[2(1+\e^{-1})^{2}e_{j,k}+y\big],
\end{eqnarray*}
and that for any $\e\in \big(0,\sqrt{2}-1\big]$
\begin{eqnarray*}
%\label{p112_local}
\frac{2(1+\e^{-1})^{2}e_{j,k}+(1+\e)^{2}\Big[
z+\tau_\e e_{j,k}+
(\e+q)\ln{\big(\delta_{j} \delta_{k}\big)}\Big]}{z+\tau_\e e_{j,k}+(\e+q)\ln{\big(\delta_{j}\vee \delta_{k}\big)}}
\leq (1+\e)^{4}(2\e)^{-1},
%\leq (1+\e)^2\left[1+4\e^{-2}(1+\e)^{-2}\tau_\e^{-1}\right],
\end{eqnarray*}
we get finally from (\ref{p11_local}) for  $j,k\geq 1$
\begin{eqnarray}
\label{p1222_local}
E_{j,k}
 \leq 2\mathrm{c}\Gamma(q+1)\big[\e^{-1}(1+\e)^{5}\big]^{q}\big[\underline{A}\vee\underline{B}\big]^q\big(\delta_{j} \delta_{k}\big)^{-\e}\exp{\left\{-z\right\}}.
\end{eqnarray}
Here we have used that $(\delta_{j}\vee \delta_{k})/(\delta_{j}\delta_{k})\leq 1$ since $\delta_{j},\delta_{k}\geq 1$.
Taking into account that
\begin{equation*}
%\label{p10_local}
\mathrm{E}\left(\sup_{\theta\in\widetilde{\Theta}}\left\{\Psi\left(\chi_{\theta}\right)-\widetilde{\mathrm{U}}^{(z,\ve,q)}_s(\theta)
\right\}\right)^q_+
%\nonumber\\*[2mm]&&
\leq
\sum_{j=1}^{\infty}\sum_{k=1}^{\infty}E_{j,k}
\end{equation*}
and that $\sum_{j=1}^{\infty}\sum_{k=1}^{\infty}\big(\delta_{j}\delta_{k}\big)^{-\e}=\big[(1+\e)^{\e}-1\big]^{-2}$,
we come to the second assertion of the proposition in view of (\ref{p1222_local}), where one can replace $(1+\e)^{5}$ par $2^{5/2}$ since
$\e\in\big(0,\sqrt{2}-1\big]$. We have also used that $\big[(1+\e)^{\e}-1\big]^{-2}\leq 2\e^{-4}$ since $\e\in\big(0,\sqrt{2}-1\big]$.

\epr

\subsection{Proof of Proposition \ref{prop_uniform_local3}}

 First we discuss the measurability issue.  Note that $\Psi\big(\chi_\bullet\big):\Omega\times\Theta\to\bR_+$ is obviously $\mathrm{P}$-a.s. continuous in $\mathrm{a}\vee\mb$ as a composition of two continuous mappings between metric spaces. Hence Lemma \ref{lem:measurability} is applicable with $\boldsymbol{\mT}=\Theta$,
$\md=\mathrm{a}\vee\mb$, $\left(\boldsymbol{\Omega},\boldsymbol{\mB},\boldsymbol{\mathrm{P}}\right)=\left(\Omega,\mB,\mathrm{P}\right)$
$\zeta(\mt,\cdot)=\Psi\big(\chi_\theta(\cdot)\big),\;\mt=\theta$, $\mZ=\mA$, $\mT_\z=\Theta_\alpha$  and $g(\mz)=\widehat{\mathrm{U}}^{(z,\e,r)}(\alpha),\;\mz=\alpha.$

\vskip0.1cm

The proof of the proposition is similar to whose of Proposition \ref{prop_uniform_local2}. Let $\e\in (0,1)$ be fixed and put $\delta_l=(1+\e)^{-l},\; l\geq 0,\;\; \delta_l=1,\; l<0$. Introduce the following sets:
\begin{eqnarray*}
\mA_{1}^{(l)}=\left\{\alpha\in\mA:\; \overline{\tau}_1\delta_{l+1}\leq \tau_1(\alpha)\leq \overline{\tau}_1\delta_{l}\right\},\quad
\mA_{2}^{(l)}=\left\{\alpha\in\mA:\; \overline{\tau}_2\delta_{l+1}\leq \tau_2(\alpha)\leq \overline{\tau}_2\delta_{l}\right\},\;\;l\in\bN.
\end{eqnarray*}
%The idea is to apply Proposition \ref{prop_uniform_local1} with $\Theta_1=\Theta_{A}^{(j)}$ and $\Theta_2=\Theta_{B}^{(k)}$ for any given $ j,k\geq 0$.
Fix  $j,k\geq 0$, $ \alpha\in \mA_{1}^{(j)}\cap\mA_{2}^{(k)}$ and  put   $u=(1+\e)\tau_1(\alpha)$, $v=(1+\e)\tau_2(\alpha)$.
Note that
\begin{eqnarray}
\label{eq1-01:proof-prop3-local}
e^{(\ra)}_{s_1(u,\cdot)}\Big(\lambda_1^{-1}g_A(u),\Theta^{\prime}_1(u)\Big)&\geq& e^{(\ra)}_{s_1(u,\cdot)}\Big(\lambda_1^{-1}g_A(u),\Theta^{\prime}_1\big( \overline{\tau}_1\delta_{j}\big)\Big)
\nonumber\\
&\geq& e^{(\ra)}_{s_1(u,\cdot)}\Big(\lambda_1^{-1}g_A\big(\overline{\tau}_1\delta_{j-1}\big),\Theta^{\prime}_1\big( \overline{\tau}_1\delta_{j}\big)\Big).
%\geq  e^{(\ra)}_{s_1(u,\cdot)}\Big(\underline{A}\delta_{j+1},\Theta_A^{(j)}\Big).
\end{eqnarray}
To get the first  inequality in (\ref{eq1-01:proof-prop3-local}) we have used that  $\Theta^{\prime}_1\big( \overline{\tau}_1\delta_{j}\big)\subseteq\Theta^{\prime}_1(u)$  in view of  $\overline{\tau}_1\delta_{j}\leq u$ since $\alpha\in \mA_{1}^{(j)}$. To get the second inequality in (\ref{eq1-01:proof-prop3-local}) we have used that $g_A\big(\overline{\tau}_1\delta_{j-1}\big)\geq g_A(u)$, since $g_A$ is increasing  and $\overline{\tau}_1\delta_{j-1}\geq u$ for  $\alpha\in \mA_{1}^{(j)}$. Moreover we have used that
 the entropy  is decreasing function of its argument.

Remembering the definition of $e^{(\ra)}_{s_1(u,\cdot)}\Big(\lambda_1^{-1}g_A\big(\overline{\tau}_1\delta_{j-1}\big),\Theta^{\prime}_1\big( \overline{\tau}_1\delta_{j}\big)\Big)$, see (\ref{eq3:def_local}), we have
\begin{eqnarray}
\label{eq1-0101:proof-prop3-local}
&&\;\;e^{(\ra)}_{s_1(u,\cdot)}\Big(\lambda_1^{-1}g_A\big(\overline{\tau}_1\delta_{j-1}\big),\Theta^{\prime}_1\big( \overline{\tau}_1\delta_{j}\big)\Big)=\sup_{\delta>0}\delta^{-2}
\mathfrak{E}_{\Theta^{\prime}_1\big( \overline{\tau}_1\delta_{j}\big),\;\mathrm{a}}\left(\lambda_1^{-1}g_A\big(\overline{\tau}_1\delta_{j-1}\big)(48\delta)^{-1}s_1(u,\delta)\right)
\nonumber\\
&&\;\;=\sup_{\delta>0}\delta^{-2}
\mathfrak{E}_{\Theta^{\prime}_1\big( \overline{\tau}_1\delta_{j}\big),\;\mathrm{a}}\left(\lambda_1^{-1}g_A\big(\overline{\tau}_1\delta_{j-1}\big)
(48\delta)^{-1}s_1\left(\overline{\tau}_1\delta_{j},\delta\right)
\left[\frac{s_1(u,\delta)}{s_1\left(\overline{\tau}_1\delta_{j},\delta\right)}\right]\right)
\nonumber\\
&&\;\;\geq\sup_{\delta>0}\delta^{-2}
\mathfrak{E}_{\Theta^{\prime}_1\big( \overline{\tau}_1\delta_{j}\big),\;\mathrm{a}}\left(g_A\big(\overline{\tau}_1\delta_{j-1}\big)(48\delta)^{-1}s_1\left(\overline{\tau}_1\delta_{j},
\delta\right)
\right)
=:e^{(\ra)}_{s_1\big(\overline{\tau}_1\delta_{j},\cdot\big)}\Big(g_A\big(\overline{\tau}_1\delta_{j-1}\big),\Theta^{\prime}_1
\big(\overline{\tau}_1\delta_{j}\big)\Big).
\end{eqnarray}
To get (\ref{eq1-0101:proof-prop3-local}) we have used that $1\leq u/\overline{\tau}_1\delta_{j}\leq 1+\e\leq\sqrt{2}$, the definition of $\lambda_1$
and, as previously, that the entropy is decreasing function of its argument. We obtain from (\ref{eq1-01:proof-prop3-local}) and (\ref{eq1-0101:proof-prop3-local})
\begin{eqnarray}
\label{eq1-0102:proof-prop3-local}
e^{(\ra)}_{s_1(u,\cdot)}\Big(\lambda_1^{-1}g_A(u),\Theta^{\prime}_1(u)\Big)\geq
e^{(\ra)}_{s_1\big(\overline{\tau}_1\delta_{j},\cdot\big)}\Big(g_A\big(\overline{\tau}_1\delta_{j-1}\big),\Theta^{\prime}_1
\big(\overline{\tau}_1\delta_{j}\big)\Big).
\end{eqnarray}
By the same reasons we have
\begin{eqnarray}
\label{eq1-0103:proof-prop3-local}
e^{(\rb)}_{s_2(v,\cdot)}\Big(\lambda_2^{-1}g_B(v),\Theta^{\prime}_2(v)\Big)\geq
e^{(\rb)}_{s_1\big(\overline{\tau}_2\delta_{k},\cdot\big)}\Big(g_B\big(\overline{\tau}_2\delta_{k-1}\big),\Theta^{\prime}_2
\big(\overline{\tau}_2\delta_{k}\big)\Big),
\end{eqnarray}
and, we get from (\ref{eq1-0102:proof-prop3-local}) and (\ref{eq1-0103:proof-prop3-local}) for any $ \alpha\in \mA_{1}^{(j)}\cap\mA_{2}^{(k)}$
\begin{eqnarray}
\label{eq1-04:proof-prop3-local}
\widehat{\cE}^{(\e)}(\alpha)&:=&\cE^\prime(u,v)\geq e^{(\ra)}_{s_1\big(\overline{\tau}_1\delta_{j},\cdot\big)}\Big(g_A\big(\overline{\tau}_1\delta_{j-1}\big),\Theta^{\prime}_1
\big(\overline{\tau}_1\delta_{j}\big)\Big)
+e^{(\rb)}_{s_1\big(\overline{\tau}_2\delta_{k},\cdot\big)}\Big(g_B\big(\overline{\tau}_2\delta_{k-1}\big),\Theta^{\prime}_2
\big(\overline{\tau}_2\delta_{k}\big)\Big)
\nonumber\\
&=&
e_{\vec{s}}\Big(\kappa,\Theta^{\prime}_1\big(\overline{\tau}_1\delta_{j}\big),\Theta^{\prime}_2\big(\overline{\tau}_2\delta_{k}\big)\Big),
\end{eqnarray}
where we have put $\kappa=\Big(g_A\big(\overline{\tau}_1\delta_{j-1}\big),g_B\big(\overline{\tau}_2\delta_{k-1}\big)\Big)$ and $\vec{s}=\Big(s_1\big(\overline{\tau}_1\delta_{j},\cdot\big),s_2\big(\overline{\tau}_2\delta_{k},\cdot\big)\Big)$.

\par

Note also that for any $ \alpha\in \mA_{1}^{(j)}\cap\mA_{2}^{(k)}$ in view of monotonicity of functions $g_A$ and $g_B$
\begin{eqnarray}
\label{eq1-0400:proof-prop3-local}
&&g_A\left([1+\e]^{2}\tau_1(\alpha)\right)\geq g_A\big(\overline{\tau}_1\delta_{j-1}\big)=:\kappa_1,\quad g_B\left([1+\e]^{2}\tau_2(\alpha)\right)\geq g_A\big(\overline{\tau}_2\delta_{k-1}\big)=:\kappa_2.
\end{eqnarray}
Moreover,  the definition of sets $\Theta^\prime_
1(\cdot)$ and $\Theta^\prime_2(\cdot)$ implies that for any $\alpha\in\mA_{1}^{(j)}\cap\mA_{2}^{(k)}$
$$
\Theta_\alpha\subseteq \Theta^{\prime}_1\big(\overline{\tau}_1\delta_{j}\big)\cap\Theta^{\prime}_2\big(\overline{\tau}_2\delta_{k}\big)
$$
and, therefore, for any $\alpha\in\mA_{1}^{(j)}\cap\mA_{2}^{(k)}$
\begin{eqnarray}
\label{p3_local3}
\sup_{\theta\in\Theta_\alpha}\Psi\left(\chi_{\theta}\right)\leq \sup_{\theta\in\Theta^{\prime}_1\big(\overline{\tau}_1\delta_{j}\big)
\cap\Theta^{\prime}_2\big(\overline{\tau}_2\delta_{k}\big)}\Psi\left(\chi_{\theta}\right).
\end{eqnarray}

\paragraph{ Probability bound}
We get from  (\ref{eq1-04:proof-prop3-local}) and (\ref{eq1-0400:proof-prop3-local}) for any $ \alpha\in \mA_{1}^{(j)}\cap\mA_{2}^{(k)}$
\begin{eqnarray*}
%\label{p2_local3}
&&\widehat{\mathrm{U}}^{(z,\e,0)}(\alpha)\geq U^{(\e)}_{\vec{s}}\Big(y,\kappa,\Theta^{\prime}_1\big(\overline{\tau}_1\delta_{j}\big),\Theta^{\prime}_2
\big(\overline{\tau}_2\delta_{k}\big)\Big),\;\; y=(1+\e)^{2}
\Big[z+R_0\big(\overline{\tau}_1\delta_{j},\overline{\tau}_2\delta_{k}\big)\Big],
\end{eqnarray*}
where we have used that $R_0$ is increasing, in both arguments, function.

It yields together with  (\ref{p3_local3}) $\forall  j,k\geq 0$
%denoting $ y=(1+\e)\Big[z+\ell\big(\delta_{j}\vee \delta_{k}\big)\Big]$,
\begin{eqnarray*}
\Psi_{j,k}^{(V)}(z)&:=&\sup_{\alpha\in\mA_{1}^{(j)}\cap\mA_{2}^{(k)}}\left(\sup_{\theta\in\Theta_\alpha}\left\{\Psi\left(\chi_{\theta}\right)-
\widehat{\mathrm{U}}^{(z,\ve,0)}(\alpha)\right\}\right)\\*[2mm]
&\leq&\sup_{\theta\in\Theta_1^{\prime}\big(\overline{\tau}_1 \delta_{j}\big)\cap\Theta_2^{\prime}\big(\overline{\tau}_2 \delta_{k}\big)\Big)}\Psi\left(\chi_{\theta}\right)-
U^{(\e)}_{\vec{s}}\Big(y,\kappa,\Theta^{\prime}_1\big(\overline{\tau}_1\delta_{j}\big),\Theta^{\prime}_2\big(\overline{\tau}_2\delta_{k}\big)\Big).
\end{eqnarray*}
Let us remark that the definition of the sets $\Theta_A^\prime(\cdot)$, $\Theta_B^\prime(\cdot)$ and the functions  $g_A$ and $g_B$ as well as their monotonicity   imply that
\begin{eqnarray}
\label{p2_local3}
&&\kappa_1:=g_A\big(\overline{\tau}_1\delta_{j-1}\big)> g_A\big(\overline{\tau}_1\delta_{j}\big)\geq g^*_A\big(\overline{\tau}_1\delta_{j}\big) =:\sup_{\theta\in\Theta_1^{\prime}\big(\overline{\tau}_1 \delta_{j}\big)}A(\theta);
\\
\label{p222_local3}
&&\kappa_2:=g_B\big(\overline{\tau}_2\delta_{k-1}\big)> g_B\big(\overline{\tau}_2\delta_{k}\big)\geq g^*_B\big(\overline{\tau}_2\delta_{k}\big) =:\sup_{\theta\in\Theta_2^{\prime}\big(\overline{\tau}_2 \delta_{k}\big)}B(\theta),
\end{eqnarray}
and, therefore, Lemma \ref{lemma_intersection-of-2-sets-local1} is applicable with $\Theta_1=\Theta_1^{\prime}\big(\overline{\tau}_1 \delta_{j}\big)$ and $\Theta_2=\Theta_2^{\prime}\big(\overline{\tau}_2 \delta_{k}\big)$.

Thus,  applying the first assertion of Lemma \ref{lemma_intersection-of-2-sets-local1}, we obtain $\forall  j,k\geq 0$  and $\forall z\geq 1$
\begin{equation*}
%\label{p1_local}
\mathrm{P}\left\{\Psi_{j,k}^{(V)}(z)\geq 0\right\}\leq 2\mathrm{c} \exp{\left\{-z-R_0\big(\overline{\tau}_1\delta_{j},\overline{\tau}_2\delta_{k}\big)\right\}}.
%\nonumber\\2\mathrm{c}\exp{\left\{-z\right\}}w_jw_k,
\end{equation*}
%where, as previously,
%$w_m=\Big[1+m\ln{(1+\ve)}\Big]^{-1}\Big[1+\ln{\left\{1+m\ln{(1+\ve)}\right\}}\Big]^{-2}.$

Noting that the union of $\Big\{\mA_{1}^{(j)}\cap\mA_{2}^{(k)},\; j=\overline{0,J},k=\overline{0,K},\Big\}$ covers
 $\mA$,    summing up the right hand side in the last inequality
%(\ref{p1_local})
over $j,k$
%and repeating the calculations made in the proof of the first assertion of Proposition \ref{prop_uniform_local2},
we come to the first statement of the proposition.

\paragraph{Moment's bound}
We get from  (\ref{eq1-04:proof-prop3-local}) and (\ref{eq1-0400:proof-prop3-local}) for any $ \alpha\in \mA_{1}^{(j)}\cap\mA_{2}^{(k)}$
\begin{eqnarray*}
%\label{p2_local3}
&&\widehat{U}^{(z,\e,q)}(\alpha)\geq U^{(\e)}_{\vec{s}}\Big(y,\kappa,\Theta^{\prime}_1\big(\overline{\tau}_1\delta_{j}\big),\Theta^{\prime}_2\big(\overline{\tau}_2\delta_{k}\big)\Big),
\end{eqnarray*}
where $\; y=(1+\e)^{2}\Big[
z+\tau_\e e_{\vec{s}}\Big(\kappa,\Theta^{\prime}_1\big(\overline{\tau}_1\delta_{j}\big),\Theta^{\prime}_2
\big(\overline{\tau}_2\delta_{k}\big)\Big)+R_q\big(\overline{\tau}_1\delta_{j},\overline{\tau}_2\delta_{k}\big)\Big],\;\tau_\e=2(2+\e)\e^{-1}$.

It yields together with  (\ref{p3_local3}) $\forall  j,k\geq 0$
%denoting $ y=(1+\e)\Big[z+\ell\big(\delta_{j}\vee \delta_{k}\big)\Big]$,
\begin{eqnarray*}
\Phi_{j,k}^{(U)}(z)&:=&\left\{\sup_{\alpha\in\mA_{A}^{(j)}\cap\mA_{B}^{(k)}}
\left(\sup_{\theta\in\Theta_\alpha}\left\{\Psi\left(\chi_{\theta}\right)-
\widehat{U}^{(z,\ve,q)}(\alpha)\right\}\right)\right\}_+\\*[2mm]
&\leq&\Bigg\{\sup_{\theta\in\Theta_1^{\prime}\big(\overline{\tau}_1 \delta_{j}\big)\cap\Theta_2^{\prime}\big(\overline{\tau}_2 \delta_{k}\big)\Big)}\Psi\left(\chi_{\theta}\right)-
U^{(\e)}_{\vec{s}}\Big(y,\kappa,\Theta^{\prime}_1\big(\overline{\tau}_1\delta_{j}\big),\Theta^{\prime}_2\big(\overline{\tau}_2
\delta_{k}\big)\Big)\Bigg\}_+.
\end{eqnarray*}
Taking into account (\ref{p2_local3}), (\ref{p222_local3})  and applying the second assertion of Lemma \ref{lemma_intersection-of-2-sets-local1}, we have, analogously to (\ref{p11_local}), $\forall  j,k\geq 0$  and $\forall z\geq 1$
\begin{eqnarray*}
&&\mathrm{E}\left(\Phi_{j,k}^{(U)}(z)\right)
\nonumber\\*[2mm]
&&\leq 2\mathrm{c}\Gamma(q+1)\left[\frac{U^{(\e)}_{\vec{s}}\Big(y,\kappa,\Theta^{\prime}_1\big(\overline{\tau}_1\delta_{j}\big),\Theta^{\prime}_2
\big(\overline{\tau}_2\delta_{k}\big)\Big)}{
z+\tau_\e e_{\vec{s}}\Big(\kappa,\Theta^{\prime}_1\big(\overline{\tau}_1\delta_{j}\big),\Theta^{\prime}_2\big(\overline{\tau}_2\delta_{k}\big)\Big)
+R_q\big(\overline{\tau}_1\delta_{j},\overline{\tau}_2\delta_{k}\big)}\right]^q
\exp{\left\{-z-R_q\big(\overline{\tau}_1\delta_{j},\overline{\tau}_2\delta_{k}\big)\right\}}.
\end{eqnarray*}
Putting for brevity $e_{j,k}=e_{\vec{s}}\left(\kappa,\Theta^{\prime}_1\big(\overline{\tau}_1\delta_{j}\big),\Theta^{\prime}_2\big(\overline{\tau}_2\delta_{k}\big)\right)$ we note that for any $j,k\geq 0$ and $z\geq 1$
\begin{eqnarray*}
&&U^{(\e)}_{\vec{s}}\Big(y,\kappa,\Theta^{\prime}_1\big(\overline{\tau}_1\delta_{j}\big),\Theta^{\prime}_2
\big(\overline{\tau}_2\delta_{k}\big)\Big)\leq 2\big[g_A\big(\overline{\tau}_1\delta_j\big)\vee g_B\big(\overline{\tau}_2\delta_k\big)\big]\left[4\e^{-2}e_{j,k}+y\right]
\end{eqnarray*}
Repeating the computation done after (\ref{p11_local}) we come to the second assertion of the proposition.

\epr

\section{Proof of Theorem \ref{th:gauss-norm}}
\label{sec:proofs-of-gauss-theorems}

 Below $c_1,c_2\ldots,$ denote the constants completely determined by $d,p,\mu$ and $\gamma$.
 We break  the proof on several steps.

\smallskip

$1^{0}$.
\iffalse
In view of Kolmogorov criterion \cite{GihScor} the trajectories of $\big\|\xi_h\big\|_p$ is stochastically continuous on $\cH$. Thus, supremum
over $\cH$ can be replaced by supremum over $\widetilde{\cH}$, where  $\widetilde{\cH}:=\left\{\tilde{h}_1,\ldots,\tilde{h}_i,\ldots\right\}$ is set of rational numbers belonging to $\cH$.
\fi
Let $\bB_q,\;1<q<\infty,$ denote the set of functions vanishing outside $\bK_\mu$ and  whose $\bL_q$-norm is less or equal to 1. Later on
the integration is always understood as the integration over $\bR^{d}$.
%the unit ball in $\bL_q\big(\bK_\mu\big)$ and for any $h\in\cH$

Consider the set of functions
$$
\Theta_h=\left\{\theta:\bR^d\to\bR:\; \theta(x)=h^{-d}\int K_h\left(t-x\right)\ell(t)\rd t,\;\;\ell\in\bB_{\frac{p}{p-1}},\;h\in\cH\right\}.
$$
Put also $\Theta=\cup_{h\in\cH}\Theta_h$ and for any $\theta\in\Theta$ introduce the gaussian process
\begin{equation}
\label{eq:def-of-chi-theta-gauss_local}
\chi_\theta=\int_{\bR^d}\theta(x)b(\rd x).
\end{equation}
In view of Young inequality \cite{folland}, Theorems 6.18 and 6.36, for any $\theta\in\Theta$
\begin{equation}
\label{eq:L2-norm-of-theta-gauss_local}
\|\theta\|_2:=\left(\int_{\bR^d}\big|\theta(x)\big|^{2}\rd x\right)^{\frac{1}{2}}\leq h^{\frac{d(2-p)}{2p}} \|K\|_{\frac{2p}{p+2}}
\leq h^{\frac{d(2-p)}{2p}}\leq \big(h^{(\min)}\big)^{\frac{d(2-p)}{2p}} <\infty,
\end{equation}
 since $2\leq p<\infty$.  Here we have also used that $\|K\|_q\leq \|K\|_\infty\leq 1 $ in view of assumption imposed on the function $K$.

Thus, the stochastic integral in (\ref{eq:def-of-chi-theta-gauss_local}) is well-defined and $\chi_\theta$ is zero-mean gaussian random function on $\Theta$ such that
\begin{equation}
\label{eq:def-of-variance-semimetric-gauss_local}
V(\theta)=\|\theta\|_2,\quad \rho(\theta_1,\theta_2)=\|\theta_1-\theta_2\|_2.
%\;\;\theta\in\Theta_{\widetilde{\cH}}.
\end{equation}
We equip $\Theta$ with the  semi-metric $\rho$ and in the next paragraph we compute the entropy of several subsets of $\Theta$. This computation allows us, in particular, to assert  that Dudley integral is finite on $\Theta$.
It yields \cite{Lif} that $\chi_\bullet$ is $\mathrm{P}$-a.s uniformly continuous on  $(\Theta,\rho)$, therefore Assumption \ref{ass:metric-case-local} holds. Moreover, we show that Assumption \ref{ass:parameter_local} is fulfilled.

We conclude that Proposition  \ref{prop_uniform_local3} are applicable to $\chi_\theta,$ on $\Theta$, with $\Theta_\alpha=\Theta_h$,
$\mA=\cH$, $\ra=\sqrt{2}\rho$, $A=\sqrt{2}V$. It remains to note that in view of duality arguments for any $h\in\cH$
$$
\big\|\xi_h\big\|_p=\sup_{\ell\in\bB_{\frac{p}{p-1}}}\int\ell(t)\xi_{h}(t)\rd t,
$$
and, therefore, for any $h\in\cH$
\begin{equation}
\label{eq:duality-gauss-local}
\|\xi_h\big\|_p=\sup_{\theta\in\Theta_h}\chi_{\theta}.
\end{equation}

\smallskip

$2^{0}$.  In order to apply Proposition  \ref{prop_uniform_local3} we need to compute several quantities. First, we have to choose the function
$\tau_1$ (since $B,\rb\equiv 0$, hence $\tau_2\equiv 0$).  Set $\tau_1(h)=h^{-1}$ and note that for any $u>0$
$$
\Theta^\prime_1(u):=\bigcup_{h:\;\tau_1(h)\leq u}\Theta_h=\bigcup_{h\geq u^{-1} }\Theta_h,
$$
We note that
computations given in (\ref{eq:L2-norm-of-theta-gauss_local}) yield
$$
g^*_A(u):=\sup_{\theta\in\Theta^{\prime}_1(u)}A(\theta)\leq \sqrt{2}u^{\frac{d(p-2)}{2p}}=:g_A(u),\;\;u\geq \left(h^{(\max)}\right)^{-1}.
$$
%where $\nu(u)=\big(\sqrt{2}u^{-1}\big)^{\frac{2p}{d(p-2)}}$and
%Define the function
%\begin{eqnarray}
%\label{eq:def-of-function-prime-e-gauss}
%&&\mathcal{E}^{\prime}(u)=e_{\ms(u,\cdot)}\Big(u,\Theta^\prime_A(u)\Big),\;\; ,
%\end{eqnarray}
%where $\ms(\cdot,\cdot)$ will be specified later.

\smallskip

$3^{0}$. Let $\mE^{(u)}(\delta),\;\delta>0,$ be the entropy of $\Theta^\prime_A(u)$ computed with respect to   semi-metric $\ra=\sqrt{2}\rho$, where, remind $\rho(\cdot)=\|\cdot\|_2$.
The following assertion is true: there exist $c_0$ completely determined by $\gamma, d$, $p$ and $\mu$ such that for any $\forall\bar{\gamma}\in\big(d/2-d/p,\gamma\big]$ and for any $u\geq\left(h^{(\max)}\right)^{-1}$
\begin{equation}
\label{eq2:proof-th-norm-gauss-local}
\mE^{(u)}(\delta)\leq c_1 u^{d}\;\delta^{-d/\bar{\gamma}},
\;\;\forall\; \delta>0.
\end{equation}
 where $c_1=c_03^{4d}d^2$.

The proof of (\ref{eq2:proof-th-norm-gauss-local}) is obtained by routine computations and it is postponed to the step  $8^{0}$.

\smallskip

$4^{0}$. Choosing $u=\left(h^{(\min)}\right)^{-1}$  (that yields $\Theta_1^\prime(u)=\Theta$) and  $\bar{\gamma}=\gamma$ we get for any $\delta>0$
$$
\mE_{\Theta,\ra}(\delta)\leq c_1\left(h^{(\min)}\right)^{-d}\left(\frac{1}{\delta}\right)^{d/\gamma}.
%\cL(v)=\left[3^{d}\Big([4\mu\sqrt{d}]\vee 1\Big)d^{\alpha/2}\right]\left[\frac{1}{V^{(\min)}}\right].
$$
In view the condition $\gamma>d/2$, Dudley integral is finite on $\Theta$ and, therefore,
 $\chi_\bullet$ is $\mathrm{P}$-a.s uniformly continuous on  $(\Theta,\rho)$.
 It complete the verification of  Assumption \ref{ass:metric-case-local}.

The last inequality shows also that there exist $\tau>0$ such that for any $s\in\bS$, satisfying $\sup_{\delta> 0}\delta^{\tau}s^{-1}(\delta)<\infty$, Assumption \ref{ass:parameter_local} is fulfilled.

\smallskip

$5^{0}$. Let us now choose the function $\ms$. Set $\delta(u)=u^{d/2-d/p}$ and for any $u\geq 1 $ let $m(u)\in\bN$ be such that $2^{m(u)}\leq \delta(u)< 2^{m(u)+1}$. Define
$$
\ms(u,\delta)=\big(3\big/4\pi^2\big)\left(1+\Big[\log_2{\big(2^{-m(u)}\delta\big)
\Big]^{2}}\right)^{-1}.
$$
We remark that
$$
\ms\big(u,2^{k/2}\big)=(3\big/4\pi^2)\left[1+\big((k/2)-m(u)\big)^{2}\right]^{-1}
$$
and, therefore, $\ms(u,\cdot)\in\bS$ for any $u\geq 1$. Moreover, we note that if $p=2$ then $\ms$ does not depend on $u$ and it is given by
$$
\ms(\delta)=\big(3\big/4\pi^2\big)\left(1+\big[\log_2{\delta}
\big]^{2}\right)^{-1}.
$$
Obviously the factor $3/4\pi^2$ can be replaced here by $6/\pi^2$.

Now, let us compute the quantity $\lambda_1$ related to the function $\ms$. Remind that
$$
\lambda_1:=\sup_{t\in\big[1,\sqrt{2}\big]}\sup_{x>\underline{\tau}_1}\sup_{\delta>0}\frac{\ms(xt,\delta)}{\ms(x,\delta)},
$$
where $\underline{\tau}_1=\left(h^{(\max)}\right)^{-1}\geq 1$.
It is evident that
\begin{eqnarray*}
\lambda_1&=&\sup_{m\geq 0}\left\{1\bigvee\sup_{\delta>0}\left[\frac{1+\Big(\log_2{\big(2^{-m}\delta\big)}\Big)^{2}}
{1+\Big(\log_2{\big(2^{-m-1}\delta\big)}\Big)^{2}}\right]\right\}=
1\bigvee\sup_{x>0}\left[\frac{1+\big(\log_2{(x)}\big)^{2}}
{1+\big(\log_2{(x/2)}\big)^{2}}\right]
\\*[2mm]
&\leq& 1+\sup_{x>0}\left|\frac{2\log_2{(x)}-1}{1+\big(\log_2{(x)}-1\big)^{2}}\right|
=1+\sup_{y\in\bR}\left|\frac{2y+1}{1+y^{2}}\right|< \;3.
\end{eqnarray*}

\smallskip

$6^{0}$.
Define
$$
\bar{\gamma}(\delta)=\left\{
\begin{array}{cc}
d/2-d/(2p),\quad &0< \delta< \delta(u);\\
\gamma,\quad &\delta\geq \delta(u),
\end{array}\right.
$$
and note that $(d/2-d/p)<d/4<\gamma(\delta)\leq \gamma$ for any $\delta>0$.

 Putting $c_2=c_1 4(144)^{4}$ we get for any $\delta>0$ in view of (\ref{eq2:proof-th-norm-gauss-local})
\begin{eqnarray}
\label{eq5:proof-th-norm-gauss-local}
&&\mE^{(u)}\Big(\lambda_1^{-1}\big[g_A(u)\ms(u,\delta)\big]/48\delta\Big)\leq
c_2 u^{d}\;\Big[u^{\frac{d(p-2)}{2p}}\ms(u,\delta)\Big]^{-d/\bar{\gamma}(\delta)}\delta^{d/\bar{\gamma}(\delta)}.
\nonumber\\
&&=c_2 u^{\frac{2d}{p}-\left(\frac{d(p-2)}{2p}\right)\left(\frac{d}{\bar{\gamma}(\delta)}-2\right)}\;
\Big[\ms(u,\delta)\Big]^{-d/\bar{\gamma}(\delta)}\delta^{d/\bar{\gamma}(\delta)}
=c_2 u^{\frac{2d}{p}}[\delta(u)]^{2-\frac{d}{\bar{\gamma}(\delta)}}\;
\Big[\ms(u,\delta)\Big]^{-4}.
\end{eqnarray}
To get the last inequality we have taken into account that $\ms(u,\delta)<1$ for any $u\geq 1,\delta>0$ and that $d/\bar{\gamma}(\delta)\leq 4$.

We obtain from (\ref{eq5:proof-th-norm-gauss-local}) for any $\delta>0$, putting $c_3=(4\pi^2/3)^4 c_2$,
\begin{eqnarray}
\label{eq6:proof-th-norm-gauss-local}
%\sup_{\delta>0}
&&\delta^{-2}\mE^{(u)}\Big(\lambda_1^{-1}\big[g_A(u)\ms(u,\delta)\big]/48\delta\Big)\leq c_2 u^{\frac{2d}{p}}\;\left(\frac{\delta}{\delta(u)}\right)^{d/\bar{\gamma}(\delta)-2}\;\Big[\ms(u,\delta)\Big]^{-4}
\nonumber\\*[2mm]
&&\leq c_2\; u^{\frac{2d}{p}}
\left[\left(\delta\Big/\delta(u)\right)^{2/(p-1)}\mathrm{1}_{(0,\delta(u))}(\delta)+\left(\delta(u)\big/\delta\right)^{2-d/\gamma}
\mathrm{1}_{[\delta(u),\infty)}(\delta)\right]\Big[\ms(u,\delta)\Big]^{-4}
\nonumber\\*[2mm]
&&\leq c_2\; u^{\frac{2d}{p}}
\left[\left(2^{-m(u)}\delta\right)^{2/(p-1)}\mathrm{1}_{(1,\delta(u))}(\delta)+4\left(2^{m(u)}\big/\delta\right)^{2-d/\gamma}
\mathrm{1}_{[\delta(u),\infty)}(\delta)\right]\Big[\ms(u,\delta)\Big]^{-4}
\nonumber\\*[2mm]
&&= 5c_3\; u^{\frac{2d}{p}}\;2^{-\alpha\left|\log_2{\left(2^{-m(u)}\delta\right)}\right|}\left(1+\Big[\log_2{\big(2^{-m(u)}\delta\big)
\Big]^{2}}\right)^{4},
\end{eqnarray}
where $\alpha=\min\big\{2/(p-1), 2-d/\gamma\big\}$.
We obtain from (\ref{eq6:proof-th-norm-gauss-local})
\begin{eqnarray}
\label{eq8:proof-th-norm-gauss-local}
\cE^\prime(u):=\sup_{\delta>0}\left[\delta^{-2}\mE^{(u)}\Big(\lambda_1^{-1}\big[g_A(u)\ms(u,\delta)\big]/48\delta\Big)\right]\leq c_4\; u^{\frac{2d}{p}},
\end{eqnarray}
where
$
c_4=5c_3\sup_{x\geq 0}\left[2^{-\alpha x}\big(1+x^2\big)^{4}\right].
$

\smallskip

$7^{0}$. Remind that $\tau_1(h)=h^{-1}$ and, in particular, $\overline{\tau}_1=\big[h^{(\min)}\big]^{-1}$. Choosing $\e=\sqrt{2}-1$
we get from (\ref{eq8:proof-th-norm-gauss-local}) putting $c_5=c_42^{d/p}$
\begin{gather*}
\widehat{\cE}^{\big(\sqrt{2}-1\big)}(h):=\cE^\prime\Big(\sqrt{2}\tau_1(h)\Big)\leq c_5h^{-\frac{2d}{p}}.
\end{gather*}
Choose also $R_r(t)=R(t)=t^{\frac{2d}{p}}$ (independent of $r$) that yields
$$
\widehat{R}^{\big(\sqrt{2}-1\big)}(h):=R\Big(\sqrt{2}\tau_1(h)\Big)=2^{d/p}h^{-\frac{2d}{p}}.
$$
Choosing finally $z=0$ and putting $\widehat{\mathrm{U}}(h)=\widehat{\mathrm{U}}^{(z,\e,r)}(h),\e=\sqrt{2}-1, z=1,$ we obtain that
\begin{eqnarray*}
\widehat{\mathrm{U}}(h)\leq \sqrt{2}g_A\Big(2\tau_1(h)\Big)\sqrt{32
c_5h^{-\frac{2d}{p}}+2^{d/p}h^{-\frac{2d}{p}}+1}\leq c_6 h^{\frac{d(2-p)}{2p}-\frac{d}{p}}=c_6h^{-\frac{d}{2}}.
\end{eqnarray*}
Let us compute now the quantities $\cR^{(\e,0)},\cR^{(\e,q)}$ defined in (\ref{eq:condition-on the functions-R}) with $\e=\sqrt{2}-1$.

Noting that $\underline{\tau}_1= \left(h^{(\max)}\right)^{-1}$ and $2^{J/2}\overline{\tau}_1\leq 2^{-1/2}\underline{\tau}_1$ we get
\begin{eqnarray*}
\cR^{(\e,0)}&:=&\sum_{j=0}^{J}\exp{\left\{-\Big(\overline{\tau}_1 2^{-(j/2)}\Big)^{\frac{2d}{p}}\right\}}\leq
c_8\exp{\left\{-2^{-3/2}\left(h^{(\max)}\right)^{-2d/p}\right\}};
\\
\cR^{(\e,q)}&:=&\sum_{j=0}^{J}\left[g_A\big(\overline{\tau}_1 2^{-j/2}\big)\right]^{q}\exp{\left\{-\Big(\overline{\tau}_1 2^{-(j/2)}\Big)^{\frac{2d}{p}}\right\}}
\\
&\leq& c_9\;\left(h^{(\max)}\right)^{\frac{qd(2-p)}{2p}}
\exp{\left\{-2^{-3/2}\left(h^{(\max)}\right)^{-2d/p}\right\}}.
\end{eqnarray*}
The assertions of the theorem follow now from Proposition \ref{prop_uniform_local3}.

\smallskip

$8^{0}$.
It remains to  prove (\ref{eq2:proof-th-norm-gauss-local}). The proof is based on the following inclusion:
for any $\bar{\gamma}\in(0,\gamma]$
\begin{equation}
\label{eq:smoothness-of-theta-gauss_local}
\Theta^\prime_1(u)\subset\bH^*_{\frac{p}{p-1}}\left(\bar{\gamma},3^{d}\sqrt{d}u^{\bar{\gamma}}\right),\;\;\forall u>0,
\end{equation}
where $\bH^{*}_{q}(\cdot,\cdot)\subset\bH_{q}(\cdot,\cdot)$ consists  of functions vanishing outside of  $\bK_{2\mu}$.

Let $\mE^*(\cdot)$ be the entropy of $\bH^{*}_{\frac{p}{p-1}}\left(\bar{\gamma},L\right),L>0$, measured in $\|\cdot\|_2$. It is well-known \cite{EdmTri}, that for any $p>1$ there exist $c_0$ completely determined by $\gamma, d$, $p$ and $\mu$ such that for any $(d/2-d/p)<\bar{\gamma}\leq\gamma$ and for any $L>0$
$$
\mE^*(\delta)\leq c_0\left(L \delta^{-1}\right)^{d/\bar{\gamma}},
\;\;\forall\; \delta>0.
$$
 Since we consider only  $(d/2-d/2p)<\bar{\gamma}\leq\gamma$, (\ref{eq2:proof-th-norm-gauss-local}) follows immediately from (\ref{eq:smoothness-of-theta-gauss_local}). Thus, we shall prove (\ref{eq:smoothness-of-theta-gauss_local}).

Fix $\theta\in\Theta^\prime_1(u)$. By its definition there exists $\ell\in\bB_{\frac{p}{p-1}}$ and $h\geq u^{-1}$ such that $\theta=K_h\ast \ell$, where $"\ast"$ stands convolution operator on $\bR^d$. First, we note that all functions belonging to $\Theta$ vanish outside the cube $\bK_{2\mu}$ in view of assumption imposed on function $K$.

Next,  for any $m=(m_1,\ldots,m_d)\in\bN^d $ put $|m|=m_1+\cdots+m_d$, and  set $\gamma=l+\alpha$ and  $\bar{\gamma}=\bar{l}+\bar{\alpha}$, where $l,\bar{l}\in\bN$ and $0<\alpha,\bar{\alpha}\leq 1$.

Then, since $K\in\bH_\infty(\gamma, 1)$ we have for any $m\in\bN^d$ such that $|m|\leq l$
$$
\sup_{x\in\bR^d}\left|\frac{\partial^{|m|}K_h(x)}{\partial^{m_1}x_1\cdots\partial^{m_d}x_d}\right|\leq h^{-|m|}\leq  \big(h^{(\min)}\big)^{-\gamma}<\infty.
$$
Above remarks allow us to assert that  all partial derivatives $\theta^{(m)}$ exist whenever $|m|\leq l$ and they are given by
$$
\theta^{(m)}(x)=\int(K_h)^{(m)}(t-x)\ell(t)\rd t,\;\; \forall x\in\bR^d,
$$
where for any function $g$ the notation $g^{(m)}$ or (if it is more convenient) $(g)^{(m)}$ is used for its partial derivative.

We obtain in view of Young inequality for any $\Delta\in\bR$ and any $m\in\bN^d$ satisfying $|m|=\bar{l}$
\begin{eqnarray*}
\left\|\theta^{(m)}(\cdot+\Delta)-\theta^{(m)}(\cdot)\right\|_{\frac{p}{p-1}}&\leq&
 \left\|(K_h)^{(m)}(\cdot+\Delta)-(K_h)^{(m)}(\cdot)\right\|_1\\*[2mm]
 &\leq& h^{-\bar{l}}\left\|K^{(m)}(\cdot+[\Delta/h])-K^{(m)}(\cdot)\right\|_1.
\end{eqnarray*}
Here we have used that $\ell\in\bB_{\frac{p}{p-1}}$.

We remark that if $h\leq |\Delta|$ then for any $u\in\bR^d$ either $K^{(m)}(u+[\Delta/h])=0$ or $K^{(m)}(u)=0$
in view of the assumption imposed on the support of $K$. Thus, if $h\leq |\Delta|$
$$
\left\|K^{(m)}(\cdot+[\Delta/h])-K^{(m)}(\cdot)\right\|_1\leq \big\|K^{(m)}\big\|_1\leq \big\|K^{(m)}\big\|_1\big(\Delta/h\big)^{\bar{\alpha}}\leq\big(\Delta/h\big)^{\bar{\alpha}},
$$
since $h\leq |\Delta|$ and $\big\|K^{(m)}\big\|_1\leq \big\|K^{(m)}\big\|_\infty \leq 1$ in view of assumption imposed on the function $K$.

If $h> |\Delta|$ then in view of the assumption imposed on the support of $K$ we have
\begin{eqnarray*}
&&\left\|K^{(m)}(\cdot+[\Delta/h])-K^{(m)}(\cdot)\right\|_1=
\int_{\left[-\frac{3}{2},\frac{3}{2}\right]^{d}}\left|K^{(m)}(u+[\Delta/h])-K^{(m)}(u)\right|\rd u
\\*[2mm]
&&\leq 3^{d}\left|\sum_{i=1}^{d}(\Delta/h)^{2}\right|^{\bar{\alpha}/2}\leq 3^{d}\sqrt{d}\big(|\Delta|/h\big)^{\bar{\alpha}}.
\end{eqnarray*}
%where $c_2=3^{d}\Big([4\mu\sqrt{d}]\vee 1\Big)d^{\bar{\alpha}/2}$. Here we have used that (\ref{eq1:proof-th-norm-gauss-local}).
Since $h\geq u^{-1}$ we conclude finally that $\forall \Delta\in\bR$
\begin{eqnarray}
\label{eq100:proof-of-theorem-gauus-norm}
\left\|\theta^{(m)}(\cdot+\Delta)-\theta^{(m)}(\cdot)\right\|_{\frac{p}{p-1}}\leq 3^{d}\sqrt{d}h^{-\bar{\gamma}}|\Delta|^{\bar{\alpha}}\leq
3^{d}\sqrt{d}u^{\bar{\gamma}}|\Delta|^{\bar{\alpha}}.
\end{eqnarray}
 It means that $\theta\in\bH_{\frac{p}{p-1}}\left(\bar{\gamma},3^{d}\sqrt{d}u^{\bar{\gamma}}\right)$.  As it was mentioned above all function belonging to $\Theta$ vanish outside the cube $\bK_{2\mu}$ that allows us to conclude that $\theta\in\bH^*_{\frac{p}{p-1}}\left(\bar{\gamma},3^{d}\sqrt{d}u^{\bar{\gamma}}\right)$ and,
 therefore, (\ref{eq:smoothness-of-theta-gauss_local}) is proved.

\epr

\section{Proof of Theorems \ref{th:empiric_totaly_bounded_case}--\ref{th:LL-nonasym}}
\label{sec:proofs-of-empiric-theorems}

\subsection{Proof of Theorem \ref{th:empiric_totaly_bounded_case}}

\subsubsection{Preliminaries}
We start the proof with several technical results used in the sequel.
Put for any $i=\overline{1,\mathbf{n_2}}$, $y\in\big[\mathbf{n_1}/\mathbf{n_2}, 1\big]$  and $\alpha=\boldsymbol{b}\big[\ln{(\mathbf{n_2})}\big]^{-1},$
$$
\cQ_i(y)=\mathrm{1}_{(i/\mathbf{n_2},1]}\big(y\big)+\big(\mathbf{n_2}y-i+1\big)^{\alpha}\mathrm{1}_{\Delta_i}(y),\quad Q_i(y)=y^{-1}\cQ_i(y).
$$
Here we have denoted $\Delta_i=\big((i-1)/\mathbf{n_2},i/\mathbf{n_2}\big],\;i=\overline{3,\mathbf{n_2}}$ and
$\Delta_{2}=\big[1/\mathbf{n_2},2/\mathbf{n_2}\big].$
%The choice of the  parameter $\beta$ will be done later.

For any $a\geq 1$ let $\lceil a\rceil$ be the smallest  integer larger or equal to  $a$.
It implies, in particular, $y\in\Delta_{\lceil \mathbf{n_2}y\rceil}$.
First we note that for any $y,\bar{y}\in[\mathbf{n_1}/\mathbf{n_2},1]$ and any $ i=\overline{1,\mathbf{n_2}}$
\begin{equation}
\label{eq1:property-of_functions-Q_i}
\cQ_i(y)\leq 1,\qquad \cQ_i(y)=0,\;\forall i> \lceil \mathbf{n_2}y\rceil,  \qquad\left|\cQ_i(y)-\cQ_i(\bar{y})\right|\leq 1\wedge\left|\mathbf{n_2}(y-\bar{y})\right|^{\alpha}.
\end{equation}
The first and  third inequalities  imply for any $i=\overline{1,\mathbf{n_2}}$ and any $ y,\bar{y}\in[\mathbf{n_1}/\mathbf{n_2},1]$
\begin{eqnarray}
\label{eq2:property-of_functions-Q_i}
&&\left|Q_i(y)-Q_i(\bar{y})\right|\leq (y\wedge\bar{y})^{-1}\left[\left|\mathbf{n_2}(y-\bar{y})\right|^{\alpha}+\left(1-\frac{y\wedge\bar{y}}{y\vee\bar{y}}\right)\right].
\end{eqnarray}
For any $z,z^\prime\in\bR_+ $ denote
$
 \mathfrak{w}\big(z,z^\prime\big)=\left(1-\sqrt{\frac{z\wedge z^\prime}{z\vee z^\prime}}\right)^{1/2},
$
and remark that $\mathfrak{w}$ is a metric on $\bR_+$. It follows from the relation
$
\sqrt{2}\mathfrak{w}\big(z,z^\prime\big)=\left[\bE\left(\frac{b(z)}{\sqrt{z}}-\frac{b(z^\prime)}{\sqrt{z^\prime}}\right)^2\right]^{1/2},
$
where $b$ is the standard Wiener process. Taking into account that $y,\bar{y}\geq 1/2$ and that $\mathfrak{w}(y\wedge\bar{y})\leq 1$ we obtain from
(\ref{eq2:property-of_functions-Q_i})
\begin{eqnarray}
\label{eq222:property-of_functions-Q_i}
&&\left|Q_i(y)-Q_i(\bar{y})\right|\leq 8e^{\boldsymbol{b}}\big[\mathfrak{w}(y,\bar{y})\big]^{\alpha}.
\end{eqnarray}
Here we have also used the definition of  $\alpha$.
Taking into account that for any $\boldsymbol{a}\leq \boldsymbol{c}$
$$
\sup_{p\in (0,1]}p^{\boldsymbol{a}}\left(1-\ln{(p)}\right)^{\boldsymbol{c}}= e^{\boldsymbol{a}-\boldsymbol{c}}\big[\boldsymbol{c}/\boldsymbol{a}\big]^{\boldsymbol{c}},
$$
we obtain from
(\ref{eq222:property-of_functions-Q_i}) for any $\boldsymbol{b}>0$,  $ y,\bar{y}\in[\mathbf{n_1}/\mathbf{n_2},1]$ and $\mathbf{n_2}\geq 3$
\begin{equation}
\label{eq3:property-of_functions-Q_i}
\sup_{i=\overline{1,\mathbf{n_2}}}\left|Q_i(y)-Q_i(\bar{y})\right|\leq 8e
\left[\frac{\ln{(\mathbf{n_2})}}{1-\ln{\big(\mathfrak{w}(y,\bar{y})\big)}}\right]^{\boldsymbol{b}}
\end{equation}
Next,  for any $y,\bar{y}\in [\mathbf{n_1}/\mathbf{n_2},1]$
%such that $\Delta_{\lceil \mathbf{n_2}t\rceil}\neq\Delta_{\lceil \mathbf{n_2}\bar{t}\rceil}$
\begin{eqnarray}
\label{eq44:property-of_functions-Q_i}
\left|\cQ_i(y)-\cQ_i(\bar{y}\right|
=0,\quad i\notin\left\{\lceil \mathbf{n_2}(y\wedge\bar{y})\rceil,\ldots,\lceil \mathbf{n_2}(y\vee\bar{y})\rceil\right\}.
\end{eqnarray}
We have for any $y,\bar{y}\in\big[\mathbf{n_1}/\mathbf{n_2}, 1\big]$ in view of the first and third bounds in (\ref{eq1:property-of_functions-Q_i}) and (\ref{eq44:property-of_functions-Q_i})
\begin{eqnarray*}
\sum_{i=1}^{\mathbf{n_2}}\left|\cQ_i(y)-\cQ_i(\bar{y}\right|^{2}\leq
\left\{\begin{array}{ll}
2\mathbf{n_2}\left|y-\bar{y}\right|,\quad& \lceil \mathbf{n_2}(y\vee\bar{y})\rceil-\lceil \mathbf{n_2}(y\wedge\bar{y})\rceil\geq 3;
\\
3\left|\mathbf{n_2}(y-\bar{y})\right|^{2\alpha} ,\quad& \lceil \mathbf{n_2}(y\vee\bar{y})\rceil-\lceil \mathbf{n_2}(y\wedge\bar{y})\rceil\leq 2.
\end{array}
\right.
\end{eqnarray*}
To get the first inequality we have also used that $\lceil \mathbf{n_2}(y\vee\bar{y})\rceil-\lceil \mathbf{n_2}(y\wedge\bar{y})\rceil\geq 3$ implies
$\mathbf{n_2}(y\vee\bar{y}-y\wedge\bar{y})>2$ and, therefore,
$
\lceil \mathbf{n_2}(y\vee\bar{y})\rceil-\lceil \mathbf{n_2}(y\wedge\bar{y})\rceil+1\leq \mathbf{n_2}(y\vee\bar{y}-y\wedge\bar{y})+2\leq
2\mathbf{n_2}(y\vee\bar{y}-y\wedge\bar{y})=2\mathbf{n_2}\left|y-\bar{y}\right|.
$
Thus,  we have for any $y,\bar{y}\in\big[\mathbf{n_1}/\mathbf{n_2}, 1\big]$
$$
\sqrt{\sum_{i=1}^{\mathbf{n_2}}\left|\cQ_i(y)-\cQ_i(\bar{y}\right|^{2}}\leq \sqrt{2\mathbf{n_2}\left|y-\bar{y}\right|}+ \sqrt{3}\left|\mathbf{n_2}(y-\bar{y})\right|^{\alpha}\leq 2\sqrt{\mathbf{n_2}}\mw\big(y,\bar{y}\big)+2\sqrt{3}e^{\boldsymbol{b}}\big[\mw\big(y,\bar{y}\big)\big]^{\alpha}.
$$
Hence we get
\begin{eqnarray*}
\sqrt{\sum_{i=1}^{\mathbf{n_2}}\left|Q_i(y)-Q_i(\bar{y})\right|^{2}}&\leq& 8\sqrt{\mathbf{n_2}}\mw\big(y,\bar{y}\big)+4\sqrt{3}e^{\boldsymbol{b}}\big[\mw\big(y,\bar{y}\big)\big]^{\alpha}
\\
&\leq& 8\sqrt{\mathbf{n_2}}\mw\big(y,\bar{y}\big)+4\sqrt{3}e
\left[\frac{\ln{(\mathbf{n_2})}}{1-\ln{\big(\mathfrak{w}(y,\bar{y})\big)}}\right]^{\boldsymbol{b}} .
\end{eqnarray*}
 Taking into account that  $\sup_{z\geq 1}z^{-1/2}\left[\ln{(2ez)}\right]^{\boldsymbol{b}}\leq (2\boldsymbol{b}/e)^{\boldsymbol{b}}$
 we obtain
\begin{equation*}
\sqrt{\sum_{i=1}^{\mathbf{n_2}}\left|Q_i(y)-Q_i(\bar{y}\right|^{2}}\leq 8\sqrt{\mathbf{n_2}}
\left[\mw\big(y,\bar{y}\big)+\sqrt{3/4}e(2\boldsymbol{b}/e)^{\boldsymbol{b}}
\left\{1-\ln{\big(\mathfrak{w}(y,\bar{y})\big)}\right\}^{-\boldsymbol{b}}\right].
\end{equation*}
Finally we get for any
$ y,\bar{y}\in[\mathbf{n_1}/\mathbf{n_2},1]$ and any $\boldsymbol{b}>1$
\begin{equation}
\label{eq4:property-of_functions-Q_i}
\sqrt{\sum_{i=1}^{\mathbf{n_2}}\left|Q_i(y)-Q_i(\bar{y}\right|^{2}}\leq 8
\left[2^{\boldsymbol{b}}+1\right](\boldsymbol{b})^{\boldsymbol{b}}\sqrt{\mathbf{n_2}}\Big[1-\ln{\big(\mathfrak{w}(y,\bar{y})\big)}
\Big]^{-\boldsymbol{b}}.
\end{equation}

\subsubsection{Constants}
\label{sec:constants-th-totally-bounded}

The following constants appeared in the description of  upper functions and inequalities
found in Theorem \ref{th:empiric_totaly_bounded_case}.
Let $\chi=0$ if $\mathbf{n_1}=\mathbf{n_2}$ and  $\chi=1$ if $\mathbf{n_1}\neq\mathbf{n_2}$.
\begin{gather*}
C_{N,R,m,k}=C^{(1)}_{N,R,m,k}+C^{(2)}_{N,R,m,k}+2\chi\mathbf{a}_{\boldsymbol{b}},\qquad \mathbf{a_{\boldsymbol{b}}}=2\delta^{-2}_*\ln(2)+2\sup_{\delta>\delta_*}(\delta^2\wedge\delta)^{-1}
\big(96\delta\big/s^*(\delta)\big)^{\frac{1}{{\boldsymbol{b}}}};
\\*[2mm]
C^{(1)}_{N,R,m,k}=\sup_{\delta>\delta_*}\delta^{-2}\left\{k\left[1+\ln{\left(\frac{9216m\delta^{2}}{[s^*(\delta)]^{2}}\right)}\right]_++N(m-k)
\left(\left[\log_2{\left\{\left(\frac{4608m R\delta^{2}}{[s^*(\delta)]^{2}}\right)\right\}}\right]_++1\right)\right\};
\\*[2mm]
C^{(2)}_{N,R,m,k}=\sup_{\delta>\delta_*}\delta^{-1}\left\{k\left[1+\ln{\left(\frac{9216m\delta}{s^*(\delta)}\right)}\right]_+
+N(m-k)\left(\left[\log_2{\left\{\left(\frac{4608m R\delta}{s^*(\delta)}\right)\right\}}\right]_++1\right)\right\}.
\end{gather*}
Put also
$
C_D:=\left[\sup_{j=0,k+1,\ldots,m}\;\sup_{z\in [0,1]}D_j^\prime(z)\right]\vee 2,
$
where $D_j^\prime$ is the first derivative of the function $D_j$. Set at last,
$\boldsymbol{c_b}=4\sqrt{2}\left[2^{\boldsymbol{b}}+1\right](\boldsymbol{b})^{\boldsymbol{b}}$ and let
$$
\lambda_1=4\sqrt{2e}\big(\sqrt{C_D}\vee[\chi\boldsymbol{c_b}]\big),\quad \lambda_2=(16/3)\big(C_D\vee 8e\big),\quad
C_{D,\boldsymbol{b}}=\big(\sqrt{2C_D}\vee[\chi\boldsymbol{c_b}]\big)\vee \big[(2/3)\big(C_D\vee 8e\big)\big].
$$

\subsubsection{Proof of the theorem}

$\mathbf{1^0.}$ Put for any $i=\overline{1,n}$
$$
\e\big(\mh,X_i\big)=G\big(\mh,X_i\big)-\bE_{\rf} G(\mh,X_i),
$$
and define for any  $y\in(\mathbf{n_1}/\mathbf{n_2},1]$ and any $\mh\in\boldsymbol{\mH}$ the random function
\begin{equation}
\label{eq:equiv_expression-for-bold-xi}
\boldsymbol{\xi}(y,\mh)=\mathbf{n_2}^{-1}\sum_{i=1}^{\mathbf{n_2}}
\e\big(\mh,X_{i}\big)Q_i(y).
\end{equation}
We remark that $
\xi_\mh(p)=\boldsymbol{\xi}\big(p/\mathbf{n_2},\mh\big)
$ for any $p\in \widetilde{\mathbf{N}}$ and any $\mh\in\boldsymbol{\mH}$. Thus,
  in order to get the assertions of the theorem it suffices to find  upper functions for
$|\boldsymbol{\xi}(\cdot,\cdot)|$ on $\big[\mathbf{n_1}/\mathbf{n_2}, 1\big]\times\mH\big(\mathbf{n}\big)$  in view of Assumption \ref{ass:dependence-on-n} and the definition of the number $\mathbf{n}$.

In view of Bernstein inequality
 Assumption \ref{ass:fixed_theta_local}   is fulfilled with
 $\theta=\mathrm{h}=:
(y,\mh)$ and $\bar{\theta}=\bar{\mathrm{h}}=:
(\bar{y},\bar{\mh})$
%and (without loss of generality we assume that $\bar{t}\leq t$)
\begin{eqnarray}
\label{eq11:def-A-B-a-b-empirical}
&&\displaystyle{A^{2}(\theta)=A^{2}_{\rf}(\mathrm{h}):=2\mathbf{n_2}^{-2}\sum_{i=1}^{\mathbf{n_2}}Q^2_i(y)\bE_fG^{2}(\mh,X_i)};
\\
\label{eq22:def-A-B-a-b-empirical}
&&\ra^2(\theta,\bar{\theta})=\ra^{2}_{\rf}(\mathrm{h},\bar{\mathrm{h}}):=2\mathbf{n_2}^{-2}\sum_{i=1}^{\mathbf{n_2}}
\bE_f\Big[Q_i(y)G(\mh,X_i)-Q_i(\bar{y})G(\bar{\mh},X_i)\Big]^{2};
\\
\label{eq33:def-A-B-a-b-empirical}
&& B(\theta)=B_{\infty}(\mathrm{h})=(4/3)\mathbf{n_2}^{-1}\Big[\sup_{i=\overline{1,\mathbf{n_2}}}Q_i(y)\Big]\sup_{x\in\cX}\big|G(\mh,x)\big|.
\\
\label{eq44:def-A-B-a-b-empirical}
&&\rb(\theta,\bar{\theta})=\rb_\infty(\mathrm{h},\bar{\mathrm{h}})
:=(2/3)\mathbf{n_2}^{-1}\sup_{i=\overline{1,n}}\sup_{x\in\cX}\left|\e(\mh,x)Q_i(y)-\e(\bar{\mh},x)Q_i(\bar{y})\right|.
%\\
%&&\rb_\infty(\mathrm{h},\bar{\mathrm{h}})=(4/3)\mathbf{n_2}^{-1}\bigg\{
%\Big[\sup_{i=\overline{1,\mathbf{n_2}}}Q_i(y)\Big]\sup_{x\in\cX}\big|G(\mh,x)-G(\bar{\mh},x)\big|
%\\
%&&\qquad\qquad\qquad\qquad\qquad\qquad\;
%%\bigg\}
%\nonumber
\end{eqnarray}
Note that $\ra_{\rf}$ and $\rb_\infty$ are semi-metrics on $\big[\mathbf{n_1}/\mathbf{n_2}, 1\big]\times\boldsymbol{\mH}$ and $\boldsymbol{\xi}(\cdot,\cdot)$ is obviously continuous on  $\big[\mathbf{n_1}/\mathbf{n_2},1\big]\times\mH\big(\mathbf{n}\big)$ in the topology generated by $\rb_\infty$. Moreover, $A_\rf$ and $B_\infty$ are bounded and, therefore, Assumption \ref{ass:metric-case-local} is fulfilled.

\noindent Later on we will use the following notation: for any  $\mathfrak{Q}:\cX\to\bR$ put $\|\mathfrak{Q}\|_\infty=\sup_{x\in\cX}|\mathfrak{Q}(x)|$.

We obtain from (\ref{eq11:def-A-B-a-b-empirical})--(\ref{eq44:def-A-B-a-b-empirical}) and (\ref{eq3:property-of_functions-Q_i}) for any $\mathrm{h},\bar{\mathrm{h}}\in[\mathbf{n_1}/\mathbf{n_2},1]\times\mH\big(\mathbf{n}\big)$
%(since  all expressions in (\ref{eq11:def-A-B-a-b-empirical})--(\ref{eq44:def-A-B-a-b-empirical}) are symmetric with respect to $y,\bar{y},$ without loss of generality we will assume that $y=y\vee\bar{y}$).
\begin{eqnarray}
\label{eq111:def-A-B-a-b-empirical}
&&A^{2}_{\rf}(\mathrm{h})\leq 2(\mathbf{n_1})^{-1}F_{\mathbf{n_2}}(\mh)G_\infty\big(\mh^{(k)}\big),
\quad B_{\infty}(\mathrm{h})\leq (4/3)(\mathbf{n_1})^{-1}G_\infty\big(\mh^{(k)}\big);
\\*[2mm]
\label{eq444:def-A-B-a-b-empirical}
&&\rb_\infty(\mathrm{h},\bar{\mathrm{h}})\leq\frac{4\ln^{\beta}{(\mathbf{n_2})}}{3\mathbf{n_1}}\bigg\{
\big\|G(\mh,\cdot)-G(\bar{\mh},\cdot)\big\|_\infty
%\\
%&&\hskip4.8cm
+\gamma 8eG_\infty\big(\bar{\mh}^{(k)}\big)
\Big[1-\ln{\big(\mathfrak{w}(y,\bar{y})\big)}\Big]^{-\boldsymbol{b}}
\bigg\},
%\nonumber
\end{eqnarray}
where, remind,
 $\gamma=0$ if $\mathbf{n_1}=\mathbf{n_2}$ and  $\gamma=1$ if $\mathbf{n_1}\neq\mathbf{n_2}$.
 Here we have used that if $\mathbf{n_1}=\mathbf{n_2}$  the second term in the last inequality disappears.

We also get  using (\ref{eq1:property-of_functions-Q_i}) and (\ref{eq4:property-of_functions-Q_i})
\begin{eqnarray}
\label{eq222:def-A-B-a-b-empirical}
\ra_{\rf}(\mathrm{h},\bar{\mathrm{h}})&\leq& \sqrt{2}\mathbf{n_2}^{-1}\Bigg\{\sqrt{\sum_{i=1}^{\mathbf{n_2}}Q^{2}_i(y)
\bE_f\Big[G(\mh,X_i)-G(\bar{\mh},X_i)\Big]^{2}}
\nonumber
\\
&&\hskip1.5cm +\sqrt{F_{\mathbf{n_2}}(\bar{\mh})G_\infty \big(\bar{\mh}^{(k)}\big)\sum_{i=1}^{\mathbf{n_2}}
\big(Q_i(y)-Q_i(\bar{y})\big)^{2}}\Bigg\};
\nonumber\\
&\leq& \sqrt{2}(\mathbf{n_1})^{-1/2}\bigg\{\sqrt{\left(F_{\mathbf{n_2}}(\mh)+F_{\mathbf{n_2}}(\bar{\mh})\right)
\big\|G(\mh,\cdot)-G(\bar{\mh},\cdot)\big\|_\infty}
\nonumber\\
&&\hskip2cm
+\chi\boldsymbol{c_b}\sqrt{2F_{\mathbf{n_2}}(\bar{\mh})G_\infty \big(\bar{\mh}^{(k)}\big)}
\Big[1-\ln{\big(\mathfrak{w}(y,\bar{y})\big)}\Big]^{-\boldsymbol{b}}
\bigg\},
%\nonumber
\end{eqnarray}
where we have put $\boldsymbol{c_b}=4\sqrt{2}\left[2^{\boldsymbol{b}}+1\right](\boldsymbol{b})^{\boldsymbol{b}}$.
 Here we have used that if $\mathbf{n_1}=\mathbf{n_2}$  the second term in the last inequality disappears.

For any  $\tau>0$ put
$
\mH\big(\mathbf{n},\tau\big)=\left\{\mh\in\mH\big(\mathbf{n}\big):\;\; F_{\mathbf{n_2}}(\bar{\mh})\leq \tau\right\}.
$
Our first step consists in establishing an upper function for $|\boldsymbol{\xi}(\cdot,\cdot)|$ on $\mathrm{H}(\tau):=\big[\mathbf{n_1}/\mathbf{n_2}, 1\big]\times\mH\big(\mathbf{n},\tau\big)$.
As always the supremum
over empty set is supposed to be zero.

\smallskip

$\mathbf{2^0.}$ Note that in view of (\ref{eq111:def-A-B-a-b-empirical})  and (\ref{eq222:def-A-B-a-b-empirical})  for any $\mathrm{h},\overline{\mathrm{h}}\in\mathrm{H}(\tau)$
\begin{eqnarray}
\label{eq1:proof-theorem4-empirical}
&&\;A^{2}_{\rf}(\mathrm{h})\leq 2\tau (\mathbf{n_1})^{-1}G_\infty\big(\mh^{(k)}\big),\quad B_{\infty}(\mathrm{h})
\leq \frac{4\ln^{\beta}{(\mathbf{n_2})}}{3\mathbf{n_1}}G_\infty\big(\mh^{(k)}\big);
\\[2mm]
\label{eq2:proof-theorem4-empirical}
&&\;\ra_{\rf}(\mathrm{h},\bar{\mathrm{h}})\leq 2\sqrt{\tau}(\mathbf{n_1})^{-1/2}\bigg\{\sqrt{
\big\|G(\mh,\cdot)-G(\bar{\mh},\cdot)\big\|_\infty}
%\\
%&&\hskip4.2cm
+
\chi\boldsymbol{c_b}\sqrt{G_\infty \big(\bar{\mh}^{(k)}\big)}
\Big[1-\ln{\big(\mathfrak{w}(y,\bar{y})\big)}\Big]^{-\boldsymbol{b}}
\bigg\}.
%\nonumber
\end{eqnarray}
%Here we have implicitly used the definition of the mappings $G_\infty(\cdot)$ and $F_n(\cdot)$, i.e (\ref{eq:def-of-function-G-infty}) and
%(\ref{eq:def-function-F_n}).

Moreover, in view of triangle inequality we obviously have for any $\mathrm{h},\overline{\mathrm{h}}\in\mathrm{H}(\tau)$
\begin{eqnarray}
\label{eq3:proof-theorem4-empirical}
&&\ra_{\rf}(\mathrm{h},\overline{\mathrm{h}})\leq A_\rf(\mathrm{h})+A_\rf(\overline{\mathrm{h}})\leq\sqrt{8\tau (\mathbf{n_1})^{-1}\left[G_\infty\Big(\mh^{(k)}\Big)\vee G_\infty\Big(\overline{\mh}^{(k)}\Big)\right]};
\\[2mm]
\label{eq4:proof-theorem4-empirical}
&&\rb_{\infty}(\mathrm{h},\overline{\mathrm{h}})\leq B_\infty(\mathrm{h})+B_\infty(\overline{\mathrm{h}})\leq \frac{8\ln^{\beta}{(\mathbf{n_2})}}{3\mathbf{n_1}}
\left[G_\infty\Big(\mh^{(k)}\Big)\vee G_\infty\Big(\overline{\mh}^{(k)}\Big)\right].
\end{eqnarray}
Set
$$
\cG\left(\mh^{(k)},\overline{\mh}^{(k)}\right)=G_\infty\Big(\mh^{(k)}\Big)\vee G_\infty\Big(\overline{\mh}^{(k)}\Big).
$$
We get for any $\mh,\overline{\mh}$, satisfying $\varrho^{(k)}\Big(\mh^{(k)},\overline{\mh}^{(k)}\Big)\vee{\displaystyle\sup_{j=\overline{k+1,m}}}\varrho_j\big(\mh_j,\mh^\prime_j\big)\leq 1$  in view of Assumption \ref{ass:bounded_case} ($\mathbf{ii}$)
\begin{equation*}
%\label{eq5005:proof-theorem4-empirical}
\big\|G(\mh,\cdot)-G(\bar{\mh},\cdot)\big\|_\infty\leq C_D \bigg\{\cG\left(\mh^{(k)},\overline{\mh}^{(k)}\right)\varrho^{(k)}\Big(\mh^{(k)},\overline{\mh}^{(k)}\Big)+
\sum_{j=k+1}^{m}L_j\left\{\cG\left(\mh^{(k)},\overline{\mh}^{(k)}\right)\right\}\varrho_j\big(\mh_j,\overline{\mh}_j\big)\bigg\}.
\end{equation*}
On the other hand, putting  $\widetilde{L}_j(y)=L_j(y)\vee y, \;j=0,k+1,\ldots m$, we have
  for any $\mh,\overline{\mh}$, satisfying $\left[\varrho^{(k)}\Big(\mh^{(k)},\overline{\mh}^{(k)}\Big)\vee\sup_{j=\overline{k+1,m}}\varrho_j\big(\mh_j,\mh^\prime_j\big)\right]> 1$
\begin{eqnarray*}
%\label{eq5005:proof-theorem4-empirical}
&&\hskip-0.9cm\big\|G(\mh,\cdot)-G(\bar{\mh},\cdot)\big\|_\infty\leq \big\|G(\mh,\cdot)\big\|_\infty+\big\|G(\bar{\mh},\cdot)\big\|_\infty\leq 2\cG\left(\mh^{(k)},\overline{\mh}^{(k)}\right)
\\
&&\hskip2.6cm \leq C_D \bigg\{\cG\left(\mh^{(k)},\overline{\mh}^{(k)}\right)\varrho^{(k)}\Big(\mh^{(k)},\overline{\mh}^{(k)}\Big)+
\sum_{j=k+1}^{m}\widetilde{L}_j\left\{\cG\left(\mh^{(k)},\overline{\mh}^{(k)}\right)\right\}\varrho_j\big(\mh_j,\overline{\mh}_j\big)\bigg\}.
\end{eqnarray*}
Here we have also used that $C_D\geq 2$.
Thus, finally we have for any $\mh,\overline{\mh}$
\begin{equation*}
%\label{eq5005:proof-theorem4-empirical}
\big\|G(\mh,\cdot)-G(\bar{\mh},\cdot)\big\|_\infty\leq C_D \bigg\{\cG\left(\mh^{(k)},\overline{\mh}^{(k)}\right)\varrho^{(k)}\Big(\mh^{(k)},\overline{\mh}^{(k)}\Big)+
\sum_{j=k+1}^{m}\widetilde{L}_j\left\{\cG\left(\mh^{(k)},\overline{\mh}^{(k)}\right)\right\}\varrho_j\big(\mh_j,\overline{\mh}_j\big)\bigg\}.
\end{equation*}
The latter inequality together with (\ref{eq444:def-A-B-a-b-empirical})  and (\ref{eq2:proof-theorem4-empirical}) yields for any $\mathrm{h},\overline{\mathrm{h}}\in\mathrm{H}(\tau)$
\begin{eqnarray}
\label{eq5:proof-theorem4-empirical}
\ra_{\rf}(\mathrm{h},\overline{\mathrm{h}})&\leq& \ma \bigg\{\bigg(\cG\left(\mh^{(k)},\overline{\mh}^{(k)}\right)\varrho^{(k)}\Big(\mh^{(k)},
\overline{\mh}^{(k)}\Big)+
\sum_{j=k+1}^{m}\widetilde{L}_j\left\{\cG\left(\mh^{(k)},\overline{\mh}^{(k)}\right)\right\}\varrho_j\big(\mh_j,\overline{\mh}_j\big)\bigg)^{1/2}
\\
&&\hskip0.5cm +\chi\sqrt{\cG\left(\mh^{(k)},\overline{\mh}^{(k)}\right)}
\Big[1-\ln{\big(\mathfrak{w}(y,\bar{y})\big)}\Big]^{-\boldsymbol{b}}
\bigg\};
\nonumber\\
\label{eq6:proof-theorem4-empirical}
\qquad\rb_{\infty}(\mathrm{h},\overline{\mathrm{h}})&\leq& \mb \bigg\{\cG\left(\mh^{(k)},\overline{\mh}^{(k)}\right)\varrho^{(k)}\Big(\mh^{(k)},\overline{\mh}^{(k)}\Big)+
\sum_{j=k+1}^{m}\widetilde{L}_j\left\{\cG\left(\mh^{(k)},\overline{\mh}^{(k)}\right)\right\}\varrho_j\big(\mh_j,\overline{\mh}_j\big)
\\
&&\hskip0.5cm
+\chi\cG\left(\mh^{(k)},\overline{\mh}^{(k)}\right)
\Big[1-\ln{\big(\mathfrak{w}(y,\bar{y})\big)}\Big]^{-\boldsymbol{b}}
\bigg\},
\nonumber
\end{eqnarray}
where we have put $\ma=2\sqrt{\tau}(\mathbf{n_1})^{-1/2}\big(\sqrt{C_D}\vee[\chi\boldsymbol{c_b}]\big),\;
\mb=\frac{4\big(C_D\vee 8e\big)\ln^{\beta}{(\mathbf{n_2})}}{3\mathbf{n_1}}$.

\smallskip

$\mathbf{3^0.}$
We note that  in view of (\ref{eq1:proof-theorem4-empirical}) ) Assumption \ref{ass:fixed_theta_local} ({\it 1})
 is verified on $\mathrm{H}(\tau)$ with
 $$
 A(\theta)=A\left(\mathrm{h}\right):=\ma\sqrt{ G_\infty\big(\mh^{(k)}\big)}, \quad B(\theta)=B\left(\mathrm{h}\right):=\mb G_\infty\big(\mh^{(k)}\big),\;\;\theta=\mh.
 $$

The idea now is to apply Proposition \ref{prop_uniform_local2} with $\Theta=\mathrm{H}(\tau)$. Put
$$
\underline{G}_\mathbf{n}[\tau]=\inf_{\mh\in\mH\big(\mathbf{n},\tau\big)}G_\infty\big(\mh^{(k)}\big),
$$
that yields $\underline{A}=\ma\sqrt{\underline{G}_\mathbf{n}[\tau]}$ and $\underline{B}=\mb \underline{G}_\mathbf{n}[\tau]$.
Choose  $s_1=s_2=s^*$.
To apply Proposition \ref{prop_uniform_local2} one has to bound from above the function
\begin{eqnarray*}
&&\mathcal{E}_{\vec{s}}(u,v)=
e^{(\ra)}_{s_1}\Big(\underline{A}u,\Theta_{A}\big(\underline{A}u\big)\Big)+ e^{(\rb)}_{s_2}\Big(\underline{B}v,\Theta_{B}\big( \underline{B}v\big)\Big),\quad u,v\geq 1,
\end{eqnarray*}
defined in (\ref{eq:def-of-function-e}). Here, in our case, $\ra=\ra_\rf$, $\rb=\rb_\infty$ and
\begin{eqnarray*}
\Theta_{A}\big(\underline{A}u\big)&=&\left\{\mh\in\mH\big(\mathbf{n},\tau\big):\;\; G_\infty\big(\mh^{(k)}\big)\leq u^2\underline{G}_\mathbf{n}[\tau]\right\}\times
\big[\mathbf{n_1}/\mathbf{n_2}, 1\big];
\\
\Theta_{B}\big(\underline{B}v\big)&=&\left\{\mh\in\mH\big(\mathbf{n},\tau\big):\;\; G_\infty\big(\mh^{(k)}\big)\leq v\underline{G}_\mathbf{n}[\tau]\right\}\times\big[\mathbf{n_1}/\mathbf{n_2}, 1\big].
\end{eqnarray*}
To compute  the  function $\widetilde{\mathcal{E}}$ let us make several remarks.

$\mathbf{3^0a.}$ First remind that
\begin{eqnarray*}
%\label{eq11:proof-theorem4-empirical}
&&\;\;e^{(\ra_{\rf})}_{s^*}
\Big(\underline{A}u,\Theta_{A}\big(\underline{A}u\big)\Big)=
\sup_{\delta>0}\delta^{-2}
\mathfrak{E}_{\Theta_{A}\big(\underline{A}u\big),\;\mathrm{a}_{\rf}}\left(\underline{A}u(48\delta)^{-1}s^*(\delta)\right);
\\
&&e^{(\rb_\infty)}_{s^*}\Big(\underline{B}v,\Theta_{B}\big(\underline{B}v\big)\Big)=
\sup_{\delta>0}\delta^{-1}
\mathfrak{E}_{\Theta_{B}\big(\underline{B}v\big),\;\mathrm{b}_{\infty}}\left(\underline{B}v(48\delta)^{-1}s^*(\delta)\right).
\end{eqnarray*}
We have in view of (\ref{eq3:proof-theorem4-empirical}) and (\ref{eq4:proof-theorem4-empirical}) that for any $\mathrm{h},\overline{\mathrm{h}}\in\mathrm{H}(\tau)$
$$
\ra_\rf\big(\mathrm{h},\overline{\mathrm{h}}\big)\leq
\ma\sqrt{\left[G_\infty\Big(\mh^{(k)}\Big)\vee
G_\infty\Big(\overline{\mh}^{(k)}\Big)\right]},\quad \rb_\infty\big(\mathrm{h},\overline{\mathrm{h}}\big)
\leq \mb \left[G_\infty\Big(\mh^{(k)}\Big)\vee G_\infty\Big(\overline{\mh}^{(k)}\Big)\right],
$$
where we have also used again that $C_D\geq 2$. Therefore,
$$
\sup_{\mathrm{h},\overline{\mathrm{h}}\in\Theta_{A}\big(\underline{A}u\big)}\ra_\rf\big(\mathrm{h},\overline{\mathrm{h}}\big)\leq \underline{A}u,\quad \sup_{\mathrm{h},\overline{\mathrm{h}}\in\Theta_{B}\big(\underline{B}v\big)}\rb_\infty\big(\mathrm{h},\overline{\mathrm{h}}\big)\leq \underline{B}v.
$$
It yields for any $\delta\leq \delta_*$, where remind
$\delta_*$ be the smallest solution of the equation  $(48\delta)^{-1}s^{*}(\delta)=1$,
$$
\mathfrak{E}_{\Theta_{A}\big(\underline{A}u\big),\;\mathrm{a}_{\rf}}\left(\underline{A}u(48\delta)^{-1}s^*(\delta)\right)=0,\quad
\mathfrak{E}_{\Theta_{B}\big(\underline{B}v\big),\;\mathrm{b}_{\infty}}\left(\underline{B}v(48\delta)^{-1}s^*(\delta)\right)=0
$$
and, therefore
\begin{eqnarray}
\label{eq110:proof-theorem4-empirical}
&&e^{(\ra_{\rf})}_{s^*}
\Big(\underline{A}u,\Theta_{A}\big(\underline{A}u\big)\Big)=
\sup_{\delta>\delta_*}\delta^{-2}
\mathfrak{E}_{\Theta_{A}\big(\underline{A}u\big),\;\mathrm{a}_{\rf}}\left(\underline{A}u(48\delta)^{-1}s^*(\delta)\right);
\\
\label{eq1111:proof-theorem4-empirical}
&&e^{(\rb_\infty)}_{s^*}\Big(\underline{B}v,\Theta_{B}\big(\underline{B}v\big)\Big)=
\sup_{\delta>\delta_*}\delta^{-1}
\mathfrak{E}_{\Theta_{B}\big(\underline{B}v\big),\;\mathrm{b}_{\infty}}\left(\underline{B}v(48\delta)^{-1}s^*(\delta)\right).
\end{eqnarray}

\quad $\mathbf{3^0b.}$
For any $t\geq 1$ put
$
\mH_{1}^k(t,\mathbf{n})=\left\{\mh^{(k)}\in\mH_1^k(\mathbf{n}):\;\;G_\infty\big(\mh^{(k)}\big)\leq \underline{G}_\mathbf{n} t \right\}
$
and note that the following obvious inclusions hold:
\begin{eqnarray}
\label{eq11:proof-theorem4-empirical}
&&\Theta_{A}\big(\underline{A}u\big)\subseteq \mH_{1}^k\left(u^2\underline{G}_\mathbf{n}[\tau]\underline{G}^{-1}_\mathbf{n},\mathbf{n}\right)\times\mH_{k+1}^m\times\big[\mathbf{n_1}/\mathbf{n_2},1\big];
\\
\label{eq11-1:proof-theorem4-empirical}
&&\Theta_{B}\big(\underline{B}v\big)\subseteq \mH_{1}^k\left(v\underline{G}_\mathbf{n}[\tau]\underline{G}^{-1}_\mathbf{n},\mathbf{n}\right)\times\mH_{k+1}^m\times\big[\mathbf{n_1}/\mathbf{n_2}, 1\big].
\end{eqnarray}
%Here we have taken into account  that $\mathbf{n_2}\leq 2\mathbf{n_1}$.
For any $\e>0$ denote by $\mathfrak{N}_t^{(k)}(\e)$ the minimal number of $\varrho_\mathbf{n}^{(k)}$-balls  of radius $\e$ needed to cover $\mH_{1}^k(t,\mathbf{n})$ , $\mathfrak{N}_j(\e),\;j=\overline{k+1,m},$ the minimal number of $\varrho_j$-balls  of radius $\delta$ needed to cover $\mH_j$ and let $\mathfrak{N}(\e)$ be the minimal number of $\mw$-balls  of radius $\e$ needed to cover $\big[\mathbf{n_1}/\mathbf{n_2}, 1\big]$.

Let $\bH$ be an arbitrary subset of $\mH_{1}^k\left(t,\mathbf{n}\right)\times\mH_{k+1}^m\times[1/2, 1]$.
It is evident  that for any given $\ve^{(k)}>0$, $\ve_j>0,\;j=\overline{k+1,m}$ and $\ve>0$ one can construct
a net $\left\{\mathrm{h}(\mathbf{i}),\;\mathbf{i}=\overline{1,\mathbf{I}\big[\bH\big]}\right\}\subset\bH$ such that $\forall\mathrm{h}=(\mh,y)\in\bH\quad\exists  \mathbf{i}\in\left\{1,\ldots, \mathbf{I}\big[\bH\big]\right\}$
\begin{eqnarray}
\label{eq111:proof-theorem4-empirical}
&& \varrho_\mathbf{n}^{(k)}\left(\mh^{(k)},\mh^{(k)}(i)\right)\leq \ve^{(k)},\quad
\varrho_j\Big(\mh_j,\mh_j(i)\Big)\leq \ve_j,\; j=\overline{k+1,m},\quad \mw(y,y(\mathbf{i}))\leq\ve;
\\
\label{eq112:proof-theorem4-empirical}
&& \mathbf{I}\big[\bH\big]\leq \mathfrak{N}(\ve/2)\;\mathfrak{N}_t^{(k)}\left(\ve^{(k)}/2\right)\prod_{j=k+1}^{m}\mathfrak{N}_j\big(\ve_j/2\big),\quad \forall \bH\subseteq \mH_{1}^k\left(t,\mathbf{n}\right)\times\mH_{k+1}^m\times\big[\mathbf{n_1}/\mathbf{n_2}, 1\big].
\end{eqnarray}
Moreover we obtain from (\ref{eq5:proof-theorem4-empirical}) and (\ref{eq6:proof-theorem4-empirical}) for any $u,v\geq 1$
\begin{eqnarray*}
%\label{eq12:proof-theorem4-empirical}
\ra_{\rf}(\mathrm{h},\overline{\mathrm{h}})&\leq& \ma \bigg\{\bigg(\underline{G}_\mathbf{n}[\tau] u^2\;\varrho_\mathbf{n}^{(k)}\Big(\mh^{(k)},\overline{\mh}^{(k)}\Big)+
\sum_{j=k+1}^{m}\widetilde{L}_j\left(\underline{G}_\mathbf{n} [\tau] u^2\right)\varrho_j\big(\mh_j,\overline{\mh}_j\big)\bigg)^{1/2}
\\
&&\hskip0.5cm +\chi u\sqrt{\underline{G}_\mathbf{n} [\tau] }
\Big[1-\ln{\big(\mathfrak{w}(y,\bar{y})\big)}\Big]^{-\boldsymbol{b}}
\bigg\},\;\quad \forall \mathrm{h},\overline{\mathrm{h}}\in \Theta_{A}\big(\underline{A}u\big);
\\
\rb_{\infty}(\mathrm{h},\overline{\mathrm{h}})&\leq& \mb \bigg\{\underline{G}_\mathbf{n} [\tau] v\;\varrho_\mathbf{n}^{(k)}\Big(\mh^{(k)},\overline{\mh}^{(k)}\Big)+
\sum_{j=k+1}^{m}\widetilde{L}_j\left(\underline{G}_\mathbf{n} [\tau] v\right)\varrho_j\big(\mh_j,\overline{\mh}_j\big)
\\
&&\hskip0.5cm
+\chi\underline{G}_\mathbf{n} [\tau] v\;
\Big[1-\ln{\big(\mathfrak{w}(y,\bar{y})\big)}\Big]^{-\boldsymbol{b}}
\bigg\},\qquad \forall \mathrm{h},\overline{\mathrm{h}}\in \Theta_{A}\big(\underline{A}u\big).
\end{eqnarray*}
Thus,  putting $t=t_1:=u^2\underline{G}_\mathbf{n}[\tau]\underline{G}^{-1}_\mathbf{n}$ and  choosing for any $\varsigma>0$
$$
\ve^{(k)}=\frac{\varsigma^2}{2\ma^{2} m\underline{G}_\mathbf{n}[\tau] u^2},\quad
\ve_j=\frac{\varsigma^2}{2\ma^{2} m\widetilde{L}_j\Big(\underline{G}_\mathbf{n}(\tau) u^2\Big)},
\quad \ve=e^{-\left(\frac{2u\ma \sqrt{\underline{G}_\mathbf{n}[\tau]}}{\varsigma}\right)^{1/\boldsymbol{b}}}.
%\quad \ve=e^{-\left(2u\ma \sqrt{\underline{G}_\infty[\tau]}\varsigma^{-1}\right)^{1/\beta}}.
$$
we obtain in view of (\ref{eq11:proof-theorem4-empirical})  and (\ref{eq111:proof-theorem4-empirical}) with $\bH=\Theta_{A}\big(\underline{A}u\big)$
\begin{eqnarray}
\label{eq113:proof-theorem4-empirical}
\forall\mathrm{h}\in\Theta_{A}\big(\underline{A}u\big)\quad\exists \mathbf{i}\in\left\{1,\ldots, \mathbf{I}\Big[\Theta_{A}\big(\underline{A}u\big)\Big]\right\}:\;\; \ra_\rf\left(\mathrm{h},\mathrm{h}(\mathbf{i})\right)\leq \varsigma.
\end{eqnarray}
Putting $t=t_2:=v\underline{G}_\mathbf{n}[\tau]\underline{G}^{-1}_\mathbf{n}$ and choosing
$$
\ve^{(k)}=\frac{\varsigma}{2\mb m\underline{G}_\mathbf{n}[\tau] v},\quad \ve_j=\frac{\varsigma}
{2\mb m\widetilde{L}_j\Big(\underline{G}_\mathbf{n}[\tau] v\Big)},\quad \ve=e^{-\left(\frac{2v\mb \underline{G}_\mathbf{n}[\tau]}{\varsigma}\right)
^{1/\boldsymbol{b}}}
$$
we obtain in view of (\ref{eq11:proof-theorem4-empirical}) and (\ref{eq111:proof-theorem4-empirical}) with $\bH=\Theta_{B}\big(\underline{B}v\big)$
\begin{eqnarray}
\label{eq114:proof-theorem4-empirical}
\forall\mathrm{h}\in\Theta_{B}\big(\underline{B}v\big)\quad\exists \mathbf{i\mathbf{}}\in\left\{1,\ldots, \mathbf{I}\Big[\Theta_{B}\big(\underline{B}v\big)\Big]\right\}:\;\; \rb_\infty\left(\mathrm{h},\mathrm{h}(\mathbf{i})\right)\leq \varsigma.
\end{eqnarray}
We get from (\ref{eq112:proof-theorem4-empirical}), (\ref{eq113:proof-theorem4-empirical}) and (\ref{eq114:proof-theorem4-empirical})
 for any $\varsigma>0$
 \begin{eqnarray}
\label{eq12:proof-theorem4-empirical}
&&\mathfrak{E}_{\Theta_{A}\big(\underline{A}u\big),\ra_{\rf}}(\varsigma)\leq\mathfrak{E}_{\mH_{1}^k(t_1,\mathbf{n}),\varrho_\mathbf{n}^{(k)}}
\left(\frac{\varsigma^{2}}{4m\ma^{2}\underline{G}_\mathbf{n}[\tau] u^{2}}\right)+\sum_{j=k+1}^{m}
\mathfrak{E}_{\mH_j,\varrho_j}\left(\frac{\varsigma^2}{4m\ma^{2}\widetilde{L}_j\Big(\underline{G}_\mathbf{n} [\tau] u^2\Big)}\right)
\\&&\hskip2.8cm
+\mathfrak{E}_{[\mathbf{n_1}/\mathbf{n_2}, 1],\mw}\bigg(2^{-1}
\exp{\bigg\{-\left(2u\ma \sqrt{\underline{G}_\mathbf{n}[\tau]}\varsigma^{-1}\right)^{1/\boldsymbol{b}}\bigg\}}\bigg);
\nonumber
\\
\label{eq13:proof-theorem4-empirical}
&&\mathfrak{E}_{\Theta_{B}\big(\underline{B}v\big),\rb_\infty}(\varsigma)\leq\mathfrak{E}_{\mH_{1}^k(t_2,\mathbf{n}),\varrho_\mathbf{n}^{(k)}}
\left(\frac{\varsigma}{4m\mb\underline{G}_\mathbf{n}[\tau] v}\right)+\sum_{j=k+1}^{m}
\mathfrak{E}_{\mH_j,\varrho_j}\left(\frac{\varsigma}{4m\mb\widetilde{L}_j\Big(\underline{G}_\mathbf{n}[\tau] v\Big)}\right)
\\&&\hskip2.9cm
+\mathfrak{E}_{[\mathbf{n_1}/\mathbf{n_2}, 1],\mw}\Big(2^{-1}
\exp{\Big\{-\left(2v\mb \underline{G}_\mathbf{n}[\tau]\varsigma^{-1}\right)^{1/\boldsymbol{b}}\Big\}}\Big).
\nonumber
\end{eqnarray}

 $\mathbf{4^0.}$
%Put $\cL_k(z)=\sum_{j=k+1}^{m}\left[\log_2{\big\{\widetilde{L}_j\left(z\right)\big\}}\right]_+$.
We get in view of Assumption \ref{ass:sec:totaly_bounded_case}
 \begin{eqnarray}
\label{eq14:proof-theorem4-empirical}
&&\sum_{j=k+1}^{m}
\mathfrak{E}_{\mH_j,\varrho_j}\left(\frac{\varsigma^2}{4m\ma^2\widetilde{L}_j\Big(\underline{G}_\mathbf{n} [\tau] u^2\Big)}\right)\leq
N\sum_{j=k+1}^{m}\left(\left[\log_2{\left\{4\ma^2 m R\widetilde{L}_j\Big(\underline{G}_\mathbf{n}[\tau] u^2\Big)\varsigma^{-2}\right\}}\right]_++1\right);
%N\cL_k\left(\underline{G}_\infty t\right)+
%N(m-k)\left(\left[\log_2{\left\{2\ma m R\varsigma^{-2}\right\}}\right]_++1\right);
\\
\label{eq15:proof-theorem4-empirical}
&&\sum_{j=k+1}^{m}
\mathfrak{E}_{\mH_j,\varrho_j}\left(\frac{\varsigma}{4m\mb\widetilde{L}_j\Big(\underline{G}_\mathbf{n} [\tau] v\Big)}\right)\leq
N\sum_{j=k+1}^{m}\left(\left[\log_2{\left\{4\mb m R\widetilde{L}_j\Big(\underline{G}_\mathbf{n} [\tau] v\Big)\varsigma^{-1}\right\}}\right]_++1\right).
%N\cL_k\left(\underline{G}_\infty t\right)+
%N(m-k)\left(\left[\log_2{\left\{2\mb m R\varsigma^{-1}\right\}}\right]_++1\right).
\end{eqnarray}
Taking into account that $\mathfrak{E}_{[\mathbf{n_1}/\mathbf{n_2}, 1],\mw}(\cdot)\equiv 0$, if $\mathbf{n_1}=\mathbf{n_2}$, and $\mathfrak{E}_{[\mathbf{n_1}/\mathbf{n_2}, 1],\mw}(\e)\leq \ln{\big(2/\e^2\big)}$ for any $\e\in (0,1]$ and any $\mathbf{n_2}\leq2\mathbf{n_1}$,
we have
 \begin{eqnarray}
\label{eq140:proof-theorem4-empirical}
&&\mathfrak{E}_{[\mathbf{n_1}/\mathbf{n_2}, 1],\mw}\bigg(2^{-1}\exp{\bigg\{-\left(2u\ma
\sqrt{\underline{G}_\mathbf{n}[\tau]}\varsigma^{-1}\right)^{1/\beta}\bigg\}}\bigg)=\chi\bigg(2\ln(2)+
2\left(2u\ma \sqrt{\underline{G}_\mathbf{n}[\tau]}\varsigma^{-1}\right)^{\frac{1}{\boldsymbol{b}}}\bigg);
\\
\label{eq150:proof-theorem4-empirical}
&&\mathfrak{E}_{[\mathbf{n_1}/\mathbf{n_2}, 1],\mw}\Big(2^{-1}\exp{\Big\{-\left(2v\mb \underline{G}_\mathbf{n}[\tau]\varsigma^{-1}\right)^{1/\beta}\Big\}}\Big)=
\chi\bigg(2\ln(2)+2\left(2v\mb \underline{G}_\mathbf{n}[\tau]\varsigma^{-1}\right)^{\frac{1}{\boldsymbol{b}}}\bigg).
\end{eqnarray}
%Here $\gamma=0$ if $\mathbf{n_1}=\mathbf{n_2}$ and  $\gamma=1$ if $\mathbf{n_1}\neq\mathbf{n_2}$.

Let us now bound from above $\mathfrak{E}_{\mH_{1}^k(t,\mathbf{n}),\varrho_\mathbf{n}^{(k)}}$. First we note that in view of Assumption \ref{ass:bounded_case} ($\mathbf{i}$)
\begin{equation}
\label{eq16:proof-theorem4-empirical}
\mH_{1}^k(t,\mathbf{n})\subseteq \left\{\mh_1\in\mH_1(\mathbf{n}):\;G_{1,\mathbf{n}}(\mh_1)\leq t\underline{G}_{1,\mathbf{n}}\right\}\times\cdots\times\left\{\mh_k\in\mH_k(\mathbf{n}):\;G_{k,\mathbf{n}}(\mh_k)\leq t\underline{G}_{k,\mathbf{n}}\right\}.
\end{equation}
Consider the hyper-rectangle $\cZ(t)=\left[\underline{G}_{1,\mathbf{n}},t\underline{G}_{1,\mathbf{n}}\right]\times
\cdots\times\left[\underline{G}_{k,\mathbf{n}},t\underline{G}_{k,\mathbf{n}}\right],\;t\geq 1,$ which we equip with the metrics
$$
\mm^{(k)}\big(z,z^\prime\big)=\max_{i=\overline{1,k}}\left|\ln\big(z_i)-\ln\big(z_i^\prime\big)\right|,\quad z,z^\prime\in\cZ(t),
$$
where $z_i,z_i^\prime,\;i=\overline{1,k}$ are the coordinates of $z,z^\prime$ respectively.
It  easily seen that for any $\varsigma>0$
$$
\mE_{\cZ(t),\mm^{(k)}}(\varsigma)\leq k\left[\ln{\ln{t}}-\ln{\ln{(1+\varsigma)}}\right]_+\leq
k\left(\ln{\big(1+\ln{t}\big)}+\big[1+\ln{(1/\varsigma)}\big]_+\right).
$$
It yields together with (\ref{eq16:proof-theorem4-empirical}) in view of obvious inequality $\mathfrak{E}_{\mH_{1}^k(t,\mathbf{n}),\varrho_\mathbf{n}^{(k)}}(\varsigma)\leq\mE_{\cZ(t),\mm^{(k)}}(\varsigma/2)$
\begin{equation}
\label{eq17:proof-theorem4-empirical}
\mathfrak{E}_{\mH_{1}^k(t,\mathbf{n}),\varrho_\mathbf{n}^{(k)}}(\varsigma)\leq
k\left(\ln{\big(1+\ln{t}\big)}+\big[1+\ln{(2/\varsigma)}\big]_+\right).
\end{equation}
We obtain from (\ref{eq17:proof-theorem4-empirical})
\begin{eqnarray}
\label{eq18:proof-theorem4-empirical}
&&
\mathfrak{E}_{\mH_{1}^k(t_1,\mathbf{n}),\varrho_\mathbf{n}^{(k)}}\left(\frac{\varsigma^{2}}{4m\ma^2\underline{G}_\mathbf{n}[\tau] u^{2}}\right)\leq
k\left(\ln{\big(1+\ln{t_1}\big)}+\left[1+\ln{\left(8m\ma^2\underline{G}_\mathbf{n}[\tau] u^{2}\varsigma^{-2}\right)}\right]_+\right);
\\
\label{eq19:proof-theorem4-empirical}
&&\mathfrak{E}_{\mH_{1}^k(t_2,\mathbf{n}),\varrho_\mathbf{n}^{(k)}}
\left(\frac{\varsigma}{4m\mb\underline{G}_\mathbf{n}[\tau] v}\right)\leq
k\left(\ln{\big(1+\ln{t_2}\big)}+\left[1+\ln{\left(8m\mb\underline{G}_\mathbf{n}[\tau] v\varsigma^{-1}\right)}\right]_+\right).
\end{eqnarray}
Putting $\widehat{L}_j(z)=z^{-1}\widetilde{L}_j(z)=\max\big\{z^{-1}L_j(z),1\big\}$, we get from (\ref{eq110:proof-theorem4-empirical}), (\ref{eq12:proof-theorem4-empirical}), (\ref{eq14:proof-theorem4-empirical}), (\ref{eq140:proof-theorem4-empirical}) and (\ref{eq18:proof-theorem4-empirical})
\begin{eqnarray}
\label{eq20:proof-theorem4-empirical}
&&e^{(\ra_{\rf})}_{s^*}
\Big(\underline{A}u,\Theta_{A}\big(\underline{A}u\big)\Big)
\leq
k\delta^{-2}_*\ln{\left(1+\ln{\left(u^2\underline{G}_\mathbf{n}[\tau]\underline{G}^{-1}_\mathbf{n}\right)}\right)}
+N\delta^{-2}_*\sum_{j=k+1}^{m}\log_2{\left\{\widehat{L}_j\left(\underline{G}_\mathbf{n}[\tau]u^2\right)\right\}}
\nonumber\\
&&\qquad +\sup_{\delta>\delta_*}\delta^{-2}\left\{k\left[1+\ln{\left(\frac{9216m\delta^{2}}{[s^*(\delta)]^{2}}\right)}\right]_++N(m-k)
\left(\left[\log_2{\left\{\left(\frac{4608m R\delta^{2}}{[s^*(\delta)]^{2}}\right)\right\}}\right]_++1\right)\right\}
\nonumber
\\
&&\qquad +\chi\bigg(2\delta^{-2}_*\ln(2)+2\sup_{\delta>\delta_*}\delta^{-2}\big(96\delta\big/s^*(\delta)\big)^{\frac{1}{\boldsymbol{b}}}\bigg)
\nonumber
\\
&&\qquad=k\delta^{-2}_*\ln{\left(1+\ln{\left(u^2\underline{G}_\mathbf{n}[\tau]\underline{G}^{-1}_\mathbf{n}\right)}\right)}+
N\delta^{-2}_*\sum_{j=k+1}^{m}\log_2{\left\{\widehat{L}_j\left(\underline{G}_\mathbf{n}[\tau]u^2\right)\right\}}+C^{(1)}_{N,R,m,k}+\chi \mathbf{a_{\boldsymbol{b}}},
\end{eqnarray}
where, remind,
$
\mathbf{a_{\boldsymbol{b}}}=2\delta^{-2}_*\ln(2)+2\sup_{\delta>\delta_*}(\delta^2\wedge\delta)^{-1}
\big(96\delta\big/s^*(\delta)\big)^{\frac{1}{{\boldsymbol{b}}}}.
$
Note that $\mathbf{a_{\boldsymbol{b}}}<\infty$ since $\boldsymbol{b}>1$.

Repeating these computations we get from (\ref{eq1111:proof-theorem4-empirical}), (\ref{eq13:proof-theorem4-empirical}), (\ref{eq15:proof-theorem4-empirical}), (\ref{eq150:proof-theorem4-empirical}) and (\ref{eq19:proof-theorem4-empirical})
\begin{eqnarray}
\label{eq21:proof-theorem4-empirical}
e^{(\rb_{\infty})}_{s^*}
\Big(\underline{B}v,\Theta_{B}\big(\underline{B}v\big)\Big)&\leq& k\delta^{-1}_*\ln{\left(1+\ln{\left(v\underline{G}_\mathbf{n}[\tau]\underline{G}^{-1}_\mathbf{n}\right)}\right)}
\nonumber\\
&&\;+
N\delta^{-1}_*\sum_{j=k+1}^{m}\log_2{\left\{\widehat{L}_j\left(\underline{G}_\mathbf{n}[\tau]u^2\right)\right\}}+C^{(2)}_{N,R,m,k}+\chi
\mathbf{a_{\boldsymbol{b}}},
\end{eqnarray}

%We note that the condition $\beta>1$ is necessary for finiteness of $\mathbf{a_\beta}$.

We deduce from (\ref{eq20:proof-theorem4-empirical}) and (\ref{eq21:proof-theorem4-empirical}) that  $\widetilde{\mathcal{E}}_{\vec{s}},\;\vec{s}=(s^*,s^*)$ is bounded from above by the function
\begin{eqnarray}
\label{eq22:proof-theorem4-empirical}
\cE(u,v)&\leq&k\delta^{-2}_*\ln{\bigg\{\left(1+\ln{\left(u^2\underline{G}_\mathbf{n}[\tau]\underline{G}^{-1}_\mathbf{n}\right)}\right)
\left(1+\ln{\left(v\underline{G}_\mathbf{n}[\tau]\underline{G}^{-1}_\infty\right)}\right)\bigg\}}
\nonumber\\
&&+N\delta^{-2}_*\sum_{j=k+1}^{m}\log_2{\bigg[\left\{\widehat{L}_j\left(\underline{G}_\mathbf{n}[\tau]u^2\right)\right\}
\left\{\widehat{L}_j\left(\underline{G}_\mathbf{n}[\tau]v\right)\right\}\bigg]}+C_{N,R,m,k}.
\end{eqnarray}
Here we have used that $\delta_*<1$. We note that (\ref{eq22:proof-theorem4-empirical}) implies in particular  Assumption  \ref{ass:parameter_local}
and, therefore, Proposition  \ref{prop_uniform_local2} is applicable  with  $\Theta=\mathrm{H}(\tau)$.

\smallskip

 $\mathbf{5^0.}$
To apply Proposition  \ref{prop_uniform_local2} on $\Theta=\mathrm{H}(\tau)$ we choose  $\e=\sqrt{2}-1$ and bound from above the quantities
\begin{eqnarray*}
%\label{eq:def-price-to-pay-proba}
 P_{\sqrt{2}-1}(\mathrm{h})&:=&4\big[\sqrt{2}-1\big]^{-2}\cE\left(\sqrt{2\underline{G}^{-1}_\mathbf{n}[\tau] G_\infty\big(\mh^{(k)}\big)},\sqrt{2}\underline{G}^{-1}_\mathbf{n}[\tau] G_\infty\big(\mh^{(k)}\big)\right)
\\
&&+
2\ell\left(\sqrt{2\underline{G}^{-1}_\mathbf{n}[\tau] G_\infty\big(\mh^{(k)}\big)}\right)+
2\ell\left(\sqrt{2}\underline{G}^{-1}_\mathbf{n}[\tau] G_\infty\big(\mh^{(k)}\big)\right);
\\*[2mm]
 M_{\sqrt{2}-1,q}(\mathrm{h})&:=&8\big[\sqrt{2}-1\big]^{-2}\cE\left(\sqrt{2\underline{G}^{-1}_\mathbf{n}[\tau] G^{-1}_\infty\big(\mh^{(k)}\big)},\sqrt{2}\underline{G}^{-1}_\mathbf{n}[\tau] G_\infty\big(\mh^{(k)}\big)\right)
\\
&&+
2\big(\sqrt{2}-1+q\big)\ln{\left(\sqrt{2\underline{G}^{-1}_\mathbf{n}[\tau] G_\infty\big(\mh^{(k)}\big)}\sqrt{2}\underline{G}^{-1}_\mathbf{n}[\tau] G_\infty\big(\mh^{(k)}\big)\right)},
\end{eqnarray*}
where remind
$
\ell(u)=\ln{\left\{1+\ln{(u)}\right\}}+2\ln{\left\{1+\ln{\left\{1+\ln{(u)}\right\}}\right\}}.
$

\smallskip

Taking into account that $
\ell(u)\leq 3\ln{\left\{1+\ln{(u)}\right\}}, u\geq 1,$ $\big[\sqrt{2}-1\big]^{-2}\leq 9$ and that
$\underline{G}_\mathbf{n}[\tau]\geq \underline{G}_\mathbf{n}$ for any $\tau$,
we obtain from (\ref{eq22:proof-theorem4-empirical})
\begin{eqnarray*}
%\label{eq:def-price-to-pay-proba}
P_{\sqrt{2}-1}(\mathrm{h})&\leq& \big[72k\delta^{-2}_*+12\big]\ln{\bigg\{1+\ln{\left(2G_\infty\big(\mh^{(k)}\big)\underline{G}^{-1}_\mathbf{n}\right)}
\bigg\}}
\nonumber\\
&&+72N\delta^{-2}_*\sum_{j=k+1}^{m}\log_2{\left\{\widehat{L}_j\left(2G_\infty\big(\mh^{(k)}\big)\right)\right\}
}+36C_{N,R,m,k}=:2P\big(\mh^{(k)}\big);
\\*[2mm]
%\label{eq:def-price-to-pay-moments}
\quad M_{\sqrt{2}-1,q}(\mathrm{h})&\leq&\big[144k\delta^{-2}_*+3(1+q)\big]\ln{\left(2G_\infty\big(\mh^{(k)}\big)\underline{G}^{-1}_\mathbf{n}\right)}
\nonumber\\
&&+144N\delta^{-2}_*\sum_{j=k+1}^{m}\log_2{\left\{\widehat{L}_j\left(2G_\infty\big(\mh^{(k)}\big)\right)\right\}
}+72C_{N,R,m,k}=2M_q\big(\mh^{(k)}\big).
\end{eqnarray*}
We remark that $P$ et $M_q$ are independent of $\tau$ and $y$.

Put for any $z\geq 0$ and any $\mathrm{h}\in\mathrm{H}(\tau)$
\begin{eqnarray*}
%\label{eq:def-of-upper-function-proba}
\check{\mathrm{V}}_\tau^{(z)}\big(\mh^{(k)}\big)&=&2\sqrt{2}\ma\sqrt{ G_\infty\big(\mh^{(k)}\big)\big[P\big(\mh^{(k)}\big)+z\big]}+
4\mb G_\infty\big(\mh^{(k)}\big)\Big[P\big(\mh^{(k)}\big)+z\Big];
\\*[2mm]
%\label{eq:def-of-upper-function-moments}
\check{\mathrm{U}}_\tau^{(z,q)}\big(\mh^{(k)}\big)&=&2\sqrt{2}\ma\sqrt{ G_\infty\big(\mh^{(k)}\big)\big[M_q\big(\mh^{(k)}\big)+z\big]}+
4\mb G_\infty\big(\mh^{(k)}\big)\Big[M_q\big(\mh^{(k)}\big)+z\Big].
 \end{eqnarray*}
where remind $\ma=2\sqrt{\tau}(\mathbf{n_1})^{-1/2}\big(\sqrt{C_D}\vee[\chi\boldsymbol{c_b}]\big),\;
\mb=\frac{4\big(C_D\vee 8e\big)\ln^{\beta}{(\mathbf{n_2})}}{3\mathbf{n_1}}$.

We conclude that Proposition \ref{prop_uniform_local2} is applicable with $\check{\mathrm{V}}_\tau^{(z)}$ and $\check{\mathrm{U}}_\tau^{(z,q)}$.
Put for any $n\in \{\mathbf{n_1},\mathbf{n_1}+1,\ldots,\mathbf{n_2}\}$
$$
\ma(n)=2\sqrt{2\tau}(n)^{-1/2}\big(\sqrt{C_D}\vee[\chi\boldsymbol{c_b}]\big),\;\;
\mb(n)=\frac{8\big(C_D\vee 8e\big)\ln^{\beta}{(2n)}}{3n}
$$
and define
\begin{eqnarray*}
%\label{eq:def-of-upper-function-proba}
\mathrm{V}_\tau^{(z)}\big(n,\mh^{(k)}\big)&=&2\sqrt{2}\ma(n)\sqrt{ G_\infty\big(\mh^{(k)}\big)\big[P\big(\mh^{(k)}\big)+z\big]}+
4\mb(n) G_\infty\big(\mh^{(k)}\big)\Big[P\big(\mh^{(k)}\big)+z\Big];
\\*[2mm]
%\label{eq:def-of-upper-function-moments}
\mathrm{U}_\tau^{(z,q)}\big(n,\mh^{(k)}\big)&=&2\sqrt{2}\ma(n)\sqrt{ G_\infty\big(\mh^{(k)}\big)\big[M_q\big(\mh^{(k)}\big)+z\big]}+
4\mb(n) G_\infty\big(\mh^{(k)}\big)\Big[M_q\big(\mh^{(k)}\big)+z\Big].
 \end{eqnarray*}
It is easily seen that $\ma(n)\geq \ma,\;\;\mb(n)\geq \mb$ for any $n\in \big\{\mathbf{n_1},\ldots, \mathbf{n_2}\big\}$ since $\mathbf{n_2}\leq 2\mathbf{n_1}$. Therefore,
$$
\mathrm{V}_\tau^{(z)}\big(n,\mh^{(k)}\big)\geq \check{\mathrm{V}}_\tau^{(z)}\big(\mh^{(k)}\big),\quad \mathrm{U}_\tau^{(z,q)}\big(n,\mh^{(k)}\big)\geq \check{\mathrm{U}}_\tau^{(z,q)}\big(\mh^{(k)}\big).
$$
 It remains to remind that  $
\xi_\mh(n)=\boldsymbol{\xi}\big(n/\mathbf{n_2},\mh\big)
$ for any $n\in \{\mathbf{n_1},\mathbf{n_1}+1,\ldots,\mathbf{n_2}\}$ and any $\mh\in\boldsymbol{\mH}$.
 All saying above allows us to assert that Proposition, \ref{prop_uniform_local2} is applicable to $|\xi_{\mh}(n)|$ on
 $\mathbf{H}(\tau):=\{\mathbf{n_1},\mathbf{n_1}+1,\ldots,\mathbf{n_2}\}\times\mH(\mathbf{n},\tau)$ for any $\tau>0 $ with
%that means in particular that
$\mathrm{V}_\tau^{(z)}(\cdot,\cdot)$ and $\mathrm{U}_\tau^{(z,q)}(\cdot,\cdot)$.

Thus, putting $\mathbf{h}=(n,\mh)$ we obtain for any $\tau>0$, any $z\geq 1$ and any $q\geq 1$
\begin{gather}
\label{eq23:proof-theorem4-empirical}
\bP_\rf\left\{\sup_{\mathbf{h}\in \mathbf{H}(\tau)}\Big[|\xi_{\mh}(n)|-\mathrm{V}_\tau^{(z)}\big(n,\mh^{(k)}\big)\Big]\geq 0\right\}
\leq 4\left[1+\Big[\ln{\left\{1+2^{-1}\ln{2}\right\}}\Big]^{-2}\right]^{2}\exp{\left\{-z\right\}};
\\*[2mm]
\label{eq24:proof-theorem4-empirical}
\bE_\rf\left\{\sup_{\mathbf{h}\in \mathbf{H}(\tau)}\Big[|\xi_{\mh}(n)|-\mathrm{U}_\tau^{(z,q)}\big(n,\mh^{(k)}\big)\Big]\right\}^q_+
%\\*[2mm]&&
\leq 2^{(5q/2)+3}3^{q+4}\Gamma(q+1)\;\;\big[\underline{A}\vee\underline{B}\big]^q\exp{\left\{-z\right\}},
\end{gather}
where, remind, $\underline{A}=\ma\sqrt{ \underline{G}_\mathbf{n}[\tau]}$ and $\underline{B}=\mb \underline{G}_\mathbf{n}[\tau]$.

To get the statements of the theorem we will have to choose   $z$.  This, in its turn, will be done for $\mathrm{V}_\tau$ and $\mathrm{U}_\tau$ differently in dependence on the  values of the parameter $\tau$.

\smallskip

 $\mathbf{6^0.}$
Let $\mathbf{r}\in\bN$ be fixed and for any  $ r\in\bN^*$ put $\tau_r=e^{r-\mathbf{r}}$. For any $r\in\bN^*$ denote
$
\widehat{\mH}(r)=\mH\big(\mathbf{n},\tau_{r}\big)\setminus\mH\big(\mathbf{n},\tau_{r-1}\big),
$
$\widehat{\mH}(0)=\mH\big(\mathbf{n},\tau_{0}\big)$ and let
$\widehat{\mathbf{H}}(r):=\{\mathbf{n_1},\mathbf{n_1}+1,\ldots,\mathbf{n_2}\}\times\widehat{\mH}(r)$.

\paragraph{Probability bound}

For any $u\geq 1$  put
$
z_r(u)=u+2\ln{\big(1+|r-\mathbf{r}|\big)}
$
and remark that
$$
z_r(u)=\left\{
\begin{array}{cc}
u+2\ln{\left(\left|\ln{(\tau_{r-1})}\right|\right)},\quad &r\leq\mathbf{r};
\\
u+2\ln{\left(1+\left|\ln{(\tau_{r})}\right|\right)},\quad &r\geq \mathbf{r}.
\end{array}
\right.
$$
We have  for any $r\in\bN$ and any $\mh\in\widehat{\mH}(r)$
\begin{eqnarray*}
 \tau_{0}=F_{\mathbf{n_2},\mathbf{r}}(\mh)&\Rightarrow& z_{0}(u)= u+2\ln{\Big\{1+\left|\ln{\left(F_{n,\mathbf{r}}(\mh)\right)}\right|\Big\}};
 \\
 \tau_{r-1}\leq F_{\mathbf{n_2}}(\mh)=F_{\mathbf{n_2},\mathbf{r}}(\mh)&\Rightarrow& z_{r}(u)\leq u+2\ln{\Big\{\left|\ln{\left(F_{n,\mathbf{r}}(\mh)\right)}\right|\Big\}},\;\; 1\leq r\leq \mathbf{r}-1;
\\
 \tau_{r}\geq F_{\mathbf{n_2}}(\mh)=F_{\mathbf{n_2},\mathbf{r}}(\mh)&\Rightarrow& z_{r}(u)\leq u+2\ln{\Big\{1+\left|\ln{\left(F_{\mathbf{n_2},\mathbf{r}}(\mh)\right)}\right|\Big\}},\;\; r\geq \mathbf{r}.
\end{eqnarray*}
Hence, we have for any $r\in\bN$
\begin{gather}
\label{eq249:proof-theorem4-empirical}
z_{r}(u)\leq u+2\ln{\Big\{1+\left|\ln{\left(F_{\mathbf{n_2},\mathbf{r}}(\mh)\right)}\right|\Big\}},\quad \forall \mh\in\widehat{\mH}(r),
\end{gather}
that yields for any $r\in\bN$
\begin{gather}
\label{eq250:proof-theorem4-empirical}
\mathrm{V}_{\tau_{r}}^{(z_{r}(u))}\big(n,\mh^{(k)}\big)\leq \cV_{\mathbf{r}}^{(u)}(n,\mh),\;\;\forall(n,\mh)\in\widehat{\mathbf{H}}(r).
\end{gather}
Here we have also taken into account that
$\tau_{r}\leq eF_{\mathbf{n_2},\mathbf{r}}(\mh),\;\forall\mh\in\widehat{\mH}(r)$ for any $r\in\bN$.

Thus, we get for any $r\in\bN$
and   $u\geq 0$, taking into account (\ref{eq250:proof-theorem4-empirical}), the inclusion  $\widehat{\mathbf{H}}(r)\subseteq\mathbf{H}(\tau_r)$ and   applying   (\ref{eq23:proof-theorem4-empirical}) with $\tau=\tau_r$,
\begin{gather}
\label{eq26:proof-theorem4-empirical}
\bP_\rf\left\{\sup_{(n,\mh)\in\widehat{\mathbf{H}}(r)}\Big[|\xi_{\mh}(n)|-\cV_{\mathbf{r}}^{(u)}(n,\mh)\Big]\geq 0\right\}
\leq \frac{4\left[1+\Big[\ln{\left\{1+2^{-1}\ln{2}\right\}}\Big]^{-2}\right]^{2}\exp{\left\{-u\right\}}}{[1+|r-\mathbf{r}|]^{2}}.
\end{gather}
Since obviously $\widetilde{\mathbf{N}}\times\mH(\mathbf{n})=\cup_{r=0}^{\infty}\widehat{\mathbf{H}}(r)$,
summing up the right hand side of (\ref{eq26:proof-theorem4-empirical})   over $r$, we come to the first assertion of the theorem.
Here we have also used that $16\left[1+\Big[\ln{\left\{1+2^{-1}\ln{2}\right\}}\Big]^{-2}\right]^{2}\leq 2419$ and the fact that
$\widetilde{\mH}(n)\subseteq \mH(\mathbf{n})$ for any $n\in\widetilde{\mathbf{N}}$ in view of Assumption \ref{ass:dependence-on-n} and the definition of the number $\mathbf{n}$.

\paragraph{Moment's bound}

For any $u\geq 1$  put
$$
z_r(u)=u+2\ln{\big(1+|r-\mathbf{r}|\big)}+q\ln{\left(\underline{G}_\mathbf{n}[\tau_r]\underline{G}^{-1}_\mathbf{n}\right)}.
$$
Similarly to (\ref{eq249:proof-theorem4-empirical}) we have for any $r\in\bN$ and any $\mh\in\widehat{\mH}(r)$
$$
z_r(u)\leq u+2\ln{\Big\{1+\left|\ln{\left(F_{\mathbf{n_2},r}(\mh)\right)}\right|\Big\}}
+q\ln{\left(\underline{G}_\mathbf{n}[\tau_r]\underline{G}^{-1}_\mathbf{n}\right)}.
$$
 Moreover, for any $r\in\bN$ by definition
$$
\underline{G}_\mathbf{n}[\tau_r]:=\inf_{\mh\in\mH(\mathbf{n},\tau_r)}G_\infty\left(\mh^{(k)}\right)
$$
and, therefore, for any $\mh\in\widehat{\mH}(r)$
$$
z_r(u)\leq  u+2\ln{\Big\{1+\left|\ln{\left(F_{\mathbf{n_2},r}(\mh)\right)}\right|\Big\}}
+q\ln{\left\{G_\infty\big(\mh^{(k)}\big)\underline{G}^{-1}_\mathbf{n}\right\}}.
$$
Similarly to (\ref{eq250:proof-theorem4-empirical}), it yields   for any $r\in\bN$
\begin{gather}
\label{eq27:proof-theorem4-empirical}
\mathrm{U}_{\tau_{r}}^{(z_r(u),q)}\big(n,\mh^{(k)}\big)\leq \cU^{(u,q)}_{\mathbf{r}}(n,\mh),\;\;\forall(n,\mh)\in\widehat{\mathbf{H}}(r).
\end{gather}
Note that for any $r\in\bN$
$$
\underline{A}\vee\underline{B}\leq 2C_{D,\boldsymbol{b}}\left[\sqrt{(\mathbf{n_1})^{-1}F_{\mathbf{n_2}}\underline{G}_\mathbf{n} }\vee\left( (\mathbf{n_1})^{-1}\ln^{\beta}{(\mathbf{n_2})}\underline{G}_\mathbf{n} \right) \right]\left[\underline{G}_\mathbf{n}[\tau_r]\underline{G}^{-1}_\mathbf{n}\right],
$$
where $C_{D,\boldsymbol{b}}=\big(\sqrt{2C_D}\vee[\gamma\boldsymbol{c_b}]\big)\vee \big[(2/3)\big(C_D\vee 8e\big)\big]$.
We get from (\ref{eq24:proof-theorem4-empirical}) and (\ref{eq27:proof-theorem4-empirical}), similarly to (\ref{eq26:proof-theorem4-empirical}),
\begin{gather*}
%\label{eq28:proof-theorem4-empirical}
\bE_\rf\left\{\sup_{(n,\mh)\in\widehat{\mathbf{H}}(r)}\Big[|\xi_{\mh}(n)|-\cU^{(u,q)}_{\mathbf{r}}(n,\mh)\Big]\right\}^q_+
%\\*[2mm]&&
\leq \frac{K_q\left[\sqrt{(\mathbf{n_1})^{-1}F_{\mathbf{n_2}}\underline{G}_\mathbf{n} }\vee\left( (\mathbf{n_1})^{-1}\ln^{\beta}{(\mathbf{n_2})}\underline{G}_\mathbf{n} \right) \right]^q e^{-u}}{[1+|r-\mathbf{r}|]^{-2}} ,
\end{gather*}
where $K_q=2^{(7q/2)+3}3^{q+4}\Gamma(q+1)(C_{D,\boldsymbol{b}})^{q}$.

Summing up the right hand side of  the last inequality   over $r$  we come to the second assertion of the theorem.

\epr

\subsection{Proof of Theorem \ref{th:LIL-nonasym}}

For any $l\in\bN^*$ set $n_l=\mathbf{j}2^l$, $\mathbf{N}_l=\big\{n_l,n_{l}+1,\ldots, n_{l+1}\big\}$ and let
$$
\zeta_{\mathbf{j}}=\sup_{n\geq \mathbf{j}}\;\sup_{\mh^{(k)}\in\overline{\mH}_1^k(n,a)}\Bigg[\frac{\sqrt{n}\;\eta_{\mh^{(k)}}(n)}
{\sqrt{G_\infty\big(\mh^{(k)}\big)\ln{\big(1+\ln{(n)}\big)}}}\Bigg].
$$
We obviously have
\begin{eqnarray*}
\bP_\rf\left\{\zeta_{\mathbf{j}}\geq \Upsilon\right\}
&\leq&\sum_{l=1}^{\infty} \bP_\rf\left\{\sup_{n\in \mathbf{N}_l}\sup_{\mh^{(k)}\in\overline{\mH}_1^k(n,a)}\Bigg[\frac{\sqrt{n}\;\eta_{\mh^{(k)}}(n)}
{\sqrt{G_\infty\big(\mh^{(k)}\big)\ln{\big(1+\ln{(n)}\big)}}}\Bigg]\geq \Upsilon\right\}
\\*[2mm]
&=&\sum_{l=1}^{\infty}\bP_\rf\left\{\sup_{n\in \mathbf{N}_l}\sup_{\mh^{(k)}\in\overline{\mH}_1^k(n,a)}
\left[\eta_{\mh^{(k)}}(n)-\Upsilon\sqrt{n^{-1}G_\infty\big(\mh^{(k)}\big)\ln{\big(1+\ln{(n)}\big)}}\right]>0\right\} .
\end{eqnarray*}
Let $l\in\bN^*$ be fixed and later on $\boldsymbol{\Upsilon}_r,\,r=1,2,3$ denote the constants independent on $l$ and $n$.

Note that in view of (\ref{eq1:LIL}), (\ref{eq2:LIL}) and (\ref{eq22:LIL}) for any $n\in\mathbf{N}_l$
\begin{eqnarray*}
\cV_{0}^{(2\ln{\left(1+\ln{(n_l)}\right)})}(n,\mh)&\leq&\lambda_1\sqrt{\Big(\mathbf{F}n^{-1}\Big) G_\infty\big(\mh^{(k)}\Big(P_n+
2\ln{\left\{1+\left|\ln{(\mathbf{F})}\right|\right\}}
+2\ln{\left(1+\ln{(n)}\right)}\Big)}
\\
&&\hskip-0.3cm +
\lambda_2\Big(n^{-1}\ln^{\boldsymbol{b}}{(n)}\Big) G_\infty\big(\mh^{(k)}\big)
\Big(P_n+2\ln{\left\{1+\left|\ln{(\mathbf{F})}\right|\right\}}+2\ln{\left(1+\ln{(n)}\right)}\Big);
\end{eqnarray*}
where we have put
$$
P_n=(36k\delta^{-2}_*+6)\ln{\left(1+\mb\ln{\left(2n\right)}\right)}+36N\delta^{-2}_*\ma\ln{\left(1+\ln{\left(2n^{\mb}\mathfrak{c}\right)}\right)}
+18C_{N,R,m,k}(\boldsymbol{b}).
$$
Hence,  for any $n\in\mathbf{N}_l$ and any $\mh\in\widetilde{\mH}(n)$
$$
\cV_{0}^{(2\ln{\left(1+\ln{(n_l)}\right)})}(n,\mh)\leq \Upsilon_1 \sqrt{\frac{G_\infty\big(\mh^{(k)}\big)\ln{\big(1+\ln{(n)}\big)}}{n}}+
\Upsilon_2 \left[\frac{G_\infty\big(\mh^{(k)}\big)\ln^{\boldsymbol{b}}{(n)}
\ln{\big(1+\ln{(n)}\big)}}{n}\right].
$$
Since $\boldsymbol{b}>1$ can be arbitrary chosen and $a>2$ let $1<\boldsymbol{b}<a/2$.  It yields  for any $n\geq 3$ and any
$\mh^{(k)}\in\overline{\mH}_1^k(n,a)$
$$
\frac{G_\infty\big(\mh^{(k)}\big)\ln^{\boldsymbol{b}}{(n)}
\ln{\big(1+\ln{(n)}\big)}}{n}\leq \Upsilon_3 \sqrt{\frac{G_\infty\big(\mh^{(k)}\big)\ln{\big(1+\ln{(n)}\big)}}{n}}
$$
and, therefore, putting $\Upsilon=\Upsilon_1+\Upsilon_2\Upsilon_3$  we get for any $n\in\mathbf{N}_l$
$$
\cV_{0}^{(2\ln{\left(1+\ln{(n_l)}\right)})}(n,\mh)\leq \Upsilon \sqrt{\frac{G_\infty\big(\mh^{(k)}\big)\ln{\big(1+\ln{(n)}\big)}}{n}}.
$$
Noting that right hand side of the latter inequality is independent of $\mh_{(k)}$ and applying the first assertion of  Theorem \ref{th:empiric_totaly_bounded_case} with $\widetilde{\mathbf{N}}=\mathbf{N}_l$, $\mathbf{r}=0$ and $u=2\ln{\left(1+\ln{(n_l)}\right)}$ we have
$$
\bP_\rf\left\{\zeta_{\mathbf{j}}\geq \Upsilon\right\}\leq 2419\sum_{l=1}^{\infty}\left(l+\ln{(\mathbf{j})}\right)^{-2}\leq \frac{2419}{\ln{(\mathbf{j})}}.
$$
\epr

\subsection{Proof of Theorem \ref{th:empiric-partially_totaly_bounded_case}}

$1^0.$ We start the proof with establishing some simple facts used in the sequel.

For any $\mathbf{i}\in\mathbf{I}$  let $\mn(\mathbf{i})\in\bN^*$ and  $\tilde{\pi}_j(\mathbf{i})\in\mathbf{I},\; j=1,\ldots,\mn(\mathbf{i}),$ be the pairwise disjoint collection which is determined  by the condition:
$
\mathrm{H}_{m,\mathbf{i}}\cap\mathrm{H}_{m,\mathbf{k}}=\emptyset,\quad\forall \mathbf{k}\notin \left\{\tilde{\pi}_1(\mathbf{i}),\ldots,\tilde{\pi}_{\mn(\mathbf{i})}(\mathbf{i})\right\}.
$
First we  have
\begin{equation*}
\label{eq0000:proof-theorem5-empirical}
1\leq \mn(\mathbf{i})\leq \mn,\quad\forall \mathbf{i}\in\mathbf{I},
\end{equation*}
and we  always put  $\tilde{\pi}_{\mn(\mathbf{i})}(\mathbf{i})=\mathbf{i}$.
It yields, in particular, that we can construct another collection of indices $\pi(\mathbf{i}):=\left\{\pi_j(\mathbf{i})\in\mathbf{I},\; j=\overline{1,\mn}\right\}$  given by
$$
\pi_j(\mathbf{i})=\left\{\begin{array}{cc}
\tilde{\pi}_{j}(\mathbf{i}), &\quad 1\leq j\leq \mn(\mathbf{i});\\
\mathbf{i}, &\quad \mn(\mathbf{i})+1\leq j\leq \mn.
\end{array}\right.
$$
%the definition of the collection $\left\{\pi_j(\mathbf{i})\in\mathbf{I},\; j=1,\ldots,\mn(\mathbf{i})\right\}$ allows us to
Note also that for any $1\leq j\leq\mn$
\begin{equation}
\label{eq00001000:proof-theorem5-empirical}
\text{card}\Big(\left\{\mathbf{i}\in\mathbf{I}:\;\;\pi_j(\mathbf{i})=\mathbf{p}\right\}\Big)\leq \mn,\quad\forall \mathbf{p}\in\mathbf{I}.
\end{equation}
Indeed, if $\text{card}\Big(\left\{\mathbf{i}\in\mathbf{I}:\;\;\pi_j(\mathbf{i})=\mathbf{p}\right\}\Big)\geq \mn+1$ for some $\mathbf{p}\in\mathbf{I}$, then $$\text{card}\Big(\left\{\mathbf{i}\in\mathbf{I}:\;\;\mathrm{H}_{m,\mathbf{p}}\cap\mathrm{H}_{m,\mathbf{i}}\neq\emptyset\right\}\Big)\geq\mn+1,
$$
that contradicts to the definition of a $\mn$-totally bounded cover.
For any $\mathbf{i}\in\mathbf{I}$ define
$$
\mathrm{H}_m(\mathbf{i})=\bigcup_{
\mathbf{k}\in\mathbf{I}:\;\mathrm{H}_{m,\mathbf{k}}\cap\mathrm{H}_{m,\mathbf{i}}\neq\emptyset}
\;\bigcup_{\mathbf{j}\in\mathbf{I}:\;\mathrm{H}_{m,\mathbf{j}}\cap\mathrm{H}_{m,\mathbf{k}}\neq\emptyset}=\bigcup_{
l=1}^\mn
\;\bigcup_{j=1}^\mn\mathrm{H}_{m,\pi_j\big(\pi_l(\mathbf{i})\big)}.
$$
First we note that the definition of the set $\mH_m(\cdot)$ implies the following inclusion:  for any $\mathbf{i}\in\mathbf{I}$
\begin{equation}
\label{eq0001:proof-theorem5-empirical}
\mH_m(\mh_m)\subseteq\mathrm{H}_m(\mathbf{i}),
\quad\forall\mh_m\in\mathrm{H}_{m,\mathbf{i}}.
\end{equation}
Next,  taking into account that
$
\sum_{\mathbf{q}\in\mathbf{I}}\mathrm{1}_{\mathrm{H}_{m,\mathbf{q}}}(\mh_m)\leq \mn$ for any $\mh_m\in\mH_m
$
in view of the definition of a $\mn$-totally bounded cover,
we obtain in view of (\ref{eq00001000:proof-theorem5-empirical})
\begin{eqnarray}
\label{eq0002:proof-theorem5-empirical}
\sum_{\mathbf{i}\in\mathbf{I}}\mathrm{1}_{\mathrm{H}_m(\mathbf{i})}(\mh_m)
&\leq&
\sum_{\mathbf{i}\in\mathbf{I}}\sum_{j=1}^\mn\sum_{l=1}^\mn\mathrm{1}_{\mathrm{H}_{m,\pi_j\big(\pi_l(\mathbf{i})\big)}}(\mh_m)=
\sum_{j=1}^\mn\sum_{l=1}^\mn\sum_{\mathbf{p}\in\mathbf{I}}\;\sum_{\mathbf{i}: \pi_l(\mathbf{i})=\mathbf{p}}\mathrm{1}_{\mathrm{H}_{m,\pi_j(\mathbf{p})}}(\mh_m)
\nonumber\\
&\leq&\mn\sum_{j=1}^\mn\sum_{l=1}^\mn\sum_{\mathbf{p}\in\mathbf{I}}\mathrm{1}_{\mathrm{H}_{m,\pi_j(\mathbf{p})}}(\mh_m)=
\sum_{j=1}^\mn\sum_{l=1}^\mn\sum_{\mathbf{q}\in\mathbf{I}}\;\sum_{\mathbf{p}: \pi_j(\mathbf{p})=\mathbf{q}}\mathrm{1}_{\mathrm{H}_{m,\mathbf{q}}}(\mh_m)
\nonumber\\
&\leq& \mn^2\sum_{j=1}^\mn\sum_{l=1}^\mn\sum_{\mathbf{q}\in\mathbf{I}}\mathrm{1}_{\mathrm{H}_{m,\mathbf{q}}}(\mh_m)\leq \mn^5,\quad \forall \mh_m\in\mH_m.
\end{eqnarray}

Define finally for any $\mathbf{i}\in\mathbf{I}$
$$
\rf_{\mathbf{i}}:=\mathbf{n_1}^{-1}\sum_{i=1}^{\mathbf{n_2}}\int_{\mathrm{H}_m(\mathbf{i})}f_{1,i}(x)\nu_1\big(\rd x\big)
$$
and let $\mathbf{I_1}=\left\{\mathbf{i}\in\mathbf{I}:\;\; \rf_{\mathbf{i}}\geq (\mathbf{n_1})^{-v}\right\}$ and $\mathbf{I_2}=\mathbf{I}\setminus\mathbf{I_1}$.

\smallskip

$2^{0}.$\;
Let us fix $\mathbf{i}\in\mathbf{I_1}$ and for any $n\geq 1$ define $\mathrm{H}_\mathbf{i}(n):=\widetilde{\mH}_1^k(n)\times\mH_{k+1}^{m-1}\times\mathrm{H}_{m,\mathbf{i}},\;\mathbf{i}\in\mathbf{I}$.
 The idea is to apply  Theorem \ref{th:empiric_totaly_bounded_case} to  $\left\{\mathrm{H}_\mathbf{i}(n),\;n\geq 1\right\}$ that is possible in view of Assumptions \ref{ass:partially-bounded-case} ($\mathbf{i}$) and \ref{ass:dependence-on-n}.
To do it we first note that the definition of $\mathbf{I_1}$ together with (\ref{eq0001:proof-theorem5-empirical}) implies for any $n\in\widetilde{\mathbf{N}}$
$$
\mL_{n,v}\big(\mh_m\big)\geq \ln{\big(1\big/\rf_{\mathbf{i}}\big)},\;\; \forall \mh_m\in\mathrm{H}_{m,\mathbf{i}}.
$$
It yields for any $n\in\widetilde{\mathbf{N}}$ and  any $\mh\in\mathrm{H}_\mathbf{i}(n)$
$$
\widetilde{\cV}^{(v,z)}_{\mathbf{r}}(n,\mh)\geq \cV^{(u)}_{\mathbf{r}}(n,\mh),\qquad \widetilde{\cU}^{(v,z,q)}_{\mathbf{r}}(n,\mh)\geq \cU^{(u,q)}_{\mathbf{r}}(n,\mh),
$$
where $u=\ln{\left(1/\rf_{\mathbf{i}}\right)}+z$. We deduce from Theorem \ref{th:empiric_totaly_bounded_case} for any $\mathbf{i}\in\mathbf{I_1}$
\begin{eqnarray}
\label{eq1:proof-theorem5-empirical}
&&\bP_\rf\left\{\sup_{n\in\widetilde{\mathbf{N}}}\;\sup_{\mh\in\mathrm{H}_\mathbf{i}(n)}
\Big[\big|\xi_\mh(n)\big|-\widetilde{\cV}_{\mathbf{r}}^{(v,z)}(n,\mh)\Big]\geq 0\right\}
\leq  2419\;\rf_{\mathbf{i}}\;e^{-z};
\\*[2mm]
\label{eq2:proof-theorem5-empirical}
&&\bE_\rf\left\{\sup_{n\in\widetilde{\mathbf{N}}}\;\sup_{\mh\in\mathrm{H}_\mathbf{i}(n)}
\Big[\big|\xi_\mh(n)\big|-\widetilde{\cU}^{(v,z,q)}_{\mathbf{r}}(n,\mh)\Big]\right\}^q_+
%\\*[2mm]&&
\leq
\rf_{\mathbf{i}}\; \Lambda_q\big(\mathbf{n_1},\mathbf{n_2}\big)e^{-z},
\end{eqnarray}
where we have put $\Lambda_q\big(\mathbf{n_1},\mathbf{n_2}\big)=c_q\left[\sqrt{(\mathbf{n_1})^{-1}F_{\mathbf{n_2}}\underline{G}_\mathbf{n} }
\vee\left( (\mathbf{n_1})^{-1}\ln^{\beta}{(\mathbf{n_2})}\underline{G}_\mathbf{n} \right) \right]^q$.

We have in view of (\ref{eq0002:proof-theorem5-empirical}), taking into account that $\mathbf{n_2}\leq 2\mathbf{n_1}$,
\begin{eqnarray}
\label{eq3:proof-theorem5-empirical}
&& \sum_{\mathbf{i}\in\mathbf{I}}\rf_{\mathbf{i}}=(\mathbf{n_1})^{-1}
\sum_{i=1}^{\mathbf{n_2}}\int f_{1,i}(x)\bigg[\sum_{\mathbf{i}\in\mathbf{I}}\mathrm{1}_{\mathrm{H}_m(\mathbf{i})}(x)\bigg]\nu_1\big(\rd x\big)\leq 2\mn^5 .
\end{eqnarray}
Putting $\widetilde{\mH}^{(1)}(n)=\bigcup_{\mathbf{i}\in\mathbf{I_1}}\mathrm{H}_{\mathbf{i}}(n),\;n\geq 1,$ we obtain from (\ref{eq1:proof-theorem5-empirical}), (\ref{eq2:proof-theorem5-empirical})  and (\ref{eq3:proof-theorem5-empirical})
\begin{eqnarray}
\label{eq4:proof-theorem5-empirical}
&&\bP_\rf\left\{\sup_{n\in\widetilde{\mathbf{N}}}\;\sup_{\mh\in\widetilde{\mH}^{(1)}(n)}
\Big[\big|\xi_\mh(n)\big|-\widetilde{\cV}_{\mathbf{r}}^{(v,z)}(n,\mh)\Big]\geq 0\right\}
\leq  4838\;\mn^{5}\;e^{-z};
\\*[2mm]
\label{eq44:proof-theorem5-empirical}
&&\bE_\rf\left\{\sup_{n\in\widetilde{\mathbf{N}}}\;\sup_{\mh\in\widetilde{\mH}^{(1)}(n)}
\Big[\big|\xi_\mh(n)\big|-\widetilde{\cU}^{(v,z,q)}_{\mathbf{r}}(n,\mh)\Big]\right\}^q_+
\leq
2\Lambda_q\big(\mathbf{n_1},\mathbf{n_2}\big)\mn^{5}\;e^{-z}.
\end{eqnarray}
To get (\ref{eq44:proof-theorem5-empirical}) we have used obvious equality: $\left[\sup_{\alpha}Q(\alpha)\right]_+^q=\sup_{\alpha}\left[Q(\alpha)\right]_+^q$.

\smallskip

$3^{0}.$\;
Fix $\mathbf{i}\in\mathbf{I_2}$ and note that in view of Assumption \ref{ass:as-assumption-totally-bounded} for any $n\geq 1$, any $\mh\in\mH(n)$
and $i\geq 1$
\begin{eqnarray}
\label{eq50000:proof-theorem5-empirical}
&&\bE_\rf\left|G\big(\mh,X_i\big)\right|= \bE_\rf\left\{\left|G\big(\mh,X_i\big)\right|\mathbf{1}_{\mH_m\big(\mh_m\big)}(X_{1,i})\right\}
+\bE_\rf\left\{\left|G\big(\mh,X_i\big)\right|\mathbf{1}_{\mH_m\setminus\mH_m\big(\mh_m\big)}(X_{1,i})\right\}
\\
&&\leq G_\infty\big(\mh^{(k)}\big)\left[\bP_\rf\Big\{X_{1,i}\in\mH_m\big(\mh_m\big)\Big\}+n^{-1}\right]
\leq G_\infty\big(\mh^{(k)}\big)\left[\bP_\rf\Big\{X_{1,i}\in\mathrm{H}_m(\mathbf{i})\Big\}+n^{-1}\right].
\nonumber
\end{eqnarray}
The last inequality follows from (\ref{eq0001:proof-theorem5-empirical}).
It yields for any $n\in\widetilde{\mathbf{N}}$ and any $\mh\in\mH(n)$
\begin{eqnarray}
\label{eq5:proof-theorem5-empirical}
n^{-1}\sum_{i=1}^{n}\bE_\rf\left|G\big(\mh,X_i\big)\right|&\leq&
 G_\infty\big(\mh^{(k)}\big)\left[\rf_{\mathbf{i}}+n^{-1}\right]\leq  2(\mathbf{n_1})^{-1}G_\infty\big(\mh^{(k)}\big),
\end{eqnarray}
since $\rf_{\mathbf{i}}\leq (\mathbf{n_1})^{-v}$ for any $\mathbf{i}\in\mathbf{I_2}$ and $v\geq 1$.

Introduce random events
$$
\cC_{\mathbf{i}}=\left\{\sum_{i=1}^{\mathbf{n_2}}\mathbf{1}_{\mathrm{H}_m(\mathbf{i})}(X_{1,i})\geq 2\right\},\;\;\mathbf{i}\in\mathbf{I_2}, \qquad \cC=\bigcup_{\mathbf{i}\in\mathbf{I_2}}\cC_{\mathbf{i}}.
$$
Note that if the random event $\bar{\cC}$ holds (where, as usual, $\bar{\cC}$ is complementary to $\cC$) then for any $n\in\widetilde{\mathbf{N}}$ and any $\mh\in\mH(n)$
in view of Assumption \ref{ass:as-assumption-totally-bounded} and (\ref{eq0001:proof-theorem5-empirical})
\begin{eqnarray}
\label{eq6:proof-theorem5-empirical}
n^{-1}\sum_{i=1}^{n}\left|G\big(\mh,X_i\big)\right|&\leq&
   2n^{-1}G_\infty\big(\mh^{(k)}\big)\leq 2(\mathbf{n_1})^{-1}G_\infty\big(\mh^{(k)}\big).
\end{eqnarray}
Taking into account that bounds found in  (\ref{eq5:proof-theorem5-empirical}) and  (\ref{eq6:proof-theorem5-empirical}) are independent of $\mathbf{i}$ we get for any $n\in\widetilde{\mathbf{N}}$ and any $\mh\in\widetilde{\mH}^{(2)}(n):=\widetilde{\mH}(n)\setminus\in\widetilde{\mH}^{(1)}(n)$
\begin{eqnarray*}
\big|\xi_\mh(n)\big|\mathbf{1}_{\bar{\cC}}&\leq&
   4(\mathbf{n_1})^{-1}G_\infty\big(\mh^{(k)}\big).
\end{eqnarray*}
Noting that  for any $\mh\in\boldsymbol{\mH}$, $z\geq 1$ and $n\in\widetilde{\mathbf{N}}$
\begin{eqnarray*}
\widetilde{\cV}_{\mathbf{r}}^{(v,z)}(n,\mh)&>& 8n^{-1}G_\infty\big(\mh^{(k)}\big)\geq 4(\mathbf{n_1})^{-1}G_\infty\big(\mh^{(k)}\big),
\\
\widetilde{\cU}_{\mathbf{r}}^{(v,z,q)}(n,\mh)&>& 8n^{-1}G_\infty\big(\mh^{(k)}\big)\geq 4(\mathbf{n_1})^{-1}G_\infty\big(\mh^{(k)}\big).
\end{eqnarray*}
and, therefore, if the random event $\bar{\cC}$ is realized we have
\begin{eqnarray*}
\label{eq7:proof-theorem5-empirical}
\sup_{n\in\widetilde{\mathbf{N}}}\sup_{\mh\in\widetilde{\mH}^{(2)}(n)}\left[\big|\xi_\mh(n)\big|-
\widetilde{\cV}_{\mathbf{r}}^{(v,z)}(n,\mh)\right]<0,\quad
\sup_{n\in\widetilde{\mathbf{N}}}\sup_{\mh\in\widetilde{\mH}^{(2)}(n)}\left[\big|\xi_\mh(n)\big|-
\widetilde{\cU}_{\mathbf{r}}^{(v,z,q)}(n,\mh)\right]<0.
\end{eqnarray*}
It yields, first,
\begin{eqnarray}
\label{eq8:proof-theorem5-empirical}
&&\quad\bP_\rf\left\{\sup_{n\in\widetilde{\mathbf{N}}}\sup_{\mh\in\widetilde{\mH}^{(2)}(n)}
\Big[\big|\xi_\mh(n)\big|-\widetilde{\cV}_{\mathbf{r}}^{(v,z)}(n,\mh)\Big]\geq 0\right\}
\leq \bP_\rf\left\{\cC\right\}\leq\sum_{\mathbf{i}\in\mathbf{I_2}}\bP_\rf\left\{\cC_{\mathbf{i}}\right\}.
\end{eqnarray}
Next, taking into account the trivial bound $\big|\xi_\mh(n)\big|\leq 2\overline{G}_\mathbf{n}$ for any $n\in\widetilde{\mathbf{N}}$ and any $\mh\in\mH(n)$, we get
%from (\ref{eq7:proof-theorem5-empirical})
\begin{eqnarray}
\label{eq9:proof-theorem5-empirical}
&&\quad\bE_\rf\left\{\sup_{n\in\widetilde{\mathbf{N}}}\sup_{\mh\in\widetilde{\mH}^{(2)}(n)}
\Big[\big|\xi_\mh(n)\big|-\widetilde{\cU}^{(v,z,q)}_{\mathbf{r}}(n,\mh)\Big]\right\}^q_+
%\\*[2mm]&&
\leq \left(2\overline{G}_\mathbf{n}\right)^q \bP_\rf\left\{\cC\right\}\leq \left(2\overline{G}_\mathbf{n}\right)^q\sum_{\mathbf{i}\in\mathbf{I_2}}\bP_\rf\left\{\cC_{\mathbf{i}}\right\}.
\end{eqnarray}
For any $\mathbf{i}\in\mathbf{I_2}$ put $\mathrm{p}_{i,\mathbf{i}}=\bP_f\Big\{X_{1,i}\in\mathrm{H}_m(\mathbf{i})\Big\}$.
%and note that
%$\mathrm{p}_{\mathbf{i}}\leq n\mathrm{f}_{\mathbf{i}}\leq n^{1-v}.$
%
Since $X_{1,i},\;i\geq 1,$ are independent random elements we have for any $\mathbf{i}\in\mathbf{I_2}$ and any $\lambda>0$
in view of exponential Markov inequality
$$
\bP_\rf\left\{\cC_{\mathbf{i}}\right\}\leq
\exp{\left\{-2\lambda+(e^{\lambda}-1)\sum_{i=1}^{\mathbf{n_2}}\mathrm{p}_{i,\mathbf{i}}\right\}}=
\exp{\left\{-2\lambda+\mathbf{n_1}(e^{\lambda}-1)\mathrm{f}_{\mathbf{i}}\right\}}.
$$
Minimizing the right hand side in $\lambda$ we obtain for any $\mathbf{i}\in\mathbf{I_2}$
$$
\bP_\rf\left\{\cC_{\mathbf{i}}\right\}\leq (e/2)^{2}(\mathbf{n_1}\mathrm{f}_{\mathbf{i}})^2\leq 2\mathrm{f}_{\mathbf{i}}\;\mathbf{n_1}^{2-v}.
$$
The last inequality follows from the definition of $\mathbf{I_2}$. We obtain finally in view of (\ref{eq3:proof-theorem5-empirical})
\begin{eqnarray}
\label{eq10:proof-theorem5-empirical}
&&\sum_{\mathbf{i}\in\mathbf{I_2}}\bP_\rf\left\{\cC_{\mathbf{i}}\right\}\leq 2\mn^5\;\mathbf{n_1}^{2-v}.
\end{eqnarray}
The assertions of the theorem follow now from (\ref{eq4:proof-theorem5-empirical}), (\ref{eq44:proof-theorem5-empirical}), (\ref{eq8:proof-theorem5-empirical}), (\ref{eq9:proof-theorem5-empirical}) and (\ref{eq10:proof-theorem5-empirical}).

\epr

\subsection{Proof of Corollary  \ref{cor:after-th:empiric-partially_totaly_bounded_case}}

To prove the assertion of the corollary it suffices to bound from above the function $\mL_{n,v}(\cdot)$. Remind that we proved, see
(\ref{eq50000:proof-theorem5-empirical}), for any $n\geq 1$, any $\mh\in\mH(n)$ and $i\geq 1$
$$
\bE_\rf\left|G\big(\mh,X_i\big)\right|
\leq G_\infty\big(\mh^{(k)}\big)\left[\bP_\rf\Big\{X_{1,i}\in\mH_m\big(\mh_m\big)\Big\}+n^{-1}\right].
$$
It yields for any $n\in \widetilde{\mathbf{N}}$ and  any $\mh\in\mH(n)$
\begin{eqnarray}
\label{eq1:cor:after-th:empiric-partially_totaly_bounded_case}
&&{F_{\mathbf{n_2}}}(\mh)\leq G_\infty\big(\mh^{(k)}\big)\left[ A_n(\mh_m)+n^{-1}\right],
\quad A_n(\mh_m)=n^{-1}\sum_{i=1}^{n}\int_{\mH_m(\mh_m)}f_{1,i}(x)\nu_1\big(\rd x\big).
\end{eqnarray}
Indeed, if $\mathbf{n_1}=\mathbf{n_2}$ then $n=\mathbf{n_2}$ and (\ref{eq1:cor:after-th:empiric-partially_totaly_bounded_case})
is obvious. If $\mathbf{n_1}\neq\mathbf{n_2}$ then $\bP_\rf\Big\{X_{1,i}\in\mH_m\big(\mh_m\big)\Big\}$ is independent of $i$
 since we supposed that $X_{1,i},\, i\geq 1$ are identically distributed.
Hence, $A_n(\cdot)$ is independent of $n$ and (\ref{eq1:cor:after-th:empiric-partially_totaly_bounded_case}) holds.
Let $n\in\widetilde{\mathbf{N}}$ be fixed and let  $\mh\in\mH(n)$ be such that $F_{\mathbf{n_2}}(\mh)\geq n^{-1/2}$.

If
$A_n(\mh_m)\leq n^{-1}$ we have
$
 G_\infty\big(\mh^{(k)}\big)\geq 2^{-1}\sqrt{n}
$
and, therefore,
$$
2v\Big|\ln{\big\{2G_\infty\big(\mh^{(k)}\big)\big\}}\Big|\geq v\ln(n)\geq \mL_{n,v}(\mh_m).
$$
If $A_n(\mh_m)> n^{-1}$ we have
$
\widehat{{F}}_{\mathbf{n_2}}(\mh)=F_{\mathbf{n_2}}(\mh)\leq 2G_\infty\big(\mh^{(k)}\big)A_n(\mh_m)
$
and, therefore,
\begin{eqnarray*}
\mL_{n,v}(\mh_m)\leq \ln{\left(A^{-1}_n(\mh_m)\right)}&\leq& \ln{\left(2G_\infty\big(\mh^{(k)}\big)\widehat{{F}}^{-1}_{\mathbf{n_2}}(\mh)\right)}
=\left|\ln{\left(2G_\infty\big(\mh^{(k)}\big)\widehat{{F}}^{-1}_{\mathbf{n_2}}(\mh)\right)}\right|
\\*[2mm]
&\leq&
\left|\ln{\left(2G_\infty\big(\mh^{(k)}\big)\right)}\right|+\left|\ln{\left(\widehat{{F}}_{\mathbf{n_2}}(\mh)\right)}\right|.
\end{eqnarray*}
Here we have also used that $A_n(\mh_m)\leq 1$. Thus, if $F_{\mathbf{n_2}}(\mh)\geq n^{-1/2}$ for any $v\geq 1$
\begin{eqnarray}
\label{eq2:cor:after-th:empiric-partially_totaly_bounded_case}
&&\mL_{n,v}(\mh_m)\leq 2v\left|\ln{\left(2G_\infty\big(\mh^{(k)}\big)\right)}\right|+\left|\ln{\left(\widehat{{F}}_{\mathbf{n_2}}(\mh)\right)}\right|.
\end{eqnarray}
If now $\mh\in\mH(n)$ be such that $F_{\mathbf{n_2}}(\mh)< n^{-1/2}$ then obviously $\widehat{{F}}_{\mathbf{n_2}}(\mh)< n^{-1/2}$ and, therefore,
$$
2v\left|\ln{\left(\widehat{{F}}_{\mathbf{n_2}}(\mh)\right)}\right|\geq v\ln(n)\geq \mL_{n,v}(\mh_m).
$$
The latter inequality together with (\ref{eq2:cor:after-th:empiric-partially_totaly_bounded_case}) yields for any $n\in\widetilde{\mathbf{N}}$,  any $\mh\in\mH(n)$ and $v\geq 1$
\begin{eqnarray}
\label{eq3:cor:after-th:empiric-partially_totaly_bounded_case}
\mL_{n,v}(\mh_m)\leq 2v\Big[\left|\ln{\left(2G_\infty\big(\mh^{(k)}\big)\right)}\right|+\left|\ln{\left(\widehat{{F}}_{\mathbf{n_2}}(\mh)\right)}\right|\Big].
\end{eqnarray}
Hence, choosing $\mathbf{r}=\ln(\mathbf{n_2})$ and  replacing $\mL_{n,v}(\cdot)$  in the expressions of
$
\widetilde{\cV}^{(v,z)}_{\mathbf{r}}(\cdot,\cdot)$ and  $\widetilde{\cU}^{(v,z,q)}_{\mathbf{r}}(\cdot,\cdot)$ by its upper bound found in (\ref{eq3:cor:after-th:empiric-partially_totaly_bounded_case})
we come to the assertion  of the corollary.

\epr

\subsection{Proof of Theorem \ref{th:LL-nonasym}}

For any $l\in\bN^*$ set $n_l=\mathbf{j}2^l$, $\mathbf{N}_l=\big\{n_l,n_{l}+1,\ldots, n_{l+1}\big\}$ and let
$$
\zeta_\mathbf{j}=\sup_{n\geq \mathbf{j}}\;\
\sup_{\mh^{(k)}\in\overline{\mH}_1^k(n,a)}\frac{\sqrt{n}\;\eta_{\mh^{(k)}}(n)}
{\sqrt{G_\infty\big(\mh^{(k)}\big)\Big[\ln{\left\{G_\infty\big(\mh^{(k)}\big)\right\}}\vee\ln{\ln{(n)}}\Big]}}.
$$
We obviously have for any $y\geq 0$
\begin{eqnarray*}
&&\bP_\rf\left\{\zeta_{\mathbf{j}}\geq \boldsymbol{\Upsilon}\right\}
\\*[2mm]
&&\leq\sum_{l=1}^{\infty}\bP_\rf\Bigg\{\sup_{n\in \mathbf{N}_l}\sup_{\mh^{(k)}\in\overline{\mH}_1^k(n,a)}
\left[\eta_{\mh^{(k)}}(n)-\boldsymbol{\Upsilon}\sqrt{n^{-1}G_\infty\big(\mh^{(k)}\big)
\Big[\ln{\left\{G_\infty\big(\mh^{(k)}\big)\right\}}\vee\ln{\ln{(n)}}\Big]}\right]\geq 0\Bigg\} .
\end{eqnarray*}
 Remind, that for any $3\leq \mathbf{n_1}\leq \mathbf{n_2}\leq 2\mathbf{n_1}$ and any  $n\in\widetilde{\mathbf{N}}$
\begin{eqnarray*}
\widehat{\cV}^{(v,z)}(n,\mh)&=&\lambda_1\sqrt{\Big(\widehat{{F}}_{\mathbf{n_2}}(\mh)n^{-1}\Big) G_\infty\big(\mh^{(k)}\big)\Big(\widehat{P}_v(\mh^{(k)}\big)+2(v+1)\big|\ln{\big\{\widehat{{F}}_{\mathbf{n_2}}(\mh)\big\}}\big|
+z\Big)}
\\
&&\hskip-0.3cm +
\lambda_2\Big(n^{-1}\ln^{\beta}{(n)}\Big) G_\infty\big(\mh^{(k)}\big)
\Big(\widehat{P}_v(\mh^{(k)}\big)+2(v+1)\big|\ln{\big\{\widehat{{F}}_{\mathbf{n_2}}(\mh)\big\}}\big|+z\Big);
%\label{eq:def-of-upper-function-moments}
\end{eqnarray*}
Let $l\in\bN^*$ be fixed and choose $v=3$ and $z=2\ln{\left(1+\ln{(n_l)}\right)}$. Later on $\boldsymbol{\Upsilon}_r,\,r=1,2,3,4$ denote the constants independent on $l$ and $n$.

We have in view of (\ref{eq1:LIL}), (\ref{eq2:LIL}) and (\ref{eq2:LL}) for any $n\in\mathbf{N}_l$ and $\mh\in\widetilde{\mH}(n)$
$$
\widehat{\cV}^{(3,\;2\ln{\left(1+\ln{(n_l)}\right)})}(n,\mh)\leq \boldsymbol{\Upsilon}_1 \sqrt{\frac{G_\infty\big(\mh^{(k)}\big)\Big[\ln{\left\{G_\infty\big(\mh^{(k)}\big)\right\}}\vee\ln{\ln{(n)}}\Big]}{n}}+
\boldsymbol{\Upsilon}_2 \left[\frac{G_\infty\big(\mh^{(k)}\big)\ln^{\boldsymbol{b}+1}{(n)}}{n}\right].
$$
To get the latter inequality we have used, first, that
$$
\widehat{{F}}_{\mathbf{n_2}}(\mh)\big|\ln{\big\{\widehat{{F}}_{\mathbf{n_2}}(\mh)\big\}}\big|\leq \sup_{x\in (0,\mathbf{F}]}
x\big|\ln(x)|=:c(\mathbf{F})<\infty, \quad \forall \mathbf{F}<\infty.
$$
Next, to get the second term, we have used that for any $n\in\mathbf{N}_l$ and  $\mh\in\widetilde{\mH}(n)$
$$
\widehat{P}_3(\mh^{(k)}\big)\leq \boldsymbol{\Upsilon}_3\ln(n),\qquad \big|\ln{\big\{\widehat{{F}}_{\mathbf{n_2}}(\mh)\big\}}\big|\leq
\max\left[\big|\ln{\big\{\mathbf{F}\big\}}\big|,\ln(n_{l+1})\right]\leq\max\left[\big|\ln{\big\{\mathbf{F}\big\}}\big|,\ln(2n)\right].
$$

Since $\boldsymbol{b}>1$ can be arbitrary chosen and $a>4$ let $1<\boldsymbol{b}<a/2-1$. It yields  for any $n\geq 3$ and any
$\mh^{(k)}\in\overline{\mH}_1^k(n,a)$
$$
\frac{G_\infty\big(\mh^{(k)}\big)\ln^{\boldsymbol{b}+1}{(n)}}{n}\leq \boldsymbol{\Upsilon}_4 \sqrt{\frac{G_\infty\big(\mh^{(k)}\big)\Big[\ln{\left\{G_\infty\big(\mh^{(k)}\big)\right\}}\vee\ln{\ln{(n)}}\Big]}{n}}
$$
and, therefore, putting $\boldsymbol{\Upsilon}=\boldsymbol{\Upsilon}_1+\boldsymbol{\Upsilon}_2\boldsymbol{\Upsilon}_4$  we get for any $n\in\mathbf{N}_l$
$$
\widehat{\cV}^{(3,\;2\ln{\left(1+\ln{(n_l)}\right)})}(n,\mh)\leq \boldsymbol{\Upsilon} \sqrt{\frac{G_\infty\big(\mh^{(k)}\big)\Big[\ln{\left\{G_\infty\big(\mh^{(k)}\big)\right\}}\vee\ln{\ln{(n)}}\Big]}{n}}.
$$
 Noting that right hand side of the latter inequality is independent of $\mh_{(k)}$ and applying the first assertion of  Corollary \ref{cor:after-th:empiric-partially_totaly_bounded_case} with $\widetilde{\mathbf{N}}=\mathbf{N}_l$ and $z=2\ln{\left(1+\ln{(n_l)}\right)}$ we obtain
$$
\bP_\rf\left\{\zeta_{\mathbf{j}}\geq \Upsilon\right\}\leq 2\mn^5\left\{2419\sum_{l=1}^{\infty}\left(l+\ln{(\mathbf{j})}\right)^{-2}
+\mathbf{j}^{-1}\sum_{l=1}^{\infty}2^{-l}\right\}\leq 2\mn^5\left\{\frac{2419}{\ln{(\mathbf{j})}}+\mathbf{j}^{-1}\right\}\leq \frac{4840\mn^5}{\ln{(\mathbf{j})}}.
$$

\epr

\section{Appendix}

\paragraph{Proof of Lemma \ref{lem:measurability}}

$1^0.$ We start the proof with the following simple fact. Let $\widetilde{\boldsymbol{\mT}}$ be an arbitrary subset of $\boldsymbol{\mT}$. Then
\begin{equation}
\label{eq1:proof-lemma-measurab}
\sup_{\mt\in\widetilde{\boldsymbol{\mT}}}\zeta(\mt,\cdot)\;\;is\;\;\boldsymbol{\mB}-measurable.
\end{equation}
Indeed, since $\boldsymbol{\mT}$ is totally bounded $\widetilde{\boldsymbol{\mT}}$ is totally bounded as well. Denote by $\widehat{\boldsymbol{\mT}}$ the union of $2^{-l}$-nets, $l\geq 0$,
in $\widetilde{\boldsymbol{\mT}}$. Let $\boldsymbol{\Omega}_0=\big\{\omega\in\boldsymbol{\Omega}:\;\; \zeta(\cdot,\omega)\;\text{is continuous}\big\}$ and let $\overline{\boldsymbol{\Omega}}_0$ be the complementary to $\boldsymbol{\Omega}_0$.
We have for any $x\in\bR$
$$
\bigg\{\omega\in\boldsymbol{\Omega}:\;\;\sup_{\mt\in\widetilde{\boldsymbol{\mT}}}\zeta(\mt,\omega)\leq x\bigg\}\cap\boldsymbol{\Omega}_0=\bigg\{\omega\in\boldsymbol{\Omega}:\;\;\sup_{\mt\in\widehat{\boldsymbol{\mT}}}\zeta(\mt,\omega)\leq x\bigg\}\cap\boldsymbol{\Omega}_0\in\boldsymbol{\mB}
$$
since $\widehat{\boldsymbol{\mT}}$ is countable dense subset of $\boldsymbol{\mT}$.
It remain to note that
 $
 \left\{\omega\in\boldsymbol{\Omega}:\;\;\sup_{\mt\in\widetilde{\boldsymbol{\mT}}}\zeta(\mt,\omega)\leq x\right\}\cap\overline{\boldsymbol{\Omega}}_0\in\boldsymbol{\mB}
 $
  since $\boldsymbol{\mathrm{P}}\big(\overline{\boldsymbol{\Omega}}_0\big)=0$ and the considered   probability space is complete.

\vskip0.1cm

$2^0.$ Set
$\mZ(n,k)=\big\{\mz\in\mZ:\;\; g(\mz)\in [k/n,(k+1)/n]\Big\}, \;
n\in\bN^*,\;k\in\bZ,$
and let $\mathbf{K}(n)\subseteq\bZ, \;n\in\bN^*,$ be defined from the relation if $k\in\mathbf{K}(n)\;\Leftrightarrow\;\mZ(n,k)\neq\emptyset$.
Put also $\boldsymbol{\mT}(k,n)=\cup_{\mz\in \mZ(k,n)}\mT_\mz$ and define
\begin{eqnarray*}
&&\xi_{k,n}(\omega)=\sup_{\mt\in\boldsymbol{\mT}(k,n)}\zeta(\mt,\omega)-(k+1)/n,\qquad \xi_{n}(\omega)=\sup_{k\in\mathbf{K}(n)}\xi_{k,n}(\omega).
\\
&&\eta_{k,n}(\omega)=\sup_{\mz\in\mZ(k,n)}\Big[\sup_{\mathfrak{t}\in\boldsymbol{\mT}_\mz}\zeta(\mt,\cdot)-g(\mz)\big],\qquad \eta(\omega)=\sup_{\mz\in\mZ}\Big[\sup_{\mathfrak{t}\in\boldsymbol{\mT}_\mz}\zeta(\mt,\omega)-g(\mz)\big].
\end{eqnarray*}
Some remarks are in order. First, the definition of $\mZ(k,n)$ implies that for any $k\in\mathbf{K}(n),\;n\in\bN^*$
\begin{equation}
\label{eq2:proof-lemma-measurab}
\xi_{k,n}(\cdot)\leq \eta_{k,n}(\cdot)\leq \xi_{k,n}(\cdot)+n^{-1}.
\end{equation}
Next, taking into account that $\mZ=\cup_{k\in\mathbf{K}(n)}\mZ(k,n)$ for any $n\in\bN^*$ we have
\begin{equation}
\label{eq3:proof-lemma-measurab}
\eta(\cdot)=\sup_{k\in\mathbf{K}(n)}\eta_{k,n}(\cdot),\quad \forall n\in\bN^*.
\end{equation}
We obtain from (\ref{eq2:proof-lemma-measurab}) and (\ref{eq3:proof-lemma-measurab}) that for any $n\in\bN^*$
$$
0\leq \eta(\cdot)-\xi_{n}(\cdot)\leq \sup_{k\in\mathbf{K}(n)}\left[\eta_{k,n}(\cdot)-\xi_{k,n}(\cdot)\right]\leq n^{-1},
$$
and, therefore,
$
\eta(\cdot)=\lim_{n\to\infty}\xi_{n}(\cdot).
$
It remains to note that $\xi_{k,n}(\cdot)$ are $\boldsymbol{\mB}$-measurable  for any $k\in\mathbf{K}(n),\;n\in\bN^*$ in view of
(\ref{eq1:proof-lemma-measurab}), that implies obviously that $\xi_{n}(\cdot)$ is $\boldsymbol{\mB}$-measurable  for any $n\in\bN^*$.
Thus, $\eta(\cdot)$ is $\boldsymbol{\mB}$-measurable as a pointwise limit of $\boldsymbol{\mB}$-measurable functions.

\epr

\paragraph{Proof of Lemma \ref{lem:loc-empiric-pointwise}}

$1^0.$ Remind that
\begin{eqnarray}
\label{eq:jjjjjjjjjj}
&& G_\infty(r)= V^{-1}_r\|g\|_\infty \|K\|_\infty,\qquad
\underline{G}_n=V^{-1}_{r^{(\max)}(n)}\|g\|_\infty \|K\|_\infty,\;\; n\geq 1.
\end{eqnarray}
Hence we have have  for  any $l=\overline{1,d}$ and any $\mh_l:=r_l\in\Big[ r_l^{(\min)}(n), r_l^{(\max)}(n)\Big]$
\begin{gather*}
\label{eq100:proof-theorem55-empirical}
G_{l,n}(r_l)=\|g\|_\infty \|K\|_\infty \left[V_{r^{(\min)}(n)}\right]^{-1}\left[r_l^{(\min)}(n)\Big/r_l\right]^{\gamma_l};
\\*[2mm]
\label{eq101:proof-theorem55-empirical}
 \underline{G}_{l,n}=\|g\|_\infty \|K\|_\infty \left[V_{r^{(\min)}(n)}\right]^{-1}\left[r_l^{(\min)}(n)\Big/r^{(\max)}_l(n)\right]^{\gamma_l}.
\end{gather*}
Thus, we get for any $n\geq 1$, $r\in\cR(n)$ and for any $j=\overline{1,d}$
$$
\frac{G_\infty(r)}{\underline{G}_n}=\prod_{l=1}^{d}\left[\frac{r^{(\max)}_l}{r_l}\right]^{\gamma_l}\geq \left[\frac{r^{(\max)}_j}{r_j}\right]^{\gamma_j}= \frac{G_{j,n}(r_j)}{\underline{G}_{j,n}}.
 $$
We conclude that Assumption \ref{ass:bounded_case} ($\mathbf{i}$) is fulfilled.

\smallskip

$2^0.$ Remind  that for any $r,r^\prime\in\cR(n)$
\begin{equation}
\label{eq333:proof-theorem55-empirical}
\varrho_n^{(d)}\Big(r,r^\prime\Big):=\max_{l=\overline{1,d}}\mm_0
\Big(G_{l,n}(r_l),G_{j,n}(r^\prime_l)\Big)=\max_{l=\overline{1,d}}\gamma_l\Big|\ln(r_l)-\ln(r^\prime_l)\Big|=:\varrho^{(d)}\Big(r,r^\prime\Big),
\end{equation}

$3^0.$ Set $\big\|K_{r}-K_{r^\prime}\big\|_\infty=\sup_{z\in\bR^d}\big|K_{r}(z)-K_{r^\prime}(z)\big|$ and note that
for any $x\in\bX_1^d\times\bX_{d+1}$ and for any $\mh=\big(r,\mz,y^{(d)}\big)$,
$\mh^{\prime}=\big(r^\prime,\mz^\prime,z^{(d)}\big)$
\begin{eqnarray*}
\big|G(\mh,x)-G(\mh^\prime,x)\big|
&\leq& \|g\|_\infty  \big\|K_{r}-K_{r^\prime}\big\|_\infty+\|K\|_\infty\left[V_r\vee V_{r^\prime}\right]^{-1}
\big|g\big(\mz,x\big)-g\big(\mz^\prime,x\big)\big|,
\\*[1mm]
&+& \|g\|_\infty V^{-1}_{r^{\prime}}\left|K\left(\vec{\rho}\big(x^{(d)},y^{(d)}\big)/r^\prime\right)-K\left(\vec{\rho}\big(x^{(d)},z^{(d)}\big)/r^\prime\right)\right|
\\
&\leq& \|g\|_\infty  \big\|K_{r}-K_{r^\prime}\big\|_\infty+
L_\alpha\|K\|_\infty \left[V_r\vee V_{r^\prime}\right]^{-1} \big[\md(\mz,\mz^\prime)\big]^{\alpha}
\\*[1mm]
&+& \|g\|_\infty V^{-1}_{r^{\prime}}\left|K\left(\vec{\rho}\big(x^{(d)},y^{(d)}\big)/r^\prime\right)-K\left(\vec{\rho}\big(x^{(d)},z^{(d)}\big)/r^\prime\right)\right|.
\end{eqnarray*}
The get the last inequality we have  used Assumption \ref{ass:on-fucnctions-g-and-K} ($\mathbf{ii}$). Using  Assumption \ref{ass:on-fucnctions-g-and-K} ($\mathbf{i}$) we have
\begin{eqnarray*}
\left|K\left(\vec{\rho}\big(x^{(d)},y^{(d)}\big)/r^\prime\right)-K\left(\vec{\rho}\big(x^{(d)},z^{(d)}\big)/r^\prime\right)\right|
&\leq& L_1\max_{l=\overline{1,d}}\left[(r_l^\prime)^{-1}\left|\rho_l\big(x_l,y_l\big)-\rho_l\big(x_l,z_l\big)\right|\right]
\\
&\leq& L_1\max_{l=\overline{1,d}}\left[(r_l^\prime)^{-1}\rho_l\big(y_l,z_l\big)\right].
\end{eqnarray*}
To get the last inequality we have taken into account that $\rho_l,\;l=\overline{1,d},$ are semi-metrics. Note also that
$(r_l^\prime)^{-1}\leq V^{-1}_{r^{\prime}}$ for any $l=\overline{1,d}$, since $r_l^\prime\leq 1$ and we obtain
\begin{eqnarray}
\label{eq1:proof-theorem6-empirical}
\left|K\left(\vec{\rho}\big(x^{(d)},y^{(d)}\big)/r^\prime\right)-K\left(\vec{\rho}\big(x^{(d)},z^{(d)}\big)/r^\prime\right)\right|
&\leq& L_1V^{-1}_{r^{\prime}}\rho^{(d)}\big(y^{(d)},z^{(d)}\big),
\end{eqnarray}
where we have put $\rho^{(d)}=\max_{l=\overline{1,d}}\rho_l$.
Obviously,
\begin{equation}
\label{eq5:proof-theorem55-empirical}
\big\|K_{r}-K_{r^\prime}\big\|_\infty\leq \|K\|_\infty\left|V_r^{-1}- V_{r^\prime}^{-1}\right|+
\left[V_r\vee V_{r^\prime}\right]^{-1}\left\|K\big(\cdot/r\big)-K\big(\cdot/r^\prime\big)\right\|_\infty.
\end{equation}
We have in view of Assumption \ref{ass:on-fucnctions-g-and-K} ($\mathbf{i}$) and  (\ref{eq333:proof-theorem55-empirical})
\begin{eqnarray*}
&&\left\|K\big(\cdot/r\big)-K\big(\cdot/r^\prime\big)\right\|_\infty\leq L_1 \sup_{u\in\bR^d} \max_{l=\overline{1,d}}\left[\frac{|u_l|\left|1/r_l-1/r_l^\prime\right|}{1+|u_l|\left(1/r_l\wedge1/r_l^\prime\right)}\right]
\leq L_1\max_{l=\overline{1,d}}\left[\frac{r_l\vee r_l^\prime}{r_l\wedge r_l^\prime}-1\right]
\\*[2mm]
&&=L_1\Big[\exp{\Big\{\max_{l=\overline{1,d}}\Big|\ln(r_l)-\ln(r^\prime_l)\Big|\Big\}}-1\Big]
\leq L_1\Big[\exp{\Big\{\gamma^{-1}\varrho^{(d)}\Big(r,r^\prime\Big)\Big\}}-1\Big],
\end{eqnarray*}
we we have put $\gamma=\min[\gamma_1,\ldots,\gamma_d]$.
 Moreover, we obviously have  for any $r,r^\prime\in (0,1]^{d}$
$$
\frac{V_r\vee V_{r^\prime}}{V_r\wedge V_{r^\prime}}\leq\frac{V_{r\vee r^\prime}}{V_{r\wedge r^\prime}} =  \exp{\left\{
\sum_{l=1}^{d}\Big|\ln{\big(r_{l}\big)}-\ln{\big(r^{\prime}_{l}\big)}\Big|
\right\}}\leq \exp{\left\{
d\varrho^{(d)}\big(r,r^\prime\big)
\right\}}.
$$
Thus, we finally obtain from (\ref{eq5:proof-theorem55-empirical})
\begin{eqnarray*}
\label{eq04:proof-theorem55-empirical}
&&\hskip-1cm\big\|K_{r}-K_{r^\prime}\big\|_\infty\leq\left[V_r\vee V_{r^\prime}\right]^{-1}
\bigg[\|K\|_\infty\left(\exp{\left\{
d\varrho^{(d)}\big(r,r^\prime\big)
\right\}}-1\right)
%\nonumber\\
%&&\hskip4cm
+L_1\Big[\exp{\Big\{\gamma^{-1}\varrho^{(d)}\Big(r,r^\prime\Big)\Big\}}-1\Big]
\bigg].
\end{eqnarray*}
This yields together with (\ref{eq1:proof-theorem6-empirical}) for any $\mh=\big(r,\mz,y^{(d)}\big)$ and
$\mh^{\prime}=\big(r^\prime,\mz^\prime,z^{(d)}\big)$
\begin{eqnarray}
\label{eq600:proof-theorem55-empirical}
&&\sup_{x\in\bX_1^d\times\bX_{d+1}}\big|G(\mh,x)-G(\mh^\prime,x)\big|
\\
&&\leq \|g\|_\infty \left[V_r\vee V_{r^\prime}\right]^{-1}
\bigg[\|K\|_\infty\left(\exp{\left\{
d\varrho^{(d)}\big(r,r^\prime\big)
\right\}}-1\right)
%\nonumber\\
%&&\hskip4cm
+L_1\Big[\exp{\Big\{\gamma^{-1}\varrho^{(d)}\Big(r,r^\prime\Big)\Big\}}-1\Big]\bigg]
\nonumber\\
&&\hskip3.6cm +L_\alpha\|K\|_\infty \left[V_r\vee V_{r^\prime}\right]^{-1} \big[\md(\mz,\mz^\prime)\big]^{\alpha}
+
L_1\|g\|_\infty V^{-2}_{r^{\prime}}\rho^{(d)}\big(y^{(d)},z^{(d)}\big)
%\nonumber\\
\nonumber\\
&&\leq \|g\|_\infty \|K\|_\infty\left[V_r\wedge V_{r^\prime}\right]^{-1}
\bigg[D_0\left(
\varrho^{(d)}\right)+D_{d+1}\big(\varrho_{d+1}\big)+\left[V_r\wedge V_{r^\prime}\right]^{-1}D_{d+2}\Big(\rho^{(d)}\big(y^{(d)},z^{(d)}\big)\Big)\bigg],
\nonumber
\end{eqnarray}
where we have put
 $\varrho_{d+1}=[\md]^{\alpha},\;\;D_{d+1}(z)=\big(L_\alpha/ \|g\|_\infty\big) z,\;\;D_{d+2}(z)=\big(L_1/ \|K\|_\infty\big) z\;$ and
$$
D_0(z)=\exp{\{dz\}}-1+(L_1/\|K\|_\infty) \Big(\exp{\left\{\gamma^{-1}z\right\}}-1\Big).
$$
Putting $L_{d+1}(z)=z$ and $L_{d+2}(z)=z^2$
we obtain from (\ref{eq:jjjjjjjjjj}) and (\ref{eq600:proof-theorem55-empirical}) for any $\mh=\big(r,\mz,y^{(d)}\big)$ and
$\mh^{\prime}=\big(r^\prime,\mz^\prime,z^{(d)}\big)$
\begin{eqnarray*}
%\label{eq600:proof-theorem6-empirical}
&&\sup_{x\in\bX_1^d\times\bX_{d+1}}\big|G(\mh,x)-G(\mh^\prime,x)\big|\leq G_\infty(r)\vee G_\infty(r^\prime)D_0\Big(\varrho^{(d)}\big(r,r^\prime\big)\Big)
\nonumber\\
&&+
L_{d+1}\left(G_\infty(r)\vee G_\infty(r^\prime)\right)
D_{d+1}\Big(\varrho_{d+1}\big(\mz,\mz^\prime\big)\Big)+L_{d+2}\left(G_\infty(r)\vee G_\infty(r^\prime)\right)D_{d+2}\left(\rho^{(d)}\big(y^{(d)},z^{(d)}\big)\right).
\end{eqnarray*}

We conclude that Assumption \ref{ass:bounded_case} ($\mathbf{ii}$) is fulfilled.  It remains to note that if $\bar{X}_1^d$ consists of a single
element then last summand in  the right hand side of the latter inequality disappears that correspond formally to $L_{d+2}\equiv 0$. This completes the proof of the lemma.

\epr

\paragraph{Proof of Lemma \ref{lem:bound-for-F_n2}}

In view of (\ref{eq:ass-boundness-of densities}) for any $r\in (0,1]^d$
\begin{eqnarray}
\label{eq1:proof-lemma-bound-for-F_n2}
&&F_{\mathbf{n_2}}\Big(r,\bar{x}^{(d)}\Big)
\leq \mathrm{f}_\infty\|g\|_\infty \int_{\bX_1^{d}}\left|K_{r}\left(\vec{\rho}\big(x^{(d)},\bar{x}^{(d)}\big)\right)\right|\mu^{(d)}(\rd x^{(d)})=:
\mathrm{f}_\infty\|g\|_\infty \cI_r.
\end{eqnarray}
Denote for any $l=\overline{1,d}$
$$
\mR_l(k_l,r_l)=\bB_l\Big(2^{k_l+1}r_l,\bar{x}_l\Big)\setminus\bB_l\Big(2^{k_l}r_l,\bar{x}_l\Big),\quad \mR_l(0,r_l)=\bB_l\Big(r_l,\bar{x}_l\Big),\;\;
k_l\in\bN.
$$
and  for any multi-index $\mathbf{k}=(k_1,\ldots,k_d)\in\bN^d$ set
$
\mR_{\mathbf{k},r}=\Pi_1(k_1,r_1)\times\cdots\times\mR_d(k_d,r_d).
$.
We get in view of Assumption \ref{ass:after-th-empirical-product-general} that
$
\bX_1^{d}=\bigcup_{\mathbf{k}\in\bN^d}\mR_{\mathbf{k},r}
$
for any $r\in (0,1]^d$ and, therefore,
$$
\cI_r=\sum_{\mathbf{k}\in\bN^d}\int_{\mR_{\mathbf{k},r}}\left|K_{r}\left(\vec{\rho}\big(x^{(d)},\bar{x}^{(d)}\big)\right)\right|\mu^{(d)}(\rd x^{(d)}).
$$
We note that for any $\mathbf{k}\in\bN^d$ that for any $x^{(d)}\in\mR_{\mathbf{k},r}$
$$
\left|K_{r}\left(\vec{\rho}\big(x^{(d)},\bar{x}^{(d)}\big)\right)\right|=V_r^{-1}
\left|K\left(\frac{\vec{\rho}\big(x^{(d)},\bar{x}^{(d)}\big)}{r}\right)\right|\leq
V_r^{-1}\sup_{|u|\notin\Pi_{\mathbf{t}(\mathbf{k})}}|K(u)|=V_r^{-1}\check{K}\big(\mathbf{t}(\mathbf{k})\big).
$$
where, we have put $\mathbf{t}(\mathbf{k})=\big(2^{k_1},\ldots,2^{k_d}\big)$ and where, remind, $\Pi_t=[0,t_1]\times\cdots\times[0,t_d]$, $\;t\in\bR^d_+$.

Thus, we obtain from (\ref{eq1:ass-normlization-using-gammas}) of Assumption \ref{ass:after-th-empirical-product-general} (remind that $\mu^{(d)}$ is a product measure)
\begin{eqnarray}
\cI_r&\leq& V_r^{-1}\sum_{\mathbf{k}\in\bN^d}\check{K}\big(\mathbf{t}(\mathbf{k})\big)\mu^{(d)}\left(\Pi_{\mathbf{t}(\mathbf{k})}\right)
\leq V_r^{-1}\sum_{\mathbf{k}\in\bN^d}\check{K}\big(\mathbf{t}(\mathbf{k})\big)
\left[\prod_{l=1}^d\mu_l\left(\bB_l\big(2^{k_l+1}r_l,\bar{x}_l\big)\right)\right]
\nonumber\\
&\leq&
\left[\prod_{l=1}^d2^{\gamma_l}L^{(l)}\right]\sum_{\mathbf{k}\in\bN^d}\check{K}\big(\mathbf{t}(\mathbf{k})\big)\left[\prod_{l=1}^d2^{\gamma_lk_l}\right].
\end{eqnarray}
We get finally from (\ref{eq2:ass-normlization-using-gammas}) of Assumption \ref{ass:after-th-empirical-product-general} that for any
$r\in(0,1]^d$
$$
\cI_r\leq 2^{d}L_2\prod_{l=1}^d2^{\gamma_l}L^{(l)}.
$$
The assertion of the lemma follows now from (\ref{eq1:proof-lemma-bound-for-F_n2}).

\epr

\paragraph{Proof of Lemma \ref{lem:verification-ass8}}

 Remind, that for the considered problem
$$
\mH\big(\mh_{d+2}\big)=\bX_1^d\big(\bar{x}^{(d)}\big):={\displaystyle\bigcup_{\mathbf{i}:\;\bar{x}^{(d)}
\in\mathrm{X}_{\mathbf{i}}}}\;\bigcup_{
\mathbf{k}:\;\mathrm{X}_{\mathbf{k}}\cap\mathrm{X}_{\mathbf{i}}\neq\emptyset}\mathrm{X}_{\mathbf{k}}.
$$
For any $\bar{x}^{(d)}\in\bX_1^d$ and any $\mathrm{r}>0$ denote
$
\bB_{\rho^{(d)}}\big(\mathrm{r},\bar{x}^{(d)}\big)=\left\{x^{(d)}\in\bX_1^d:\;\;\rho^{(d)}\big(x^{(d)},\bar{x}^{(d)}\big)\leq \mathrm{r} \right\}
$
where, remind, $\rho^{(d)}=\max[\rho_1,\ldots,\rho_d]$.
The following inclusion holds in view of Assumption \ref{ass:sec-sup-norm-case-general}
\begin{equation}
\label{eq2:proof-theorem6-empirical}
\bB_{\rho^{(d)}}\big(\mt,\bar{x}^{(d)}\big)\subseteq \bX_1^d\big(\bar{x}^{(d)}\big),\;\;\forall \bar{x}^{(d)}\in\bX_1^d.
\end{equation}
Indeed, suppose that $\exists y^{(d)}\in \bB_{\rho^{(d)}}\big(\mt,\bar{x}^{(d)}\big)$ such that $ y^{(d)}\notin\bX_1^d\big(\bar{x}^{(d)}\big)$. Then, the definition of $\bX_1^d\big(\bar{x}^{(d)}\big)$ implies that for any $\mathbf{p},\mathbf{q}\in\mathbf{I}$ such that $\bar{x}^{(d)}\in\mathrm{X}_{\mathbf{p}},\;\;y^{(d)}\in\mathrm{X}_{\mathbf{q}}$ necessarily
$$
\mathrm{X}_{\mathbf{p}}\cap\mathrm{X}_{\mathbf{q}}=\emptyset.
$$
Hence, in view of Assumption \ref{ass:sec-sup-norm-case-general}, $\rho^{(d)}\big(y^{(d)},\bar{x}^{(d)}\big)>\mt$ and, therefore,
$ y^{(d)}\notin \bB_{\rho^{(d)}}\big(\mt,\bar{x}^{(d)}\big)$. The obtained contradiction proves (\ref{eq2:proof-theorem6-empirical}).

%\smallskip

Note that in view of Assumption \ref{ass:on-fucnctions-g-and-K} $(\mathbf{ii})$ for any $x\in\bX_1^d\times\bX_{d+1}$ and any $\mh=\big(r,\mz,\bar{x}^{(d)}\big)$
$$
|G(\mh,x)|\leq \|g\|_\infty V_r^{-1}\left|K\left(\vec{\rho}\big(x^{(d)},\bar{x}^{(d)}\big)/r\right)\right|
$$
and, therefore, we get from (\ref{eq2:proof-theorem6-empirical}) and (\ref{eq:condition-on-r-max})
\begin{eqnarray*}
\sup_{x\in\bX_1^d\times\bX_{d+1}:\; x^{(d)}\notin\bX_1^d\big(\bar{x}^{(d)}\big)}|G(\mh,x)|&\leq& \|g\|_\infty V_r^{-1}
\sup_{ x^{(d)}\notin \bB_{\rho^{(d)}}\big(\mt,\bar{x}^{(d)}\big)}\left|K\left(\vec{\rho}\big(x^{(d)},\bar{x}^{(d)}\big)/r\right)\right|
\\
&\leq &\|g\|_\infty V_r^{-1}\sup_{r\in\cR(n)}\sup_{u\notin [0,\mt]^d}|K(u/r)|\leq \|g\|_\infty \|K\|_\infty V_r^{-1}n^{-1}
\nonumber
\\*[2mm]
&=:& n^{-1}G_\infty(r)=n^{-1}G_\infty\big(\mh^{(d)}\big).
\nonumber
\end{eqnarray*}

\epr

\bibliographystyle{agsm}

\end{document}